\documentclass[nospthms,envcountsect]{svmult} 
 \makeatletter
 \def\title{\@ifstar\s@title\s@title}
 \makeatother
\usepackage{makeidx}     
\makeindex
\usepackage{multicol}    
\usepackage{graphicx}    
\usepackage{latexsym}           
\usepackage{amsfonts}           
\usepackage[all]{xy}            
\usepackage{amsmath}            
\usepackage{amssymb}            
\usepackage{amscd}              
\usepackage[mathscr]{euscript}  
\usepackage{enumerate}          
\usepackage{verbatim}           
\def\href#1#2{#1}
\usepackage{epsf}           

\DeclareMathAlphabet\EuR{U}{eur}{m}{n}
\SetMathAlphabet\EuR{bold}{U}{eur}{b}{n}

\usepackage{amsthm}             
  \theoremstyle{plain}
    \newtheorem{theorem}                    {Theorem}       [section]
    \newtheorem{lemma}      [theorem]       {Lemma}
    \newtheorem{corollary}  [theorem]       {Corollary}
    
    \newtheorem{propositionnew}[theorem]       {Proposition} 
    
    \theoremstyle{definition}
    \newtheorem{definition} [theorem]       {Definition}
     \newtheorem{remarknew}     [theorem]       {Remark}
    \newtheorem{examplenew}    [theorem]       {Example}
    \newtheorem{conjecture} [theorem]       {Conjecture}
    \theoremstyle{remark}

\begin{document}

\setcounter{minitocdepth}{2}


\newcommand{\abs}[1]{\lvert#1\rvert} 
\newcommand{\version}[1]{\begin{center} last edited on #1 \\ last compiled on \today \end{center}} 
\newcommand{\indexnotation}[1]{\label{#1}}

\hyphenation{equi-variant}


\newcommand{\calall}{{\mathcal A} {\mathcal L}{\mathcal L}}
\newcommand{\calfcyc}{{\mathcal F}{\mathcal C}{\mathcal Y}}
\newcommand{\calcyc}{{\mathcal C}{\mathcal Y}{\mathcal C}}
\newcommand{\calfin}{{\mathcal F}{\mathcal I}{\mathcal N}}
\newcommand{\calko}{{\mathcal K}{\mathcal O}}
\newcommand{\calmfin}{{\mathcal M}{\mathcal F}{\mathcal I}}
\newcommand{\calvcyc}{{\mathcal V}{\mathcal C}{\mathcal Y}}
\newcommand{\caltr}{ \{ \! 1 \! \} }
\newcommand{\trivial}{ \{ \! 1 \! \} }

\newcommand{\cala}{{\cal A}}
\newcommand{\calb}{{\cal B}}
\newcommand{\calc}{{\cal C}}
\newcommand{\cald}{{\cal D}}
\newcommand{\cale}{{\cal E}}
\newcommand{\calf}{{\cal F}}
\newcommand{\calg}{{\cal G}}
\newcommand{\calh}{{\cal H}}
\newcommand{\cali}{{\cal I}}
\newcommand{\calj}{{\cal J}}
\newcommand{\calk}{{\cal K}}
\newcommand{\call}{{\cal L}}
\newcommand{\calm}{{\cal M}}
\newcommand{\caln}{{\cal N}}
\newcommand{\calo}{{\cal O}}
\newcommand{\calp}{{\cal P}}
\newcommand{\calq}{{\cal Q}}
\newcommand{\calr}{{\cal R}}
\newcommand{\cals}{{\cal S}}
\newcommand{\calt}{{\cal T}}
\newcommand{\calu}{{\cal U}}
\newcommand{\calv}{{\cal V}}
\newcommand{\calw}{{\cal W}}
\newcommand{\calx}{{\cal X}}
\newcommand{\caly}{{\cal Y}}
\newcommand{\calz}{{\cal Z}}


%
\newcommand{\IA}{{\mathbb A}} 
\newcommand{\IB}{{\mathbb B}} 
\newcommand{\IC}{{\mathbb C}}
\newcommand{\ID}{{\mathbb D}}  
\newcommand{\IE}{{\mathbb E}} 
\newcommand{\IF}{{\mathbb F}}
\newcommand{\IG}{{\mathbb G}}  
\newcommand{\IH}{{\mathbb H}}
\newcommand{\II}{{\mathbb I}}
\newcommand{\IK}{{\mathbb K}} 
\newcommand{\IL}{{\mathbb L}} 
\newcommand{\IM}{{\mathbb M}} 
\newcommand{\IN}{{\mathbb N}}
\newcommand{\IO}{{\mathbb O}}
\newcommand{\IP}{{\mathbb P}} 
\newcommand{\IQ}{{\mathbb Q}} 
\newcommand{\IR}{{\mathbb R}} 
\newcommand{\IS}{{\mathbb S}}
\newcommand{\IT}{{\mathbb T}} 
\newcommand{\IU}{{\mathbb U}}
\newcommand{\IV}{{\mathbb V}} 
\newcommand{\IW}{{\mathbb W}} 
\newcommand{\IX}{{\mathbb X}} 
\newcommand{\IY}{{\mathbb Y}}
\newcommand{\IZ}{{\mathbb Z}}



\newcommand{\bfA}{{\mathbf A}} 
\newcommand{\bfB}{{\mathbf B}} 
\newcommand{\bfC}{{\mathbf C}}
\newcommand{\bfD}{{\mathbf D}}  
\newcommand{\bfE}{{\mathbf E}}
\newcommand{\bff}{{\mathbf f}} 
\newcommand{\bfF}{{\mathbf F}}
\newcommand{\bfG}{{\mathbf G}}  
\newcommand{\bfH}{{\mathbf H}}
\newcommand{\bfI}{{\mathbf I}}
\newcommand{\bfJ}{{\mathbf J}} 
\newcommand{\bfK}{{\mathbf K}} 
\newcommand{\bfL}{{\mathbf L}} 
\newcommand{\bfM}{{\mathbf M}} 
\newcommand{\bfN}{{\mathbf N}}
\newcommand{\bfO}{{\mathbf O}}
\newcommand{\bfP}{{\mathbf P}} 
\newcommand{\bfQ}{{\mathbf Q}} 
\newcommand{\bfR}{{\mathbf R}} 
\newcommand{\bfS}{{\mathbf S}}
\newcommand{\bfT}{{\mathbf T}} 
\newcommand{\bfU}{{\mathbf U}}
\newcommand{\bfV}{{\mathbf V}} 
\newcommand{\bfW}{{\mathbf W}} 
\newcommand{\bfX}{{\mathbf X}} 
\newcommand{\bfY}{{\mathbf Y}}
\newcommand{\bfZ}{{\mathbf Z}} 
\newcommand{\bfWh}{{\mathbf W}{\mathbf h}}


\let\sect=\S
\newcommand{\curs}{\EuR}
\newcommand{\ALGEBRAS}{\curs{ALGEBRAS}}
\newcommand{\CHAINCOMPLEXES}{\curs{CHCOM}}
\newcommand{\COMPLEXES}{\curs{COMPLEXES}}
\newcommand{\GROUPOIDS}{\curs{GROUPOIDS}}
\newcommand{\GROUPS}{\curs{GROUPS}}
\newcommand{\PAIRS}{\curs{PAIRS}}
\newcommand{\FGINJ}{\curs{FGINJ}}
\newcommand{\MODULES}{\curs{MODULES}}
\newcommand{\Or}{\curs{Or}}
\newcommand{\RINGS}{\curs{RINGS}}
\newcommand{\SPACES}{\curs{SPACES}}
\newcommand{\SPECTRA}{\curs{SPECTRA}}
\newcommand{\Sub}{\curs{Sub}}




\newcommand{\alg}{\operatorname{alg}}
\newcommand{\asmb}{\operatorname{asmb}}
\newcommand{\aut}{\operatorname{aut}}
\newcommand{\ch}{\operatorname{ch}}
\newcommand{\class}{\operatorname{class}}
\newcommand{\Cliff}{\operatorname{Cliff}}
\newcommand{\cok}{\operatorname{coker}}
\newcommand{\colim}{\operatorname{colim}}
\newcommand{\conhom}{\operatorname{conhom}}
\newcommand{\cone}{\operatorname{cone}}
\newcommand{\con}{\operatorname{con}}
\newcommand{\Con}{\operatorname{Con}}
\newcommand{\conn}{\operatorname{con}}
\newcommand{\cyl}{\operatorname{cyl}}
\newcommand{\cyclic}{\operatorname{cyclic}}
\newcommand{\diff}{\operatorname{diff}}
\newcommand{\Dtr}{\operatorname{Dtr}}
\newcommand{\finker}{\operatorname{finker}}
\newcommand{\FPMOD}{\operatorname{FPMOD}}
\newcommand{\Gen}{\operatorname{Gen}}
\newcommand{\hocolim}{\operatorname{hocolim}}
\newcommand{\HS}{\operatorname{HS}}
\newcommand{\id}{\operatorname{id}}
\newcommand{\ind}{\operatorname{ind}}
\newcommand{\inj}{\operatorname{inj}}
\newcommand{\inn}{\operatorname{inn}}
\newcommand{\inv}{\operatorname{inv}}
\newcommand{\Is}{\operatorname{is}}
\renewcommand{\Im}{\operatorname{im}}
\newcommand{\Ker}{\operatorname{ker}}
\newcommand{\map}{\operatorname{map}}
\newcommand{\mor}{\operatorname{mor}} 
\newcommand{\Ob}{\operatorname{ob}}
\newcommand{\pr}{\operatorname{pr}}
\newcommand{\Rep}{\operatorname{Rep}}
\newcommand{\res}{\operatorname{res}}
\newcommand{\sign}{\operatorname{sign}}
\newcommand{\Spin}{\operatorname{Spin}}
\newcommand{\Sw}{\operatorname{Sw}}
\newcommand{\topo}{\operatorname{top}}
\newcommand{\weak}{\operatorname{weak}}
\newcommand{\Wh}{\operatorname{Wh}}


\newcommand{\ab}{\operatorname{ab}}
\newcommand{\Aut}{\operatorname{Aut}}
\newcommand{\zentrum}{\operatorname{center}}
\newcommand{\coker}{\operatorname{coker}}
\newcommand{\Diff}{\operatorname{Diff}}
\newcommand{\dist}{\operatorname{dist}}
\newcommand{\dom}{\operatorname{dom}} 
\newcommand{\End}{\operatorname{End}} 
\newcommand{\GL}{\operatorname{GL}}    
\newcommand{\Idem}{\operatorname{Idem}}
\newcommand{\im}{\operatorname{im}}
\newcommand{\inc}{\operatorname{inc}}
\newcommand{\iso}{\operatorname{iso}}
\newcommand{\M}{\operatorname{M}}      
\newcommand{\MOD}{\operatorname{MOD}} 
\newcommand{\Nil}{\operatorname{Nil}}  
\newcommand{\obj}{\operatorname{obj}}
\newcommand{\pt}{\operatorname{pt}}
\newcommand{\St}{\operatorname{St}}    
\newcommand{\supp}{\operatorname{supp}} 
\newcommand{\Top}{\operatorname{Top}}
\newcommand{\Tor}{\operatorname{Tor}}
\newcommand{\tr}{\operatorname{tr}}
\newcommand{\Z}{\operatorname{Z}}      

\newcommand{\sma}{{\wedge}} 
\newcommand{\normal}{\triangleleft}
\newcommand{\semidirect}{\rtimes}
\newcommand{\nnt}{{\rm NZD}}
\newcommand{\dtr}{{\rm dtr}}
\newcommand{\bfdtr}{\ensuremath{\mathbf{dtr}}}
\newcommand{\bfpt}{\ensuremath{\mathbf{pt}}}


\newcommand{\entry}[2]{#1, ~ \pageref{#2}}   
\newcommand{\EGF}[2]{E_{#2}(#1)}               
\newcommand{\OrGF}[2]{\Or_{#2}(#1)}               
\newcommand{\SubGF}[2]{\Sub_{#2}(#1)}               
\newcommand{\comsquare}[8]                   
{\begin{CD}
#1 @>#2>> #3\\
@V{#4}VV @VV{#5}V\\
#6 @>>#7> #8
\end{CD}
}

\newcommand{\indextheorem}[1]{\index{#1}}         


\theoremstyle{plain}
\newtheorem{metaconjecture}[theorem]{Meta-Conjecture}
\newtheorem{addendum}[theorem]{Addendum}
\newtheorem{consequence}[theorem]{Consequence}

\theoremstyle{definition}
\newtheorem{notation}[theorem]{Notation}
\newtheorem{convention}[theorem]{Convention}
\newtheorem{situation}[theorem]{Situation}

\theoremstyle{remark}
\newtheorem{digression}[theorem]{Digression}

\renewcommand{\labelenumi}{(\roman{enumi})}
\renewcommand{\theenumi}{(\roman{enumi})}




\typeout{----------------------------  bcsfinal.tex  ----------------------------}

\title{The Baum-Connes and the Farrell-Jones Conjectures in K- and L-Theory}
\author{
Wolfgang L\"uck\thanks{\noindent email:
lueck@math.uni-muenster.de\protect\\
www: ~http://www.math.uni-muenster.de/u/lueck/\protect\\
fax: +49 251 8338370\protect}\; 
and Holger Reich\thanks{\noindent email: 
holger.reich@math.uni-muenster.de\protect\\
www: ~http://www.math.uni-muenster.de/u/reichh/\protect\\
}}

\institute{
Fachbereich Mathematik\\ Universit\"at M\"unster\\
Einsteinstr.~62\\ 48149 M\"unster\\Germany}

\maketitle

\typeout{-----------------------  Abstract  ------------------------}

\begin{abstract}
We give a survey of the meaning, status and applications of the
Baum-Connes Conjecture about the topological $K$-theory of the reduced
group $C^*$-algebra and the Farrell-Jones Conjecture about the algebraic
$K$- and $L$-theory of the group ring of a (discrete) group $G$. 

\smallskip
\noindent
Key words:  $K$- and $L$-groups of group rings and group $C^*$-algebras,
Baum-Connes Conjecture, Farrell-Jones Conjecture.

\smallskip
\noindent
Mathematics subject classification 2000: 19A31, 19B28, 19D99, 19G24,
19K99, 46L80. 
\end{abstract}

\typeout{-----------------------  Introduction  ------------------------}
\section*{Introduction}

This survey article is devoted to the Baum-Connes Conjecture about the
topological $K$-theory of the reduced group $C^*$-algebra and the
Farrell-Jones Conjecture about the algebraic $K$- and $L$-theory of
the group ring of a discrete group $G$. We will present a unified
approach to these conjectures hoping that it will stimulate further interactions and
exchange of methods and ideas between algebraic and geometric topology
on the one side and non-commutative geometry on the other.



Each of the above mentioned conjectures has already been proven for astonishingly large classes of groups
using a variety of different methods coming
from operator theory,  controlled topology and homotopy theory. Methods
have been developed for this purpose which turned out to be fruitful
in other contexts. The conjectures imply many other well-known and important 
conjectures. Examples are the
Borel Conjecture about the topological rigidity of closed aspherical
manifolds, the Novikov Conjecture about the homotopy invariance of
higher signatures, the stable Gromov-Lawson-Rosenberg Conjecture about the existence of 
Riemannian metrics with positive scalar curvature 
and the Kadison Conjecture about idempotents in the
reduced $C^*$-algebra of a torsionfree discrete group $G$. 


\subsection*{Formulation of the Conjectures}

The Baum-Connes and 
Farrell-Jones Conjectures predict that for every discrete group $G$ the following so called ``assembly maps''
are isomorphisms.
\begin{eqnarray*}
K_n^G(\EGF{G}{\calfin}) & \to & K_n(C^*_r(G));\\
H_n^G( \EGF{G}{\calvcyc} ; \bfK_R ) & \to &  K_n( RG ); \\
H_n^G( \EGF{G}{\calvcyc} ; \bfL_R^{\langle - \infty \rangle}
) & \to &  L_n^{\langle -\infty \rangle} ( RG ).
\end{eqnarray*}
Here the targets are the groups one would like to understand, namely the
topological $K$-groups of the reduced group $C^*$-algebra in the
Baum-Connes case and the algebraic  $K$- or $L$-groups of the group
ring $RG$ for $R$ an associative ring with unit. 
In each case the source is a
$G$-homology theory evaluated on a certain classifying space.
In the Baum-Connes Conjecture the $G$-homology theory is equivariant
topological $K$-theory and the classifying space $\EGF{G}{\calfin}$ is the classifying
space of the family of finite subgroups, which is often called
the classifying space for proper $G$-actions and denoted 
$\underline{E}G$ in the literature. In the Farrell-Jones
Conjecture the $G$-homology theory is given by a certain $K$- or
$L$-theory spectrum over the orbit category, and the classifying space
$\EGF{G}{\calvcyc}$ is the one associated to the family of virtually cyclic subgroups. 
The conjectures say that these assembly maps are isomorphisms.

These conjectures were stated in 
\cite[Conjecture 3.15 on page 254]{Baum-Connes-Higson(1994)} and
\cite[1.6 on page 257]{Farrell-Jones(1993a)}.  
Our formulations differ from the original ones, but are equivalent. 
In the case of the Farrell-Jones Conjecture
we slightly generalize the original conjecture by allowing arbitrary coefficient rings instead of $\IZ$. 
At the time of writing no counterexample to the Baum-Connes Conjecture~\ref{con: Baum-Connes Conjecture}
or the Farrell-Jones Conjecture~\ref{con: Farrell-Jones Conjecture} is known to the authors.



One can apply methods from algebraic topology such as spectral
sequences and Chern characters to the sources of the assembly maps. 
In this sense the sources
are  much more accessible than the targets. The conjectures hence lead to very concrete calculations.
Probably even more important is the structural insight: to what extent do the target groups 
show a homological behaviour.
These aspects can be treated completely analogously in the Baum-Connes and the Farrell-Jones setting.

However, the conjectures are not merely computational tools. Their importance comes from 
the fact that the assembly maps have geometric 
interpretations in terms of indices in the Baum-Connes case and in
terms of surgery theory in the Farrell-Jones case. These
interpretations are the key ingredient in applications 
and the reason that the Baum-Connes and Farrell-Jones Conjectures imply so many other conjectures
in non-commutative geometry, geometric topology and algebra.


\subsection*{A User's Guide}

A reader who wants to get specific information or focus on a certain
topic should consult the detailed table of contents, the index and the index of
notation in order to find the right place in the paper.
We have tried to write the text in a way such
that one can read small units independently from the rest.
Moreover, a reader who may only be interested in the
Baum-Connes Conjecture or only in the Farrell-Jones Conjecture for
$K$-theory or for $L$-theory can ignore the other parts. But we
emphasize again that one basic idea of this paper is to explain the
parallel treatment of these conjectures.

A reader without much prior knowledge about the Baum-Connes
Conjecture or the Farrell-Jones Conjecture should begin with 
Chapter~\ref{chap: torsion free}. There, the special case of a
torsionfree group is treated, since the formulation of the
conjectures is less technical in this case and  there are already many interesting
applications. The applications themselves however, are not needed later.
A more experienced reader may pass directly to
Chapter \ref{chap: general formulation}.

Other (survey) articles on the Farrell-Jones Conjecture and the 
Baum-Connes Conjecture are \cite{Farrell-Jones(1993a)}, \cite{Ferry-Ranicki-Rosenberg(1995a)},
\cite{Higson(1998a)}, \cite{Mislin-Valette(2003)}, \cite{Valette(2002)}.


\subsection*{Notations and Conventions}

Here is a briefing on our main notational conventions.
Details are of course discussed in the text. 
The columns in the following table contain our notation for:
the spectra, their  associated homology theory, the right hand side of the corresponding assembly maps,
the functor from groupoids to spectra and finally the $G$-homology theory associated to these spectra valued functors. 
\begin{quote}
\begin{tabular}{|| p{4em}| p{6.5em}|p{4em}|p{2em}| p{5.5em}||}
\hline
$\bfB \bfU$ & 
$K_n ( X ) $  & 
$K_n (C_r^{\ast} G)$  &   
$\bfK^{\topo}$  &   
$H_n^G(X ; \bfK^{\topo} )$ \\
\hline
$\bfK\! ( R ) $  &  
$ H_n ( X ; \bfK ( R ))$  & 
$K_n ( RG )$              &
$\bfK_R$ &
$H_n^G( X ; \bfK_R )$  \\
\hline
$\bfL^{\langle j\rangle} \! ( R )  $  & 
$H_n ( X ; \bfL^{\langle j \rangle} ( R ))$ &
$L_n^{\langle j \rangle}  (RG )$  & 
$\bfL^{\langle j \rangle}_R$ &
$H_n^G( X ; \bfL^{\langle j \rangle}_R )$  \\
\hline
\end{tabular}
\end{quote}
We would like to stress that $\bfK$ without any further decoration will always refer to the non-connective $K$-theory spectrum.
$\bfL^{\langle j \rangle}$ will always refer to quadratic $L$-theory with decoration $j$.
For a $C^{\ast}$- or Banach algebra $A$ the symbol  $K_n ( A)$ has two possible interpretations
but we will mean the topological $K$-theory.

A ring is always an associative ring with unit, and ring homomorphisms are always unital. Modules are left modules.
We will always work in the category of compactly generated spaces, compare~\cite{Steenrod(1967)} and \cite[I.4]{Whitehead(1978)}.
For our conventions concerning spectra see Section~\ref{sec: Spectra over the Orbit Category}.
Spectra are denoted with boldface letters such as $\bfE$.


\subsection*{Acknowledgements}

We would like to thank 
Arthur Bartels,
Indira Chatterji,
Tom Farrell, 
Joachim Grunewald,
Nigel Higson,  
Lowell Jones, 
Michel Matthey, 
Guido Mislin,
Erik Pedersen, 
David Rosenthal,
Marco Schmidt,
Alain Valette, 
Marco Varisco, 
Julia Weber,
Shmuel Weinberger, 
Bruce Williams,
Guoliang Yu and the referee
for useful comments and discussions. The first author 
thanks the Max-Planck Institute for Mathematics in Bonn for its
hospitality during his visit in November and December 2002, and both
authors thank the Department of Mathematics of the University
of Chicago,  in particular Shmuel Weinberger, for the pleasant stay in
March 2003, when parts of this paper were written.
We particularly enjoyed the ``Windy City''-Martini respectively the ``Earl Grey'' 
on the 96th floor of the Hancock building, where this introduction was drafted.\\

\hfill M{\"u}nster, August 2003

\bigskip


\setcounter{minitocdepth}{3} 
\setcounter{secnumdepth}{3}  
\makeatletter\tocsubsecnum=27\p@\makeatother 
\dominitoc


\typeout{---  Formulations and Applications of the Conjectures in the torsion free case -------}

\section{The Conjectures in the Torsion Free Case}
\label{chap: torsion free}

In this chapter we discuss the Baum-Connes and Farrell-Jones Conjectures
in the case of a torsion free group. Their formulation is less technical 
than in the general case, but  already in the torsion free case there are many interesting and illuminating conclusions. 
In fact some of the most important consequences of the conjectures, like for example the Borel Conjecture 
(see Conjecture~\ref{con: Borel Conjecture})
or the Kadison Conjecture (see Conjecture~\ref{con: Kadison Conjecture}), refer exclusively to the torsion free case.
On the other hand in the long run the general case, involving groups with torsion,
seems to be unavoidable. The general formulation yields a clearer and more complete picture, and furthermore
there are proofs of the conjectures for torsion free groups, where in
intermediate steps of the proof it is essential to have the general 
formulation available (compare Section~\ref{sec: Methods to improve control}).

The statement of the general case and further applications will
be presented  in the next chapter. 
The reader may actually skip this chapter and start immediately with 
Chapter \ref{chap: general formulation}.

We have put some effort into dealing with coefficient rings $R$ other than the integers.
A topologist may a priori be interested only in the case $R= \IZ$ but other cases are interesting
for algebraists and also do occur in computations for integral group rings.



\subsection{Algebraic {$K$}-Theory - Low Dimensions}
\label{sec: algebraic K low}

A ring $R$ is always understood to be associative with unit.
We denote by $K_n(R)$%
\indexnotation{K_n(R)} 
\emph{the algebraic $K$-group}\/ of $R$%
\index{K-groups@$K$-groups!algebraic $K$-groups of a ring}
for $n \in \IZ$. 
In particular $K_0(R)$ is the Grothendieck group of finitely generated projective $R$-modules
and elements in $K_1(R)$ can be represented by automorphisms of such modules.
In this section we are mostly interested in the $K$-groups $K_n ( R)$ with $n \leq 1$.
For definitions of these groups we refer to 
\cite{Milnor(1971)}, 
\cite{Rosenberg(2004)},
\cite{Silvester(1981)},
\cite{Swan(1970)},
\cite{Weibel(2003)}
 for $n=0$, $1$ and to \cite{Bass(1968)} and
\cite{Rosenberg(1994)} for $n \leq 1$.

For a ring $R$ and a group $G$ we denote by 
\[
A_0 = K_0(i) \colon K_0 (R) \to K_0 ( RG )
\]
the map induced by the natural inclusion $i\colon R \to RG$.
Sending $(g , [P]) \in G \times K_0( R )$ to the class of the $RG$-automorphism
\[
R[G] \otimes_R P \to R[G] \otimes_R P, \quad  u \otimes x \mapsto ug^{-1} \otimes x
\] 
defines a map $\Phi\colon  G_{\ab} \otimes_{\IZ}  K_0(R) \to K_1(RG)$, where $G_{\ab}$
denotes the abelianized group. We set
\[
A_1= \Phi \oplus K_1(i) \colon G_{\ab} \otimes_{\IZ} K_0(R) \oplus K_1(R) \to K_1(RG).
\]

We recall the notion of a regular ring.
We think of modules as left modules unless stated explicitly differently.
Recall that $R$ is \emph{Noetherian}%
\index{ring!Noetherian}
if any submodule of a finitely generated $R$-module is again finitely generated.
It is called \emph{regular}%
\index{ring!regular}
if it is Noetherian and any $R$-module has a finite-dimensional projective resolution.
Any principal ideal domain such as $\IZ$ or a field is regular.

The Farrell-Jones Conjecture about algebraic $K$-theory implies 
for a torsion free group the following conjecture 
about the low dimensional $K$-theory groups.

\begin{conjecture}[The Farrell-Jones Conjecture for Low Dimensional $K$-Theory and Torsion Free Groups]
\label{con: FJK lower}
\index{Conjecture!Farrell-Jones Conjecture!for Low Dimensional $K$-Theory and Torsion Free Groups}
Let $G$ be a torsion free group and let $R$ be a regular ring. Then 
\[
K_n(RG) =  0 \quad \text{ for } \quad n \le -1
\]
and the maps
\begin{eqnarray*}
K_0(R) & \xrightarrow{A_0} &   K_0(RG) \quad \mbox{ and } \\
G_{\ab} \otimes_{\IZ} K_0(R) \oplus K_1(R) & \xrightarrow{A_1} & K_1(RG)
\end{eqnarray*}
are both isomorphisms.
\end{conjecture} 
Every regular ring satisfies $K_n(R) = 0$ for $n \le -1$
\cite[5.3.30 on page 295]{Rosenberg(1994)} and hence the first statement 
is equivalent to $K_n (i) \colon K_n ( R ) \to K_n ( RG )$ being an isomorphism for $n \leq -1$.
In Remark~\ref{rem: Bass-Heller-Swan decomposition} 
below we explain why we impose the regularity assumption on the ring $R$.

For a regular ring $R$ and a group $G$ we define 
$\Wh_1^R(G)$
\index{Whitehead group!generalized}
\indexnotation{Wh_1^R(G)}
as the cokernel of the map $A_1$ and 
$\Wh_0^R (G)$%
\indexnotation{Wh_0^R(G)} as the cokernel of the map $A_0$.
In the important case where $R = \IZ$ the group $\Wh^{\IZ}_1 (G)$ coincides with the classical 
\emph{Whitehead group}%
\index{Whitehead group}
$\Wh(G)$%
\indexnotation{Wh(G)} 
which is the quotient of $K_1(\IZ G)$ by the subgroup consisting  of 
the classes of the units $\pm g \in (\IZ G)^{\inv}$ for $g \in G$. Moreover for every ring $R$ we define
the \emph{reduced algebraic $K$-groups $\widetilde{K}_n (R)$}%
\index{K-groups@$K$-groups!reduced algebraic $K$-groups of a ring}
\indexnotation{widetildeK_n(R)} 
as the cokernel of the natural map $K_n (\IZ) \to K_n( R )$. 
Obviously $\Wh_0^{\IZ} (G)= \widetilde{K}_0 ( \IZ G )$.

\begin{lemma}
The map $A_0$ is always injective. If $R$ is commutative and the natural map $\IZ \to K_0(R)$, $1 \mapsto [R]$ is an 
isomorphism, then the map $A_1$ is injective. 
\end{lemma}
\begin{proof}
The augmentation $\epsilon\colon  R G \to R$, which maps each group element $g$ to $1$, yields a retraction for the
inclusion $i \colon R \to RG$ and hence induces a retraction for $A_0$. If the map $\IZ \to K_0( R )$, $1 \mapsto [R]$
induces an isomorphism and 
$R$ is commutative, then we have the commutative diagram
\[
\xymatrix{
G_{\ab}  \otimes_{\IZ} K_0 ( R ) \oplus K_1 ( R ) \ar[dd]_-{ \cong } \ar[r]^-{A_1}  & K_1 ( R G ) \ar[d] \\
 & K_1 ( RG_{\ab} ) \ar[d]^-{( \det , K_1 ( \epsilon ) )} \\
G_{\ab}  \oplus K_1 ( R ) \ar[r] & RG_{\ab}^{\inv} \oplus K_1 ( R ),
         }
\]
where the upper vertical arrow on the right is induced from the map $G \to G_{\ab}$ to the abelianization.
Since $RG_{\ab}$ is a commutative ring we have the determinant 
$\det \colon K_1 ( R G_{\ab} ) \to (RG_{\ab})^{\inv}$. 
The lower horizontal arrow is induced from the obvious inclusion of $G_{\ab}$ into the invertible elements of the 
group ring $RG_{\ab}$ and in particular injective.
\end{proof}

In the special case $R= \IZ$  Conjecture~\ref{con: FJK lower}
above is equivalent to the following conjecture.
\begin{conjecture}[Vanishing of Low Dimensional $K$-Theory for Torsionfree Groups and Integral Coefficients]
\index{Conjecture!Farrell-Jones Conjecture!Vanishing of Low
Dimensional $K$-Theory for Torsionfree Groups and Integral Coefficients}
\label{con: vanishing of lower K}
For every torsion free group $G$ we have
\[
K_n(\IZ G)  = 0 \text{ for } n \le -1, \quad  \widetilde{K}_0 ( \IZ G ) = 0 \quad \mbox{and} \quad \Wh(G)=0.
\]
\end{conjecture}
\begin{remarknew}[Torsionfree is Necessary] \label{rem: torsion free is necessary} 
In general $\widetilde{K}_0(\IZ G)$ and $\Wh(G)$ do not vanish for finite groups. 
For example $\widetilde{K}_0 ( \IZ [\IZ / 23 ] ) \cong \IZ / 3 $ \cite[page 29,~30]{Milnor(1971)} and 
$\Wh ( \IZ / p ) \cong \IZ^{\frac{p-3}{2}}$ for $p$ an odd prime \cite[11.5 on page 45]{Cohen(1973)}.
This shows that the assumption that $G$ is torsion free is crucial in the formulation 
of Conjecture~\ref{con: FJK lower} above.

For more information on $\widetilde{K}_0( \IZ G )$ and Whitehead groups of finite groups see for instance
\cite[Chapter XI]{Bass(1968)},
\cite{Curtis-Reiner(1987)},
\cite{Milnor(1966)},
\cite{Oliver(1985a)},
\cite{Oliver(1989)} and
\cite{Swan(1970)}.

\end{remarknew}


\subsection{Applications I}
\label{sec: Applications I}
We will now explain 
the geometric relevance of the groups whose vanishing is predicted by 
Conjecture~\ref{con: vanishing of lower K}.


\subsubsection{The s-Cobordism Theorem and the Poincar\'e Conjecture}
\label{subsec: The s-Cobordism Theorem and the Poincare Conjecture}

The Whitehead group $\Wh(G)$ plays a key role if one studies manifolds because of the so called
s-Cobordism Theorem. In order to state it, we explain the notion of an h-cobordism.

Manifold always means smooth manifold unless it is explicitly stated differently.
We say that $W$ or more precisely $(W;M^-,f^-,M^+,f^+)$
is an $n$-dimensional \emph{cobordism}
\index{cobordism}%
over $M^-$ if $W$ is a compact $n$-di\-men\-sio\-nal
manifold together with the following: a disjoint decomposition of its boundary $\partial W$ into 
two closed $(n-1)$-dimensional manifolds $\partial^- W$ and $\partial^+ W$, two closed $(n-1)$-dimensional
manifolds $M^-$ and $M^+$ and diffeomorphisms $f^-\colon M^- \to \partial^- W$ and $f^+\colon M^+ \to \partial^+ W$.
The cobordism is called an \emph{h-cobordism}%
\index{h-cobordism}
if the inclusions $i^-\colon \partial^- W \to W$ and 
$i^+\colon \partial^+ W \to W$ are both homotopy equivalences.
Two cobordisms $(W;M^-,f^-,M^+,f^+)$ and $(W';M^-,f'^-,M'^+,f'^+)$
over $M^-$ are \emph{diffeomorphic relative $M^-$}
\index{cobordism!diffeomorphic relative $M^-$} if there is a
diffeomorphism $F \colon  W \to W'$
with $F \circ f^- = f'^-$. We call
a cobordism over $M^-$ \emph{trivial}%
\index{cobordism!trivial}, if it is diffeomorphic relative $M^-$
to the trivial h-cobordism given by the cylinder
$M^- \times[0,1]$ together with the obvious inclusions of $M^- \times \{0\}$ and $M^- \times \{1\}$.
Note that ``trivial'' implies in particular that $M^-$ and $M^+$ are diffeomorphic.

The question whether a given h-cobordism is trivial is decided by the Whitehead torsion $\tau(W; M^-) \in \Wh ( G )$
where $G=\pi_1(M^-)$.
For the details of the definition of $\tau(W;M^-)$ the reader should consult \cite{Cohen(1973)}, \cite{Milnor(1966)} or 
Chapter~2 in \cite{Lueck(2002c)}. Compare also \cite{Rosenberg(2004)}.

\begin{theorem}[s-Cobordism Theorem] \label{the: s-cobordism theorem}
\indextheorem{s-Cobordism Theorem}
Let $M^-$ be a closed connected oriented manifold of dimension $n \ge 5$
with fundamental group $G = \pi_1(M^-)$. Then
\begin{enumerate}

\item \label{the: s-cobordism theorem: triviality}
An h-cobordism  $W$ over $M^-$ is
trivial if and only if its Whitehead torsion  
$\tau(W,M^-) \in \Wh(G)$ vanishes.

\item \label{the: s-cobordism theorem: bijection}
Assigning to an h-cobordism over
$M^-$ its Whitehead torsion yields a bijection from
the diffeomorphism classes relative $M^-$ of h-cobordisms over $M^-$
to the Whitehead group $\Wh(G)$. 
\end{enumerate}
\end{theorem}

The s-Cobordism Theorem is due to Barden, Mazur and Stallings. There are also topological and PL-versions.
Proofs can be found for instance in 
\cite{Kervaire(1965)},
\cite[Essay~III]{Kirby-Siebenmann(1977)},
\cite{Lueck(2002c)} and 
\cite[page~87-90]{Rourke-Sanderson(1982)}.

The s-Cobordism Theorem tells us that the vanishing of the Whitehead group (as predicted in
Conjecture~\ref{con: vanishing of lower K}
for torsion free groups) has the following geometric interpretation.

\begin{consequence} \label{s-cobordisms and Whitehead group}
For a finitely presented group $G$ the vanishing of
the Whitehead group $\Wh(G)$ 
is equivalent to the statement that each 
h-cobordism over a closed connected manifold $M^-$ of dimension $\dim(M^-) \ge 5$
with fundamental group $\pi_1(M^-)\cong G$ is trivial. 
\end{consequence}

Knowing that all h-cobordisms over a given manifold are trivial is a strong and useful statement.
In order to illustrate this we would like to discuss the case where the fundamental group is trivial.

Since the ring $\IZ$ has a Gaussian algorithm, the determinant induces an
isomorphism $K_1(\IZ) \xrightarrow{\cong} \{\pm 1\}$ (compare~\cite[Theorem~2.3.2]{Rosenberg(1994)})
and the Whitehead group $\Wh(\{1\})$ of the trivial group vanishes. 
Hence any h-cobordism over a simply connected closed manifold of dimension $\ge 5$ is trivial.
As a consequence one obtains the Poincar\'e Conjecture for high dimensional manifolds.

\begin{theorem}[Poincar\'e Conjecture] \label{the: Poincare conjecture}
\indextheorem{Poincar\'e Conjecture}
Suppose $n \geq 5$. If the closed manifold $M$ is homotopy equivalent to the sphere $S^n$, then it is homeomorphic
to $S^n$.
\end{theorem}
\begin{proof} We only give the proof for $\dim(M) \ge 6$.
Let $f\colon  M \to S^n$ be a homotopy equivalence.
Let $D^n_- \subset M$ and $D^n_+ \subset M$ 
be two disjoint embedded disks. Let $W$ be the complement of the 
interior of the two disks in $M$. 
Then $W$ turns out to be a simply connected
h-cobordism over $\partial D^n_-$. Hence we can find a diffeomorphism
\[
F\colon  (\partial D^n_- \times [0,1];
\partial D^n_- \times \{0\}, \partial D^n_- \times \{1\})
\to (W;\partial D^n_-,\partial D^n_+)
\]
which is the identity on $\partial D^n_- = \partial D^n_- \times \{0\}$
and induces some (unknown) diffeomorphism
$f^+\colon  \partial D^n_- \times \{1\} \to \partial D^n_+$.
By the \emph{Alexander trick}%
\index{Alexander trick} one can extend
$f^+\colon  \partial D^n_- = \partial D^n_- \times \{1\} \to \partial D^n_+$
to a homeomorphism
$\overline{f^+}\colon  D^n_-  \to D^n_+$. Namely, any homeomorphism
$f\colon  S^{n-1} \to S^{n-1}$ extends to a homeomorphism
$\overline{f}\colon  D^n  \to D^n$ by sending
$t\cdot x$ for $t \in [0,1]$ and $x \in S^{n-1}$ to
$t \cdot f(x)$. Now define a homeomorphism
$h\colon  D^n_- \times \{0\} \cup_{i_-} \partial D^n_- \times [0,1]
\cup_{i_+} D^n_- \times \{1\} \to M$ for the canonical inclusions
$i_k\colon  \partial D^n_- \times \{k\} \to \partial D^n_- \times [0,1]$
for $k = 0,1$ by
$h|_{D^n_- \times \{0\}} = \id$, $h|_{\partial D^n_- \times [0,1]} = F$ and
$h|_{D^n_- \times \{1\}} = \overline{f^+}$. Since the source of $h$
is obviously homeomorphic to $S^n$, Theorem~\ref{the: Poincare conjecture} follows. 
\end{proof}

The Poincar\'e Conjecture
(see Theorem~\ref{the: Poincare conjecture})
is at the time of writing known in all dimensions except
dimension $3$.  It is essential in its formulation that one concludes $M$ to 
be homeomorphic (as opposed to diffeomorphic) to $S^n$. The Alexander trick does not work differentiably.
There are \emph{exotic spheres},%
\index{exotic sphere}
i.e.\ smooth manifolds which are homeomorphic but not
diffeomorphic to $S^n$ \cite{Milnor(1956)}.

More information about the Poincar\'e Conjecture, the Whitehead torsion and the
s-Cobordism Theorem can be found for instance in
\cite{Cappell-Shaneson(1985)},
\cite{Cohen(1973)},  
\cite{Donaldson(1987a)}, 
\cite{Freedman(1982)}, 
\cite{Freedman(1983)},
\cite{Hambleton(2002)},
\cite{Kervaire(1965)}, 
\cite{Lueck(2002c)}, 
\cite{Milnor(1965b)},
\cite{Milnor(1966)},  
\cite{Rosenberg(2004)} and 
\cite{Rourke-Sanderson(1982)}.


\subsubsection{Finiteness Obstructions}
\label{subsec: Finiteness Obstructions}

We now discuss the geometric relevance of $\widetilde{K}_0(\IZ G)$. 

Let $X$ be a $CW$-complex. It is called \emph{finite}
\index{$CW$-complex!finite} if it consists of finitely many cells. It is called
\emph{finitely dominated}%
\index{$CW$-complex!finitely dominated}
if there is a finite $CW$-complex $Y$ together with maps $i\colon X \to Y$ and $r\colon Y \to X$
such that $r \circ i$ is homotopic to the identity on $X$. The fundamental group of a finitely dominated
$CW$-complex is always finitely presented.

While studying existence problems for spaces with prescribed properties (like for example group actions),
it happens occasionally that it is relatively easy to construct a finitely dominated $CW$-complex within a 
given homotopy type, whereas it is not at all clear whether one can also 
find a homotopy equivalent \emph{finite} $CW$-complex.
\emph{Wall's finiteness obstruction}%
\index{finiteness obstruction},
a certain obstruction element 
$\widetilde{o}(X)%
\indexnotation{widetilde o(X)}%
\in \widetilde{K}_0(\IZ \pi_1(X))$, decides the question.

\begin{theorem}[Properties of the Finiteness Obstruction] \label{the: finiteness obstruction}
\indextheorem{Properties of the Finiteness Obstruction}
Let $X$ be a finitely dominated $CW$-complex with fundamental group $\pi=\pi_1(X)$. 
\begin{enumerate}
\item
The space $X$ is homotopy equivalent to a finite $CW$-complex if and only if 
$\widetilde{o}(X)=0 \in \widetilde{K}_0(\IZ \pi )$.
\item
Every element in $\widetilde{K}_0(\IZ G)$ can be realized
as the finiteness obstruction $\widetilde{o}(X)$ of a finitely dominated $CW$-complex $X$ with $G = \pi_1(X)$, 
provided that $G$ is finitely presented. 
\item
Let $Z$ be a space such that $G=\pi_1( Z )$ is finitely presented. Then there is a bijection between
$\widetilde{K}_0 ( \IZ G )$ and  the set of equivalence classes of maps $f \colon X \to Z$ with $X$ 
finitely dominated under the equivalence relation 
explained below.
\end{enumerate}
\end{theorem}
The equivalence relation in (iii) is defined as follows:
Two maps $f \colon X \to Z$ and $f' \colon X' \to Z$
with $X$ and $X'$ finitely dominated 
are equivalent if there exists a commutative diagram
\[
\xymatrix{
X \ar[drr]_-{f} \ar[r]^-{j} & X_1 \ar[dr]^-{f_1} \ar[r]^-{h} & X_2 \ar[d]_-{f_2} & 
X_3 \ar[dl]_-{f_3} \ar[l]_-{h'} &  X' \ar[dll]^-{f'} \ar[l]_-{j'} \\
& & Z & &
         },
\]
where $h$ and $h'$ are homotopy equivalences and $j$ and $j'$ are 
inclusions of subcomplexes for which $X_1$, respectively $X_3$,
is obtained from $X$, respectively $X'$, by attaching a finite number of cells.

The vanishing of $\widetilde{K}_0 ( \IZ G )$ as predicted in 
Conjecture~\ref{con: vanishing of lower K}
for torsion free groups hence has the following interpretation.
\begin{consequence} \label{lem: finiteness obstructions and FJC}
For a finitely presented group $G$ the vanishing of
$\widetilde{K}_0(\IZ G)$ is equivalent to the statement that any finitely dominated
$CW$-complex $X$ with $G \cong \pi_1(X)$ is homotopy equivalent to a finite $CW$-complex.
\end{consequence}

For more information about the finiteness obstruction we refer for instance to
\cite{Ferry(1981a)},
\cite{Ferry-Ranicki(2001)}, 
\cite{Lueck(1987b)}, 
\cite{Mislin(1995)}, 
\cite{Ranicki(1985)},
\cite{Rosenberg(2004)},
\cite{Varadarajan(1989a)},
\cite{Wall(1965a)} and 
\cite{Wall(1966)}.


\subsubsection{Negative $K$-Groups and Bounded h-Cobordisms} 
\label{subsec: bounded h-cobordisms}

One possible geometric interpretation of negative $K$-groups is in terms of bounded h-cobordisms.
Another interpretation will be explained in Subsection~\ref{subsec: negative K pseudo} below.

We consider \emph{manifolds $W$ parametrized over} $\IR^k$,%
\index{manifold parametrized over $\IR^k$} 
i.e.\ manifolds which are equipped with a 
surjective proper map $p \colon W \to \IR^k$. 
We will always assume that the fundamental group(oid) is bounded, compare~\cite[Definition~1.3]{Pedersen(1986)}.
A map $f\colon W \to W'$ between two manifolds parametrized over $\IR^k$ is
\emph{bounded}%
\index{bounded map}
if $\{p' \circ f(x) - p(x) \; | \; x \in W \}$ is a bounded subset of $\IR^k$.

A \emph{bounded cobordism}%
\index{cobordism!bounded}
$(W;M^-,f^-,M^+,f^+)$ is defined just as in 
Subsection~\ref{subsec: The s-Cobordism Theorem and the Poincare Conjecture}
but compact manifolds are replaced by manifolds parametrized over $\IR^k$ and the
parametrization for $M^{\pm}$ is given by $p_W \circ f^{\pm}$.
If we assume that the inclusions $i^{\pm}\colon  \partial^{\pm} W \to W$ are 
homotopy equivalences, then there exist deformations $r^{\pm}\colon W \times I \to W$, $(x,t) \mapsto r^{\pm}_t(x)$
such that $r_0^{\pm}= \id_W$ and $r_1^{\pm}(W) \subset \partial^{\pm} W$. 

A bounded cobordism is called a \emph{bounded h-cobordism}%
\index{h-cobordism!bounded} 
if the inclusions $i^{\pm}$ are homotopy equivalences and  additionally the deformations can be chosen
such that the two sets 
\[
S^{\pm} = \{p_W \circ r^{\pm}_t ( x ) - p_W \circ r^{\pm}_1 (x) \; | \; x \in W , t \in [0,1] \}
\]
are bounded subsets of $\IR^k$. 

The following theorem (compare \cite{Pedersen(1986)} and \cite[Appendix]{Weiss-Williams(1988)})
contains the s-Cobordism Theorem~\ref{the: s-cobordism theorem}
as a special case, gives another  interpretation of elements in
$\widetilde{K}_0 ( \IZ \pi )$ and explains one aspect of the geometric relevance of negative
$K$-groups.

\begin{theorem}[Bounded h-Cobordism Theorem]
 \label{the: bounded h-cobordism}
\indextheorem{Bounded h-Cobordism Theorem}
Suppose that $M^-$ is parametrized over $\IR^k$ and satisfies $\dim M^- \geq 5$. Let
$\pi$ be its fundamental group(oid). 
Equivalence classes of bounded h-cobordisms over $M^-$ modulo bounded diffeomorphism relative $M^-$ correspond
bijectively to elements in $\kappa_{1-k} (\pi)$, where
\[
\kappa_{1-k} ( \pi ) = \left\{ \begin{array}{lll} \Wh (\pi ) & \quad & \mbox{ if } k=0 ,\\
                                                  \widetilde{K}_0 ( \IZ \pi ) & \quad & \mbox{ if } k=1 , \\
                                                  K_{1-k} ( \IZ \pi ) & \quad & \mbox{ if } k \geq 2.
                               \end{array}
                       \right.
\]
\end{theorem}

More information about negative $K$-groups can be found for instance in
\cite{Anderson-Hsiang(1977)}, 
\cite{Bass(1968)}, 
\cite{Carter(1980b)}, 
\cite{Carter(1980)},
\cite{Farrell-Jones(1995)},
\cite{Madsen-Rothenberg(1989a)}, 
\cite{Pedersen(1984)},
\cite{Pedersen(1986)}, 
\cite{Quinn(1982a)},
\cite{Ranicki(1992a)},
\cite{Rosenberg(1994)} 
and \cite[Appendix]{Weiss-Williams(1988)}.


\subsection{Algebraic {$K$}-Theory - All Dimensions}
\label{sec: algebraic K all}

So far we only considered  the $K$-theory groups in dimensions $\leq 1$. 
We now want to explain how Conjecture~\ref{con: FJK lower}
generalizes to higher algebraic $K$-theory. For the definition of higher algebraic $K$-theory groups and the (connective) 
$K$-theory spectrum see 
\cite{Berrick(1982)},
\cite{Carlsson(2004)},
\cite{Inassaridze(1995)}, 
\cite{Quillen(1973)}, 
\cite{Rosenberg(1994)}, 
\cite{Srinivas(1991)},
\cite{Waldhausen(1985)} and
\cite{Weibel(2003)}. 
We would like to stress that for us 
$\bfK(R)$%
\indexnotation{bfK(R)}
will always denote 
the \emph{non-connective algebraic $K$-theory spectrum}%
\index{K-theory spectrum@$K$-theory spectrum!algebraic $K$-theory spectrum of a ring} 
for which $K_n (R) = \pi_n( \bfK (R) )$ holds for all $n \in \IZ$. 
For its definition see \cite{Carlsson(2004)}, \cite{Loday(1976)} and \cite{Pedersen-Weibel(1985)}.

The Farrell-Jones Conjecture for algebraic $K$-theory
reduces for a torsion free group to the following conjecture.

\begin{conjecture}[Farrell-Jones Conjecture for Torsion Free Groups and K-Theory]
\label{con: FJK torsion free all}
\index{Conjecture!Farrell-Jones Conjecture!for Torsion Free Groups and $K$-Theory}
Let $G$ be a torsion free group. Let $R$ be a regular ring.
Then the assembly map
\[
H_n(BG;\bfK (R)) \to K_n(RG)
\]
is an isomorphism for $n \in \IZ$.
\end{conjecture} 

Here $H_n(-;\bfK (R))$ denotes the homology theory which is associated to the spectrum $\bfK (R)$.
It has the property that
$H_n(\pt;\bfK(R))$ is $K_n(R)$ for  $n \in \IZ$, where here and elsewhere 
$\pt$%
\indexnotation{pt}
denotes the space consisting of one point. The space $BG$%
\indexnotation{BG} 
is the \emph{classifying space of the group $G$},%
\index{classifying space!of a group}
which up to homotopy is characterized by the property that it is a 
$CW$-complex with $\pi_1(BG) \cong G$ whose universal covering
is contractible. The technical details of the construction of
$H_n(-;\bfK (R))$ and the assembly map 
will be explained in a more general setting
in Section~\ref{sec: Formulation of the Conjectures}.

The point of Conjecture~\ref{con: FJK torsion free all}
is that on the right-hand side of the assembly map we have the group $K_n(RG)$ we are interested in, whereas
the left-hand side is a homology theory and hence much easier to compute. For every homology theory associated to 
a spectrum we have the Atiyah-Hirzebruch spectral sequence, which in our case has 
$E^2_{p,q} = H_p(BG;K_q(R))$ and converges to $H_{p+q}(BG;\bfK (R))$.

If $R$ is regular, then the negative $K$-groups of $R$ vanish and the spectral sequence lives in the first 
quadrant. Evaluating the spectral sequence for $n=p+q \leq 1$ shows that 
Conjecture~\ref{con: FJK torsion free all} above implies 
Conjecture~\ref{con: FJK lower}.


\begin{remarknew}[Rational Computation] \label{rem: rational computation K}
Rationally an Atiyah-Hirzebruch spectral sequence 
collapses always 
and the homological Chern character gives an isomorphism
\[
\ch\colon \bigoplus_{p+q = n} H_p(BG;\IQ) \otimes_{\IQ} \left( K_q(R)  \otimes_{\IZ} \IQ \right)
\xrightarrow{\cong}
H_n(BG;\bfK(R)) \otimes_{\IZ} \IQ.
\]
\end{remarknew}
The Atiyah-Hirzebruch spectral sequence and the Chern character will be discussed in a much more general setting in 
Chapter \ref{chap: Computations}. 


\begin{remarknew}[Separation of Variables]
 \label{rem: separation of variables for K-theory in the torsion free case}
We see that the left-hand side of the isomorphism 
in the previous remark
consists of a group homology part and a part which is the rationalized $K$-theory of $R$.
(Something similar happens before we rationalize at the level of spectra: The left
hand side of Conjecture~\ref{con: FJK torsion free all} 
can be interpreted as the homotopy groups of the spectrum $BG_+ \sma \bfK (R)$.)
So essentially Conjecture~\ref{con: FJK torsion free all} 
predicts that the $K$-theory of $RG$ is built up out of two independent parts:
the $K$-theory of $R$ and the group homology of $G$.
We call this principle \emph{separation of variables}.
\index{principle!separation of variables}
This principle also applies to other theories such as algebraic $L$-theory
or topological $K$-theory. See also Remark \ref{rem: Separation of Variables for calfin}.
\end{remarknew}


\begin{remarknew}[$K$-Theory of the Coefficients] \label{rem: coefficients K}
Note that Conjecture~\ref{con: FJK torsion free all} 
can only help us to explicitly compute the $K$-groups of $RG$ in cases where we know enough
about the $K$-groups of $R$. We obtain no new information about the $K$-theory of $R$ itself.
However, already for very simple rings the computation of their algebraic $K$-theory groups is an
extremely hard problem.

It is known that the groups $K_n( \IZ )$ are finitely generated abelian groups \cite{Quillen(1973a)}. 
Due to Borel \cite{Borel(1974)} we know that
\[
K_n(\IZ) \otimes_{\IZ} \IQ \cong \left\{ \begin{array}{ll}
                                      \IQ & \mbox{ if } n = 0; \\
                                      \IQ & \mbox{ if } n = 4k+1 \mbox{ with  } k \ge 1; \\
                                        0 & \mbox{ otherwise. }
                                    \end{array}
                              \right.
\]
Since $\IZ$ is regular we know that $K_n(\IZ)$ vanishes for $n \le -1$. 
Moreover, $K_0(\IZ) \cong \IZ$ and $K_1(\IZ) \cong \{\pm 1\}$, where the isomorphisms
are given by the rank and the determinant. One also knows that $K_2(\IZ) \cong \IZ/2$,
$K_3(\IZ) \cong \IZ/48$  \cite{Lee-Szczarba(1976)} and $K_4(\IZ) \cong 0$ \cite{Rognes(2000)}.
Finite fields belong to the few rings where one has a complete and 
explicit knowledge of all $K$-theory groups \cite{Quillen(1972)}.
We refer the reader for example to
\cite{Kolster(2004)},
\cite{Mitchell(1994)}, 
\cite{Rognes-Weibel(2000)}, 
\cite{Weibel(2004)} and Soul{\'e}'s article in 
\cite{Lluis-Puebla(1992)} for more information about the algebraic $K$-theory of the integers or more generally of 
rings of integers in number fields.

Because of Borel's calculation 
the left hand side of the isomorphism described in Remark~\ref{rem: rational computation K}
specializes for $R= \IZ$ to 
\begin{eqnarray} 
\label{Chern character for H_n(BG;bfK(IZ))}
H_n(BG;\IQ) \oplus \bigoplus_{k=1}^{\infty}  H_{n - (4k +1 )} (BG;\IQ) 
\end{eqnarray}
and Conjecture~\ref{con: FJK torsion free all} predicts that this group
is isomorphic to $K_n ( \IZ G ) \otimes_{\IZ} \IQ$.
\end{remarknew}


Next we discuss the case where the group $G$ is infinite cyclic.
\begin{remarknew}[Bass-Heller-Swan Decomposition]
 \label{rem: Bass-Heller-Swan decomposition} 
The so called \emph{Bass-Heller-Swan-decomposition},%
\index{Bass-Heller-Swan-decomposition} also known as the
\emph{Fundamental Theorem of algebraic $K$-theory},%
\indextheorem{Fundamental Theorem of algebraic $K$-theory} computes the algebraic $K$-groups of
$R[\IZ]$ in terms of the algebraic $K$-groups and Nil-groups of $R$:
\[
K_n(R[\IZ]) ~ \cong ~ K_{n-1}(R) \oplus K_n(R) \oplus N\!K_n(R) \oplus N\!K_n(R).
\]
Here the group $N\!K_n(R)$%
\indexnotation{NK_n(R)} 
is defined as the cokernel of the split injection $K_n(R) \to K_n(R[t])$.
It can be identified with the cokernel of the split injection $K_{n-1}(R) \to K_{n-1} ( \caln \! il (R) )$. 
Here $K_n(\caln \! il (R))$%
\indexnotation{K_n(calnil(R))}
denotes the $K$-theory of the exact category of nilpotent
endomorphisms  of finitely generated projective $R$-modules. 
For negative $n$ it is defined with the help of Bass' notion of a contracting functor 
\cite{Bass(1968)} (see also \cite{Carter(1980b)}). The groups are known as 
\emph{Nil-groups}%
\index{Nil-groups}
and often denoted $\Nil_{n-1} (R)$.%
\indexnotation{Nil_n(R)}

For proofs of these facts and more information the reader should consult
\cite[Chapter~XII]{Bass(1968)},
\cite{Bass-Heller-Swan(1964)}, 
\cite[Theorem on page~236]{Grayson(1976)}, 
\cite[Corollary in \S 6 on page~38]{Quillen(1973)},
\cite[Theorems~3.3.3 and 5.3.30]{Rosenberg(1994)},
\cite[Theorem~9.8]{Srinivas(1991)}
and \cite[Theorem~10.1]{Swan(1995)}.

If we iterate and use $R[\IZ^n]=R[\IZ^{n-1}][\IZ]$ we see that a computation of $K_n(RG)$ must in general take into account 
information about  $K_i(R)$ for all $i \le n$. In particular we see that it is important 
to formulate Conjecture~\ref{con: FJK torsion free all}
with the non-connective $K$-theory spectrum.

Since $S^1$ is a model for $B\IZ$, we get an isomorphism
$$H_n(B\IZ;\bfK (R)) ~ \cong ~  K_{n-1}(R) \oplus K_n(R)$$
and hence Conjecture~\ref{con: FJK torsion free all} predicts
$$K_n(R[\IZ]) ~ \cong ~ K_{n-1}(R) \oplus K_n(R).$$
This explains why in the formulation of Conjecture~\ref{con: FJK torsion free all} the condition
that $R$ is regular appears. It guarantees that  $N\!K_n(R) =
0$~\cite[Theorem 5.3.30 on page 295]{Rosenberg(1994)}. 
There are weaker conditions which imply that $N\!K_n ( R )=0$ but ``regular'' has the advantage that
$R$ regular implies that $R[t]$ and $R[\IZ]=R[t^{\pm 1}]$ are again regular, 
compare the discussion in Section~2 in~\cite{Bass(1973)}.

The Nil-terms $N\!K_n(R)$ seem to be hard to compute. For instance
$N\!K_1(R)$ either vanishes or is infinitely generated as an abelian group
\cite{Farrell(1977)}. 
In Subsection~\ref{subsec: The Isomorphism Conjecture for NK-groups}
we will discuss the Isomorphism Conjecture for $N\!K$-groups.
For more information about Nil-groups see for instance
 \cite{Connolly-da-Silva(1995)}, \cite{Connolly-Kozniewski(1995)}, \cite{Hesselholt-Madsen(2001)},
\cite{Weibel(1980)} and \cite{Weibel(1981)}.
\end{remarknew}


\subsection{Applications II}

\subsubsection{The Relation to Pseudo-Isotopy Theory} 
\label{subsec: pseudoisotopies}

Let $I$ denote the unit interval $[0,1]$.
A topological \emph{pseudoisotopy}%
\index{pseudoisotopy}
of a compact  manifold $M$ is a 
homeomorphism $h \colon M \times I \to M \times I$, which restricted
to $M \times \{ 0 \} \cup \partial M \times I$ is the obvious inclusion. 
The space $P ( M )$%
\indexnotation{P(M)}
 of pseudoisotopies is the (simplicial) group
of all such homeomorphisms.
Pseudoisotopies play an important role if one tries to understand the homotopy type  of the space $\Top ( M )$
of self-homeomorphisms of a manifold $M$. We will see below in Subsection~\ref{subsec: automorphisms of manifolds}
how the results about pseudoisotopies discussed in this section combined with surgery theory lead to quite explicit
results about the homotopy groups of $\Top(M)$.

There is a stabilization map $P( M ) \to P ( M \times I )$ given by
crossing a pseudoisotopy  with the identity on the interval $I$
and the stable pseudoisotopy space is defined as 
$\calp ( M )%
\indexnotation{calp(M)}
 = \colim_k P ( M \times I^k )$. In fact $\calp(-)$ can be extended to a functor
on all spaces \cite{Hatcher(1978)}. The natural inclusion $P(M) \to \calp ( M )$ induces an isomorphism
on the $i$-th homotopy group if the dimension of $M$ is large compared to $i$, 
see~\cite{Burghelea-Lashof(1977)} and \cite{Igusa(1988)}.

Waldhausen \cite{Waldhausen(1978)}, \cite{Waldhausen(1985)} defines the algebraic $K$-theory of spaces functor 
$\bfA (X)$%
\indexnotation{bfA(X)}
and the functor 
$\bfWh^{PL}( X )$%
\indexnotation{bfWh^PL(X)}
from spaces to spectra (or infinite loop spaces) and a fibration sequence
\[
X_+ \sma \bfA ( \pt ) \to \bfA( X )  \to \bfWh^{PL} ( X ).
\]
Here $X_+ \sma \bfA ( \pt ) \to \bfA( X) $ is an assembly map, which can be compared to the
algebraic $K$-theory assembly map that appears in 
Conjecture~\ref{con: FJK torsion free all} via a commutative diagram
\[
\xymatrix{
H_{n} (X ; \bfA ( \pt ) ) \ar[d] \ar[r] \ar[d] & \pi_n ( \bfA ( X )) \ar[d] \\
H_n ( B \pi_1 ( X ) ; \bfK ( \IZ )) \ar[r] & K_n ( \IZ \pi_1 ( X ) ).
         }
\]
In the case where $X \simeq BG$ is aspherical the vertical maps induce 
isomorphisms after rationalization for $n \geq 1$, 
compare~\cite[Proposition~2.2]{Waldhausen(1978)}.
Since $\Omega^2 \Wh^{PL}(X) \simeq \calp (X)$ (a guided tour through the 
literature concerning this and related results can be found
in \cite[Section~9]{Dwyer-Weiss-Williams(2001)}),
Conjecture~\ref{con: FJK torsion free all} implies  rational vanishing results
for the groups $\pi_n ( \calp (M ))$ if $M$ is an aspherical manifold. 
Compare also Remark~\ref{rem: relating K and P}.
\begin{consequence} \label{conseq: rational vanishing of pseudoisotopies} 
Suppose $M$ is a closed aspherical manifold and Conjecture~\ref{con: FJK torsion free all}
holds for $R=\IZ$ and $G = \pi_1 (M)$, then for all $n \geq 0$ 
\[
\pi_n ( \calp ( M ) ) \otimes_{\IZ} \IQ = 0.
\]
\end{consequence}

Similarly as above one defines smooth pseudoisotopies and the space of stable smooth
pseudoisotopies 
$\calp^{\Diff} ( M )$.%
\indexnotation{calp^Diff(M)}
There is also a smooth version of the Whitehead space 
$\Wh^{\Diff}( X )$%
\indexnotation{Wh^Diff(X)}
and $\Omega^2 \Wh^{\Diff}( M ) \simeq \calp^{\Diff} ( M )$.
Again there is a close relation to $A$-theory via the natural splitting 
$\bfA(X) \simeq \Sigma^{\infty} (X_+) \vee \bfWh^{\Diff}( X )$, see
\cite{Waldhausen(1987a)}. 
Here $\Sigma^{\infty} (X_+)$ denotes the suspension spectrum associated to $X_+$.
Using this one
can split off an assembly map $H_n( X ; \bfWh^{\Diff} (\pt ) ) \to \pi_n (\bfWh^{\Diff} (X))$ 
from the $A$-theory assembly map. Since
for every space $\pi_n ( \Sigma^{\infty} ( X_+ ) ) \otimes_{\IZ} \IQ \cong H_n ( X ; \IQ )$ 
Conjecture~\ref{con: FJK torsion free all} combined with the rational computation in 
\eqref{Chern character for H_n(BG;bfK(IZ))} yields the following result.
\begin{consequence} \label{conseq: computation of diff pseudos}
Suppose $M$ is a closed aspherical manifold and Conjecture~\ref{con: FJK torsion free all} 
holds for $R = \IZ$ and $G= \pi_1(M)$. Then for $n \geq 0$ we have 
\[
\pi_n ( \calp^{\Diff} ( M ) )\otimes_{\IZ} \IQ = \bigoplus_{k = 1}^{\infty} H_{n-4k+1} (M ; \IQ ).
\]
\end{consequence}

Observe that the  fundamental difference between the smooth and the topological case occurs already when
$G$ is the trivial group. 


\subsubsection{Negative $K$-Groups and Bounded Pseudo-Isotopies} 
\label{subsec: negative K pseudo}

We briefly explain a further geometric interpretation of negative $K$-groups, which  
parallels the discussion of bounded h-cobordisms in Subsection~\ref{subsec: bounded h-cobordisms}.

Let $p\colon M \times \IR^k \to \IR^k$ denote the natural projection.
The space $P_b( M ; \IR^k )$%
\indexnotation{P_b(M;R^k)}
of bounded pseudoisotopies is the space of all self-homeomorphisms
$h\colon M \times \IR^k \times I \to M \times \IR^k \times I$ 
such that restricted to $M \times \IR^k \times \{ 0 \}$ the map 
$h$ is the inclusion and such that $h$ is bounded, i.e.\ 
the set $\{ p \circ h ( y ) - p( y ) \; | \; y \in M \times \IR^k \times I \}$
is a bounded subset of $\IR^k$.
There is again a stabilization map $P_b( M ; \IR^k ) \to P_b ( M \times I ; \IR^k )$ 
and a stable bounded pseudoisotopy
space $\calp_b ( M ; \IR^k )= \colim_j P_b ( M \times I^j ; \IR^k )$.
There is a homotopy equivalence $\calp_b ( M ; \IR^k ) \to \Omega \calp_b ( M ; \IR^{k+1} )$ 
\cite[Appendix~II]{Hatcher(1978)} and hence the sequence of 
spaces $\calp_b ( M ; \IR^k )$ for $k= 0 ,1 , \dots $ is an $\Omega$-spectrum 
$\bfP (M)$.%
\indexnotation{bfP(M)} 
Analogously one defines the differentiable bounded pseudoisotopies $\calp_b^{\diff}(M ; \IR^k)$ and an $\Omega$-spectrum 
$\bfP^{\diff}(M)$.
The negative homotopy groups of these spectra have an 
interpretation in terms of low and negative dimensional $K$-groups.
In terms of unstable homotopy groups this is explained in the following theorem which 
is closely related to Theorem~\ref{the: bounded h-cobordism} about bounded h-cobordisms. 

\begin{theorem}[Negative Homotopy Groups of Pseudoisotopies]
 \label{the: negative homotopy of pseudoisotopy}
\indextheorem{Negative Homotopy Groups of Pseudoisotopies}
Let $G=\pi_1(M)$.
Suppose $n$ and $k$ are such that $n+k \geq 0$, then
for $k \geq 1$ there are isomorphisms 
\[
\pi_{n+k} ( \calp_b ( M ; \IR^k ) ) =  \left\{ \begin{array}{lll} \Wh (G ) & \quad & \mbox{ if } n=-1 ,\\
                                                  \widetilde{K}_0 ( \IZ G ) & \quad & \mbox{ if } n=-2 , \\
                                                  K_{n +2} ( \IZ G ) & \quad & \mbox{ if } n < -2
                                              \end{array}
                                     \right.
\]
The same result holds in the differentiable case.
\end{theorem}

Note that Conjecture~\ref{con: FJK torsion free all} predicts that these groups vanish if $G$ is torsionfree.
The result above is due to Anderson and Hsiang \cite{Anderson-Hsiang(1977)} and 
is also discussed in \cite[Appendix]{Weiss-Williams(1988)}.


\subsection{$L$-Theory}
\label{sec: L-Theory in the Torsion Free Case}

We now move on to the $L$-theoretic version of the Farrell-Jones Conjecture. We will still stick to the 
case where the group is torsion free. 
The conjecture is obtained by replacing  $K$-theory and 
the $K$-theory spectrum in Conjecture~\ref{con: FJK torsion free all}
by $4$-periodic 
$L$-theory and the $L$-theory spectrum $\bfL^{\langle - \infty \rangle} ( R )$. Explanations will follow below.

\begin{conjecture}[Farrell-Jones Conjecture for Torsion Free Groups and L-Theory]
\label{con: FJL torsion free}
\index{Conjecture!Farrell-Jones Conjecture!for Torsion Free Groups and $L$-Theory}
Let $G$ be a torsion free group and let $R$ be a ring with involution.
Then the assembly map 
\[
H_n(BG;\bfL^{\langle - \infty \rangle}(R)) \to L_n^{\langle - \infty \rangle}(RG)
\]
is an isomorphism for $n \in \IZ$.
\end{conjecture}

To a ring with involution one can associate (decorated) symmetric or
quadratic algebraic $L$-groups,%
\index{L-groups@$L$-groups}
compare 
\cite{Cappell-Ranicki-Rosenberg(2000)},
\cite{Cappell-Ranicki-Rosenberg(2001)},
\cite{Ranicki(1981)},
\cite{Ranicki(1992a)} and \cite{Williams(2004)}.
We will exclusively deal with the quadratic algebraic $L$-groups and denote them by $L_n^{\langle j \rangle}(R)$%
\indexnotation{L_n^j(R)}.
Here $n \in \IZ$ and $j \in \{-\infty\} \amalg \{j \in \IZ\mid j \le 2\}$ is the so called \emph{decoration}.%
\index{decoration}
The decorations $j = 0,1$ correspond to the decorations $p,h$ and $j=2$ is related to the decoration $s$ appearing in the literature.
Decorations will be discussed in Remark~\ref{rem: decorations} below.
The $L$-groups $L_n^{\langle j \rangle}(R)$ are $4$-periodic, i.e.\
$L_n^{\langle j \rangle}(R) \cong L_{n+4}^{\langle j \rangle}(R)$ for $n \in \IZ$.

If we are given an involution $r \mapsto \overline{r}$ on a ring $R$, we will always equip
$RG$ with the involution that extends the given one and satisfies $\overline{g}=g^{-1}$.
On $\IZ$, $\IQ$ and $\IR$ we always use the trivial involution and on 
$\IC$ the complex conjugation.

One can construct an \emph{$L$-theory spectrum}%
\index{L-theory spectrum@$L$-theory spectrum!algebraic $L$-theory spectrum of a ring}
$\bfL^{\langle j  \rangle}(R)$%
\indexnotation{bfL^j(R)} such that $\pi_n ( \bfL^{\langle j \rangle}(R))= L_n^{\langle j \rangle}(R)$,
compare \cite[\S~13]{Ranicki(1992)}. 
Above and in the sequel $H_n(-;\bfL^{\langle j\rangle}(R))$ denotes the 
homology theory which is associated to this spectrum. In particular we have 
$H_n(\pt;\bfL^{\langle j \rangle}(R))= L_n^{\langle j \rangle}(R)$. We
postpone the discussion  of the assembly map to
 Section~\ref{sec: Formulation of the Conjectures}
where we will construct it in greater generality.


\begin{remarknew}[The Coefficients in the $L$-Theory Case]
\label{rem: coefficients L}
In contrast to $K$-theory (compare Remark~\ref{rem: coefficients K}) the  $L$-theory of the most interesting
coefficient ring $R= \IZ$ is well known. 
The groups $L_n^{ \langle j \rangle } ( \IZ )$ for fixed $n$ and varying 
$j \in \{- \infty\} \amalg \{j \in \IZ\mid j \le 2\}$ are all naturally isomorphic 
(compare Proposition~\ref{pro: interplay K L} below) and 
we have $L_0^{\langle j \rangle}(\IZ) \cong \IZ$ and 
$L_2^{\langle j \rangle}(\IZ) \cong \IZ/2$, 
where the isomorphisms are given by the signature divided by $8$ 
and the Arf invariant,  and $L_1^{\langle j \rangle}(\IZ ) = 
L_3^{\langle j \rangle}(\IZ) =0$, see 
\cite[Chapter III]{Browder(1972)}, 
\cite[Proposition 4.3.1 on page 419]{Ranicki(1981)}.
\end{remarknew}


\begin{remarknew}[Decorations]  \label{rem: decorations}
$L$-groups are designed as obstruction groups for surgery problems. The decoration reflects what kind of surgery problem
one is interested in. All $L$-groups can be described as cobordism 
classes of suitable quadratic Poincar{\'e} chain complexes.
If one works with chain complexes of finitely generated free based $R$-modules
and requires that the torsion of the Poincar{\'e} chain homotopy equivalence vanishes in
$\widetilde{K}_1 ( R )$, then  the corresponding $L$-groups are denoted $L^{\langle 2 \rangle}_n (R)$.
If one drops the torsion condition, one obtains $L^{\langle 1 \rangle}_n (R)$, which is usually denoted $L^h (R)$.
If one works with finitely generated projective modules, one obtains $L^{\langle 0 \rangle} ( R)$, which is also known as $L^p (R)$.

The L-groups with negative decorations can be defined inductively via the Shaneson splitting,
compare Remark~\ref{rem: Shaneson splitting} below. Assuming that the $L$-groups with decorations $j$ have already been defined
one sets
\[
L^{<j-1>}_{n-1}(R) = 
\coker ( L_{n}^{<j>}(R) \to L_{n}^{<j>}(R[\IZ])).
\]
Compare \cite[Definition 17.1 on page 145]{Ranicki(1992a)}.
Alternatively these groups can be obtained via a process which is in the 
spirit of Subsection~\ref{subsec: bounded h-cobordisms} and Subsection~\ref{subsec: negative K pseudo}. 
One can define them as L-theory groups 
of suitable categories of modules parametrized over $\IR^k$. For details
the reader could consult \cite[Section 4]{Carlsson-Pedersen(1995a)}.
There are forgetful maps $L^{\langle j+1 \rangle}_n (R) \to L^{\langle j \rangle}_n (R)$.
The group $L_n^{\langle - \infty \rangle} ( R )$ is defined as the colimit over these maps.
For more information see \cite{Ranicki(1973b)}, \cite{Ranicki(1992a)}.

For group rings we also define $L^s_n ( RG )$ similar to $L^{\langle 2 \rangle}_n ( RG )$ but we require
the torsion to lie in $\im A_1 \subset \widetilde{K}_1 ( RG )$, where $A_1$ is the map defined in
Section~\ref{sec: algebraic K low}.
Observe that $L^s_n ( RG )$ really depends on the pair $(R , G)$ and differs in general from $L^{\langle 2 \rangle}_n ( RG)$.
\end{remarknew}


\begin{remarknew}[The Interplay of $K$- and $L$-Theory] \label{rem: Rothenberg sequence}
For $j \leq 1$ 
there are forgetful maps
$L_n^{\langle j +1 \rangle}(R) \to L_n^{\langle j \rangle}(R)$ which sit inside the following sequence, which
is known as the 
\emph{Rothenberg sequence}%
\index{Rothenberg sequence} 
\cite[Proposition~1.10.1 on page~104]{Ranicki(1981)},
\cite[17.2]{Ranicki(1992a)}. 
\begin{multline}
\ldots \to L_n^{\langle j+1 \rangle}(R  ) \to L_n^{\langle j \rangle}(R  ) \to
\widehat{H}^n(\IZ/ 2; \widetilde{K}_{j}(R)) 
\\ \to   L_{n-1}^{\langle j+1 \rangle}(R ) \to L_{n-1}^{\langle j \rangle}(R) \to \ldots.
\label{Rothenberg sequence}
\end{multline}
Here $\widehat{H}^n(\IZ/ 2;\widetilde{K}_{j}(R))$ is the Tate-cohomology of the group
$\IZ/ 2$ with coefficients in the $\IZ[\IZ/ 2]$-module $\widetilde{K}_{j}(R)$. The involution on 
$\widetilde{K}_{j}(R)$ comes from the involution on $R$. 
There is a similar sequence relating $L^s_n ( R G)$ and $L^h_n ( R G )$, where the
third term is the $\IZ/ 2$-Tate-cohomology of $\Wh^R_1 (G)$. 
Note that Tate-cohomology groups of the group
$\IZ/ 2 $ are always annihilated by multiplication with $2$.
In particular $L_n^{\langle j \rangle}(R)[\frac{1}{2}] = L_n^{\langle j \rangle} (R) \otimes_{\IZ} \IZ [ \frac{1}{2} ]$ 
is always independent of $j$. 

Let us formulate explicitly what we obtain from the above sequences for 
$R = \IZ G$.
\end{remarknew}

\begin{propositionnew}  \label{pro: interplay K L}
Let $G$ be a torsion free group, then
Conjecture~\ref{con: vanishing of lower K} about the vanishing of $\Wh(G)$, $\widetilde{K_0} ( \IZ G )$ and 
$K_{-i} ( \IZ G )$ for $i \geq 1$ implies that 
for fixed $n$ and varying $j \in \{ - \infty \} \amalg \{ j \in \IZ \; | \; j \leq 1 \}$ 
the $L$-groups
$L_n^{ \langle j \rangle } ( \IZ G )$ are all naturally isomorphic and 
moreover $L^{\langle 1 \rangle}_n ( \IZ G ) = L^h_n ( \IZ G ) \cong L^s_n ( \IZ G )$.
\end{propositionnew}


\begin{remarknew}[Rational Computation] \label{rem: rational computation L}
As in the $K$-theory case we have an 
Atiyah-Hirzebruch spectral sequence:
\[
E^2_{p,q} = H_p(BG;L_q^{\langle -\infty \rangle}(R)) \quad \Rightarrow  \quad H_{p+q}(BG;\bfL^{\langle -\infty \rangle}(R)).
\]
Rationally this spectral sequence collapses 
and the homological Chern character gives for $n \in \IZ$ an isomorphism
\begin{multline}
\ch\colon \bigoplus_{p+q = n} H_p(BG;\IQ) \otimes_{\IQ} 
\left( L_q^{\langle -\infty \rangle} (R) \otimes_{\IZ} \IQ \right)
\\
 \xrightarrow{\cong} ~ 
H_n(BG;\bfL^{\langle -\infty \rangle}(R)) \otimes_{\IZ} \IQ.
\label{Chern character for  H_n(BG;bfL^{langle -infty rangle})}
\end{multline}
In particular we obtain 
in the case $R = \IZ$ 
from Remark~\ref{rem: coefficients L} 
for all $n \in \IZ$ and all decorations $j$ an isomorphism
\begin{eqnarray} \label{Chern hcaracter for H_n(BG;bfL^{langle j rangle}(IZ))}
\xymatrix{
\ch \colon \bigoplus_{k=0}^{\infty}  H_{n-4k}(BG;\IQ) 
\ar[r]^-{\cong} &
H_n(BG;\bfL^{\langle j \rangle}(\IZ)) \otimes_{\IZ} \IQ.
         }
\end{eqnarray}
This spectral sequence and the 
Chern character above will be discussed in a much more general setting in 
Chapter \ref{chap: Computations}. 
\end{remarknew}


\begin{remarknew}[Torsion Free is Necessary] \label{rem: L for finite groups}
If $G$ is finite, $R = \IZ$ and $n=0$, then 
the rationalized left hand side of the assembly equals $\IQ$, whereas the right hand side is isomorphic 
to the rationalization of the real representation ring. 
Since the group homology of a finite group vanishes rationally except in dimension $0$,  the previous
remark shows that we need to assume the group to be torsion free in Conjecture~\ref{con: FJL torsion free}
\end{remarknew}


\begin{remarknew} [Shaneson splitting] \label{rem: Shaneson splitting} 
The Bass-Heller-Swan decomposition in $K$-theory 
(see Remark \ref{rem: Bass-Heller-Swan decomposition}) 
has the following analogue for the 
algebraic $L$-groups, which is known as the 
\emph{Shaneson splitting}%
\index{Shaneson-splitting}
\cite{Shaneson(1969)}
\begin{equation} \label{equ: shaneson splitting}
L^{\langle j \rangle}_n(R[\IZ]) ~ \cong ~ L^{\langle j-1 \rangle}_{n-1}(R) 
\oplus L^{\langle j \rangle}_n(R).
\end{equation}
Here for the decoration $j = - \infty$ one has to interpret $j-1$ as $-\infty$.
Since $S^1$ is a model for $B\IZ$, we get an isomorphisms
\[
H_n (B\IZ;\bfL^{\langle j \rangle}(R)) ~ \cong ~ L^{\langle j \rangle}_{n-1} (R) \oplus
L^{\langle j \rangle}_n(R).
\]
This explains why in the formulation of 
the $L$-theoretic Farrell-Jones Conjecture for torsion free groups 
(see Conjecture~\ref{con: FJL torsion free})
we use the decoration
$j = - \infty$. 

As long as one deals with torsion free groups and one believes in the low dimensional part of the 
$K$-theoretic 
Farrell-Jones Conjecture 
(predicting the vanishing of $\Wh(G)$, $\widetilde{K}_0( \IZ G)$ and of the negative $K$-groups, see 
Conjecture~\ref{con: vanishing of lower K})
there is no difference between the various decorations $j$, compare Proposition~\ref{pro: interplay K L}.
But as soon as one allows
torsion in $G$, the decorations make a difference and it indeed turns out that if one replaces the decoration $j = - \infty$
by $j = s,h$ or $p$ there are counterexamples for the $L$-theoretic version of Conjecture~\ref{con: Farrell-Jones Conjecture}
even for $R = \IZ$ \cite{Farrell-Jones-Lueck(2002)}.

Even though in the above Shaneson splitting (\ref{equ: shaneson splitting}) there are no terms analogous to 
the Nil-terms in Remark~\ref{rem: Bass-Heller-Swan decomposition} such 
Nil-phenomena do also occur in $L$-theory,  as soon as one considers 
amalgamated free products. The corresponding groups are the UNil-groups. 
They vanish if one inverts $2$ \cite{Cappell(1974c)}.
For more information about the UNil-groups we refer to 
\cite{Banagl-Ranicki(2003)}
\cite{Cappell(1973a)},
\cite{Cappell(1973)}, 
\cite{Connolly-Kozniewski(1995)}, 
\cite{Connolly-Ranicki(2003)},
\cite{Farrell(1979)},
\cite{Ranicki(1995b)}.
\end{remarknew}


\subsection{Applications III}
\label{sec: Applications III}


\subsubsection{The Borel Conjecture} 
\label{subsec: Borel}

One of the driving forces for the development of the Farrell-Jones Conjectures 
is still the following topological rigidity conjecture
about closed aspherical manifolds, compare~\cite{Farrell-Jones(1989)}. Recall that a manifold, or more generally a $CW$-complex, 
is called \emph{aspherical}%
\index{aspherical}
if its universal covering is contractible.
An aspherical $CW$-complex $X$ with $\pi_1 (X)=G$ is a model for the classifying space $BG$. If $X$ is an aspherical manifold
and hence finite dimensional, then $G$ is necessarily torsionfree.

\begin{conjecture}[Borel Conjecture] \label{con: Borel Conjecture}
\index{Conjecture!Borel Conjecture}
Let $f \colon M \to N$ be a homotopy equivalence of aspherical closed topological
manifolds. Then $f$ is homotopic to a homeomorphism. In particular two closed
aspherical manifolds with isomorphic fundamental groups are homeomorphic.
\end{conjecture}

Closely related to the Borel Conjecture is the conjecture that each 
aspherical finitely dominated Poincar\'e complex is homotopy
equivalent to a closed topological manifold. The Borel Conjecture \ref{con: Borel Conjecture} 
is false in the smooth category, i.e.\ if one replaces topological manifold by smooth
manifold and homeomorphism by diffeomorphism \cite {Farrell-Jones(1989b)}.

Using surgery theory one can show that in dimensions $\geq 5$ 
the Borel Conjecture is implied by the $K$-theoretic vanishing 
Conjecture~\ref{con: vanishing of lower K}
combined with the $L$-theoretic Farrell-Jones Conjecture.

\begin{theorem}[The Farrell-Jones Conjecture Implies the Borel Conjecture]
\label{the: The Farrell-Jones Conjecture Implies the Borel Conjecture}
\indextheorem{The Farrell-Jones Conjecture Implies the Borel Conjecture}
Let $G$ be a torsion free group.
If $\Wh(G)$, $\widetilde{K}_0( \IZ G )$ and all the groups $K_{-i} ( \IZ G )$ with $i \geq 1$ vanish
and if the assembly map 
\[
H_n ( BG ; \bfL^{\langle - \infty \rangle}( \IZ ) ) \to L_n^{\langle - \infty \rangle} ( \IZ G )
\]
is an isomorphism for all $n$, then the Borel Conjecture holds for all orientable manifolds of dimension $\geq 5$ whose 
fundamental group is $G$.
\end{theorem}

The Borel Conjecture \ref{con: Borel Conjecture} can be reformulated in the language of
surgery theory to the statement that the topological structure set $\cals^{\topo}(M)$
of an aspherical closed topological manifold $M$ consists of a single point. This set is the set of equivalence
classes of homotopy equivalences $f\colon M' \to M$ with a topological closed manifold as source
and $M$ as target under the equivalence relation, for which  $f_0 \colon M_0 \to M$ and
$f_1 \colon M_1 \to M$ are equivalent if there is a homeomorphism $g\colon M_0 \to M_1$
such that $f_1 \circ g$ and $f_0$ are homotopic.

The surgery sequence of a closed
orientable topological manifold $M$ of dimension $n \ge 5$ is the 
exact sequence
\begin{multline*}
\ldots \to \caln_{n+1}(M\times [0,1],M \times \{0,1\}) \xrightarrow{\sigma}  L^s_{n+1}(\IZ\pi_1(M)) 
\xrightarrow{\partial} \cals^{\topo}(M) 
\\
\xrightarrow{\eta} \caln_n(M) \xrightarrow{\sigma} L_n^s(\IZ\pi_1(M)),
\end{multline*}
which extends infinitely to the left.
It is the basic tool for the classification of topological manifolds. 
(There is also a smooth version of it.)
The map $\sigma$ appearing in the sequence
sends a normal map of degree one to its surgery obstruction. 
This map can be identified with the version of the $L$-theory assembly map where one works with the 
$1$-connected cover $\bfL^s ( \IZ ) \langle 1 \rangle$ of $\bfL^s( \IZ )$. 
The map  $H_k(M;\bfL^s (\IZ)\langle 1 \rangle ) \to H_k(M;\bfL^s (\IZ))$
is injective for $k=n$ and an isomorphism for $k >n$. Because of the $K$-theoretic assumptions we can replace the 
$s$-decoration with the $\langle - \infty \rangle$-decoration, compare Proposition~\ref{pro: interplay K L}.
Therefore the Farrell-Jones 
Conjecture~\ref{con: FJL torsion free} implies
that the maps $\sigma\colon \caln_n(M) \to L_n^s(\IZ\pi_1(M))$ and
$\caln_{n+1}(M\times [0,1],M \times \{0,1\}) \xrightarrow{\sigma}  L^s_{n+1}(\IZ\pi_1(M))$ are
injective  respectively bijective and thus by the surgery
sequence that $\cals^{\topo}(M)$ is a point and hence the Borel Conjecture 
\ref{con: Borel Conjecture} holds for $M$. More details can be found e.g.\ in
\cite[pages 17,18,28]{Ferry-Ranicki-Rosenberg(1995)}, \cite[Chapter 18]{Ranicki(1992)}.

For more information about surgery theory we refer for instance to
\cite{Browder(1972)}, 
\cite{Cappell-Ranicki-Rosenberg(2000)}, 
\cite{Cappell-Ranicki-Rosenberg(2001)}, 
\cite{Farrell-Goettsche-Lueck(2002a)}, 
\cite{Farrell-Goettsche-Lueck(2002b)},
\cite{Karoubi(2004)},
\cite{Kreck(1999)}, 
\cite{Ranicki(2002a)},
\cite{Stark(2000)},
\cite{Stark(2002)},
and 
\cite{Wall(1999)}.


\subsubsection{Automorphisms of Manifolds} 
\label{subsec: automorphisms of manifolds}

If one additionally also assumes the Farrell-Jones Conjectures for higher $K$-theory,
one can combine the surgery theoretic
results with the results about pseudoisotopies from 
Subsection~\ref{subsec: pseudoisotopies} to obtain the following results.

\begin{theorem}[Homotopy Groups of $\Top(M)$]
\label{the: auto top}
\indextheorem{Homotopy Groups of $\Top(M)$}
Let $M$ be an orientable closed aspherical manifold of dimension $> 10$ with fundamental group $G$.
Suppose the $L$-theory assembly map 
\[
H_n ( B G ; \bfL^{\langle - \infty \rangle} (\IZ) ) \to L_n^{ \langle - \infty \rangle } ( \IZ G )
\]
is an isomorphism for all $n$ and suppose the $K$-theory assembly map
\[
H_n ( BG ; \bfK  ( \IZ ) ) \to K_n ( \IZ G )
\]
is an isomorphism for $n \leq 1$ and a rational isomorphism for $n \geq 2$.
Then for $1 \leq i \leq ( \dim M -7 ) / 3$
one has
\[
\pi_i ( \Top( M ) ) \otimes_{\IZ} \IQ 
              = \left\{ \begin{array}{lll} \zentrum ( G ) \otimes_{\IZ} \IQ  & \quad & \mbox{ if } i=1 ,\\
                                                  0  & \quad & \mbox{ if } i > 1
                         \end{array}    
                \right.
\]
\end{theorem}

In the differentiable case one additionally needs to study involutions on the higher $K$-theory groups.
The corresponding result reads:
\begin{theorem}[Homotopy Groups of $\Diff(M)$]
\label{the: auto diff}
\indextheorem{Homotopy Groups of $\Diff(M)$}
Let $M$ be an orientable closed aspherical differentiable manifold of dimension $>10$ with fundamental group $G$.
Then under the same assumptions as in Theorem~\ref{the: auto top} we have for $1 \leq i \leq (\dim M -7)/3 $
\[
\pi_i ( \Diff ( M ) ) \otimes_{\IZ} \IQ =
                \left\{ \begin{array}{lll} \zentrum ( G ) \otimes_{\IZ} \IQ  &  & \mbox{ if } i=1 ;\\
     \bigoplus_{j=1}^{\infty} H_{(i +1) - 4j} ( M ; \IQ ) &  & \mbox{ if } i > 1 \mbox{ and } \dim M \mbox{ odd } ;\\
                                                      0 & & \mbox{ if } i > 1 \mbox{ and } \dim M \mbox{ even }.
                         \end{array}    
                \right.
\]
\end{theorem}

See for instance 
\cite{Farrell-Hsiang(1978)}, 
\cite[Section~2]{Farrell-Jones(1990b)}
and \cite[Lecture~5]{Farrell(2002)}. 
For a modern survey on automorphisms of manifolds
we refer to \cite{Weiss-Williams(2001)}.


\subsection{The Baum-Connes Conjecture in the Torsion Free Case}
\label{sec: The Baum Connes Conjecture in the Torsion Free Case}

We denote by $K_{\ast} (Y)$
the complex $K$-homology of a topological space $Y$ and by
$K_{\ast} ( C_r^{\ast} (G) )$ the (topological) $K$-theory of the reduced group $C^{\ast}$-algebra.
More explanations will follow below.
\begin{conjecture}[Baum-Connes Conjecture for Torsion Free Groups]
\label{con: BC torsion free}
\index{Conjecture!Baum-Connes Conjecture!for Torsion Free Groups}
Let $G$ be a torsion free group. Then the Baum-Connes assembly map
\[
K_n(BG) \to K_n(C^*_r(G))
\]
is bijective for all $n \in \IZ$.
\end{conjecture}

Complex $K$-homology 
\index{K-groups@$K$-groups!topological $K$-homology groups of a space}
$K_*(Y)$%
\indexnotation{K_*(Y)}
is the homology theory associated to 
the topological (complex) $K$-theory spectrum $\bfK^{\topo}$ (which is 
is often denoted $\bfB\bfU$) and could also be written as $K_{\ast}(Y)=H_{\ast}( Y ; \bfK^{\topo})$. 
The \emph{co}homology theory  associated to the spectrum $\bfK^{\topo}$ is the well known complex $K$-theory 
defined in terms of complex vector bundles. 
Complex $K$-homology is a $2$-periodic theory, i.e.\ $K_{n}(Y) \cong K_{n+2}(Y)$. 

Also the \emph{topological $K$-groups}%
\index{K-groups@$K$-groups!topological $K$-groups of a $C^*$-algebra}
$K_n(A)$%
\indexnotation{K_n(A)}
of a $C^*$-algebra $A$ are $2$-periodic. 
Whereas $K_0( A )$ coincides with the algebraically defined  $K_0$-group,
the other groups $K_n( A)$ take the topology of the $C^{\ast}$-algebra $A$ into account, for instance
$K_n ( A ) = \pi_{n-1} ( GL ( A ) )$ for $n \geq 1$.

Let  $\calb( l^2(G))$ denote the bounded linear operators on the Hilbert space $l^2(G)$ whose orthonormal basis is $G$.
The 
\emph{reduced complex group $C^*$-algebra}%
\index{C^*-algebra@$C^*$-algebra!reduced complex group $C^*$-algebra}
$C^*_r(G)$%
\indexnotation{C^*_r(G)}
is the closure in the norm topology
of the image of the regular representation $\IC G \to \calb(l^2(G))$, which sends an
element $u \in \IC G$ to the (left) $G$-equivariant bounded operator $l^2(G) \to l^2(G)$
given by right multiplication with $u$. In particular one has natural inclusions
\[
\IC G \subseteq C^*_r(G) \subseteq \calb(l^2(G))^G \subseteq \calb(l^2(G)).
\]
It is essential to use the reduced
group $C^*$-algebra in the Baum-Connes Conjecture, there are counterexamples for the version with 
the maximal group $C^*$-algebra, compare
Subsection~\ref{subsec: The Baum-Connes Conjecture for Maximal Group $C^*$-Algebras}.
For information about 
$C^*$-algebras and their topological $K$-theory
we refer for instance to 
\cite{Blackadar(1986)}, 
\cite{Connes(1994)}, 
\cite{Davidson(1996)},
\cite{Higson-Roe(2000)}, 
\cite{Lance(1995)}, 
\cite{Murphy(1990)}, 
\cite{Schick(2002)} and 
\cite{Wegge-Olsen(1993)}.

\begin{remarknew}[The Coefficients in the Case of Topological $K$-Theory] \label{rem: coefficients Ktop}
If we specialize to the trivial group $G=\{ 1\}$, then the complex reduced group $C^{\ast}$-algebra reduces to 
$C_r^{\ast}(G)=\IC$ and the topological $K$-theory is well known: by periodicity it suffices to know that
$K_0 ( \IC ) \cong \IZ$, where the homomorphism is given by the dimension, and $K_1 ( \IC )=0$.
Correspondingly we have $K_q(\pt) = \IZ$ for $q$ even and $K_q(\pt) = 0$ for odd $q$. 
\end{remarknew}

\begin{remarknew}[Rational Computation] \label{rem: rational computation Ktop}
There is an Atiyah-Hirzebruch spectral sequence
which converges to $K_{p+q}(BG)$ and whose $E^2$-term is
$E^2_{p,q} = H_p(BG;K_q(\pt))$. 
Rationally this spectral sequence collapses 
and the homological Chern character gives 
an isomorphism for $n \in \IZ$ 
\begin{multline}
\ch \colon \bigoplus_{k \in \IZ}  H_{n-2k}(BG;\IQ) 
~ = ~ \bigoplus_{p + q = n} H_p(BG;\IQ) \otimes_{\IQ} \left( K_q (\IC)  \otimes_{\IZ} \IQ \right)
\\
\xrightarrow{\cong}   K_n(BG) \otimes_{\IZ} \IQ.
\label{Chern character for K_n(BG)}
\end{multline}
\end{remarknew}

\begin{remarknew}[Torsionfree is Necessary] \label{rem: BC for finite groups}
In the case where $G$ is a finite group the reduced group $C^{\ast}$-algebra $C_r^{\ast}(G)$
coincides with the complex group ring $\IC G$ and $K_0 ( C_r^{\ast}(G))$ coincides with the complex representation
ring of $G$. 
Since the group homology of a finite group vanishes rationally except in dimension $0$, the previous
remark shows that we need to assume the group to be torsion free in Conjecture~\ref{con: BC torsion free}.
\end{remarknew}

\begin{remarknew}{\bf (Bass-Heller-Swan-Decomposition for Topological $K$-Theory)}
\label{Computing K(C^*_r(G times Z)} 
There is an analogue of the Bass-Heller-Swan decomposition in algebraic $K$-theory
(see Remark \ref{rem: Bass-Heller-Swan decomposition}) or of the Shaneson
splitting in $L$-theory (see Remark \ref{rem: Shaneson splitting}) for
topological $K$-theory. Namely  we have 
\[
K_n(C^*_r(G \times \IZ)) ~ \cong ~ K_n(C^*_r(G)) \oplus K_{n-1}(C^*_r(G)),
\]
see \cite[Theorem 3.1 on page 151]{Pimsner-Voiculescu(1982)} or more generally \cite[Theorem 18 on page 632]{Pimsner(1986)}.
This is consistent with the obvious isomorphism
$$K_n(B(G \times \IZ)) = K_n(BG \times S^1) \cong  K_{n-1}(BG) \oplus K_n(BG).$$
Notice that here in contrast to the algebraic $K$-theory no Nil-terms occur 
(see Remark \ref{rem: Bass-Heller-Swan decomposition}) and that there is no analogue of 
the complications in algebraic $L$-theory coming from the varying decorations
(see Remark \ref{rem: Shaneson splitting}). 
This absence of Nil-terms or decorations is the reason why  
in the final formulation of the Baum-Connes Conjecture it suffices to deal with the
family of finite subgroups, whereas in the algebraic $K$- and $L$-theory case one must
consider the larger and harder to handle family of virtually cyclic subgroups.
This in some sense makes the computation of topological $K$-theory of
reduced group $C^*$-algebras easier than the computation of $K_n(\IZ G)$ or  $L_n(\IZ G)$.
\end{remarknew}

\begin{remarknew}[Real Version of the Baum-Connes Conjecture]
\label{rem: real version of Baum-Connes} 
There is an obvious \emph{real version of the Baum-Connes Conjecture}.%
\index{Conjecture!Baum-Connes Conjecture!real version}
It says  that for a torsion free group the real assembly map 
\[
KO_n(BG) \to KO_n(C^*_r(G;\IR))
\]
is bijective for $n \in \IZ$.
We will discuss in Subsection~\ref{subsec: The Real Version of BC}
below that this real version of the Baum-Connes Conjecture 
is implied by the complex version Conjecture~\ref{con: BC torsion free}.

Here $KO_n(C^*_r(G;\IR))$%
\indexnotation{KO_n(C^*_r(G;R)} is the topological $K$-theory of the 
\emph{real reduced group $C^*$-algebra}
\index{C^*-algebra@$C^*$-algebra!reduced real group $C^*$-algebra} $C^*_r(G;\IR)$.%
\indexnotation{C^*_r(G;IR)} 
We use $KO$ instead of $K$ as a reminder that we work with real $C^{\ast}$-algebras.
The topological real $K$-theory $KO_*(Y)$%
\indexnotation{KO_n(Y)} is the homology theory associated to the spectrum 
$\bfB\bfO$, whose associated \emph{co}homology theory is given in terms of real vector bundles.
Both, topological $K$-theory of a real $C^{\ast}$-algebra  and $KO$-homology of a space are 
$8$-periodic and
$KO_n(\pt) = K_n(\IR)$ is $\IZ$, if $n = 0,4 \;(8)$, is $\IZ/ 2 $ if 
$n = 1,2\; (8)$ and is $0$ if $n = 3,5,6,7\; (8)$.

More information about the $K$-theory of real $C^{\ast}$-algebras can be found in \cite{Schroeder(1993)}.
\end{remarknew}


\subsection{Applications  IV}
\label{sec: Applications  IV}


We now discuss some consequences of the Baum-Connes Conjecture for Torsion Free Groups
\ref{con: BC torsion free}.

\subsubsection{The Trace Conjecture in the Torsion Free Case}
\label{subsec: The Trace Conjecture in the Torsion Free Case}

The assembly map appearing in the Baum-Connes Conjecture has an interpretation in terms of
index theory. This is important for geometric applications. It is of the same significance as
the interpretation of the $L$-theoretic assembly map as the map $\sigma$
appearing in the exact surgery sequence discussed in Section~\ref{sec: L-Theory in the Torsion Free Case}. 
We proceed to explain this.

An element $\eta \in K_0(BG)$ can be
represented by a pair $(M,P^*)$  consisting of a cocompact free proper smooth $G$-manifold $M$ with
Riemannian metric together with an elliptic $G$-complex $P^*$ of differential operators of
order $1$ on $M$ \cite{Baum-Douglas(1982)}. 
To such a pair one can assign an index $\ind_{C^*_r(G)}(M,P^*)$ in
$K_0(C^*_r(G))$ \cite{Mishchenko-Fomenko(1979)} which is the image of $\eta$ under the assembly map 
$K_0(BG) \to K_0(C^*_r(G))$ appearing in  
Conjecture~\ref{con: BC torsion free}.
With this interpretation the surjectivity of the assembly map for a torsion free group says that 
any element in $K_0(C^*_r(G))$ can be realized as an index. This allows to apply
index theorems to get interesting information. 

Here is a prototype of such an argument.
The \emph{standard trace}
\index{trace!standard trace of $C^*_r(G)$} 
\begin{eqnarray}
\tr_{C^*_r(G)}%
\indexnotation{tr_{C^*_r(H)}}
 \colon C^*_r(G) &  \to & \IC \label{standard trace tr_{C^*_r(G)}}
\end{eqnarray}
sends an element $f \in C^*_r(G) \subseteq \calb(l^2(G))$ to $\langle f(1),1\rangle_{l^2(G)}$. 
Applying the trace to idempotent matrices yields a homomorphism
\[
\tr_{C^*_r(G)} \colon K_0(C^*_r(G)) \to \IR.
\]
Let $\pr\colon BG \to \pt$ be the projection.
For a group $G$
the following diagram commutes 
\begin{eqnarray}
\label{square for BCC implies TC for torsion free G}
\xymatrix{
K_0(BG) \ar[d]_-{K_0 ( \pr )} \ar[r]^-A & K_0(C^*_r(G)) \ar[r]^-{\tr_{C^*_r(G)}} & \IR  \\
K_0(\pt) \ar[r]^-{\cong}                &  K_0(\IC) \ar[r]_-{\tr_{\IC}}^-{\cong} & \IZ \ar[u]_-{i}.
         }
\end{eqnarray}
Here $i \colon \IZ \to \IR$ is the inclusion and $A$ is the assembly map. 
This non-trivial statement follows from
Atiyah's $L^2$-index theorem \cite{Atiyah(1976)}. Atiyah's theorem says that
the $L^2$-index $\tr_{C^*_r(G)} \circ A (\eta)$ of an element 
$\eta \in K_0(BG)$, which is represented by a pair $(M,P^*)$, 
agrees with the ordinary index of 
$(G\backslash M;G\backslash P^*)$, which is $\tr_{\IC} \circ K_0(\pr)(\eta) \in \IZ$. 

The following conjecture is taken from \cite[page 21]{Baum-Connes(1982)}.
\begin{conjecture}[Trace Conjecture for Torsion Free Groups]
\label{con: Trace Conjecture for Torsion Free Groups}
\index{Conjecture!Trace Conjecture!for Torsion Free Groups}
For a torsion free group $G$ the image of
$$\tr_{C^*_r(G)} \colon K_0(C^*_r(G)) \to \IR$$
consists of the integers.
\end{conjecture}

The commutativity of diagram \eqref{square for BCC implies TC for torsion free G} above implies

\begin{consequence} \label{cons: BCC implies Trace}
The surjectivity of the Baum-Connes assembly map
\[
K_0(BG) \to K_0(C^*_r(G))
\]
implies Conjecture~\ref{con: Trace Conjecture for Torsion Free Groups},
the Trace Conjecture for Torsion Free Groups.
\end{consequence}


\subsubsection{The Kadison Conjecture}
\label{subsec: The Kadison Conjecture and Related Conjectures}

\begin{conjecture}[Kadison Conjecture]
\label{con: Kadison Conjecture}
\index{Conjecture!Kadison Conjecture}
If $G$ is a torsion free  group, 
then the only idempotent elements in $C^*_r(G)$ are $0$ and $1$.
\end{conjecture}

\begin{lemma} \label{lem: TC implies Kadison}
The Trace Conjecture for Torsion Free Groups 
\ref{con: Trace Conjecture for Torsion Free Groups} implies the
Kadison Conjecture \ref{con: Kadison Conjecture}.
\end{lemma}
\begin{proof} Assume that $p \in C^*_r(G)$ is an idempotent different from
$0$ or $1$. 
From $p$ one can construct a non-trivial projection $q\in C^*_r(G)$, i.e.\ $q^2 = q$,
$q^* = q$, with $\im(p) = \im(q)$  and hence with $0 < q < 1$. 
Since the standard trace $\tr_{C^*_r(G)}$ is faithful, we
conclude $\tr_{C^*_r(G)}(q) \in \IR$ with $0 < \tr_{C^*_r(G)}(q) < 1$. Since
$\tr_{C^*_r(G)}(q)$ is by definition the image of the element $[\im(q)] \in K_0(C^*_r(G))$
under $\tr_{C^*_r(G)} \colon K_0(C^*_r(G)) \to \IR$, we get a
contradiction to the assumption $\im(\tr_{C^*_r(G)}) \subseteq \IZ$. 
\end{proof}

Recall that a ring $R$ is called an \emph{integral domain}
\index{ring!integral domain}
if it has no non-trivial zero-divisors, i.e.\ if $r,s \in R$ satisfy
$rs = 0$, then $r$ or $s$ is $0$.
Obviously the Kadison Conjecture \ref{con: Kadison Conjecture} 
implies for $R \subseteq \IC$ the following.

\begin{conjecture}[Idempotent Conjecture] 
\label{con: Idempotent Conjecture}
\index{Conjecture!Idempotent Conjecture}
Let $R$ be an integral domain and let $G$ be a torsion free group. Then the only idempotents
in $RG$ are $0$ and $1$.
\end{conjecture}
The statement in the conjecture above is a purely algebraic statement. If $R = \IC$,
it is by the arguments above related to
questions about operator algebras, and thus methods from operator algebras can be used
to attack it.


\subsubsection{Other Related Conjectures}
\label{subsec: Other Related Conjectures}

We would now like to mention  several conjectures which are not directly implied by 
the Baum-Connes or Farrell-Jones Conjectures, but which are closely related to the 
Kadison Conjecture and the Idempotent Conjecture mentioned above.

The next conjecture is also called the \emph{Kaplansky Conjecture}.%
\index{Conjecture!Kaplansky Conjecture}

\begin{conjecture}[Zero-Divisor-Conjecture] 
\label{con: Zero-Divisor-Conjecture}
\index{Conjecture!Zero-Divisor-Conjecture}
Let $R$ be an integral domain and $G$ be a torsion free group. Then $RG$ is an integral domain.
\end{conjecture}

Obviously the Zero-Divisor-Conjecture \ref{con: Zero-Divisor-Conjecture} implies the
Idempotent Conjecture \ref{con: Idempotent Conjecture}. The 
Zero-Divisor-Conjecture for $R = \IQ$ is implied by 
the following version of the Atiyah Conjecture
(see \cite[Lemma~10.5 and Lemma~10.15]{Lueck(2002)}).

\begin{conjecture}[Atiyah-Conjecture for Torsion Free Groups] 
\label{con: Atiyah-Conjecture for Torsion Free Groups}
\index{Conjecture!Atiyah-Conjecture for Torsion Free Groups}
Let $G$ be a torsion free group and let $M$ be a closed Riemannian manifold. Let
$\overline{M} \to M$ be a regular covering with $G$ as group of deck transformations.
Then all $L^2$-Betti numbers $b_p^{(2)}(\overline{M};\caln(G))$ are integers.
\end{conjecture}

For the precise definition and more information about $L^2$-Betti numbers and the group von Neumann algebra
$\caln(G)$ we refer for instance to
\cite{Lueck(2002)}, \cite{Lueck(2003h)}.

This more geometric formulation of the Atiyah Conjecture is in fact implied by the 
following more operator theoretic version. 
(The two would be equivalent if one would work with rational instead of complex coefficients below.)

\begin{conjecture}[Strong Atiyah-Conjecture for Torsion Free Groups] 
\label{con: Strong Atiyah-Conjecture for Torsion Free Groups}
\index{Conjecture!Atiyah-Conjecture for Torsion Free Groups!Strong}
Let $G$ be a torsion free group. Then for all $(m,n)$-matrices $A$ over $\IC G$ the von Neumann dimension of the 
kernel of the induced $G$-equivariant bounded operator 
\[
r_A^{(2)}\colon l^2(G)^m \to l^2(G)^n
\]
is an integer.
\end{conjecture}

The Strong Atiyah-Conjecture for Torsion Free Groups implies
both the Atiyah-Conjecture for Torsion Free Groups 
\ref{con: Atiyah-Conjecture for Torsion Free Groups} 
\cite[Lemma~10.5 on page 371]{Lueck(2002)} and the
Zero-Divisor-Conjecture \ref{con: Zero-Divisor-Conjecture} for $R = \IC$
\cite[Lemma~10.15 on page 376]{Lueck(2002)}.

\begin{conjecture}[Embedding Conjecture] 
\label{con: Embedding Conjecture}
\index{Conjecture!Embedding Conjecture}
Let $G$ be a torsion free group. Then $\IC G$ admits an embedding into a skewfield.
\end{conjecture}

Obviously the Embedding Conjecture implies
the  Zero-Divisor-Conjecture \ref{con: Zero-Divisor-Conjecture} for $R = \IC$.
If $G$ is a torsion free amenable group, then the
Strong Atiyah-Conjecture for Torsion Free Groups 
\ref{con: Strong Atiyah-Conjecture for Torsion Free Groups} and
the Zero-Divisor-Conjecture \ref{con: Zero-Divisor-Conjecture} for $R = \IC$
are equivalent \cite[Lemma~10.16 on page 376]{Lueck(2002)}.
For more information about the Atiyah Conjecture we refer for instance to
\cite{Linnell(1993)}, 
\cite[Chapter 10]{Lueck(2002)} and 
\cite{Reich(2003a)}. 

Finally we would like to mention the Unit Conjecture.
\begin{conjecture}[Unit-Conjecture] 
\label{con: Unit-Conjecture}
\index{Conjecture!Unit-Conjecture}
Let $R$ be an integral domain and $G$ be a torsion free group. Then every unit in $RG$ is trivial, i.e.\
of the form $r \cdot g$ for some unit $r \in R^{\inv}$ and $g \in G$.
\end{conjecture}

The Unit Conjecture \ref{con: Unit-Conjecture} implies the 
Zero-Divisor-Conjecture \ref{con: Zero-Divisor-Conjecture}. 
For a proof of this fact and for more information
we refer to \cite[Proposition 6.21 on page 95]{Lam(1991)}.

\subsubsection{$L^2$-Rho-Invariants and $L^2$-Signatures}
\label{subsec: L2-rho and signatures}
Let $M$ be a closed connected orientable Riemannian manifold. Denote by  $\eta(M) \in \IR$
\indexnotation{eta(M)}
the \emph{eta-invariant}%
\index{eta-invariant}
of $M$ and by $\eta^{(2)}(\widetilde{M}) \in \IR$%
\indexnotation{eta^(2)(widetilde M)}
the \emph{$L^2$-eta-invariant}%
\index{L2-eta-invariant@$L^2$-eta-invariant}
of the $\pi_1(M)$-covering given by the universal
covering $\widetilde{M} \to M$. Let
$\rho^{(2)}(M) \in \IR$%
\indexnotation{rho^(2)(M)}
be the \emph{$L^2$-rho-invariant}%
\index{L2-Rho-invariant@$L^2$-Rho-invariant}
which is defined to be the difference $\eta^{(2)}(\widetilde{M}) - \eta(M)$. These invariants
were studied by Cheeger and Gromov \cite{Cheeger-Gromov(1985)},
\cite{Cheeger-Gromov(1985a)}. 
They show that $\rho^{(2)}(M)$ depends only on the diffeomorphism type of  $M$
and is in contrast to $\eta(M)$ and $\eta^{(2)}(\widetilde{M})$ independent of the choice
of Riemannian metric on $M$. The following conjecture is taken from Mathai
\cite{Mathai(1992b)}. 

\begin{conjecture}[Homotopy Invariance of the $L^2$-Rho-Invariant for Torsionfree Groups] 
\label{con: Homotopy invariance of the L2-rho-invariant}
\index{Conjecture!Homotopy Invariance of the $L^2$-Rho-Invariant for Torsionfree Groups}
If $\pi_1(M)$ is torsionfree, then $\rho^{(2)}(M)$ is a homotopy invariant.
\end{conjecture}

Chang-Weinberger \cite{Chang-Weinberger(2003b)} assign to a closed connected oriented 
$(4k-1)$-dimensional manifold $M$ a Hirzebruch-type invariant $\tau^{(2)}(M) \in \IR$%
\indexnotation{tau^(2)} 
as follows. By a result of Hausmann \cite{Hausmann(1981)} 
there is a closed connected oriented $4k$-dimensional manifold $W$ with $M = \partial W$ such that
the inclusion $\partial W \to W$ induces an injection on the fundamental groups. 
Define $\tau^{(2)}(M)$ as the difference $\sign^{(2)}(\widetilde{W}) - \sign(W)$ of the
$L^2$-signature of the $\pi_1(W)$-covering given by the universal
covering $\widetilde{W} \to W$ and the signature of $W$. This is indeed independent
of the choice of $W$. It is reasonable to believe that $\rho^{(2)}(M) =
\tau^{(2)}(M)$ is always true.
Chang-Weinberger \cite{Chang-Weinberger(2003b)} 
use $\tau^{(2)}$ to prove that if $\pi_1(M)$ is not torsionfree there are infinitely many diffeomorphically distinct
manifolds of dimension $4k+3$ with $k \geq 1$, which are tangentially simple homotopy equivalent to $M$.

\begin{theorem}[Homotopy Invariance of $\tau^{(2)}(M)$ and  $\rho^{(2)}(M)$]
\label{the: Homotopy invariance of tau(2)(M) and  rho(2)(M)}
\indextheorem{Homotopy invariance of $\tau^{(2)}(M)$ and  $\rho^{(2)}(M)$}
 
\label{the: homotopy invariance of the rho-invariant}
Let $M$ be a closed connected oriented 
$(4k-1)$-dimensional manifold $M$ such that $G = \pi_1(M)$ is torsionfree. 

\begin{enumerate}
\item \label{the: homotopy invariance of the rho-invariant: Keswani}
If the assembly map $K_0 (BG) \to K_0(C^*_{\max}(G))$ for the maximal group
$C^*$-algebra is surjective 
(see Subsection \ref{subsec: The Baum-Connes Conjecture for Maximal Group $C^*$-Algebras}), 
then $\rho^{(2)}(M)$ is a homotopy invariant.

\item \label{the: homotopy invariance of the rho-invariant: Chang}
Suppose that the Farrell-Jones Conjecture for $L$-theory \ref{con: FJL torsion free} is
rationally true for $R = \IZ$, i.e.\ the rationalized assembly map
$$H_n(BG;\bfL^{\langle - \infty \rangle}(\IZ)) \otimes_{\IZ} \IQ \to 
L_n^{\langle - \infty \rangle}(\IZ G)\otimes_{\IZ} \IQ $$
is an isomorphism for $n \in \IZ$. Then $\tau^{(2)}(M)$ is a homotopy invariant.
If furthermore $G$ is residually finite, then $\rho^{(2)}(M)$ is a homotopy invariant.

\end{enumerate}
\end{theorem}
\begin{proof}
\ref{the: homotopy invariance of the rho-invariant: Keswani} This is proved by Keswani
\cite{Keswani(1998a)}, \cite{Keswani(2000)}. 
\\[2mm]
\ref{the: homotopy invariance of the rho-invariant: Chang}
This is proved by Chang \cite{Chang(2002)} and Chang-Weinberger 
\cite{Chang-Weinberger(2003b)} using \cite{Lueck-Schick(2001)}. 
\end{proof}

\begin{remarknew} \label{rem: L2-signature theorem}
Let $X$ be a $4n$-dimensional  Poincar{\'e} space over $\IQ$.
Let $\overline{X}\to X$ be a normal covering with \emph{torsion-free}
covering group $G$. 
Suppose that the assembly map $K_0 (BG) \to K_0(C^*_{\max}(G))$ for the maximal group
$C^*$-algebra is surjective 
(see Subsection \ref{subsec: The Baum-Connes Conjecture for Maximal Group $C^*$-Algebras})
or suppose that the rationalized assembly map for $L$-theory
$$H_{4n}(BG;\bfL^{\langle - \infty \rangle}(\IZ)) \otimes_{\IZ} \IQ \to 
L_{4n}^{\langle - \infty \rangle}(\IZ G)\otimes_{\IZ} \IQ $$
is an isomorphism.
Then the following $L^2$-signature
theorem is proved in L\"uck-Schick \cite[Theorem 2.3]{Lueck-Schick(2003)} 
\begin{equation*}
     \sign^{(2)}(\overline{X}) = \sign(X).
\end{equation*}

If one drops the condition that $G$ is torsionfree this equality
becomes false. Namely, Wall has constructed a finite Poincar\'e space $X$
with a finite $G$ covering $\overline{X} \to X$ for which 
$\sign(\overline{X}) \not= |G| \cdot \sign(X)$ holds (see \cite[Example 22.28]{Ranicki(1992)}, \cite[Corollary
5.4.1]{Wall(1967)}).
\end{remarknew}

\begin{remarknew} \label{rem: Cochran-Orr-Teichner}
Cochran-Orr-Teichner give in \cite{Cochran-Orr-Teichner(2003)} new obstructions
for a knot to be slice which are sharper than the Casson-Gordon
invariants. They use $L^2$-signatures and the Baum-Connes Conjecture 
\ref{con: Baum-Connes Conjecture}. We also refer to the survey article
\cite{Cochran(2001)} about non-commutative geometry and knot theory.
\end{remarknew}


\subsection{Applications V} 
\label{sec: Applications V}

\subsubsection{Novikov Conjectures}
\label{subsec: Novikov torsion free}

The Baum-Connes and Farrell-Jones Conjectures discussed so far 
imply obviously that for torsion free groups the rationalized assembly maps
\begin{eqnarray*}
H_{\ast} ( BG ; \bfK ( \IZ ) ) \otimes_{\IZ} \IQ  & \to & K_{\ast} ( \IZ G ) \otimes_{\IZ} \IQ \\
H_{\ast} (BG ; \bfL^{\langle - \infty \rangle} ( \IZ ) ) \otimes_{\IZ} \IQ & \to & 
L_{\ast}^{\langle - \infty \rangle} ( \IZ G ) \otimes_{\IZ} \IQ \\
K_{\ast} ( BG ) \otimes_{\IZ} \IQ & \to & K_{\ast} ( C_r^{\ast} ( G ) ) \otimes_{\IZ} \IQ
\end{eqnarray*}
are injective. For reasons that will be explained below these ``rational injectivity conjectures'' are known as 
``Novikov Conjectures''.
In fact one expects these injectivity results also when the groups contain torsion. So there are the following
conjectures.
\begin{conjecture}[$K$- and $L$-theoretic Novikov Conjectures] 
\label{con: K and L Novikov}
\index{Conjecture!Novikov Conjecture!for $K$ and $L$-theory}
For every group $G$ the assembly maps
\begin{eqnarray*}
H_{\ast} ( BG ; \bfK ( \IZ ) ) \otimes_{\IZ} \IQ  & \to & K_{\ast} ( \IZ G ) \otimes_{\IZ} \IQ \\
H_{\ast} (BG ; \bfL^p ( \IZ ) ) \otimes_{\IZ} \IQ & \to & 
L_{\ast}^p ( \IZ G ) \otimes_{\IZ} \IQ \\
K_{\ast} ( BG ) \otimes_{\IZ} \IQ & \to & K_{\ast} ( C_r^{\ast} ( G ) ) \otimes_{\IZ} \IQ
\end{eqnarray*}
are injective.
\end{conjecture}

Observe that, since the $\IZ / 2 $-Tate cohomology groups vanish rationally, there is no difference between the various
decorations in $L$-theory  because of the Rothenberg sequence. We have chosen the $p$-decoration above.


\subsubsection{The Original Novikov Conjecture}
\label{subsec: The Original Novikov Conjecture}

We now explain the Novikov Conjecture in its original formulation.

Let $G$ be a (not necessarily torsion free)  group and $u \colon M \to BG$ be a map from a closed oriented 
smooth manifold $M$ to $BG$. Let $\call(M)%
\index{L-class@$L$-class}
\indexnotation{call(M)}
\in \prod_{k \ge 0} H^k(M;\IQ)$ 
be the \emph{$L$-class of $M$},
which is a certain polynomial in the Pontrjagin classes 
and hence depends a priori on the tangent bundle and hence on the differentiable structure of $M$.
For 
$x \in  \prod_{k \ge 0} H^k(BG;\IQ)$ define the 
\emph{higher signature of $M$ associated to $x$ and $u$}%
\index{signature!higher} 
to be
\begin{eqnarray}
\sign_x(M,u)%
\indexnotation{sign_x(M)} & := & \langle \call(M) \cup u^*x,[M]\rangle \hspace{5mm} \in \IQ.
\label{higher signature}
\end{eqnarray}
The Hirzebruch signature formula says that for $x=1$ the signature $\sign_1(M,u)$ coincides with the ordinary signature
$\sign(M)$ of $M$, if $\dim(M) = 4n$, and is zero, if 
$\dim(M)$ is not divisible by four. Recall that for $\dim(M) = 4n$ the \emph{signature}%
\index{signature}
$\sign(M)$%
\indexnotation{sign(M)}
of $M$ is the signature of the non-degenerate bilinear symmetric pairing
on the middle cohomology $H^{2n}(M;\IR)$ given by the intersection pairing
$(a,b) \mapsto \langle a\cup b,[M]\rangle$.
Obviously $\sign(M)$ depends only on the oriented homotopy type of $M$.
We say that $\sign_x$ for $x \in H^*(BG;\IQ)$ is \emph{homotopy invariant} if
for two closed oriented smooth manifolds $M$ and $N$ with reference maps
$u\colon M \to BG$ and $v \colon N \to BG$ we have 
\[
\sign_x(M,u) = \sign_x(N,v)
\]
if there is an orientation preserving homotopy equivalence $f \colon M \to N$ such that
$v \circ f$ and $u$ are homotopic.

\begin{conjecture}[Novikov Conjecture]%
\index{Conjecture!Novikov Conjecture}
\label{con: Novikov Conjecture}
Let $G$ be a group. 
Then $\sign_x$ is homotopy invariant for all $x \in  \prod_{k \ge 0} H^k(BG;\IQ)$.
\end{conjecture}

By Hirzebruch's signature formula
the Novikov Conjecture \ref{con: Novikov Conjecture} is true for $x = 1$. 


\subsubsection{Relations between the Novikov Conjectures}
\label{subsec: Relations between the Novikov Conjectures}

Using surgery theory one can show \cite[Proposition 6.6 on page 300]{Ranicki(1995b)} the following.

\begin{propositionnew}\label{lem: Novikov = Strong Novikov for L-theory}

For a group $G$ the original Novikov Conjecture~\ref{con: Novikov Conjecture} is equivalent to 
the $L$-theoretic Novikov Conjecture, i.e.\ the 
injectivity of the assembly map
\[
H_{\ast} (BG ; \bfL^{p} ( \IZ ) ) \otimes_{\IZ} \IQ  \to 
L_{\ast}^{p} ( \IZ G ) \otimes_{\IZ} \IQ .\\
\]
\end{propositionnew}

In particular for torsion free groups the $L$-theoretic Farrell-Jones Conjecture~\ref{con: FJL torsion free}
implies the Novikov Conjecture \ref{con: Novikov Conjecture}.
Later in Proposition~\ref{pro: relating L and K}
we will prove in particular the following statement.

\begin{propositionnew} \label{pro: Novi for to K-theory implies Novikov}
The Novikov Conjecture for topological $K$-theory, i.e.\ the injectivity of the assembly map
\[
K_{\ast} ( BG ) \otimes_{\IZ} \IQ  \to  K_{\ast} ( C_r^{\ast} ( G ) ) \otimes_{\IZ} \IQ
\]
implies the $L$-theoretic Novikov Conjecture and hence the original Novikov Conjecture.
\end{propositionnew}

For more information about the Novikov Conjectures we refer for instance to
\cite{Boekstedt-Hsiang-Madsen(1993)},
\cite{Carlsson(2004)}, 
\cite{Carlsson-Pedersen(1995a)},
\cite{Davis(2000)},
\cite{Farrell(2002)}, 
\cite{Ferry-Ranicki-Rosenberg(1995)}, 
\cite{Kreck-Lueck(2004)},
\cite{Ranicki(1992)} and 
\cite{Rosenberg(1995)}.


\subsubsection{The Zero-in-the-Spectrum Conjecture}
\label{subsec: Zero-in-the-Spectrum Conjecture}

The following Conjecture
is due to Gromov \cite[page 120]{Gromov(1986a)}. 

\begin{conjecture}[Zero-in-the-spectrum Conjecture]%
\index{Conjecture!Zero-in-the-spectrum Conjecture}
\label{con: zero-in-the-spectrum Conjecture}
Suppose that $\widetilde{M}$ is the universal covering of an
aspherical closed Riemannian manifold $M$ (equipped with the lifted Riemannian metric).
Then  zero is  in the
spectrum of the minimal closure
\[
(\Delta_p)_{\min} \colon  
L^2 \Omega^p ( \widetilde{M} ) 
\supset
 \dom  ( \Delta_p )_{\min}
\to L^2\Omega^p(\widetilde{M}),
\]
for some $p \in \{  0, 1, \dots , \dim M \}$,
where $\Delta_p$ denotes the Laplacian acting on smooth $p$-forms on $\widetilde{M}$.
\end{conjecture}

\begin{propositionnew} \label{the: Novikov implies zero}
Suppose that $M$ is an aspherical closed Riemannian manifold with fundamental group $G$, then
the injectivity of the assembly map
\[
K_{\ast} ( BG ) \otimes_{\IZ} \IQ  \to  K_{\ast} ( C_r^{\ast} ( G ) ) \otimes_{\IZ} \IQ
\]
implies the Zero-in-the-spectrum Conjecture for $\widetilde{M}$.
\end{propositionnew}
\begin{proof} 
We give a sketch of the proof. More details can be found
in \cite[Corollary 4]{Lott(1996b)}. We only explain
that the assumption that in every dimension zero is not in the spectrum of
the Laplacian on $\widetilde{M}$, yields
a contradiction in the case that $n = \dim(M)$ is even. Namely,
this assumption implies that the $C_r^*(G)$-valued index of the
signature operator twisted with the flat bundle
$\widetilde{M} \times_{G} C_r^*(G) \to M$
in $K_0(C_r^*(G))$ is zero, where $G = \pi_1(M)$.  This index is the image of the class $[S]$
defined by the signature operator in $K_0(BG)$ under the assembly
map $K_0(BG) \to K_0(C^*_r(G))$. Since by assumption the assembly
map is rationally injective, this implies $[S] = 0$ in $K_0(BG) \otimes_{\IZ} \IQ$. Notice that
$M$ is aspherical by assumption and hence $M = BG$. The
homological Chern character defines an isomorphism
\[
K_0(BG) \otimes_{\IZ} \IQ = K_0(M) \otimes_{\IZ} \IQ \xrightarrow{\cong}
\bigoplus_{p \ge 0} H_{2p}(M;\IQ)
\]
which sends $[S]$ to the Poincar\'e dual $\call (M) \cap [M]$
of the Hirzebruch $L$-class  $\call (M) \in \bigoplus_{p \ge 0} H^{2p}(M;\IQ)$. This implies that
$\call (M) \cap [M] = 0$ and hence $\call(M) = 0$.
This contradicts the fact that the component of $\call(M)$ in $H^0(M;\IQ)$ is
$1$. \end{proof}

More information about the Zero-in-the-spectrum Conjecture 
\ref{con: zero-in-the-spectrum Conjecture} can be found for instance in
\cite{Lott(1996b)} and 
\cite[Section 12]{Lueck(2002)}.

\typeout{---  Formulation of the Conjectures in the general  case -------}

\section{The Conjectures in the General Case}
\label{chap: general formulation}

In this chapter we will formulate the Baum-Connes and Farrell-Jones Conjectures.
We try to emphasize the unifying principle that underlies these conjectures.
The point of view taken in this chapter is that all three conjectures are conjectures
about specific equivariant homology theories. 
Some of the technical details concerning the actual construction of these homology theories
are deferred to Chapter~\ref{chap: Equivariant Homology Theories}.

\subsection{Formulation of the Conjectures}
\label{sec: Formulation of the Conjectures}

Suppose we are given 
\begin{itemize}
\item
A discrete group $G$;
\item
A  family $\calf$ of subgroups of $G$, i.e.\ a set of subgroups 
which is closed under conjugation with elements of $G$ and under taking finite intersections;
\item
A $G$-homology theory $\calh^G_{\ast} ( - )$.
\end{itemize}
Then one can formulate the following Meta-Conjecture. 
\begin{metaconjecture} \label{metaconjecture}
\index{Meta-Conjecture} \index{Conjecture!Meta-Conjecture}
The assembly map%
\index{assembly map}
\indexnotation{A_calf}
\[
A_{\calf} \colon \calh^G_n ( \EGF{G}{\calf}) \to \calh^G_n( \pt )
\]
which is the map induced by the projection $\EGF{G}{\calf} \to \pt$,
is an isomorphism for $n \in \IZ$.
\end{metaconjecture}
Here $\EGF{G}{\calf }$ is the classifying space of the family $\calf$, 
a certain $G$-space which specializes to the universal free $G$-space
 $EG$ if the family contains only the trivial subgroup.
A $G$-homology theory is the ``obvious'' $G$-equivariant generalization of the concept of 
a homology theory to a suitable category of $G$-spaces, in particular it is a functor on such spaces
and the map $A_{\calf}$ is simply the map induced by the projection $\EGF{G}{\calf } \to \pt$.
We devote the Subsections~\ref{subsec: G-CW complexes} to \ref{subsec: G-Homology Theories} below to a 
discussion of $G$-homology theories, classifying spaces for families 
of subgroups and related things. The reader who never encountered
these concepts should maybe first take a look at these subsections.

Of course the conjecture above is not true for arbitrary $G$, $\calf$ and $\calh^G_{\ast}( - )$, but 
the Farrell-Jones and Baum-Connes Conjectures state that for specific $G$-homology theories
there is a natural choice of a family $\calf=\calf(G)$ of subgroups for every group $G$ such that $A_{\calf(G)}$
becomes an isomorphism for all groups $G$.

Let $R$ be a ring (with involution).
In Proposition \ref{pro: Or(G)-spectra yield a G-homology theory} 
we will describe the construction of $G$-homology theories
which will be denoted 
\[
H_n^G ( - ; \bfK_R ),  \quad H_n^G ( - ; \bfL_R^{\langle -\infty
  \rangle}  ) 
\quad \mbox{ and } \quad H_n^G ( - ; \bfK^{\topo} ).
\]
If $G$ is the trivial group, these homology theories specialize 
to the (non-equivariant) homology theories with similar names that
appeared in Chapter~\ref{chap: torsion free}, namely to 
\[
H_n ( - ; \bfK (R) ) , \quad H_n ( - ; \bfL^{\langle -\infty \rangle} (R)) \quad \mbox{ and } \quad K_n( - ).
\] 
Another main feature of these $G$-homology theories is that evaluated on the one  
point space $\pt$ (considered as a trivial $G$-space)
we obtain the $K$- and $L$-theory of the group ring $RG$, respectively 
the topological $K$-theory of the reduced $C^{\ast}$-algebra 
(see Proposition \ref{pro: Or(G)-spectra yield a G-homology theory} and
Theorem~\ref{the: K- and L-Theory Spectra over Groupoids}
\ref{the: K- and L-Theory Spectra over Groupoids: values at groups})
\begin{eqnarray*}
K_n( RG ) & \cong & H_n^G ( \pt  ; \bfK_R ), \\
L_n^{\langle - \infty \rangle} ( RG ) & \cong & H_n^G ( \pt ;
\bfL_R^{\langle -\infty \rangle}  )  
\quad  \mbox{ and }  \\
K_n ( C^{\ast}_r (G ) ) & \cong & H_n^G ( \pt ; \bfK^{\topo} ).
\end{eqnarray*}


We are now prepared to formulate the conjectures around which this article is centered.
Let $\calfin$ be the family of finite subgroups and let $\calvcyc$ be the family 
of virtually cyclic subgroups.

\begin{conjecture}[Farrell-Jones Conjecture for $K$- and $L$-theory] \label{con: Farrell-Jones Conjecture}
Let $R$ be a ring (with involution) and let $G$ be a group.
Then for all $n \in \IZ$ the maps
\begin{eqnarray*}
A_{\calvcyc} \colon H_n^G( \EGF{G}{\calvcyc} ; \bfK_R ) & \to & 
H_n^G ( \pt ; \bfK_R ) \cong K_n( RG ); \\
A_{\calvcyc} \colon H_n^G( \EGF{G}{\calvcyc} ; \bfL_R^{\langle - \infty \rangle} ) & \to & 
H_n^G ( \pt ; \bfL_R^{\langle - \infty \rangle} )  \cong L_n^{\langle -\infty \rangle} ( RG ),
\end{eqnarray*}
which are induced by the projection $\EGF{G}{\calvcyc} \to \pt$, are isomorphisms.
\end{conjecture}

The conjecture for the topological $K$-theory of $C^{\ast}$-algebras is known as the Baum-Connes Conjecture
and reads as follows.
\begin{conjecture}[Baum-Connes Conjecture] \label{con: Baum-Connes Conjecture}
Let $G$ be a group. Then for all $n \in \IZ$ the map
\[
A_{\calfin} \colon H_n^G( \EGF {G}{\calfin} ; \bfK^{\topo}  ) \to 
H_n^G ( \pt ; \bfK^{\topo} ) \cong K_n( C^{\ast}_r (G) ) \\
\]
induced by the projection $\EGF{G}{\calfin} \to \pt$ is an isomorphism.
\end{conjecture}

We will explain the analytic assembly map 
$\ind_G \colon K_n^G(X) \to K_n(C^*_r(G))$, which can be identified with the
assembly map appearing in the Baum-Connes Conjecture \ref{con: Baum-Connes Conjecture} in
Section \ref{sec: The Analytic Assembly Map}.

\begin{remarknew} Of course the conjectures really come to life only if the abstract point of view taken in this
chapter is connected up with more concrete descriptions of the assembly maps.
We have already discussed a surgery theoretic description in Theorem~\ref{the: The Farrell-Jones Conjecture Implies the Borel Conjecture}
and an interpretation in terms index theory in Subsection \ref{subsec: The Trace Conjecture in the Torsion Free Case}.
More information about alternative interpretations of assembly maps can be found
in Section \ref{sec: The Analytic Assembly Map} and \ref{sec:Assembly as forget control}.
These concrete interpretations of the assembly maps lead to applications.
We already discussed many such applications in Chapter~\ref{chap: torsion free} and encourage the reader to go 
ahead and browse through Chapter~\ref{chap: more applications} in order to get further ideas about these 
more concrete aspects.
\end{remarknew}

\begin{remarknew}[Relation to the ``classical'' assembly maps]
The value of an equivariant homology theory $\calh_{\ast}^G ( - )$
on the universal free $G$-space $EG=\EGF{G}{\trivial}$ (a free $G$-$CW$-complex whose quotient $EG/G$ is a model for $BG$)
can be identified with the corresponding non-equivariant homology
theory evaluated on $BG$, if we assume that $\calh^G_*$ is the special value of
an equivariant homology theory $\calh^?_*$ at $?=G$.
This means that there exists an induction structure 
(a mild condition satisfied in our examples, compare Section
\ref{sec: The Definition of an Equivariant Homology Theory}), which yields
an identification
\[
\calh^G_n ( EG ) \cong \calh^{\trivial}_n ( BG )= \calh_n ( BG ).
\]
Using these identifications the ``classical'' assembly maps, which
appeared in Chapter~\ref{chap: torsion free} in the versions of the Farrell-Jones and Baum-Connes Conjectures 
for torsion free groups  (see Conjecture~\ref{con: FJK torsion free
  all},~\ref{con: FJL torsion free} 
and~\ref{con: BC torsion free}), 
\begin{eqnarray*}
H_n ( BG ; \bfK (R ))   \cong   H_n^G ( EG ; \bfK_R ) & \to &
H_n^G ( \pt  ; \bfK_R ) \cong K_n( RG );  \\
H_n ( BG ; \bfL^{\langle -\infty \rangle } (R ))  \cong   H_n^G ( EG ; \bfL_R^{\langle - \infty \rangle} ) & \to &
H_n^G ( \pt  ; \bfL_R^{ \langle -\infty \rangle}  ) \cong L_n^{\langle -\infty \rangle } ( RG ); \\
\mbox{ and } \quad K_n ( BG ) \cong H_n^G ( EG ; \bfK^{\topo} ) 
& \to & H_n^G ( \pt ; \bfK^{\topo} ) \cong K_n ( C_r^{\ast} ( G ) ),
\end{eqnarray*}
correspond to the assembly maps for the family $\calf=\trivial$ consisting only of the trivial group and 
are simply the maps induced by the projection $EG \to \pt$. 
\end{remarknew}

\begin{remarknew}[The choice of the right family]
As explained above the Farrell-Jones and Baum-Connes Conjectures~\ref{con: Farrell-Jones Conjecture} and 
\ref{con: Baum-Connes Conjecture}
can be considered as 
special cases of the Meta-Conjecture~\ref{metaconjecture}.
In all three cases we are interested in a computation of the right hand side $\calh_n^G ( \pt )$ of the 
assembly map, which can be identified with 
$K_n ( RG )$, $L_n^{\langle -\infty \rangle} ( RG )$ or $K_n ( C_r^{\ast} ( G ) )$.
The left hand side $\calh^G_n ( \EGF{G}{\calf})$ of such an assembly map is much more accessible
and the smaller $\calf$ is, the easier it is to compute $\calh^G_n(\EGF{G}{\calf})$  using  homological 
methods like spectral sequences, Mayer-Vietoris arguments and Chern characters.

In the extreme case where $\calf= \calall$ is  the family of all subgroups the assembly map 
$A_{\calall} \colon \calh_n^G ( \EGF{G}{\calall} ) \to \calh_n^G ( \pt )$ is always an isomorphism for the 
trivial reason that  the one point space $\pt$ is a model for
$\EGF{G}{\calall}$ (compare 
Subsection~\ref{subsec: EGFs})
and hence the assembly map is the identity. The goal however is to
have an isomorphism for 
a family which is as small as possible.

We have already seen in  Remark~\ref{rem: torsion free is necessary}, Remark~\ref{rem: L for finite groups}
and Remark~\ref{rem: BC for finite groups} that in all three cases the classical assembly map, which 
corresponds to the trivial family, is not surjective for finite groups. 
This forces one to include at least the family $\calfin$ of finite groups.
The $K$- or $L$-theory of the finite subgroups of the given group 
$G$ will then enter in a computation of the left hand side of 
the assembly map similar as the $K$- and $L$-theory of the trivial subgroup appeared on the left hand side in the 
classical case, compare e.g.\ Remark~\ref{rem: coefficients K}.
In the Baum-Connes case the family $\calfin$ seems to suffice.
However in the case of algebraic $K$-theory we saw in Remark~\ref{rem: Bass-Heller-Swan decomposition} 
that already the simplest torsion free group, the infinite cyclic group, causes problems because of 
the Nil-terms that appear in the Bass-Heller-Swan formula. 
The infinite dihedral group is a ``minimal counterexample'' which shows that the family $\calfin$ is not 
sufficient in the $L_{\IZ}^{\langle - \infty \rangle}$-case. There are non-vanishing UNil-terms, 
compare~\ref{subsec: L invert 2} and \ref{subsec: L second kind}.  Also the version of the $L$-theoretic 
Farrell-Jones Conjecture with the decoration $s$, $h= \langle 1 \rangle$ or $p = \langle 0 \rangle$  instead of 
$\langle -\infty \rangle$ is definitely false.
Counterexamples are given in \cite{Farrell-Jones-Lueck(2002)}.
Recall that there were no Nil-terms in the topological $K$-theory
context, compare Remark~\ref{Computing K(C^*_r(G times Z)}.

The choice of the family
$\calvcyc$ of virtually cyclic subgroups 
in the Farrell-Jones conjectures pushes all the Nil-problems appearing in algebraic $K$- and
$L$-theory into the source of the assembly map so that they do not occur if  one tries to 
prove the Farrell-Jones Conjecture~\ref{con: Farrell-Jones Conjecture}.
Of course they do come up again when one wants to compute the source of the assembly map. 
\end{remarknew}


We now take up the promised detailed discussion of some notions like equivariant homology theories and 
classifying spaces for families we used above.
The reader who is familiar with these concepts may of course skip the following subsections.

\subsubsection{$G$-CW-Complexes} 
\label{subsec: G-CW complexes}

A \emph{$G$-$CW$-complex $X$}%
\index{G-CW-complex@$G$-$CW$-complex}
is a $G$-space $X$ together with a filtration
$X_{-1} = \emptyset \subseteq X_0 \subseteq X_1 \subseteq \ldots  \subseteq X$ such that
$X = \colim_{n \to \infty} X_n$ and for each $n$ there is a $G$-pushout
\[
\comsquare{\coprod_{i \in I_n} G/H_i \times S^{n-1}}{\coprod_{i \in I_n} q^n_i}{X_{n-1}}
{}{}{\coprod_{i \in I_n} G/H_i \times D^n}{\coprod_{i \in I_n} Q^n_i}{X_n}.
\]
This definition makes also sense for topological groups. The following alternative definition only applies to discrete groups.
A $G$-$CW$-complex is a $CW$-complex 
with a $G$-action by cellular maps such that for each open cell
$e$ and each $g \in G$ with $ge \cap e \not= \emptyset$ we have $gx = x$ for all $x \in e$.
There is an obvious notion of a $G$-$CW$-pair.

A $G$-$CW$-complex $X$ is called \emph{finite}%
\index{G-CW-complex@$G$-$CW$-complex!finite} if it is built out of finitely many $G$-cells $G/H_i \times D^n$.
This is the case if and only if it is \emph{cocompact},%
\index{cocompact}
i.e.\ the quotient space $G \backslash X$ is compact.
More information about $G$-$CW$-complexes can be found for instance in
\cite[Sections 1 and 2]{Lueck(1989)}, \cite[Sections II.1 and II.2]{Dieck(1987)}.


\subsubsection{Families of Subgroups}
\label{subsec: families}

A \emph{family $\calf$ of subgroups of $G$}%
\index{family of subgroups}
is a set of subgroups of $G$ closed under conjugation, i.e.\ $H \in \calf, g \in G$
implies $g^{-1}Hg \in \calf$, and finite intersections, i.e.\ $H,K \in
\calf$ implies $H \cap K \in \calf$. 
Throughout the text we will use the notations 
$$
\caltr,%
\indexnotation{trivial family}  
\quad  
\calfcyc,%
\indexnotation{familycalfcyc}
 \quad \calfin,%
\indexnotation{familycalfin}
 \quad 
\calcyc,%
\indexnotation{calcyc}
\quad 
\calvcyc_{I},%
\indexnotation{calvcyc_I}  
\quad 
\calvcyc%
\indexnotation{familycalvcyc} 
\quad
\quad \mbox{ and } \quad 
\calall%
\indexnotation{familycalall}
$$
for the families consisting of the trivial, all finite cyclic, all finite, all (possibly infinite) cyclic,  
all virtually cyclic of the first kind, all virtually cyclic, respectively all subgroups of a given group $G$. 
Recall that a group is called \emph{virtually cyclic}%
\index{group!virtually cyclic} \index{virtually cyclic}
if it is finite or contains an infinite cyclic subgroup of finite index.
A group is \emph{virtually cyclic of the first kind}%
\index{group!virtually cyclic of the first kind} if it admits a surjection 
onto an infinite cyclic group with finite kernel, compare
Lemma~\ref{lem: virtually cyclic}.

\subsubsection{Classifying Spaces for Families}
\label{subsec: EGFs}

Let $\calf$ be a family of subgroups of $G$. A $G$-CW-complex, all whose isotropy groups belong to $\calf$ and
whose $H$-fixed point sets are contractible for all $H \in \calf$, is called a 
\emph{classifying space for the family $\calf$}%
\index{classifying space!for a family} 
and will be denoted $\EGF{G}{\calf}$.%
\indexnotation{E(G;calf)}
Such a space is unique up to $G$-homotopy because it is characterized
by the property that for any $G$-$CW$-complex $X$, all whose isotropy groups
belong to $\calf$, there is  up to $G$-homotopy precisely one $G$-map
from $X$ to $\EGF{G}{\calf}$. These spaces  were introduced by tom Dieck
\cite{Dieck(1972)}, \cite[I.6]{Dieck(1987)}.

A functorial ``bar-type'' construction is given in
\cite[section 7]{Davis-Lueck(1998)}. 

If $\calf \subset \calg$ are families of subgroups for $G$, then by
the universal property there is up to $G$-homotopy precisely one
$G$-map $\EGF{G}{\calf} \to \EGF{G}{\calg}$.

The space $\EGF{G}{\caltr}$ is the same as the space $EG$%
\indexnotation{EG}
which is by definition the total space of the universal $G$-principal bundle $G \to EG
\to BG$, or, equivalently, the universal covering of $BG$.
A model for $\EGF{G}{\calall}$ is given by the space $G/G =\pt$  consisting of
 one point.  

The space $\EGF{G}{\calfin}$ is also known as
\emph{the classifying space for proper $G$-actions}%
\index{classifying space!for proper G-actions@for proper $G$-actions}
and denoted by $\underline{E}G$%
\indexnotation{underline{E}G}
in the literature. Recall that a $G$-$CW$-complex $X$ is
proper%
\index{G-CW-complex@$G$-$CW$-complex!proper}
 if and only if all its isotropy groups are finite
(see for instance \cite[Theorem 1.23 on page 18]{Lueck(1989)}). There are often nice models
for $\EGF{G}{\calfin}$. If $G$ is word hyperbolic in the sense of Gromov, then the
Rips-complex is a finite model \cite{Meintrup(2000)}, \cite{Meintrup-Schick(2002)}. 

If $G$ is a discrete subgroup of 
a Lie group $L$ with finitely many path components, then for any maximal compact subgroup $K \subseteq L$
the space $L/K$ with its left $G$-action is a model for $\EGF{G}{\calfin}$
\cite[Corollary 4.14]{Abels(1978)}. More information about $\EGF{G}{\calfin}$ can be found
for instance in 
\cite[section 2]{Baum-Connes-Higson(1994)},
\cite{Kropholler-Mislin(1998)}, 
\cite{Lueck(2000a)}, 
\cite{Lueck(2004a)}, 
\cite{Lueck-Meintrup(2000)} and 
\cite{Serre(1971)}.

\subsubsection{$G$-Homology Theories}
\label{subsec: G-Homology Theories}

Fix a group $G$ and an associative commutative ring $\Lambda$ with unit.
A \emph{$G$-homology theory $\calh_*^G$ with values
in $\Lambda$-modules}%
\indexnotation{calh_*^G}
\index{homology theory!G-homology theory@$G$-homology theory}
is a collection of
covariant functors $\calh^G_n$ from the category of
$G$-$CW$-pairs to the category of
$\Lambda$-modules indexed by $n \in \IZ$ together with natural transformations
$$\partial_n^G(X,A)\colon \calh_n^G(X,A) \to
\calh_{n-1}^G(A):= \calh_{n-1}^G(A,\emptyset)$$
for $n \in \IZ$ such that the following axioms are satisfied:
\begin{enumerate}

\item $G$-homotopy invariance \\[1mm]
If $f_0$ and $f_1$ are $G$-homotopic maps $(X,A) \to (Y,B)$
of $G$-$CW$-pairs, then $\calh_n^G(f_0) = \calh^G_n(f_1)$ for $n \in \IZ$.

\item Long exact sequence of a pair \\[1mm]
Given a pair $(X,A)$ of $G$-$CW$-complexes,
there is a long exact sequence
\begin{multline*}\ldots \xrightarrow{\calh^G_{n+1}(j)}
\calh_{n+1}^G(X,A) \xrightarrow{\partial_{n+1}^G}
\calh_n^G(A) \xrightarrow{\calh^G_{n}(i)} \calh_n^G(X)
\\
\xrightarrow{\calh^G_{n}(j)} \calh_n^G(X,A)
\xrightarrow{\partial_n^G} \calh_{n-1}^G(A) \xrightarrow{\calh^G_{n-1}(i)} \ldots,
\end{multline*}
where $i \colon A \to X$ and $j\colon  X \to (X,A)$ are the inclusions.

\item Excision \\[1mm]
Let $(X,A)$ be a $G$-$CW$-pair and let
$f\colon  A \to B$ be a cellular $G$-map of
$G$-$CW$-complexes. Equip $(X\cup_f B,B)$ with the induced structure
of a $G$-$CW$-pair. Then the canonical map
$(F,f)\colon  (X,A) \to (X\cup_f B,B)$ induces for each $n \in \IZ$ an isomorphism
\[
\calh_n^G(F,f)\colon  \calh_n^G(X,A) \xrightarrow{\cong}
\calh_n^G(X\cup_f B,B).
\]

\item Disjoint union axiom%
\index{disjoint union axiom}
\\[1mm]
Let $\{X_i\mid i \in I\}$ be a family of
$G$-$CW$-complexes. Denote by
$j_i\colon  X_i \to \coprod_{i \in I} X_i$ the canonical inclusion.
Then the map
$$\bigoplus_{i \in I} \calh^G_{n}(j_i)\colon  \bigoplus_{i \in I} \calh_n^G(X_i)
\xrightarrow{\cong} \calh_n^G\left(\coprod_{i \in I} X_i\right)$$
is bijective for each $n \in \IZ$.

\end{enumerate}

Of course a $G$-homology theory for the trivial group $G=\{ 1 \}$ 
is a homology theory 
(satisfying the disjoint union axiom) in the classical non-equivariant sense.


The disjoint union axiom ensures that we can pass from 
finite $G$-CW-complexes to arbitrary ones using the following lemma.

\begin{lemma} \label{lem: G-homology theory and colimit}
Let $\calh^G_*$ be a $G$-homology theory. 
Let $X$ be a $G$-$CW$-complex and $\{X_i \mid i \in I\}$ be a directed system of
$G$-$CW$-subcomplexes directed by inclusion such that $X = \cup_{i \in I} X_i$.
Then for all $n \in \IZ$ the natural map
$$\colim_{i \in I} \calh_n^G(X_i) \xrightarrow{\cong}   \calh^G_n(X)$$
is bijective.
\end{lemma}
\begin{proof} Compare for example with \cite[Proposition 7.53 on page 121]{Switzer(1975)},
where the non-equivariant case for $I=\IN$ is treated.
\end{proof}

\begin{examplenew}[Bredon Homology]  \label{Bredon homology as G-homology theory}   
The most basic $G$-homology theory is \emph{Bredon homology}.%
\index{Bredon homology}
The \emph{orbit category}%
\index{orbit category}
$\Or (G)$%
\indexnotation{Or(G)} has as objects the homogeneous spaces $G/H$ and as
morphisms $G$-maps. Let $X$ be a $G$-$CW$-complex. 
It defines a contravariant functor 
from the orbit category $\Or(G)$ to the category of $CW$-complexes by sending $G/H$ to
$\map_G(G/H,X) = X^H$. Composing it with the functor
cellular chain complex  yields a contravariant functor
\[
C_*^c(X) \colon \Or (G) \to \IZ\text{-}\CHAINCOMPLEXES
\]
into the category of $\IZ$-chain complexes. Let $\Lambda$ be a commutative ring and let
$$M \colon \Or(G) \to \Lambda\text{-}\MODULES$$
be a covariant functor. Then one can form the tensor product over the orbit category
(see for instance \cite[9.12 on page 166]{Lueck(1989)})
and obtains the  $\Lambda$-chain complex $C_*^c(X) \otimes_{\IZ\Or(G)} M$.
Its homology is the Bredon homology  of $X$ with coefficients in $M$
\[
H_*^G(X;M)%
\indexnotation{H_*^G(X;M)}
 = H_* ( C_*^c(X) \otimes_{\IZ\Or(G)} M ).
\]
Thus we get a $G$-homology theory $H^G_*$ with values in $\Lambda$-modules.
For a trivial group $G$ this reduces to the cellular homology of
$X$ with coefficients in the $\Lambda$-module $M$. 
\end{examplenew}

More information about equivariant homology theories will be given in
Section \ref{sec: The Definition of an Equivariant Homology Theory}.


\subsection{Varying the Family of Subgroups}
\label{sec: Varying the Family of Subgroups}

Suppose we are given a family of subgroups $\calf^{\prime}$ and a subfamily $\calf \subset \calf^{\prime}$.
Since all isotropy groups of $\EGF{G}{\calf}$ lie in $\calf^{\prime}$ we know from the universal property of 
$\EGF{G}{\calf^{\prime}}$ (compare Subsection~\ref{subsec: EGFs}) 
that there is a $G$-map $\EGF{G}{\calf} \to \EGF{G}{\calf^{\prime}}$ which is unique up to 
$G$-homotopy. For every $G$-homology theory $\calh_*^G$ we hence obtain a \emph{relative assembly map}%
\index{assembly map!relative} \index{relative assembly map}
\[
A_{\calf \to \calf^{\prime}} \colon \calh_n^G ( \EGF{G}{\calf} ) \to \calh_n^G ( \EGF{G}{\calf^{\prime}} ) .
\]
If $\calf^{\prime}= \calall$, then $\EGF{G}{\calf^{\prime}} = \pt$ and $A_{\calf \to \calf^{\prime}}$ specializes
to the assembly map $A_{\calf}$ we discussed in the previous section.
If we  now gradually increase the family, we obtain a 
factorization of the classical assembly $A=A_{\caltr \to \calall}$
into several relative assembly maps. 
We obtain for example from the inclusions 
\[
\caltr \subset \calfcyc \subset \calfin \subset \calvcyc \subset \calall
\]
for every $G$-homology theory $\calh_n^G ( - )$ the following commutative diagram.
\begin{eqnarray} \label{dia: varying families}
\xymatrix{
\calh_n^G ( EG ) \ar[rr]^-{A} \ar[d]   & & \calh_n^G ( \pt ) \\
\calh_n^G ( \EGF{G}{\calfcyc} ) \ar[r]  & 
\calh_n^G ( \EGF{G}{\calfin} )  \ar[r] \ar[ur]^-{A_{\calfin}}&  
\calh_n^G ( \EGF{G}{\calvcyc} ) \ar[u]_-{A_{\calvcyc}}.
         }
\end{eqnarray}
Here $A$ is the ``classical'' assembly map and $A_{\calfin}$ and $A_{\calvcyc}$ are the assembly maps
that for specific $G$-homology theories appear in the Baum-Connes and Farrell-Jones Conjectures.

Such a factorization is extremely useful because one can study the relative assembly map $A_{\calf \to \calf^{\prime}}$ 
in terms of absolute assembly maps corresponding to groups in the bigger family. For example the relative assembly map 
\[
A_{\calfin \to \calvcyc} \colon \calh^G_n ( \EGF{G}{\calfin} ) \to \calh^G_n ( \EGF{G}{\calvcyc} )
\]
is an isomorphism
if for all virtually cyclic subgroups $V$ of $G$ the assembly map
\[
A_{\calfin} = A_{\calfin \to \calall} \colon \calh^V_n ( \EGF{V}{\calfin} ) \to \calh^V_n ( \pt )
\]
is an isomorphism. Of course here we need to assume that the $G$-homology theory $\calh^G_*$ 
and the $V$-homology theory $\calh^V_*$ are somehow related. 
In fact all the $G$-homology theories $\calh^G_*$ we care about 
are defined simultaneously for all groups $G$ and for varying $G$ these $G$-homology theories are related via
a so called ``induction structure''. Induction structures will be
discussed in detail in 
Section \ref{sec: The Definition of an Equivariant Homology Theory}.

For a family $\calf$ of subgroups of $G$ and a subgroup $H \subset G$ we define
a family of subgroups of $H$
\[
\calf \cap H = \{ K \cap H \; | \; K \in \calf \}.
\]
The general statement about relative assembly maps reads now as
follows.
\begin{theorem}[Transitivity Principle] \label{the: transitivity}
\indextheorem{Transitivity Principle}
Let $\calh_{\ast}^? ( - )$ be an equivariant homology theory in the sense
of Section \ref{sec: The Definition of an Equivariant Homology Theory}.
Suppose $\calf \subset \calf^{\prime}$ are two families of subgroups of $G$.
Suppose that $K \cap H \in \calf$ for each $K \in \calf$ and $H
\in \calf^{\prime}$ (this is automatic if $\calf$ is closed under taking subgroups).
Let $N$ be an integer. If for every $H \in \calf^{\prime}$ and every $n \leq N$ the assembly map 
\[
A_{\calf \cap H \to \calall} \colon \calh_n^H ( \EGF{H}{\calf \cap H} ) \to \calh_n^H ( \pt )
\]
is an isomorphism, then for every $n \leq N$ the relative assembly map
\[
A_{\calf \to \calf^{\prime}} \colon \calh_n^G ( \EGF{G}{\calf} ) \to \calh_n^G ( \EGF{G}{\calf^{\prime}} )
\]
is an isomorphism.
\end{theorem}
\begin{proof} 
If we equip $\EGF{G}{\calf} \times \EGF{G}{\calf^{\prime}}$ with the
diagonal $G$-action, it is a model for $\EGF{G}{\calf}$. 
Now apply Lemma~\ref{lem: EGF{G}{calf} times Z to Z}
in the special case $Z = \EGF{G}{\calf^{\prime}}$. 
\end{proof}

This principle can be used in many ways.
For example we will derive from it that the general versions 
of the Baum-Connes and Farrell-Jones Conjectures specialize to the conjectures
we discussed in Chapter~\ref{chap: torsion free} in the case where the group is torsion free.
If we are willing to make compromises, e.g.\ to invert $2$, to 
rationalize the theories or to restrict ourselves to small dimensions
or special classes of groups,
then it is often possible to get away with a smaller family, 
i.e.\ to conclude from the Baum-Connes or Farrell-Jones Conjectures
that an assembly map with respect to a family smaller than the family of finite or virtually cyclic subgroups 
is an isomorphism. The left hand side becomes more computable and this 
leads to new corollaries of the Baum-Connes and Farrell-Jones Conjectures.


\subsubsection{The General Versions Specialize to the Torsion Free Versions}
\label{subsec: The general versions specialize to the torsion free versions}

If $G$ is a torsion free group, then the family $\calfin$ obviously coincides with the trivial family $\caltr$.
Since a nontrivial torsion free virtually cyclic group is infinite cyclic we also know that
the family $\calvcyc$ reduces to the family of all cyclic subgroups, denoted $\calcyc$.

\begin{propositionnew} \label{pro: relative assembly torsion free}
Let $G$ be a torsion free group. 
\begin{enumerate}
\item \label{rel ass tfree i}
If $R$ is a regular ring, then the relative assembly map
\[
A_{\caltr \to \calcyc} \colon H_n^G ( \EGF{G}{\caltr} ; \bfK_R ) 
\to H_n^G ( \EGF{G}{\calcyc} ; \bfK_R ) 
\]
is an isomorphism.
\item \label{rel ass tfree ii}
For every ring $R$ the relative assembly map
\[
A_{\caltr \to \calcyc} \colon H_n^G ( \EGF{G}{\caltr} ; \bfL_R^{\langle - \infty \rangle} ) 
\to H_n^G ( \EGF{G}{\calcyc} ; \bfL_R^{\langle - \infty \rangle} ) 
\]
is an isomorphism.
\end{enumerate}
\end{propositionnew}
\begin{proof}
Because of the Transitivity Principle~\ref{the: transitivity} 
it suffices in both cases to prove that the classical assembly 
map $A=A_{\caltr \to \calall}$ is an isomorphism in the case 
where $G$ is an infinite cyclic group. For regular rings in the $K$-theory case
and with the $- \infty$-decoration in the $L$-theory case this is true as we discussed in 
Remark~\ref{rem: Bass-Heller-Swan decomposition} respectively Remark~\ref{rem: Shaneson splitting}.
\end{proof}

As an immediate consequence we obtain.
\begin{corollary} \label{cor: general implies torsion free}
\begin{enumerate}
\item
For a torsion free group the Baum-Connes Conjecture~\ref{con: Baum-Connes Conjecture}
is equivalent to its ``torsion free version'' Conjecture~\ref{con: BC torsion free}.
\item
For a torsion free group the Farrell-Jones Conjecture~\ref{con: Farrell-Jones Conjecture} for algebraic $K$-
is equivalent to the ``torsion free version'' Conjecture~\ref{con: FJK torsion free all}, provided
$R$ is regular.
\item
For a torsion free group the Farrell-Jones Conjecture~\ref{con: Farrell-Jones Conjecture} for algebraic 
$L$-theory is equivalent to the ``torsion free version'' Conjecture~\ref{con: FJL torsion free}.
\end{enumerate}
\end{corollary}


\subsubsection{The Baum-Connes Conjecture and the Family $\calvcyc$}
\label{subsec: The Baum-Connes Conjecture and the family calvcyc}

Replacing the family $\calfin$ of finite subgroups by the family $\calvcyc$ of virtually cyclic subgroups would not
make any difference in the Baum-Connes Conjecture~\ref{con: Baum-Connes Conjecture}.
The Transitivity Principle \ref{the: transitivity}
and the fact that the Baum-Connes Conjecture
\ref{con: Baum-Connes Conjecture} is known for virtually cyclic groups implies the following.
\begin{propositionnew} \label{pro: vcyc and BC}
For every group $G$ and every $n \in \IZ$ the relative assembly map for topological $K$-theory
\[
A_{\calfin \to \calvcyc} \colon H_n^G ( \EGF{G}{\calfin} ; \bfK^{\topo} ) \to H_n^G ( \EGF{G}{\calvcyc} ; \bfK^{\topo} )
\]
is an isomorphism.
\end{propositionnew}


\subsubsection{The Baum-Connes Conjecture and the Family $\calfcyc$}
\label{subsec: The Baum-Connes Conjecture and the family calfcyc}

The following result is proven in \cite{Matthey-Mislin(2003)}.
\begin{propositionnew} \label{pro:Matthey-Mislin}
For every group $G$ and every $n \in \IZ$ the relative assembly map for topological $K$-theory
\[
A_{\calfcyc \to \calfin} \colon H_n^G ( \EGF{G}{\calfcyc} ; \bfK^{\topo} ) \to H_n^G ( \EGF{G}{\calfin} ; \bfK^{\topo} )
\]
is an isomorphism.
\end{propositionnew}
In particular the Baum-Connes Conjecture predicts that the $\calfcyc$-assembly map
\[
A_{\calfcyc} \colon H_n^G ( \EGF{G}{\calfcyc} ; \bfK^{\topo} )  \to
K_n ( C_r^{\ast} ( G ) )
\]
is always an isomorphism.


\subsubsection{Algebraic $K$-Theory for Special Coefficient Rings} 
\label{subsec: K for special coefficients}

In the algebraic $K$-theory case we can reduce to the family of finite subgroups if we assume special coefficient
rings.

\begin{propositionnew}  \label{pro: K for special coefficients}
Suppose $R$ is a regular ring in which the orders of all finite subgroups of $G$ are invertible.
Then for every $n \in \IZ$ the relative 
assembly map for algebraic $K$-theory
\[
A_{\calfin \to \calvcyc} \colon H_n^G ( \EGF{G}{\calfin} ; \bfK_R )
\to H_n^G ( \EGF{G}{\calvcyc} ; \bfK_R)
\]
is an isomorphism. 
In particular if $R$ is a regular ring which is a $\IQ$-algebra (for example a field of characteristic $0$)
the above applies to all groups $G$.
 \end{propositionnew}
\begin{proof} We first show that $RH$ is regular for a finite group $H$. 
Since $R$ is Noetherian and $H$ is finite, $RH$ is Noetherian.
It remains to show that every $RH$-module $M$ has a finite dimensional projective resolution.
By assumption $M$ considered as an $R$-module has a finite dimensional projective resolution.
If one applies $RH \otimes_R -$ this yields a finite dimensional $RH$-resolution  of $RH \otimes_R \res M$.
Since $|H|$ is invertible, the $RH$-module $M$ is a direct summand of 
$RH \otimes_R \res M$ and hence has a finite dimensional projective 
resolution.

Because of the Transitivity Principle~\ref{the: transitivity}
we need to prove that the $\calfin$-assembly map $A_{\calfin}$ 
is an isomorphism for virtually cyclic groups $V$.
Because of Lemma~\ref{lem: virtually cyclic} we can assume that 
either $V\cong H \semidirect \IZ$ or $V \cong K_1 \ast_{H} K_2$
with finite groups $H$, $K_1$ and $K_2$. From \cite{Waldhausen(1978a)} 
we obtain in both cases long exact sequences involving 
the algebraic $K$-theory of the constituents, 
the algebraic $K$-theory of $V$ and also additional Nil-terms. 
However, in both cases the Nil-terms vanish if $RH$ is a 
regular ring (compare Theorem~4 on page~138 and the Remark on page~216
 in \cite{Waldhausen(1978a)}). Thus we get long exact sequences
$$\ldots \to K_n(RH) \to K_n(RH) \to K_n(RV) \to K_{n-1}(RH) \to K_{n-1}(RH) \to \ldots$$
and
\begin{multline*}
\ldots \to K_n(RH) \to K_n(RK_1) \oplus K_n(RK_2) \to K_n(RV) \\
\to  K_{n-1}(RH) \to K_{n-1}(RK_1) \oplus K_{n-1}(RK_2) \to \ldots
\end{multline*} 
One obtains analogous exact sequences for the sources
of the various assembly maps from the fact that the sources are equivariant homology
theories and one can find specific models for $\EGF{V}{\calfin}$. 
These sequences are compatible with the assembly maps. The assembly maps
for the finite groups $H$, $K_1$ and $K_2$ are bijective.  Now a Five-Lemma argument shows
that also the one for $V$ is bijective.
\end{proof}

In particular for regular coefficient rings $R$ which are $\IQ$-algebras 
the $K$-theoretic Farrell-Jones Conjecture 
specializes to the conjecture that the assembly map
\[
A_{\calfin} \colon H_n^G ( \EGF{G}{\calfin} ; \bfK_R) \to H_n^G ( \pt ; \bfK_R) \cong K_n ( RG )
\]
is an isomorphism.

In the proof above we used the following important fact about virtually cyclic groups.
\begin{lemma} \label{lem: virtually cyclic}
If $G$ is an infinite  virtually cyclic group then  we have the following dichotomy.
\begin{enumerate}
\item[ (I) ]
Either $G$ admits a surjection with finite kernel onto the infinite cyclic group $\IZ$, or
\item[ (II) ]
$G$ admits a surjection with finite kernel onto the infinite dihedral group $\IZ/ 2 \ast \IZ / 2$.
\end{enumerate}
\end{lemma}
\begin{proof} This is not difficult and proven as Lemma~2.5 in \cite{Farrell-Jones(1995)}. \end{proof}


\subsubsection{Splitting off Nil-Terms and Rationalized Algebraic $K$-Theory}  
\label{subsec: split off nil}

Recall that the Nil-terms, which prohibit the classical assembly map
from being an isomorphism, are direct summands of the $K$-theory of the 
infinite cyclic group (see Remark~\ref{rem: Bass-Heller-Swan decomposition}).
Something similar remains true in general \cite{Bartels(2003h)}.
\begin{propositionnew} \label{pro: splitting of calfin in K- and L-theory}
\begin{enumerate}
\item
For every group $G$, every ring $R$ and every $n \in \IZ$ the relative assembly map 
\[
A_{\calfin \to \calvcyc} \colon H_n^G ( \EGF{G}{\calfin} ; \bfK_R ) \to
H_n^G ( \EGF{G}{\calvcyc} ; \bfK_R ) 
\]
is split-injective. 
\item
Suppose $R$ is such that 
$K_{-i} ( RV )= 0$ for all virtually cyclic subgroups $V$ of $G$ and for sufficiently large $i$
(for example $R= \IZ$ will do, compare Proposition~\ref{pro: vanishing lower vcyc}). Then the relative assembly map
\[
A_{\calfin \to \calvcyc} \colon H_n^G ( \EGF{G}{\calfin} ; \bfL_R^{\langle - \infty \rangle} ) 
\to
H_n^G ( \EGF{G}{\calvcyc} ; \bfL_R^{\langle - \infty \rangle} ) 
\]
is split-injective.
\end{enumerate}
\end{propositionnew}

Combined with the Farrell-Jones Conjectures we obtain that the homology group 
$H_n^G ( \EGF{G}{\calfin} ; \bfK_R )$ is a direct summand
in $K_n ( RG )$. It is much better understood 
(compare Chapter~\ref{chap: Computations})
than the remaining summand which is isomorphic to 
$H_n^G ( \EGF{G}{\calvcyc}  , \EGF{G}{\calfin} ; \bfK_R )$. This remaining summand is the one which plays the role 
of the Nil-terms for a general group. It is known that for $R=\IZ$ 
the negative  dimensional Nil-groups which are responsible
for virtually cyclic groups vanish 
\cite{Farrell-Jones(1995)}.
They vanish rationally, in dimension $0$ by
\cite{Connolly-Prassidis(2002a)} and in higher dimensions by \cite{Kuku-Tang(2003)}.
For more information see also
\cite{Connolly-Prassidis(2002b)}.
Analogously to the proof of Proposition~\ref{pro: K for special coefficients} 
we obtain the following proposition.
\begin{propositionnew} \label{pro: rational vanishing of lower nil}
We have 
\begin{eqnarray*}
H_n^G ( \EGF{G}{\calvcyc} , \EGF{G}{\calfin} ; \bfK_{\IZ} ) = 0 \quad \mbox{ for } n < 0 \quad \mbox{ and }\\
H_n^G ( \EGF{G}{\calvcyc} , \EGF{G}{\calfin} ; \bfK_{\IZ} ) \otimes_{\IZ} \IQ = 0 \quad \mbox{ for all } n \in \IZ .
\end{eqnarray*}
\end{propositionnew}
In particular the Farrell-Jones Conjecture for the  algebraic $K$-theory of the integral group ring predicts that
the map
\[
A_{\calfin} \colon H_n^G ( \EGF{G}{\calfin} ; \bfK_{\IZ} )  \otimes_{\IZ} \IQ  \to
K_n ( \IZ G  ) \otimes_{\IZ} \IQ
\]
is always an isomorphism.


\subsubsection{Inverting $2$ in $L$-Theory} 
\label{subsec: L invert 2}

\begin{propositionnew} \label{pro: L invert 2}
For every group $G$, every ring $R$ with involution, 
every decoration $j$ and all $n \in \IZ$ the relative assembly map
\[
A_{\calfin \to \calvcyc} \colon H_n^G ( \EGF{G}{\calfin} ; \bfL_R^{\langle j \rangle} ) [\frac{1}{2}]
\to
H_n^G ( \EGF{G}{\calvcyc} ; \bfL_R^{\langle j \rangle} ) [\frac{1}{2}]
\]
is an isomorphism.
\end{propositionnew}
\begin{proof} 
According to the Transitivity Principle it suffices to prove the claim for a virtually
cyclic group. Now argue analogously  to the proof of
Proposition \ref{pro: K for special coefficients} using the exact sequences
in \cite{Cappell(1974b)} 
and the fact that the UNil-terms appearing there vanish after
inverting two \cite{Cappell(1974b)}. 
Also recall from Remark~\ref{rem: Rothenberg sequence} 
that after inverting $2$ there are no differences between the  decorations.
\end{proof}

In particular the $L$-theoretic Farrell-Jones Conjecture implies  that for every decoration $j$ the assembly map
\[
A_{\calfin} \colon H_n^G ( \EGF{G}{\calfin} ; \bfL_R^{\langle j \rangle} ) [\frac{1}{2}]
\to
L_n^{\langle j \rangle} ( RG ) [\frac{1}{2}]
\]
is an isomorphism.


\subsubsection{$L$-theory and Virtually Cyclic Subgroups of the First Kind} 
\label{subsec: L second kind}

Recall that a group is virtually cyclic of the first kind if it admits a surjection with finite kernel onto the
infinite cyclic group. The family of these groups is denoted $\calvcyc_{I}$.
\begin{propositionnew}
\label{pro: first kind of vcyc and L}
For all groups $G$, all rings $R$ and all $n \in \IZ$ the relative assembly map
\[
A_{\calfin \to \calvcyc_I} \colon
H_n^G ( \EGF{G}{\calfin} ; \bfL_R^{\langle - \infty \rangle} ) 
\to H_n^G ( \EGF{G}{\calvcyc_{I}} ; \bfL_R^{\langle - \infty \rangle} ) 
\]
is an isomorphism.
\end{propositionnew}
\begin{proof}
The point is that there are no  UNil-terms for infinite virtually cyclic
groups of the first kind.  This follows essentially from
\cite{Ranicki(1973b)} and \cite{Ranicki(1973c)} as carried out in
\cite{Lueck(2003a)}.
\end{proof}


\subsubsection{Rationally $\calfin$ Reduces to $\calfcyc$} 
\label{subsec: fcyc to fin}

We will see later 
(compare Theorem~\ref{the: Rational Computation of Topological K-Theory for Infinite Groups}, 
\ref{the: Rational Computation of Algebraic K-Theory for Infinite Groups} and
\ref{the: Rational Computation of Algebraic L-Theory for Infinite Groups})
that in all three cases, topological $K$-theory, algebraic $K$-theory
and $L$-theory, the rationalized left hand side of the $\calfin$-assembly map can be computed very explicitly using 
the equivariant Chern-Character.
As a by-product these  computations yield that after rationalizing 
the family $\calfin$ can be reduced to the family $\calfcyc$ of finite cyclic groups and that the rationalized
relative assembly maps $A_{\caltr \to \calfcyc}$ are injective.

\begin{propositionnew}
\label{pro: Rational fcyc to fin}
For every ring $R$, every group $G$ and all $n \in \IZ$ the relative assembly maps 
\begin{eqnarray*}
A_{\calfcyc \to \calfin} \colon H_n^G ( \EGF{G}{\calfcyc} ; \bfK_R ) \otimes_{\IZ} \IQ
& \to &  H_n^G ( \EGF{G}{\calfin} ; \bfK_R ) \otimes_{\IZ} \IQ \\
A_{\calfcyc \to \calfin} \colon H_n^G ( \EGF{G}{\calfcyc} ; \bfL_R^{\langle - \infty \rangle}) \otimes_{\IZ} \IQ
& \to &  H_n^G ( \EGF{G}{\calfin} ; \bfL_R^{\langle - \infty \rangle} ) \otimes_{\IZ} \IQ \\
A_{\calfcyc \to \calfin} \colon H_n^G ( \EGF{G}{\calfcyc} ; \bfK^{\topo} ) \otimes_{\IZ} \IQ
& \to &  H_n^G ( \EGF{G}{\calfin} ; \bfK^{\topo} )  \otimes_{\IZ} \IQ
\end{eqnarray*}
are isomorphisms and the corresponding relative assembly maps $A_{\caltr \to \calfcyc}$ are all
rationally injective.
\end{propositionnew}
Recall that the statement for topological $K$-theory is even known integrally, compare Proposition~\ref{pro:Matthey-Mislin}.
Combining the above with Proposition~\ref{pro: rational vanishing of lower nil} and Proposition~\ref{pro: L invert 2} we
see that the Farrell-Jones Conjecture predicts in particular that the $\calfcyc$-assembly maps
\begin{eqnarray*}
A_{\calfcyc} \colon H_n^G ( \EGF{G}{\calfcyc} ; \bfL_R^{\langle - \infty \rangle} ) \otimes_{\IZ} \IQ & \to &
L_n^{\langle -\infty \rangle} ( R G  ) \otimes_{\IZ} \IQ \\
A_{\calfcyc} \colon H_n^G ( \EGF{G}{\calfcyc} ; \bfK_{\IZ} ) \otimes_{\IZ} \IQ & \to &
K_n ( \IZ G ) \otimes_{\IZ} \IQ \\
\end{eqnarray*}
are always isomorphisms.


\typeout{----------------------------- More Applications--------------------------}
\section{More Applications} \label{chap: more applications}


\subsection{Applications VI}
\label{sec: Applications VI}

\subsubsection{Low Dimensional Algebraic $K$-Theory}
\label{subsec: Low dimensional algebraic K-theory}
As opposed to topological $K$-theory and $L$-theory, which are periodic,
the algebraic $K$-theory groups of coefficient rings such as $\IZ$, $\IQ$ or $\IC$ are known to be bounded below. 
Using the spectral sequences for the left hand side of an assembly map that will be discussed
in Subsection~\ref{subsec: Equivariant Atiyah-Hirzebruch Spectral Sequence},
this leads to vanishing results in negative dimensions and a concrete 
description of the groups in the first non-vanishing dimension.

The following conjecture is a consequence of the $K$-theoretic Farrell-Jones Conjecture in the case $R = \IZ$.
Note that by the results discussed in  Subsection~\ref{subsec: split off nil} we know that in negative dimensions 
we can reduce to the family of finite subgroups.
Explanations about the colimit that appears follow below.

\begin{conjecture}[The Farrell-Jones Conjecture for $K_n( \IZ G )$ for $n \leq -1$] \label{con: colim K-1}
\index{Conjecture!Farrell-Jones Conjecture!for $K_n( \IZ G )$ for $n \leq -1$}
For every group $G$ we have
\[
K_{-n} ( \IZ G ) = 0 \quad \mbox{ for } n \geq 2,
\]
and the map
\[
\xymatrix{
\colim_{H \in \SubGF{G}{\calfin}} K_{-1}(\IZ H) \ar[r]^-{\cong} & K_{-1}(\IZ G)
         }
\]
is an isomorphism. 
\end{conjecture}

We can consider a family $\calf$ of subgroups of $G$ as a category
$\SubGF{G}{\calf}$%
\indexnotation{sub_FG}
as follows. The objects
are the subgroups $H$ with $H \in \calf$. For $H,K \in \calf$ let $\conhom_G(H,K)$ be
the set of all group homomorphisms $f \colon  H \to K$, for which there exists a group element $g \in G$
such that $f$ is given by conjugation with $g$. The group of inner
automorphism $\inn(K)$%
\indexnotation{inn(K)}
consists of those automorphisms $K \to K$,
which are given by conjugation with an element $k \in K$. It acts on
$\conhom(H,K)$ from the left by composition. Define the set of
morphisms in $\SubGF{G}{\calf}$ from $H$ to $K$ to be
$\inn(K)\backslash\conhom(H,K)$. Composition of group homomorphisms defines the
composition of morphisms in $\SubGF{G}{\calf}$. We mention that
$\SubGF{G}{\calf}$ is a quotient category of the orbit category
$\OrGF{G}{\calf}$ which we will introduce in 
Section \ref{sec: Homotopy-Theoretic Versions of the Conjectures}.
Note that there is a morphism from $H$ to $K$ only if $H$ is conjugate
to a subgroup of $K$.
Clearly $K_n ( R ( - ) )$ yields a functor from $\SubGF{G}{\calf}$ to abelian
groups since inner automorphisms on a group $G$ induce the identity on
$K_n(RG)$. Using the inclusions into $G$, one obtains a map
\[
\colim_{H \in \SubGF{G}{\calf}} K_n ( R H ) \to K_n ( RG ).
\]
The colimit can be interpreted as the $0$-th Bredon homology group
\[
H^G_0 ( \EGF{G}{\calf} ; K_n ( R ( ? ) ))
\]
(compare 
Example~\ref{Bredon homology as G-homology theory}) and 
the map is the edge homomorphism in the equivariant Atiyah-Hirzebruch
spectral sequence discussed in Subsection~\ref{subsec: Equivariant Atiyah-Hirzebruch Spectral Sequence}.
In Conjecture~\ref{con: colim K-1} we consider the first non-vanishing entry in the lower left hand corner of the $E_2$-term 
because of the following vanishing result 
\cite[Theorem~2.1]{Farrell-Jones(1995)}
which generalizes vanishing results for finite groups from \cite{Carter(1980b)}.
\begin{propositionnew} \label{pro: vanishing lower vcyc}
If $V$ is a virtually cyclic group, then $K_{-n} ( \IZ V ) = 0$ for $n \geq 2$.
\end{propositionnew}

If our coefficient ring $R$ is a regular ring in which the orders of all finite subgroups of $G$ are invertible, then
we know already from Subsection~\ref{subsec: K for special coefficients}
that we can reduce to the family of finite subgroups. In the proof of 
Proposition~\ref{pro: K for special coefficients}
we have seen that then $RH$ is again regular if $H \subset G$ is finite. 
Since negative $K$-groups vanish for regular rings
\cite[5.3.30 on page~295]{Rosenberg(1994)}, 
the following is implied by the Farrell-Jones Conjecture~\ref{con: Farrell-Jones Conjecture}.

\begin{conjecture}[Farrell-Jones Conjecture for $K_0( \IQ G )$] \label{con: colim over finite for special coefficients}
\index{Conjecture!Farrell-Jones Conjecture!for $K_0( \IQ G )$}
Suppose $R$ is a regular ring in which the orders of all finite subgroups of $G$ are invertible
(for example a field of characteristic $0$), then
\[
K_{-n} ( RG ) = 0 \quad \mbox{ for } n  \geq 1
\]
and the map
\[
\xymatrix{
\colim_{H \in \SubGF{G}{\calfin}} K_0 (R H) \ar[r]^-{\cong} & K_0 (R G)
         }
\]
is an isomorphism.
\end{conjecture}

The conjecture above holds if $G$ is virtually poly-cyclic. 
Surjectivity is proven in \cite{Moody(1989)} (see also \cite{Cliff-Weiss(1988)}
and Chapter~8 in \cite{Passman(1989)}), injectivity in \cite{Rosenthal(2002)}.
We will show in 
Lemma~\ref{pro: Farrell-Jones and Strong Bass}~\ref{pro: Farrell-Jones and Strong Bass: FJC for IC implies SBC for IC}
that the map appearing in the 
conjecture is always rationally injective for $R = \IC$. 

The conjectures above describe the first non-vanishing term in the equivariant Atiyah-Hirzebruch
spectral sequence. Already the next step is much harder to analyze in general because there are potentially non-vanishing 
differentials. We know however that after rationalizing the 
equivariant Atiyah-Hirzebruch spectral sequence for the left hand side of 
the $\calfin$-assembly map
collapses. As a consequence we obtain that the following conjecture follows from the $K$-theoretic 
Farrell-Jones Conjecture~\ref{con: Farrell-Jones Conjecture}.
\begin{conjecture} \label{con: colim injects rationally}
For every group $G$, every ring $R$ and every $n \in \IZ$ the map
\[
\colim_{H \in \SubGF{G}{\calfin}}  K_n( R H ) \otimes_{\IZ} \IQ  \to K_n ( R G ) \otimes_{\IZ} \IQ
\]
is injective. 
\end{conjecture}

Note that for $K_0 ( \IZ G ) \otimes_{\IZ} \IQ$ the conjecture above is always true but not very interesting, because
for a finite group $H$ it is known that $\widetilde{K}_0 ( \IZ H ) \otimes_{\IZ} \IQ = 0$, 
compare~\cite[Proposition 9.1]{Swan(1960a)}, and hence 
the left hand side reduces to $K_0 ( \IZ ) \otimes_{\IZ} \IQ$. 
However, the full answer for $K_0 ( \IZ G )$ should involve the negative $K$-groups, compare
Example~\ref{exa: Formula for K_0(IZ G)}.

Analogously to Conjecture~\ref{con: colim injects rationally} the following 
can be derived from the $K$-theoretic Farrell-Jones Conjecture~\ref{con: Farrell-Jones Conjecture}, 
compare~\cite{Lueck-Reich-Rognes-Varisco(2003)}. 

\begin{conjecture} \label{con: colim over rationalized Whitehead}
The map
\[
\colim_{H \in \SubGF{G}{\calfin}}  \Wh ( H ) \otimes_{\IZ} \IQ  \to \Wh ( G ) \otimes_{\IZ} \IQ
\]
is always injective.
\end{conjecture}

In general one does not expect this map to be an isomorphism. There should be additional contributions coming
from negative $K$-groups.
Conjecture~\ref{con: colim over rationalized Whitehead}  
is true for groups satisfying a mild homological finiteness condition, see Theorem~\ref{the: LRRM for middle and lower K-theory}.


\begin{remarknew}[The Conjectures as Generalized Induction Theorems] 
\label{rem: isomorphisms conjectures as induction theorem}  
The above discussion shows that one
may think of the Farrell-Jones Conjectures \ref{con: Farrell-Jones Conjecture} and
the Baum-Connes Conjecture \ref{con: Baum-Connes Conjecture} as ``generalized induction
theorems''. 
The prototype of an induction theorem is Artin's theorem
about the complex representation ring $R_{\IC}(G)$ of a finite group $G$. Let us recall
Artin's theorem.

For finite groups $H$ the  complex representation ring 
$R_{\IC} ( H )$ coincides with $K_0 ( \IC H )$. 
Artin's Theorem%
\indextheorem{Artin's}
\cite[Theorem~17 in 9.2 on page 70]{Serre(1977)}
implies that the obvious induction homomorphism
\[
\colim_{H \in \SubGF{G}{\calcyc}} R_{\IC}(H) \otimes_{\IZ} \IQ  \xrightarrow{\cong} R_{\IC}(G) \otimes_{\IZ} \IQ
\]
is an isomorphism. Note that this is a very special case of
Theorem~\ref{the: Rational Computation of Topological K-Theory for Infinite Groups} or 
\ref{the: Rational Computation of Algebraic K-Theory for Infinite Groups}, 
compare Remark~\ref{rem: already interesting for finite}.

Artin's theorem says that rationally one can compute
$R_{\IC}(G)$ if one knows all the values
$R_{\IC}(C)$ (including all maps coming from induction with group
homomorphisms induced by  conjugation with elements in $G$) for all
cyclic subgroups $C \subseteq G$. 
The idea behind the Farrell-Jones Conjectures \ref{con: Farrell-Jones Conjecture}
and the Baum-Connes Conjecture \ref{con: Baum-Connes Conjecture} is 
analogous. We want to compute the functors $K_n(RG)$, $L_n(RG)$ and $K_n(C^*_r(G))$ 
from their values (including their functorial properties under induction) 
on elements of the family $\calfin$ or $\calvcyc$. 

The situation in the Farrell Jones and Baum-Connes Conjectures is more complicated than 
in Artin's Theorem, since we have already seen in Remarks 
\ref{rem: Bass-Heller-Swan decomposition}, \ref{rem: Shaneson splitting} and
\ref{Computing K(C^*_r(G times Z)} that a computation of $K_n(RG)$, 
$L^{\langle -\infty  \rangle}_n(RG)$ and $K_n(C^*_r(G))$ does involve also the values 
$K_p(RH)$, $L^{\langle -\infty  \rangle}_p(RH)$ and $K_p(C^*_r(H))$ for $p \le n$.  
A degree mixing occurs.
\end{remarknew}


\subsubsection{$G$-Theory}
\label{subsec: G-Theory}

Instead of considering finitely generated projective modules one may apply the standard $K$-theory 
machinery to the category of finitely generated modules.
This leads to the definition of the groups
$G_n(R)$%
\indexnotation{G_n(R)}
for $n \ge 0$. For instance $G_0(R)$ is the abelian group whose generators are isomorphism
classes $[M]$ of finitely generated $R$-modules and whose relations are given by
$[M_0] - [M_1] + [M_2]$ for any exact sequence $0 \to M_0 \to M_1\to M_2 \to 0$ 
of finitely generated  modules. 
One may ask whether versions of the Farrell-Jones Conjectures for $G$-theory instead of
$K$-theory might be true. The answer is negative as the following discussion explains.

For a finite group $H$ the ring $\IC H$ is semisimple.
Hence any finitely generated $\IC H$-module is automatically projective and $K_0 ( \IC H )
= G_0 ( \IC H )$. Recall that a group $G$ is called \emph{virtually poly-cyclic}%
\index{group!virtually poly-cyclic}
if there exists a subgroup of finite index $H \subseteq G$
together with a filtration $\{1\} = H_0 \subseteq H_1 \subseteq H_2 \subseteq \ldots
\subseteq H_r = H$ such that $H_{i-1}$ is normal in $H_i$ and the quotient
$H_i/H_{i-1}$ is cyclic.
More generally for all $n \in \IZ$ the forgetful map 
\[
f \colon K_n(\IC G) \to G_n(\IC G)
\]
is an isomorphism if $G$ is virtually poly-cyclic, since then
$\IC G$ is regular \cite[Theorem~8.2.2 and Theorem~8.2.20]{Rowen(1988b)} and 
the forgetful map $f$ is an isomorphism for regular rings,
compare~\cite[Corollary~53.26 on page~293]{Rosenberg(1994)}.
In particular this applies to virtually cyclic groups and so the left hand side of the Farrell-Jones assembly map
does not see the difference between $K$- and $G$-theory if we work with complex coefficients. We obtain a commutative 
diagram
\begin{eqnarray} \label{dia: compare K to G}
\xymatrix{
\colim_{H \in \SubGF{G}{\calfin}} K_0(\IC H) \ar[d]_-{\cong} \ar[r] &   K_0(\IC G)  \ar[d]^-f  \\
\colim_{H \in \SubGF{G}{\calfin}} G_0(\IC H) \ar[r] &   G_0(\IC G)
         }
\end{eqnarray}
where, as indicated, the left hand vertical map is an isomorphism.
Conjecture~\ref{con: colim over finite for special coefficients}, which is implied by the Farrell-Jones Conjecture,
says that the upper horizontal arrow is an isomorphism. A $G$-theoretic analogue of the Farrell-Jones Conjecture would 
say that the lower horizontal map is an isomorphism.
There are however cases where the upper horizontal arrow is known to be an isomorphism, but the forgetful map 
$f$ on the right is not injective or not surjective:

If $G$ contains a non-abelian free subgroup, then the class $[\IC G] \in G_0(\IC G)$ vanishes
\cite[Theorem~9.66 on page 364]{Lueck(2002)} and hence the map $f \colon K_0 ( \IC G ) \to G_0 ( \IC G )$
has an infinite kernel ($[\IC G]$ generates an infinite cyclic subgroup in $K_0(\IC G)$).
The Farrell-Jones Conjecture for $K_0 ( \IC G )$ is known for non-abelian free groups.

The Farrell-Jones Conjecture  is also known for 
$A = \bigoplus_{n \in \IZ} \IZ/2$ and hence $K_0(\IC A)$ is countable, whereas
$G_0(\IC A)$ is not countable
\cite[Example 10.13 on page 375]{Lueck(2002)}. 
Hence the map $f$ cannot be surjective.

At the time of writing we do not know a counterexample to the statement
that for an amenable group $G$, for which there is an upper bound on the orders 
of its finite subgroups, the forgetful map $f \colon K_0 ( \IC G ) \to G_0 ( \IC G )$ is an isomorphism.
We do not know a counterexample to the statement
that for a group $G$, which is not amenable,
$G_0(\IC G) = \{0\}$. We also do not know whether $G_0(\IC G) = \{0\}$ is true for
$G = \IZ \ast \IZ$.

For more information about $G_0(\IC G)$ we refer for instance to
\cite[Subsection 9.5.3]{Lueck(2002)}.



\subsubsection{Bass Conjectures}
\label{subsec: The Bass Conjecture}

Complex representations of a finite group can be studied using characters.
We now want to define the Hattori-Stallings 
rank of a finitely generated projective $\IC G$-module which 
should  be seen as a generalization of characters to infinite groups. 

Let $\con(G)$
\indexnotation{con(G)}
be the set of conjugacy classes $(g)$%
\indexnotation{(g)}
of elements $g \in G$. Denote by 
$\con(G)_f$
\indexnotation{con(G)_f}
the subset of $\con(G)$ consisting of those conjugacy classes $(g)$ for which each
representative $g$ has finite order. Let $\class_0(G)$%
\indexnotation{class_0(G)}
and $\class_0(G)_f$%
\indexnotation{class_0(G)_f}
be the $\IC$-vector space with the set $\con(G)$ and $\con(G)_f$ as basis. 
This is the same as the $\IC$-vector space of $\IC$-valued functions on $\con(G)$ and $\con(G)_f$ with finite support.
Define the \emph{universal $\IC$-trace}%
\index{trace!universal $\IC$-trace} \indexnotation{tr^u_IC G} as
\begin{eqnarray}  \tr_{\IC G}^u \colon \IC G \to \class_0(G), \quad
\sum_{g \in G} \lambda_g \cdot g ~ \mapsto ~ \sum_{g \in G} \lambda_g \cdot (g). 
\label{universal CG-trace}
\end{eqnarray}
It extends to a function $\tr_{\IC G}^u \colon M_n(\IC G) \to \class_0(G)$ 
on $(n,n)$-matrices over $\IC G$ by taking the sum of the traces of the diagonal entries.
Let $P$ be a finitely generated projective $\IC G$-module. Choose a matrix $A \in M_n(\IC G)$ such that
$A^2 = A$ and the image of the $\IC G$-map $r_A \colon \IC G^n \to \IC G^n$ given by right 
multiplication with $A$ is $\IC G$-isomorphic to $P$. Define the \emph{Hattori-Stallings rank}%
\index{Hattori-Stallings-rank}
of $P$ as 
\begin{eqnarray} \HS_{\IC G}(P) = \tr_{\IC G}^u(A)    \in \class_0(G).
\label{Hattori-Stallings rank} \indexnotation{HS_IC G}
\end{eqnarray}
The Hattori-Stallings rank depends only on the isomorphism class of the $\IC G$-module $P$ and induces a homomorphism
$\HS_{\IC G} \colon K_0(\IC G) \to \class_0(G)$.

\begin{conjecture}[Strong Bass Conjecture for $K_0 ( \IC G )$]
\label{con: Strong Bass Conjecture for K_0(IC G)}
\index{Conjecture!Bass Conjecture!Strong Bass Conjecture for $K_0( \IC G)$}
The $\IC$-vector space spanned by the image of the map 
\[
\HS_{\IC G} \colon K_0(\IC G) \to \class_0(G)
\]
is $\class_0(G)_f$.
\end{conjecture}

This conjecture is implied by the surjectivity of the map
\begin{eqnarray} \label{dia: colim over fin tensor C}
\colim_{H \in \SubGF{G}{\calfin}} K_0 ( \IC H ) \otimes_{\IZ} \IC \to K_0 ( \IC G ) \otimes_{\IZ} \IC,
\end{eqnarray}
(compare Conjecture~\ref{con: colim over finite for special coefficients}) and hence by the $K$-theoretic 
Farrell-Jones Conjecture for $K_0( \IC G )$. 
We will see below 
that the surjectivity of the map~\eqref{dia: colim over fin tensor C} also implies that the map
$K_0 ( \IC G ) \otimes_{\IZ } \IC \to \class_0(G)$,
which is induced by the Hattori-Stallings rank, is injective. Hence 
we do expect that the Hattori-Stallings rank induces an isomorphism
\[
K_0 ( \IC G ) \otimes_{\IZ} \IC \cong \class_0(G)_f.
\]

There are also versions of the Bass conjecture for other coefficients than $\IC$. 
It follows from results of Linnell~\cite[Theorem~4.1 on page 96]{Linnell(1983b)} 
that the following version is implied by the Strong Bass Conjecture for $K_0 ( \IC G )$.

\begin{conjecture}[The Strong Bass Conjecture for $K_0(\IZ G)$] 
\label{con: The Strong Bass Conjecture for K_0(bbZ G)}
\index{Conjecture!Bass Conjecture!Strong Bass Conjecture for $K_0(\IZ G)$}
The image of the composition
\[
K_0(\IZ G) \to  K_0(\IC G)
\xrightarrow{\HS_{\IC G}} \class_0(G)
\]
is contained in the $\IC$-vector space of those functions $f \colon \con(G) \to \IC$ which vanish for
$(g) \in \con(g)$ with $g \not= 1$.
\end{conjecture} 

The conjecture says that for every finitely generated projective $\IZ G$-module $P$ the Hattori-Stallings rank of 
$\IC G \otimes_{\IZ G } P$ looks like the Hattori-Stallings rank of a free $\IC G$-module.
A natural explanation for this behaviour is the following conjecture which clearly implies the Strong Bass Conjecture for
$K_0 ( \IZ G )$.

\begin{conjecture}[Rational $\widetilde{K}_0(\IZ G)$-to-$\widetilde{K}_0(\IQ G)$-Conjecture]
\label{con: Rational K_0(bbZ G)-to-K_0(bbQ G)-Conjecture}
\index{Conjecture!Rational K_0(ZG)-to-K_0(QG)-Conjecture@Rational $K_0(\IZ G)$-to-$K_0(\IQ G)$-Conjecture}
For every group $G$ the map 
\[
\widetilde{K}_0(\IZ G) \otimes_{\IZ} \IQ  \to \widetilde{K}_0(\IQ G) \otimes_{\IZ} \IQ
\]
induced by the change of coefficients is trivial. 
\end{conjecture}

Finally we mention the following variant of the Bass Conjecture.

\begin{conjecture}[The Weak Bass Conjecture] 
\label{con: The Weak Bass Conjecture}
\index{Conjecture!Bass Conjecture!Weak Bass Conjecture}
Let $P$ be a finitely generated projective $\IZ G$-module.
The value 
of the Hattori-Stallings rank of $\IC G \otimes_{\IZ G } P$ at the conjugacy class of the identity element 
is given by
\[
\HS_{\IC G}(\IC G  \otimes_{\IZ G} P)((1)) ~ = ~ \dim_{\IZ}(\IZ \otimes_{\IZ G} P).
\]
Here $\IZ$ is considered as a $\IZ G$-module via the augmentation.
\end{conjecture} 

The $K$-theoretic Farrell-Jones Conjecture implies all four
conjectures above. 
More precisely we have the following proposition.

\begin{propositionnew} \label{pro: Farrell-Jones and Strong Bass}
\begin{enumerate} 
\item \label{pro: Farrell-Jones and Strong Bass: FJC for IC implies SBC for IC}
The map
\[
\colim_{H \in \SubGF{G}{\calfin}} K_0 ( \IC H ) \otimes_{\IZ} \IQ  \to K_0 ( \IC G ) \otimes_{\IZ} \IQ
\]
is always injective. If the map is also surjective 
(compare Conjecture~\ref{con: colim over finite for special coefficients})
then the Hattori-Stallings rank induces an isomorphism
\[
K_0 ( \IC G ) \otimes_{\IZ} \IC \cong \class_0 ( G )_f
\]
and in particular the 
Strong Bass Conjecture for $K_0 ( \IC G )$ and hence also the Strong Bass Conjecture for $K_0 ( \IZ G )$ hold.
\item \label{pro: Farrell-Jones and Strong Bass: FJC for IZ implies SBC for IZ}
The surjectivity of the map
\[
A_{\calvcyc} \colon H_0^G ( \EGF{G}{\calvcyc} ; \bfK_{\IZ} ) \otimes_{\IZ} \IQ \to K_0 ( \IZ G ) \otimes_{\IZ} \IQ
\]
implies the Rational $\widetilde{K}_0 ( \IZ G )$-to-$\widetilde{K}_0 (
\IQ G )$~Conjecture 
and hence also the Strong Bass Conjecture for
$K_0 ( \IZ G )$.
\item \label{pro: Farrell-Jones and Strong Bass: the different Bass Conjectures}
The Strong Bass Conjecture for $K_0 ( \IC G )$ implies the Strong Bass Conjecture for $K_0 ( \IZ G )$.
The Strong Bass Conjecture for $K_0 ( \IZ G )$ implies the Weak Bass Conjecture.
\end{enumerate}
\end{propositionnew}

\begin{proof}
\ref{pro: Farrell-Jones and Strong Bass: FJC for IC implies SBC for IC} follows from
the following commutative diagram, compare~\cite[Lemma~2.15 on page 220]{Lueck(1998b)}.
\[
\xymatrix{
\colim_{H \in \SubGF{G}{\calfin}} K_0 ( \IC H ) \otimes_{\IZ} \IC
\ar[d]_-{\cong} \ar[rr] 
& & K_0 ( \IC G ) \otimes_{\IZ} \IC \ar[d] \\
\colim_{H \in \SubGF{G}{\calfin}} \class_0 ( H ) \ar[r]^-{\cong} & \class_0 ( G )_f \ar[r]^-{i} & \class_0 ( G ).
         }
\]
Here the vertical maps are induced by the Hattori-Stallings rank, the map $i$ is the natural inclusion and in 
particular injective and  we have the indicated isomorphisms.

\ref{pro: Farrell-Jones and Strong Bass: FJC for IZ implies SBC for IZ} According to 
Proposition~\ref{pro: rational vanishing of lower nil} the surjectivity of the map $A_{\calvcyc}$ 
appearing in \ref{pro: Farrell-Jones and Strong Bass: FJC for IZ
  implies SBC for IZ} 
implies the surjectivity of the corresponding 
assembly map $A_{\calfin}$ (rationalized and with $\IZ$ as coefficient ring) for the family of finite subgroups.
The map $A_{\calfin}$ is natural with respect to the change of the coefficient ring from $\IZ$ to $\IQ$.
By Theorem~\ref{the: Rational Computation of Algebraic K-Theory for Infinite Groups} 
we know that for every coefficient ring $R$ there is an isomorphism
from
\[
\bigoplus_{p,q, p+q=0} \bigoplus_{( C ) \in ( \calfcyc )} 
H_p ( B Z_G C ; \IQ ) \otimes_{\IQ [ W_G C ]} 
\Theta_C \cdot K_q ( R C ) \otimes_{\IZ} \IQ
\]
to the $0$-dimensional part of the left hand side of the rationalized 
$\calfin$-assembly map $A_{\calfin}$. The isomorphism is natural with respect to a change of coefficient rings.
To see that the Rational $\widetilde{K}_0 ( \IZ G )$-to-$\widetilde{K}_0 ( \IQ G )$ Conjecture follows, it hence suffices
to show that the summand corresponding to $C= \{ 1 \}$ and $p=q=0$ is the only one where the map induced from
$\IZ \to \IQ$ is possibly non-trivial. But $K_q ( \IQ C )=0$ if $q <0$, 
because $\IQ C$ is semisimple and hence regular, 
and for a finite cyclic group $C \neq \{ 1 \}$ we have by \cite[Lemma~7.4]{Lueck(1998b)}
\[
\Theta_C \cdot K_0 ( \IZ C ) \otimes_{\IZ} \IQ  = 
\coker \left( \bigoplus_{D \subsetneq C} K_0 ( \IZ D ) \otimes_{\IZ} \IQ 
~ \to ~ K_0 ( \IZ C ) \otimes_{\IZ } \IQ \right) = 0,
\]
since by a result of Swan $K_0 ( \IZ ) \otimes_{\IZ} \IQ \to K_0 ( \IZ H )\otimes_{\IZ} \IQ$ is an isomorphism
for a finite group $H$, see~\cite[Proposition 9.1]{Swan(1960a)}. 

\ref{pro: Farrell-Jones and Strong Bass: the different Bass Conjectures}
As already mentioned the first statement follows from \cite[Theorem~4.1 on page 96]{Linnell(1983b)}. The second statement 
follows from the formula
\[
\sum_{(g) \in \con(G)} 
\HS_{\IC G}(\IC \otimes_{\IZ} P)(g) =  \dim_{\IZ}(\IZ \otimes_{\IZ G} P).
\]
\end{proof}

The next result is due to Berrick, Chatterji and Mislin 
\cite[Theorem~5.2]{Berrick-Chatterji-Mislin(2004)}. The Bost
Conjecture is a variant 
of the Baum-Connes Conjecture and is explained
in Subsection~\ref{subsec: The Bost Conjecture}.

\begin{theorem} \label{the: Bost implies Bass}
If the assembly map appearing in the Bost Conjecture~\ref{con: Bost Conjecture} is rationally surjective, then
the Strong Bass Conjecture for $K_0( \IC G )$ (compare~\ref{con: Strong Bass Conjecture for K_0(IC G)}) is true. 
\end{theorem}

We now discuss some further questions and facts that seem to be relevant in the context of the Bass Conjectures.

\begin{remarknew}[Integral $\widetilde{K}_0(\IZ G)$-to-$\widetilde{K}_0(\IQ G)$-Conjecture]
\label{rem: Integral version of K_0(ZG)-toK_0(QG) Conjecture}
\index{Conjecture!Integral K_0(ZG)-to-K_0(QG)-Conjecture@Integral $K_0(\IZ G)$-to-$K_0(\IQ G)$-Conjecture}
We do not know a counterexample to the Integral 
$\widetilde{K}_0 ( \IZ G )$-to-$\widetilde{K}_0( \IQ G )$ Conjecture, i.e.\
to the statement that the map 
\[
\widetilde{K}_0( \IZ G ) \to \widetilde{K}_0 ( \IQ G )
\]
itself is trivial. But we also do not know a proof which shows that the
$K$-theoretic Farrell-Jones Conjecture implies this integral version.
Note that the Integral $\widetilde{K}_0 ( \IZ G )$-to-$\widetilde{K}_0( \IQ G )$ Conjecture
would imply that the 
following diagram commutes.
\[
\begin{CD}
K_0(\IZ G) @>  >> K_0(\IQ G)
\\
@Vp_*VV @AiAA
\\
K_0(\IZ) @> \dim_{\IZ} > \cong > \IZ .
\end{CD}
\]
Here $p_*$ is induced by the projection $G \to \{1\}$ and $i$ sends $1 \in \IZ$ to the class
of $\IQ G$. 
\end{remarknew}

\begin{remarknew}[The passage from $\widetilde{K}_0(\IZ G)$ to $\widetilde{K}_0(\caln(G))$]
 \label{rem: K_0(ZG) to K_0(N(G)}
Let $\caln(G)$ denote the group von Neumann algebra of $G$.
It is known that for every group $G$ the composition
\[
\widetilde{K}_0(\IZ G) \to \widetilde{K}_0(\IQ G) \to \widetilde{K}_0(\IC G) 
\to \widetilde{K}_0(C^*_r(G)) \to \widetilde{K}_0(\caln(G))
\]
is the zero-map (see for instance \cite[Theorem~9.62 on page 362]{Lueck(2002)}). 
Since the group von Neumann 
algebra $\caln(G)$ is not functorial under arbitrary group homomorphisms such as $G \to \{1\}$,
this does \emph{not} imply that the diagram 
\[
\begin{CD}
K_0(\IZ G) @>  >> K_0(\caln(G))
\\
@Vp_*VV @AiAA
\\
K_0(\IZ) @> \dim_{\IZ} > \cong > \IZ
\end{CD}
\]
commutes. However, commutativity would follow from the Weak Bass Conjecture~\ref{con: The Weak Bass Conjecture}.
For a discussion of these questions see \cite{Eckmann(1996b)}.
\end{remarknew}

More information and further references 
about the Bass Conjecture can be found for instance in
\cite{Bass(1976)},
\cite[Section 7]{Berrick-Chatterji-Mislin(2004)},
\cite{Burger-Valette(1998)},
\cite{Eckmann(1986)},
\cite{Eckmann(1996b)},
\cite{Farrell-Linnell(2003b)},
\cite{Linnell(1983b)} 
\cite[Subsection 9.5.2]{Lueck(2002)}, 
and 
\cite[page 66ff]{Mislin-Valette(2003)}.



\subsection{Applications VII}
\label{sec: Applications VII}

\subsubsection{Novikov Conjectures}
\label{subsec: Novikov Conjectures}

In Subsection~\ref{subsec: Novikov torsion free} we discussed the Novikov
Conjectures. 
Recall that one possible reformulation of 
the original Novikov Conjecture 
says that for every group $G$ the rationalized classical assembly map in $L$-theory
\[
A  \colon H_n ( BG ; \bfL^p ( \IZ )  ) \otimes_{\IZ} \IQ  \to L_n^p ( \IZ G ) \otimes_{\IZ} \IQ
\]
is injective. Since $A$ can be identified with $A_{\caltr \to
  \calall}$ and 
we know from Subsection~\ref{subsec: fcyc to fin}
that the relative assembly map 
\[
A_{\caltr \to \calfin} \colon H_n^G ( \EGF{G}{\caltr} ; \bfL_{\IZ}^{p} ) \otimes_{\IZ} \IQ \to 
H_n^G ( \EGF{G}{\calfin} ; \bfL_{\IZ}^{p} ) \otimes_{\IZ} \IQ 
\]
is injective we obtain the following proposition.

\begin{propositionnew} \label{pro: FJC for L implies Novikov} 
The rational injectivity of the assembly map appearing in $L$-theoretic Farrell-Jones Conjecture 
(Conjecture~\ref{con: Farrell-Jones Conjecture}) implies the 
$L$-theoretic Novikov Conjecture (Conjecture~\ref{con: K and L Novikov})
and hence the original Novikov Conjecture~\ref{con: Novikov Conjecture}.
\end{propositionnew}

Similarly the Baum-Connes Conjecture~\ref{con: Baum-Connes Conjecture} 
implies the injectivity of the rationalized classical assembly map
\[
A \colon H_n ( BG ; \bfK^{\topo} ) \otimes_{\IZ} \IQ  \to K_n ( C_r^{\ast} (G)) \otimes_{\IZ} \IQ.
\]
In the next subsection we discuss how one can relate assembly maps 
for topological $K$-theory with $L$-theoretic assembly maps.
The results imply in particular the following proposition.

\begin{propositionnew} \label{BCC implies Novikov}
The rational injectivity of the assembly map appearing in the 
Baum-Connes Conjecture (Conjecture~\ref{con: Baum-Connes Conjecture})
implies the Novikov Conjecture (Conjecture~\ref{con: Novikov Conjecture}).
\end{propositionnew}

Finally we would like to mention that by combining the results 
about the rationalization of $A_{\caltr \to \calfin}$ from
Subsection~\ref{subsec: fcyc to fin} with the splitting result about $A_{\calfin \to \calvcyc}$ from 
Subsection~\ref{subsec: split off nil} we obtain the following result

\begin{propositionnew} \label{pro: FJC for K implies K-Novikov}
The rational injectivity of the assembly map appearing in 
the Farrell-Jones Conjecture for algebraic $K$-theory (Conjecture~\ref{con: Farrell-Jones Conjecture})
implies the $K$-theoretic Novikov Conjecture, i.e.\
the injectivity of 
\[
A \colon H_n ( BG ; \bfK (R) ) \otimes_{\IZ} \IQ  \to K_n ( RG ) \otimes_{\IZ} \IQ.
\]
\end{propositionnew}

\begin{remarknew}[Integral Injectivity Fails]
\label{rem: Integral Injectivity Fails Novikov}
In general the classical assembly maps $A=A_{\caltr}$ themselves, i.e.\ without rationalizing,
are  not injective.
For example one can use the Atiyah-Hirzebruch spectral sequence to see that for $G= \IZ / 5$
\[
H_1(B G; \bfK^{\topo}) \quad \mbox{ and } \quad
H_1(B G ; \bfL^{\langle - \infty \rangle}( \IZ ) )
\]
contain $5$-torsion, whereas for every finite group $G$ the topological $K$-theory of $\IC G$ 
is torsionfree and the torsion in the $L$-theory of $\IZ G$ is always $2$-torsion, 
compare 
Proposition~\ref{pro: Topological K-theory for finite groups}~\ref{pro: Topological K-theory for finite groups: complex} 
and
Proposition~\ref{pro: Algebraic L-theory for finite groups}~\ref{pro: Algebraic L-theory for finite groups: fin. gen 2-tors}.
\end{remarknew}

\subsubsection{Relating Topological $K$-Theory and $L$-Theory}
\label{subsec: Relating the Farrell-Jones and Baum-Connes Conjecture}

For every  real $C^*$-algebra $A$ there is an isomorphism 
$L_n^p(A)[1/2] \xrightarrow{\cong} K_n(A)[1/2]$ \cite{Rosenberg(1995)}.
This can be used to compare 
$L$-theory to topological $K$-theory and leads to the following result.
\begin{propositionnew}\label{pro: relating L and K}
Let $\calf \subseteq \calfin$ be a family of finite subgroups of $G$.
If the topological $K$-theory assembly map
\[
A_{\calf} \colon H_n^G ( \EGF{G}{\calf} ; \bfK^{\topo} ) [\frac{1}{2}] 
 \to    K_n ( C_r^{\ast} ( G ) ) [\frac{1}{2}] 
\]
is injective, then for an arbitrary decoration $j$ also the map
\[
A_{\calf} \colon H_n^G ( \EGF{G}{\calf} ; \bfL_{\IZ}^{\langle j \rangle} ) [\frac{1}{2}] 
 \to    L^{\langle j \rangle}_n ( \IZ G  ) [\frac{1}{2}] 
\]
is injective.
\end{propositionnew}

\begin{proof}
First recall from Remark~\ref{rem: Rothenberg sequence} that after inverting $2$
there is no difference between the different  decorations and we can hence work with the $p$-decoration.
One can construct for any subfamily $\calf \subseteq \calfin$ 
the following commutative diagram \cite[Section~7.5]{Lueck(2002c)}
$$
\begin{CD}
H_n^G(\EGF{G}{\calf};\bfL^p_{\IZ} [1/2]) @> A_{\calf}^1 >> L^p_n(\IZ G)[1/2]
\\
@Vi_1 V\cong V @Vj_1 V \cong V 
\\
H_n^G(\EGF{G}{\calf};\bfL^p_{\IQ} [1/2]) @> A_{\calf}^2 >> L^p_n(\IQ G)[1/2]
\\
@Vi_2 V\cong V @Vj_2VV 
\\
H_n^G(\EGF{G}{\calf};\bfL^p_{\IR} [1/2]) @> A_{\calf}^3 >> L^p_n(\IR G)[1/2]
\\
@Vi_3 V\cong V @Vj_3VV 
\\
H_n^G(\EGF{G}{\calf};\bfL^p_{C_r^*(?;\IR )} [1/2]) @> A_{\calf}^4 >> L^p_n(C_r^*(G;\IR))[1/2]
\\
@Vi_4 V\cong V @Vj_4V \cong V 
\\
H_n^G(\EGF{G}{\calf};\bfK^{\topo}_{\IR} [1/2]) 
@> A_{\calf}^5 >> K_n (C_r^*(G;\IR))[1/2]
\\
@Vi_5 VV @Vj_5VV 
\\
H_n^G(\EGF{G}{\calf};\bfK^{\topo}_{\IC} [1/2]) 
@> A_{\calf}^6 >> K_n (C_r^*(G))[1/2]
\end{CD}
$$
Here 
\begin{eqnarray*}
& \bfL^p_{\IZ}[1/2], \quad   \bfL^p_{\IQ}[1/2],  \quad  \bfL^p_{\IR} [1/2], \quad 
\bfL_{C_r^*(?;\IR)} [1/2], & \\
& \quad  \bfK^{\topo}_{\IR} [1/2] \quad \mbox{ and } \quad \bfK^{\topo}_{\IC} [1/2] &
\end{eqnarray*}
are covariant $\Or(G)$-spectra 
(compare Section~\ref{sec: Spectra over the Orbit Category} and in particular 
Proposition~\ref{pro: Or(G)-spectra yield a G-homology theory})
such that the $n$-th homotopy group of their evaluations  at $G/H$ are given by 
\begin{eqnarray*}
& L^p_n(\IZ H)[1/2], \quad L^p_n(\IQ H)[1/2], \quad  L^p_n(\IR H)[1/2], \quad 
L^p_n(C_r^*(H;\IR))[1/2], & \\ 
& K_n(C_r^*(H;\IR))[1/2] \quad \mbox{ respectively } \quad K_n(C_r^*(H)[1/2]. &
\end{eqnarray*} 
All horizontal
maps are assembly maps induced by the projection 
$\pr \colon \EGF{G}{\calf} \to \pt$. 
The maps $i_k$ and $j_k$ for $k = 1,2,3$ are induced from a change of rings.
The isomorphisms  $i_4$ and $j_4$ come from the general isomorphism
for any real $C^*$-algebra $A$
$$L_n^p(A)[1/2] \xrightarrow{\cong} K_n (A)[1/2]$$
and its spectrum version \cite[Theorem~1.11 on page 350]{Rosenberg(1995)}. 
The maps $i_1$, $j_1$, $i_2$ are isomorphisms by 
\cite[page 376]{Ranicki(1981)} and \cite[Proposition 22.34 on page 252]{Ranicki(1992)}.
The map $i_3$ is bijective since for a finite group
$H$ we have $\IR H = C_r^*(H;\IR)$. The maps $i_5$ and $j_5$ are given by extending 
the scalars from $\IR$ to $\IC$ by induction. 
For every real $C^{\ast}$-algebra $A$ the composition 
\[
K_n (A) [ 1/2 ] \to
K_n (A \otimes_{\IR} \IC ) [ 1/2 ] \to
K_n (M_2 ( A ) ) [1/2 ] 
\]
is an isomorphism and hence $j_5$ is split injective. An $\Or (G)$-spectrum version of this argument
yields that also $i_5$ is split injective.
\end{proof}

\begin{remarknew} \label{rem: FJC for L_*(ZG) equivalent to real BCC after inv. 2} 
One may conjecture that the 
right vertical maps $j_2$ and $j_3$ are isomorphisms and try to prove this directly. Then if we invert $2$ everywhere
the Baum-Connes Conjecture \ref{con: Baum-Connes Conjecture} 
for the real reduced group $C^*$-algebra, would be equivalent to the 
Farrell-Jones Isomorphism Conjecture  for $L_*(\IZ G)[1/2]$.
\end{remarknew}


\subsection{Applications VIII} 
\label{sec: Applications VIII}

\subsubsection{The Modified Trace Conjecture}
\label{subsec: The Modified TraceConjecture}

Denote by $\Lambda^G$ the subring of $\IQ$ which is obtained from $\IZ$ by inverting all orders $| H |$ of finite
subgroups $H$ of $G$, i.e.\ 
\begin{eqnarray}
\Lambda^G   =  \IZ\left[  |H|^{-1} \mid H \subset G, \; |H| < \infty \right].
\label{def of Lambda}
\indexnotation{Lambda^G}
\end{eqnarray}
The following conjecture generalizes Conjecture~\ref{con: Trace Conjecture for Torsion Free Groups} to the case where
the group need no longer be torsionfree. For the standard trace compare~\eqref{standard trace tr_{C^*_r(G)}}.

\begin{conjecture}[Modified Trace Conjecture for a group $G$]
\label{con: Modified Trace Conjecture}
\index{Conjecture!Trace Conjecture!Modified Trace Conjecture}
Let $G$ be a  group. Then the image of the homomorphism induced by the standard trace
\begin{eqnarray}
\tr_{C^*_r(G)} \colon K_0(C^*_r(G)) \to \IR
\end{eqnarray}
is contained in $\Lambda^G$.
\end{conjecture}

The following result is proved in \cite[Theorem~0.3]{Lueck(2002d)}.

\begin{theorem} \label{the: Baum Connes Conjecture implies modified Trace Conjecture}
Let $G$ be a  group. Then the image of the composition
\[
K_0^G(\EGF{G}{\calfin}) \otimes_{\IZ} \Lambda^G 
\xrightarrow{A_{\calfin} \otimes_{\IZ} \id}
K_0(C^*_r(G)) \otimes_{\IZ} \Lambda^G 
 \xrightarrow{\tr_{C^*_r(G)}} \IR
\]
is $\Lambda^G$. Here $A_{\calfin}$ is the map appearing in the 
Baum-Connes Conjecture \ref{con: Baum-Connes Conjecture}.
In particular the Baum-Connes Conjecture \ref{con: Baum-Connes Conjecture} implies the
Modified Trace Conjecture.
\end{theorem}

The original version of the Trace Conjecture due to Baum and Connes \cite[page 21]{Baum-Connes(1982)}
makes the stronger statement that the image of 
$\tr_{C^*_r(G)} \colon K_0(C^*_r(G)) \to \IR$ 
is the additive subgroup of $\IQ$  generated by all numbers
$\frac{1}{|H|}$, where $H \subset G$ runs though all finite subgroups of $G$.
Roy has constructed a counterexample to this version in
\cite{Roy(1999)} based on her article \cite{Roy(2000)}. The examples of Roy do \emph{not}
contradict the Modified Trace Conjecture
\ref{con: Modified Trace Conjecture} or the 
Baum-Connes Conjecture \ref{con: Baum-Connes Conjecture}.



\subsubsection{The Stable Gromov-Lawson-Rosenberg Conjecture}
\label{subsec: The Stable Gromov-Lawson-Rosenberg Conjecture}

The Stable Gromov-Lawson-Rosenberg Conjecture is a typical conjecture relating Riemannian geometry to topology.
It is concerned with the question when a given manifold admits a metric of positive scalar curvature.
To discuss its relation with the Baum-Connes Conjecture 
we will need the real version of the Baum-Connes Conjecture, compare
Subsection~\ref{subsec: The Real Version of BC}.

Let $\Omega^{\Spin}_n(BG)$ be the bordism group of closed
$\Spin$-manifolds $M$ of dimension n with a reference map to
$BG$. Let $C^*_r(G;\IR)$ be the real reduced group $C^*$-algebra 
and let $KO_n(C^*_r(G;\IR)) = K_n(C^*_r(G;\IR))$ be its topological $K$-theory.
We use $KO$ instead of $K$ as a reminder that we here use the real
reduced group $C^*$-algebra. Given an element $[u \colon M \to BG] \in \Omega^{\Spin}_n(BG)$,
we can take the $C^*_r(G;\IR)$-valued index of the equivariant Dirac operator 
associated to the $G$-covering $\overline{M} \to M$ determined by $u$. Thus we get
a homomorphism
\begin{eqnarray} \ind_{C^*_r(G;\IR)} \colon \Omega^{\Spin}_n(BG) & \to & KO_n(C^*_r(G;\IR)). 
\label{index colon Omega^{Spin}_n to  KO_n(C^*_r(G;bbR))}
\end{eqnarray}
A \emph{Bott manifold}%
\index{Bott manifold}
is any simply connected closed $\Spin$-manifold $B$ of dimension $8$ whose
$\widehat{A}$-genus $\widehat{A}(B)$ is $8$. 
We fix such a choice, the particular choice
does not matter for the sequel. Notice that 
$\ind_{C^*_r(\{1\};\IR)}(B) \in KO_8(\IR) \cong \IZ$ is a generator
and the product with this element induces the Bott periodicity isomorphisms
$KO_n(C^*_r(G;\IR)) \xrightarrow{\cong}  KO_{n+8}(C^*_r(G;\IR))$. 
In particular 
\begin{eqnarray}
\ind_{C^*_r(G;\IR)}(M) & = & \ind_{C^*_r(G;\IR)}(M \times B),
\label{ind(M) = ind(M times B)}
\end{eqnarray}
if we identify $KO_n(C^*_r(G;\IR)) = KO_{n+8}(C^*_r(G;\IR))$ via Bott periodicity.

\begin{conjecture}[Stable Gromov-Lawson-Rosenberg Conjecture]
\label{con: Stable Gromov-Lawson-Rosenberg Conjecture}
\index{Conjecture!Gromov-Lawson-Rosenberg Conjecture!Stable Gromov-Lawson-Rosenberg Conjecture}
Let $M$ be a closed connected $\Spin$-manifold of dimension $n \ge 5$. 
Let $u_M \colon M \to B\pi_1(M)$ be the classifying
map of its universal covering.
Then $M \times B^k$ carries for some integer $k \ge
0$ a Riemannian metric with positive scalar curvature if and only if 
$$\ind_{C^*_r(\pi_1(M);\IR)}([M,u_M]) ~ = ~ 0 \hspace{5mm} \in KO_n(C^*_r(\pi_1(M);\IR)).$$
\end{conjecture}

If $M$ carries a Riemannian metric
with positive scalar curvature, then the index of the Dirac operator must vanish
by the Bochner-Lichnerowicz formula \cite {Rosenberg(1986b)}. 
The converse statement that the vanishing of the index
implies the existence of a Riemannian metric with positive scalar curvature is the hard part
of the conjecture. The following result is due to Stolz. 
A sketch of the proof can be found in
\cite[Section 3]{Stolz(2002)}, 
details are announced to appear in a different paper.

\begin{theorem}\label{BCC implies SGLR}
If the assembly map for the real version of the Baum-Connes Conjecture
(compare Subsection~\ref{subsec: The Real Version of BC})
is injective for the group $G$, then
the Stable Gromov-Lawson-Rosenberg Conjecture 
\ref{con: Stable Gromov-Lawson-Rosenberg Conjecture} is true for all closed
$\Spin$-manifolds of dimension $\ge 5$ with $\pi_1(M) \cong G$.
\end{theorem}

The requirement $\dim(M) \ge 5$ is essential in the 
Stable Gromov-Lawson-Rosenberg Conjecture,
since in dimension four 
new obstructions, the Seiberg-Witten invariants, occur.
The unstable version of this conjecture
says that $M$ carries a Riemannian metric with positive scalar curvature if and only if 
$\ind_{C^*_r(\pi_1(M);\IR)}([M,u_M]) = 0$. Schick \cite{Schick(1998e)}
constructs counterexamples to the unstable version using
 minimal hypersurface methods due to Schoen and Yau (see also 
\cite{Dwyer-Schick-Stolz(2002)}).
It is not known at the time of writing whether the unstable version is true for
finite fundamental groups. Since the Baum-Connes Conjecture 
\ref{con: Baum-Connes Conjecture} is true for finite groups (for the trivial reason
that $\EGF{G}{\calfin} = \pt$ for  finite groups $G$), Theorem
\ref{BCC implies SGLR} implies  that the Stable Gromov-Lawson Conjecture
\ref{con: Stable Gromov-Lawson-Rosenberg Conjecture} holds for finite fundamental groups
(see also \cite{Rosenberg-Stolz(1995)}).

The index map appearing in \eqref{index colon Omega^{Spin}_n to  KO_n(C^*_r(G;bbR))}
can be factorized as a composition
\begin{eqnarray} &
\ind_{C^*_r(G;\IR)} \colon \Omega^{\Spin}_n(BG) \xrightarrow{D} 
KO_n(BG) \xrightarrow{A} KO_n(C^*_r(G;\IR))
\label{Omega^{Spin}_n xrightarrow{D} KO_n(BG) xrightarrow{asmb}  KO_n(C^*_r(G;bbR))},
\end{eqnarray}
where $D$ sends $[M,u]$ to the class of the $G$-equivariant Dirac operator 
of the $G$-manifold $\overline{M}$
given by $u$ and $A=A_{\caltr}$ is the real version of the classical assembly map.
The homological Chern character defines an isomorphism
\[
KO_n (BG) \otimes_{\IZ} \IQ \xrightarrow{\cong} \bigoplus_{p \in \IZ}
H_{n + 4p}(BG;\IQ).
\]
Recall that associated to $M$ there is  the \emph{$\widehat{A}$-class}%
\index{A-class@$\widehat{\cala}$-class}
\begin{eqnarray} \widehat{\cala}(M)%
\indexnotation{widehat cala(M)}
 & \in & \prod_{p \ge 0} H^p(M;\IQ)
\label{widehat{cala}-genus}
\end{eqnarray}
which is a certain polynomial in the Pontrjagin classes.
The map $D$ appearing in \eqref{Omega^{Spin}_n xrightarrow{D} 
KO_n(BG) xrightarrow{asmb}  KO_n(C^*_r(G;bbR))}
sends the class of $u \colon M \to BG$ to $u_*(\widehat{\cala}(M) \cap [M])$, i.e.\
the image of the Poincar\'e dual of $\widehat{\cala}(M)$
under the map induced by $u$ in rational homology.
Hence $D([M,u]) = 0$ if and only if $u_*(\widehat{\cala}(M) \cap [M])$ vanishes. 
For  $x \in  \prod_{k \ge 0} H^k(BG;\IQ)$ define the 
\emph{higher $\widehat{A}$-genus of $(M,u)$ associated to $x$}%
\index{A-genus@$\widehat{A}$-genus!higher} 
to be
\begin{eqnarray}
\widehat{A}_x (M,u)%
\indexnotation{widehat{A}_x(M)} =  \langle \widehat{\cala}(M) \cup u^*x,[M]\rangle 
= 
\langle x,u_*(\widehat{\cala}(M) \cap [M]) \rangle \hspace{1mm} \in \IQ.
\label{higher widehat{A}-genus}
\end{eqnarray}
The vanishing of $\widehat{\cala}(M)$ is equivalent to the vanishing of all
higher $\widehat{A}$-genera $\widehat{A}_x(M,u)$.
The following conjecture is a weak version of the 
Stable Gromov-Lawson-Rosenberg Conjecture.

\begin{conjecture}[Homological Gromov-Lawson-Rosenberg Conjecture]
 \label{con: Homological Gromov-Lawson-Rosenberg Conjecture}
\index{Conjecture!Gromov-Lawson-Rosenberg Conjecture!Homological Gromov-Lawson-Rosenberg Conjecture}
Let $G$ be a group. Then for any closed $\Spin$-manifold $M$,
which admits a Riemannian metric with positive scalar curvature,
the $\widehat{A}$-genus $\widehat{A}_x(M,u)$ vanishes for all maps $u \colon M \to BG$ and elements
$x \in  \prod_{k \ge 0} H^k(BG;\IQ)$.
\end{conjecture}

From the discussion above we obtain the following result.

\begin{propositionnew} \label{pro: SNC for K_*(C^*_r(G)) implies HGLR}
If the assembly map
\[
KO_n ( BG ) \otimes_{\IZ} \IQ \to KO_n ( C_r^{\ast} ( G ; \IR )) \otimes_{\IZ} \IQ
\]
is injective for all $n \in \IZ$, then the Homological Gromov-Lawson-Rosenberg Conjecture holds for $G$.
\end{propositionnew}



\typeout{------------------------  Generalizations and Related Conjectures  --------------------------}


\section{Generalizations and Related Conjectures}
\label{chap: Generalizations and Related Conjectures}


\subsection{Variants of the Baum-Connes Conjecture}
\label{sec: Variants of the Baum-Connes Conjecture}

\subsubsection{The Real Version}
\label{subsec: The Real Version of BC}

There is an obvious real version of the Baum-Connes Conjecture,
 which predicts that for all $n \in \IZ$ and groups $G$ 
the assembly map
\[
A_{\calfin}^{\IR} \colon H_n^G ( \EGF{G}{\calf} ; \bfK^{\topo}_{\IR} ) \to KO_n ( C_r^{\ast} ( G ; \IR ))
\]
is an isomorphism. Here $H_n^G ( - ; \bfK^{\topo}_{\IR} )$ is an 
equivariant homology theory whose distinctive feature is that
$H_n^G( G/H ; \bfK^{\topo}_{\IR} ) \cong KO_n ( C_r^{\ast} ( H ; \IR ))$. 
Recall that we write $KO_n ( - )$ only to remind ourselves that
the $C^{\ast}$-algebra we apply it to is a real $C^{\ast}$-algebra, 
like for example the real reduced group $C^{\ast}$-algebra 
$C_r^{\ast} ( G ; \IR )$. The following result appears in \cite{Baum-Karoubi-Roe(2002)}.
\begin{propositionnew} \label{BC implies BC real}
The Baum-Connes Conjecture~\ref{con: Baum-Connes Conjecture} implies the real version of the Baum-Connes Conjecture.
\end{propositionnew}

In the proof of Proposition~\ref{pro: relating L and K} we have already seen 
that after inverting $2$ the ``real assembly map'' is a retract of the 
complex assembly map. In particular with $2$-inverted or after rationalizing also 
injectivity results or surjectivity results about the 
complex Baum-Connes assembly map yield the corresponding results for the real
Baum-Connes assembly map.

\subsubsection{The Version for Maximal Group $C^*$-Algebras}
\label{subsec: The Baum-Connes Conjecture for Maximal Group $C^*$-Algebras}

For a group $G$ let 
$C^*_{\max}(G)$%
\indexnotation{C^*_{max}(G)}
be its \emph{maximal group $C^*$-algebra},%
\index{C^*-algebra@$C^*$-algebra!maximal complex group $C^*$-algebra} 
compare~\cite[7.1.5 on page~229]{Pedersen(1979)}.
The maximal group $C^*$-algebra has the advantage  that every homomorphism of groups 
$\phi\colon G \to H$ induces a homomorphism $C^*_{\max}(G) \to C^*_{\max}(H)$ of
$C^*$-algebras. 
This is not true for the reduced
group $C^*$-algebra $C^*_r(G)$. Here is a counterexample: 
since  $C^*_r(F)$ is a simple algebra if $F$ is a non-abelian free group
\cite{Powers(1975)}, there is no unital algebra homomorphism $C^*_r(F) \to C^{\ast}_r ( \{ 1 \})= \IC $.

One can construct a version of the Baum-Connes assembly map using an equivariant homology theory 
$H_n^G ( - ; \bfK^{\topo}_{\max} )$ which evaluated on $G/H$ yields the $K$-theory of $C_{\max}^{\ast} ( H )$
(use Proposition~\ref{pro: Or(G)-spectra yield a G-homology theory} and a suitable modification of $\bfK^{\topo}$, 
compare Section~\ref{sec: K and L-Theory Spectra over Groupoids}).

Since on the left hand side of a $\calfin$-assembly map only the maximal group $C^{\ast}$-algebras 
for finite groups $H$ matter, and clearly 
$C^{\ast}_{\max} (H) = \IC H = C^{\ast}_r ( H )$ for such $H$, this left hand side coincides 
with the left hand side of the usual Baum-Connes Conjecture. 
There is always a $C^*$-homomorphism $p \colon C^*_{\max}(G) \to C^*_r(G)$ 
(it is an isomorphism if and only if $G$ is amenable \cite[Theorem~7.3.9 on page 243]{Pedersen(1979)})
and hence we obtain the following 
factorization of the usual Baum-Connes assembly map
\begin{eqnarray}
&
\xymatrix{                           
 & &   K_n(C^*_{\max}(G)) \ar[d]^-{K_n (p)} \\
H^G_n(\EGF{G}{\calfin};\bfK^{\topo}) \ar[rr]_-{A_{\calfin}} \ar[urr]^-{A_{\calfin}^{\max}} 
& & K_n ( C^{\ast}_{r} ( G ) )
         }
&
\label{assembly map for maximal C^*-algebra}
\end{eqnarray}
It is known that the  map $A_{\calfin}^{\max}$ is in general not surjective. 
The Baum-Connes Conjecture would imply that the map is $A_{\calfin}^{\max}$ is always injective, and that it is surjective
if and only if the vertical map $K_n ( p )$ is injective.

A countable group $G$ is called \emph{$K$-amenable}%
\index{group!K-amenable@$K$-amenable}
if the map $p \colon C_{\max}^{\ast} ( G ) \to C_r^{\ast} (G)$ induces a $KK$-equivalence (compare \cite{Cuntz(1983)}).
This implies in particular that the 
vertical map $K_n ( p )$ is an isomorphism for all $n \in \IZ$.
Note that for $K$-amenable groups the Baum-Connes Conjecture holds if and only if the
``maximal'' version of the assembly map $A_{\calfin}^{\max}$ is an isomorphism for all $n \in \IZ$.
A-T-menable groups are $K$-amenable, compare 
Theorem~\ref{the: BCC with coefficients for a-Tmenable groups}.
But $K_0 ( p )$ is not injective for every  infinite group which has property (T) 
such as for example $SL_n(\IZ)$ for $n \ge 3$,
compare for instance the discussion in \cite{Julg(1997)}.
There are groups with property (T) for which the Baum-Connes Conjecture is known 
(compare Subsection~\ref{subsec: Status of the Baum-Connes Conjecture})
and hence there are counterexamples to 
the conjecture that $A_{\calfin}^{\max}$ is an isomorphism. 

In Theorem~\ref{the: Homotopy invariance of tau(2)(M) and  rho(2)(M)} and Remark~\ref{rem: L2-signature theorem}
we discussed applications of the maximal $C^{\ast}$-algebra version of the Baum-Connes Conjecture.

\subsubsection{The Bost Conjecture}
\label{subsec: The Bost Conjecture}

Some of the strongest results about the Baum-Connes Conjecture are proven using the so called 
Bost Conjecture (see \cite{Lafforgue(2002)}).
The Bost Conjecture is the version of the Baum-Connes Conjecture, where one replaces 
the reduced group $C^*$-algebra $C^*_r(G)$  by the Banach algebra $l^1(G)$ of absolutely summable functions on $G$.
Again one can use the spectra approach 
(compare Subsection~\ref{sec: Spectra over the Orbit Category} and \ref{sec: K and L-Theory Spectra over Groupoids} and in particular 
Proposition~\ref{pro: Or(G)-spectra yield a G-homology theory}) to 
produce a variant of equivariant $K$-homology denoted $H_n^G ( - ; \bfK^{\topo}_{l^1} )$ which this time
evaluated on $G/H$ yields $K_n ( l^1 (H))$, the topological $K$-theory of the Banach algebra $l^1 ( H )$.

As explained in the beginning of Chapter~\ref{chap: general formulation},
we obtain an associated assembly map and we believe that 
it coincides with the one defined using a Banach-algebra version of $KK$-theory in \cite{Lafforgue(2002)}.

\begin{conjecture}[Bost Conjecture]
\label{con: Bost Conjecture}
\index{Conjecture!Bost Conjecture}
Let $G$ be a countable group.
Then the assembly map 
\begin{eqnarray*}
A_{\calfin}^{l^1}  \colon
H^G_n(\EGF{G}{\calfin};\bfK^{\topo}_{l^1}) \to K_n(l^1(G))
\end{eqnarray*}
is an isomorphism.
\end{conjecture}

Again the left hand side coincides with the left hand side of the Baum-Connes assembly map 
because for finite groups $H$ we have
$l^1 ( H ) =  \IC H = C_r^{\ast} ( H )$. There is always a homomorphism of Banach algebras 
$q \colon l^1 ( G ) \to C_r^{\ast} ( G ) $ and one obtains a factorization
of the usual Baum-Connes assembly map
\[
\xymatrix{
& & K_n ( l^1 (G) ) \ar[d]^-{K_n (q)} \\
H_n^G ( \EGF{G}{\calfin} ; \bfK^{\topo} ) \ar[urr]^-{A_{\calfin}^{l^1}} \ar[rr]^-{A_{\calfin}} & & K_n ( C_r^{\ast} ( G ) ).
         }
\]
Every group homomorphism $G \to H$ induces a
homomorphism of Banach algebras $l^1(G) \to l^1(H)$.
So similar as in the maximal group $C^{\ast}$-algebra case this approach repairs the lack of functoriality
for the reduced group $C^{\ast}$-algebra.

The disadvantage
of $l^1(G)$ is however that indices of operators 
tend to take values in the topological $K$-theory of the group $C^*$-algebras, not in $K_n( l^1(G))$.
Moreover the representation theory of $G$ is closely related to the group $C^{\ast}$-algebra, whereas the relation to
$l^1 ( G)$ is not well understood.

For more information 
about the Bost Conjecture \ref{con: Bost Conjecture} see 
\cite{Lafforgue(2002)},
\cite{Skandalis(1999)}. 

\subsubsection{The Baum-Connes Conjecture with Coefficients}
\label{subsec: The Baum-Connes Conjecture with Coefficients}

The Baum-Connes Conjecture~\ref{con: Baum-Connes Conjecture} can be generalized to the Baum-Connes Conjecture with 
Coefficients.
Let $A$ be a separable $C^*$-algebra with an action of the countable
group $G$. Then there is an assembly map
\begin{eqnarray}
KK^G_n(\EGF{G}{\calfin};A)  \to  K_n(A\rtimes G)
\label{assembly map for the version with coefficients}
\end{eqnarray}
defined in terms of equivariant $KK$-theory, compare Sections \ref{sec: Equivariant KK-theory} and
\ref{sec: The Dirac-Dual Dirac Method}.

\begin{conjecture}[Baum-Connes Conjecture with Coefficients]
\label{con: Baum-Connes Conjecture with coefficients}
\index{Conjecture!Baum-Connes Conjecture!with Coefficients}
For every separable $C^{\ast}$-algebra $A$ with an action of a countable group $G$ and every $n \in \IZ$  
the assembly map~\eqref{assembly map for the version with coefficients} is an isomorphism.
\end{conjecture}

There are counterexamples to the 
Baum-Connes Conjecture with Coefficients,
compare Remark~\ref{rem: status of BC with coefficients}.
If we take $A = \IC$ with the trivial action, the map~\eqref{assembly map for the version with coefficients} 
can be identified with the assembly map appearing in the ordinary 
Baum-Connes Conjecture~\ref{con: Baum-Connes Conjecture}.


\begin{remarknew}[A Spectrum Level Description]
There is a formulation of the Baum-Connes Conjecture with Coefficients 
in the framework explained in Section~\ref{sec: Spectra over the Orbit Category}.
Namely, construct an appropriate  covariant functor $\bfK^{\topo}( A \rtimes \calg^G ( - )) \colon \Or(G) \to \SPECTRA$ 
such that 
\[
\pi_n(\bfK^{\topo} ( A \rtimes \calg^G  (G/H)) \cong K_n(A \rtimes H)
\]
holds for all subgroups $H \subseteq G$ and all $n \in \IZ$, and consider the associated $G$-homology theory
$H^G_*(-;\bfK^{\topo} ( A \rtimes \calg^G( - )))$. Then the map 
\eqref{assembly map for the version with coefficients} can be
identified with the map which the projection $\pr\colon
\EGF{G}{\calfin} \to \pt$ induces for this homology theory.
\end{remarknew}


\begin{remarknew}[Farrell-Jones Conjectures with Coefficients]
\label{rem: FJ-Conjecture with Coefficients} 
One can also formulate a ``Farrell-Jones Conjecture with Coefficients''. 
(This should not be confused with the Fibered Farrell-Jones Conjecture discussed in 
Subsection~\ref{subsec: Fibered Version of the Farrell-Jones Conjecture}.)
Fix a ring $S$ and an action of $G$ on it by isomorphisms
of rings. Construct an appropriate  
covariant functor $\bfK (S \rtimes \calg^G ( - )) \colon \Or(G) \to \SPECTRA$ 
such that 
\[
\pi_n(\bfK (S \rtimes \calg^G (G/H))) \cong K_n(S \rtimes H)
\]
holds for all subgroups $H \subseteq G$ and $n \in \IZ$, where $S \rtimes H$ 
is the associated twisted group ring. Now consider the associated $G$-homology theory
$H^G_*(-;\bfK (S \rtimes \calg^G ( - )))$. There is an analogous construction for $L$-theory.
A \emph{Farrell-Jones Conjecture with Coefficients}%
\index{Conjecture!Farrell-Jones Conjecture!with Coefficients}
would say that the map induced on these homology theories by the projection
$\pr \colon \EGF{G}{\calvcyc} \to \pt$ is always an isomorphism.
We do not know whether there are counterexamples to the Farrell-Jones Conjectures with 
Coefficients, compare Remark~\ref{rem: status of BC with coefficients}.
\end{remarknew}

\subsubsection{The Coarse Baum Connes Conjecture}
\label{subsec: The Coarse Baum Connes Conjecture}

We briefly explain the Coarse Baum-Connes Conjecture, a variant of the Baum-Connes Conjecture,
which applies to metric spaces. Its importance lies in the fact that isomorphism results about the Coarse Baum-Connes Conjecture
can be used to prove injectivity results about the classical assembly map for topological $K$-theory.
Compare also  Section~\ref{sec: The descent principle}.

Let $X$ be a proper (closed balls are compact) metric space and $H_X$ a separable Hilbert space with a faithful 
nondegenerate $\ast$-representation of $C_0(X)$, the algebra of complex valued continuous functions which vanish 
at infinity. A bounded linear operator $T$ has a support
$\supp T \subset X \times X$, which is defined as the complement of the set of all pairs $(x, x')$, for 
which there exist functions $\phi$ and $\phi' \in C_0(X)$ such that $\phi(x) \neq 0$, $\phi'(x') \neq 0$ and $\phi' T \phi = 0$.
The operator $T$ is said to be a finite propagation operator if there exists a constant $\alpha$ such that 
$d(x,x') \leq \alpha$ for all pairs in the support of $T$. The operator is said to 
be \emph{locally compact} if $\phi T$ and $T\phi$ are
compact for every $\phi \in C_0(X)$. An operator is called \emph{pseudolocal} if $\phi T \psi$ is a compact operator 
for all pairs of continuous functions $\phi$ and $\psi$ with compact and disjoint supports.

The Roe-algebra $C^{\ast}( X )=C(X , H_X)$ is the operator-norm closure of the $\ast$-algebra of all 
locally compact finite propagation operators on $H_X$. The algebra $D^{\ast}(X) =D^{\ast}(X , H_X)$ is the 
operator-norm closure of the pseudolocal finite propagation operators. One can show that the topological $K$-theory 
of the quotient algebra $D^{\ast}(X) / C^{\ast}(X)$ coincides up to an index shift with the analytically defined 
(non-equivariant)
$K$-homology $K_{\ast}(X)$, compare Section~\ref{sec: Analytic Equivariant $K$-Homology}. For a uniformly contractible proper
metric space the coarse assembly map $K_{n} ( X ) \to K_n ( C^{\ast}(X) )$
is the boundary map in the long exact sequence associated to the short exact sequence of $C^{\ast}$-algebras
\[
0 \to C^{\ast}(X) \to D^{\ast} ( X ) \to D^{\ast} (X) / C^{\ast} (X) \to 0.
\]
For general metric spaces one first approximates the metric space by spaces with nice local behaviour, compare~\cite{Roe(1996)}.
For simplicity we only explain the case, where $X$ is a discrete metric space.
Let $P_d (X)$ the Rips complex for a fixed distance $d$, i.e.\ the simplicial complex with vertex set $X$, where a simplex is spanned 
by every collection of points in which every two points are a distance less than $d$ apart. Equip $P_d ( X)$ with the spherical metric,
compare \cite{Yu(1997)}.

A discrete metric space has \emph{bounded geometry}%
\index{bounded geometry} if for each  $r > 0$ there exists a
$N(r)$ such that for all $x$ the ball of radius $r$ centered at  $x \in X$ 
contains at most $N(r)$ elements. 

\begin{conjecture}[Coarse Baum-Connes Conjecture]
\label{con:CBC}
Let $X$ be a proper discrete metric space of bounded geometry. Then for $n=0$, $1$ the 
coarse assembly map
\[
\colim_{d} K_n (P_d(X)) \to \colim_d K_n (C^*(P_d (X))) \cong K_n ( C^{\ast} ( X) )
\]
is an isomorphism.
\end{conjecture}

The conjecture is false if one drops the bounded geometry hypothesis. A counterexample can be found in \cite[Section~8]{Yu(1998)}.
Our interest in the conjecture stems from the following fact, compare~\cite[Chapter~8]{Roe(1996)}.
\begin{propositionnew}
\label{pro:descentCBC}
Suppose the finitely generated group $G$ admits a classifying space $BG$ of finite type. If $G$ considered as a metric
space via a word length metric satisfies the Coarse Baum-Connes Conjecture~\ref{con:CBC} then the classical assembly map
$A \colon K_{\ast} ( B G ) \to K_{\ast} ( C_r^{\ast} G )$ which appears in Conjecture~\ref{con: BC torsion free} is injective.
\end{propositionnew}
The Coarse Baum-Connes Conjecture for a discrete group $G$ (considered as a metric space) can be interpreted as a case of the Baum-Connes 
Conjecture with Coefficients~\ref{con: Baum-Connes Conjecture with coefficients} for the group $G$ with 
a certain specific choice of coefficients, compare \cite{Yu(1995)}.

Further information about the coarse Baum-Connes Conjecture can be
found for instance in 
\cite{Higson-Pedersen-Roe(1997)}, 
\cite{Higson-Roe(1995)}, 
\cite{Higson-Roe(2000)},  
\cite{Roe(1996)}, 
\cite{Yu(1995a)},
\cite{Yu(1995b)},
\cite{Yu(1997)},
\cite{Yu(1998a)},
and 
\cite{Yu(2000)}.

\subsubsection{The Baum-Connes Conjecture for Non-Discrete Groups}
\label{subsec: The Baum-Connes Conjecture for Non-Discrete Groups}

Throughout this subsection let $T$ be a 
locally compact second countable topological Hausdorff group.
There is a notion of a classifying space for proper $T$-actions 
$\underline{E}T$ (see 
\cite[Section 1 and 2]{Baum-Connes-Higson(1994)}
\cite[Section I.6]{Dieck(1987)}, \cite[Section 1]{Lueck-Meintrup(2000)})
and one can define its equivariant topological
$K$-theory $K^{T}_n(\underline{E}T)$. The definition of a reduced 
$C^*$-algebra $C_r^*(T)$ and its topological
$K$-theory $K_n(C^*_r(T))$ makes sense also for $T$. 
There is an 
assembly map defined in terms of equivariant index theory
\begin{eqnarray}
A_{\calk} \colon K_n^{T}(\underline{E}T) & \to & K_n(C^*_r(T)).
\label{assembly for topological groups}
\end{eqnarray}
The Baum-Connes Conjecture for $T$ says that this map is
bijective for all $n \in \IZ$ \cite[Conjecture 3.15 on page 254]{Baum-Connes-Higson(1994)}.

Now consider the special case where $T$ is a connected Lie group. Let $\calk$ be the family
of compact subgroups of $T$. There is a notion of a $T$-$CW$-complex and of
a classifying space $\EGF{T}{\calk}$ defined as in Subsection~\ref{subsec: G-CW complexes} and \ref{subsec: EGFs}.
The classifying space 
$\EGF{T}{\calk}$ yields  a model for $\underline{E}T$. 
Let $K\subset T$ be a maximal compact subgroup. It 
is unique up to conjugation. The space $T/K$ is contractible and in fact a
model for $\underline{E}T$ 
(see \cite[Appendix, Theorem~A.5]{Abels(1974)}, \cite[Corollary 4.14]{Abels(1978)},
\cite[Section 1]{Lueck-Meintrup(2000)}).
One knows (see \cite[Proposition 4.22]{Baum-Connes-Higson(1994)}, \cite{Kasparov(1988)})
\begin{eqnarray*}
K_n^T ( \underline{E}T ) = K_n^{T}(T/K) & = & \left\{\begin{array}{lll}
R_{\IC}(K) & & n  = \dim(T/K) \mod 2,\\
0 & & n  = 1 + \dim(T/K) \mod 2, 
\end{array}\right.
\end{eqnarray*}
where $R_{\IC}(K)$ is the complex representation ring of $K$.

Next we consider the special case where $T$ is a totally disconnected group.
Let $\calko$ be the family of compact-open subgroups of $T$.
A $T$-$CW$-complex and 
a classifying space $\EGF{T}{\calko}$ for $T$ and $\calko$ are defined similar as in 
Subsection~\ref{subsec: G-CW complexes} and \ref{subsec: EGFs}.
Then $\EGF{T}{\calko}$ is a model for $\underline{E}T$ 
since any compact subgroup is contained in a compact-open subgroup,
and the Baum-Connes Conjecture says that the assembly map yields for $n \in \IZ$ an isomorphism
\begin{eqnarray}
A_{\calko} \colon K_n^{T}(\EGF{T}{\calko}) & \to & K_n(C^*_r(T)).
\label{assembly for totally disconnected groups}
\end{eqnarray}
For more information see \cite{Baum-Higson-Plymen(2000a)}.


\subsection{Variants of the Farrell-Jones Conjecture}
\label{sec: Variants of the Farrell-Jones Conjectures}

\subsubsection{Pseudoisotopy Theory}
\label{subsec: Pseudoisotopy Theory}

An important variant of the Farrell-Jones Conjecture deals with 
the pseudoisotopy spectrum functor $\bfP$, which we already discussed briefly in Section~\ref{subsec: negative K pseudo}.
In fact it is this variant of the Farrell-Jones Conjecture 
(and its fibered version which will be explained in the next subsection)
for which the strongest results are known at the time of writing.

In Proposition~\ref{pro: GROUPOID-spectra and equivariant homology theories}
we will explain that every functor 
$\bfE \colon \GROUPOIDS \to \SPECTRA,$
which sends equivalences of groupoids to stable weak equivalences of spectra, 
yields a corresponding equivariant homology theory $H_n^G ( - ; \bfE )$.
Now whenever we have a functor $\bfF \colon \SPACES \to \SPECTRA$, we can 
precompose it with the functor ``classifying space'' which sends a groupoid $\calg$ to its classifying space $B\calg$.
(Here $B \calg$ is simply the realization of the nerve of $\calg$ considered as a category.) 
In particular this applies to the pseudoisotopy functor $\bfP$. Thus we obtain 
a homology theory $H_n^G ( - ; \bfP \circ  B )$ whose essential feature is that
\[
H_n^G ( G/H ; \bfP \circ B ) \cong \pi_n ( \bfP ( B H )),
\]
i.e.\ evaluated at $G/H$ one obtains the homotopy groups of the 
pseudoisotopy spectrum of the classifying space $BH$ of the 
group $H$. As the reader may guess there is the following conjecture.
\begin{conjecture}[Farrell-Jones Conjecture for Pseudoisotopies of Aspherical Spaces] 
\label{con: FJ for pseudoisotopies unfibered}
\index{Conjecture!Farrell-Jones Conjecture!for Pseudoisotopies of Aspherical Spaces}
For every group $G$ and all $n \in \IZ$ the assembly map
\[
H_n^G ( \EGF{G}{\calvcyc} ; \bfP \circ B ) \to H_n^G ( \pt ; \bfP \circ B ) \cong \pi_n ( \bfP ( BG ) )
\]
is an isomorphism. Similarly for $\bfP^{\diff}$, the pseudoisotopy functor which is defined using differentiable pseudoisotopies.
\end{conjecture}
A formulation of a conjecture
for spaces which are not necessarily aspherical will be given in the next subsection, see in particular 
Remark~\ref{rem: special case universal covering}.


\begin{remarknew}[Relating $K$-Theory and Pseudoisotopy Theory]
\label{rem: relating K and P}

We already outlined in Subsection~\ref{subsec: pseudoisotopies} 
the relationship between $K$-theory and pseudoisotopies. The comparison in positive dimensions 
described there can be extended to all dimensions. 
Vogell constructs in \cite{Vogell(1990)} a version of $A$-theory using retractive spaces that are bounded 
over $\IR^k$ (compare Subsection~\ref{subsec: bounded h-cobordisms} and \ref{subsec: negative K pseudo}). 
This leads to a functor $\bfA^{-\infty}$ from spaces to non-connective spectra.
Compare also 
\cite{Carlsson-Pedersen-Vogell(1998)},
\cite{Vogell(1991)}, 
\cite{Vogell(1995)} and \cite{Weiss(2002)}. 
We define ${\bf Wh}_{PL}^{-\infty}$ via the fibration sequence
\[
X_+ \sma \bfA^{-\infty} ( \pt ) \to \bfA^{-\infty} ( X ) \to {\bf Wh}_{PL}^{-\infty} ( X ),
\]
where the first map is the assembly map.
The natural equivalence  
\[
\Omega^2 {\bf Wh}^{-\infty}_{PL}(X)  \simeq \bfP (X)
\]
seems to be hard to trace down in the
literature but should be true. We will assume it in the following discussion.

Precompose the functors above with the classifying space functor $B$ to obtain functors from groupoids to 
spectra.
The pseudoisotopy assembly map which 
appears in 
Conjecture~\ref{con: FJ for pseudoisotopies unfibered} is an isomorphism if and only if
the $A$-theory assembly map
\[
H_{n+2}^G ( \EGF{G}{\calvcyc} ; \bfA^{-\infty} \circ B) \to 
H_{n+2}^G ( \pt ; \bfA^{-\infty} \circ B) \cong \pi_{n+2} ( \bfA^{-\infty}( BG) )
\]
is an isomorphism. This uses a $5$-lemma argument and the fact that for a fixed spectrum $\bfE$ the assembly map
\[
H_n^G(\EGF{G}{\calf}; B \calg^G ( - )_+ \sma \bfE ) \to H_n^G ( \pt ; B \calg^G ( - )_+ \sma \bfE )
\]
is always bijective.
There is a linearization map $\bfA^{-\infty} ( X ) \to \bfK ( \IZ \Pi ( X )_{\oplus} )$ (see the next subsection for the notation) 
which is always $2$-connected and a rational equivalence if $X$ is 
aspherical (recall that $\bfK$ denotes the non-connective $K$-theory spectrum). 
For finer statements about the linearization map, compare also~\cite{Nicas(1985)}.

The above discussion yields in particular the following, compare \cite[1.6.7 on page 261]{Farrell-Jones(1993a)}. 

\begin{propositionnew} \label{pro: relating K and pseudoisotopies}
The rational version of the $K$-theoretic Farrell-Jones 
Conjecture~\ref{con: Farrell-Jones Conjecture} 
is equivalent  to the rational version of the 
Farrell-Jones Conjecture  for Pseudoisotopies of Aspherical Spaces~\ref{con: FJ for pseudoisotopies unfibered}.
If the assembly map in the conjecture for pseudoisotopies is (integrally) an isomorphism for $n \leq -1$, then
so is the assembly map in the $K$-theoretic Farrell-Jones Conjecture for $n \leq 1$.
\end{propositionnew}

\end{remarknew}

\subsubsection{Fibered Versions}
\label{subsec: Fibered Version of the Farrell-Jones Conjecture}

Next we present the more general fibered versions of the Farrell-Jones Conjectures.
These fibered versions have better inheritance properties, compare Section~\ref{sec: Inheritance Properties}.

In the previous section we considered functors 
$\bfF \colon \SPACES \to \SPECTRA$, like $\bfP$, $\bfP^{\diff}$ and $\bfA^{-\infty}$, and the associated equivariant homology theories
$H_n^G ( - ; \bfF \circ B )$ (compare Proposition~\ref{pro: GROUPOID-spectra and equivariant homology theories}). 
Here $B$ denotes the classifying space functor, which sends a groupoid $\calg$ to its classifying space $B \calg$.
In fact all equivariant homology theories we considered so far can be obtained in this fashion for special choices of $\bfF$.
Namely, let $\bfF$ be one of the functors
\begin{eqnarray*}
\bfK ( R \Pi ( - )_{\oplus}) , \quad 
\bfL^{\langle -\infty \rangle} (  R \Pi (- )_{\oplus} )  \quad \mbox{ or } \quad  
\bfK^{\topo} (  C_r^{\ast}  \Pi ( - )_{\oplus} ) ,
\end{eqnarray*}
where $\Pi(X)$ denotes the fundamental groupoid of a space, $R \calg_{\oplus}$ respectively $C_r^{\ast} \calg_{\oplus}$
is the $R$-linear respectively the $C^{\ast}$-category associated to a groupoid $\calg$ and $\bfK$, $\bfL^{\langle -\infty \rangle}$
and $\bfK^{\topo}$ are suitable functors which send additive respectively $C^{\ast}$-categories to spectra, compare the proof of 
Theorem~\ref{the: K- and L-Theory Spectra over Groupoids}.
There is a natural equivalence
$\calg \to \Pi B \calg$. Hence, if we precompose the functors above 
with the classifying space functor $B$ we obtain functors which are equivalent 
to the functors we have so far been calling
\[
\bfK_R, \quad 
\bfL_R^{\langle - \infty \rangle} \quad \mbox{ and } \quad 
\bfK^{\topo},
\]
compare Theorem~\ref{the: K- and L-Theory Spectra over Groupoids}.
Note that in contrast to these three cases
the pseudoisotopy functor $\bfP$ depends on more than just the fundamental groupoid. However
Conjecture~\ref{con: FJ for pseudoisotopies unfibered} above only deals with aspherical spaces.

Given a $G$-CW-complex $Z$ and  a functor $\bfF $ from spaces to spectra we obtain a functor
$X \mapsto \bfF( Z \times_{G} X)$ which digests $G$-CW-complexes. In particular we can restrict it to the 
orbit category to obtain a functor
\[
\bfF ( Z \times_G - ) \colon \Or (G) \to \SPECTRA.
\]
According to Proposition~\ref{pro: Or(G)-spectra yield a G-homology theory} we obtain 
a corresponding $G$-homology theory 
\[
H_n^G ( - ; \bfF ( Z \times_G - ) )
\]
and associated assembly maps. Note that restricted to the orbit category 
the functor $EG \times_G - $ is equivalent to the classifying space functor $B$ and so
$H_n^G ( - ; \bfF \circ B )$ can be considered as a special case of this construction.

\begin{conjecture}[Fibered Farrell-Jones Conjectures] \label{con: Fibered Farrell-Jones Conjecture}
\index{Conjecture!Farrell-Jones Conjecture!fibered version}
Let $R$ be a ring (with involution).
Let $\bfF \colon \SPACES \to \SPECTRA$ be one of the functors 
\[
\bfK( R \Pi ( - )_{\oplus} ), \quad  \bfL^{\langle -\infty \rangle} ( R \Pi ( - )_{\oplus} ), \quad \bfP ( - ) , \quad \bfP^{\diff} ( - )
\quad \mbox{ or } \quad \bfA^{-\infty} ( - ).
\]
Then for every free $G$-CW-complex $Z$ and all $n \in \IZ$ the 
associated assembly map 
\[
H_n^G ( \EGF{G}{\calvcyc} ; \bfF( Z \times_G - ) ) \to H_n^G ( \pt ; \bfF ( Z \times_G - )) \cong \pi_n ( \bfF ( Z/G ) )
\]
is an isomorphism.
\end{conjecture}

\begin{remarknew}[A Fibered Baum-Connes Conjecture]
With the family $\calfin$ instead of $\calvcyc$
and the functor $\bfF=\bfK^{\topo} ( C_r^{\ast} \Pi ( - )_{\oplus} )$ one obtains
a \emph{Fibered Baum-Connes Conjecture}.%
\index{Conjecture!Baum-Connes Conjecture!fibered version}
\end{remarknew}

\begin{remarknew}[The Special Case $Z = \widetilde{X}$]
\label{rem: special case universal covering}
Suppose $Z = \widetilde{X}$ is the universal covering of a space $X$ 
equipped with the action of its fundamental group $G= \pi_1 ( X )$. Then in the algebraic $K$- and $L$-theory case 
the conjecture above specializes to the ``ordinary'' Farrell-Jones Conjecture~\ref{con: Farrell-Jones Conjecture}.
In the pseudoisotopy and $A$-theory case one obtains a formulation of an (unfibered) conjecture about 
$\pi_n ( \bfP ( X ))$ or $\pi_n ( \bfA^{-\infty} ( X ))$ for spaces $X$ which are not necessarily aspherical.
\end{remarknew}

\begin{remarknew}[Relation to the Original Formulation]
\label{rem: relation to original}
In \cite{Farrell-Jones(1993a)} Farrell and Jones formulate a fibered version of their conjectures for every
(Serre) fibration $Y \to X$ over a connected CW-complex $X$. 
In our set-up this corresponds to choosing $Z$ to be the total space of the fibration obtained from $Y \to X$
by pulling back along the universal covering projection $\widetilde{X} \to X$. This space is a free $G$-space
for $G= \pi_1 ( X )$. 
Note that an arbitrary free $G$-$CW$-complex  $Z$ can always be obtained in this fashion from a map
$Z/G \to BG$, compare~\cite[Corollary 2.2.1 on page 264]{Farrell-Jones(1993a)}.
\end{remarknew}

\begin{remarknew}[Relating $K$-Theory and Pseudoisotopy Theory in the Fibered Case]
\label{rem: relating K and P fibered}
The linearization map $\pi_n ( \bfA^{-\infty} ( X ) ) \to K_n ( \IZ \Pi ( X ) )$ 
is always $2$-connected, but for spaces which are not aspherical it need not be a rational equivalence.
Hence the comparison results discussed in Remark~\ref{rem: relating K and P} apply for the fibered versions only in
dimensions $n \le 1$.
\end{remarknew}

\subsubsection{The Isomorphism Conjecture for $N\! K$-groups}
\label{subsec: The Isomorphism Conjecture for NK-groups}

In Remark~\ref{rem: Bass-Heller-Swan decomposition} we defined the groups $N\! K_n ( R )$ for a ring $R$.
They are the simplest kind of Nil-groups responsible for the infinite cyclic group.
Since the functor $\bfK_R$ is natural with respect to ring homomorphism we can define $\bfN \! \bfK_{R}$ 
as the (objectwise) homotopy cofiber of $\bfK_R \to \bfK_{R[t]}$. There is an associated assembly map.
\begin{conjecture}[Isomorphism Conjecture for $N \! K$-groups]
\label{con: Isomorphism Conjecture for NK}
\index{Conjecture!Isomorphism Conjecture for NK-group@Isomorphism Conjecture for $NK$-groups}
The assembly map
\[
H_n^G ( \EGF{G}{\calvcyc} ; \bfN \! \bfK_R ) \to H_n^G ( \pt ; \bfN \! \bfK_R ) \cong N\! K_n ( RG )
\]
is always an isomorphism. 
\end{conjecture}
There is a weak equivalence $\bfK_{R[t]} \simeq \bfK_R \vee \bfN \! \bfK_{R}$ of functors  from $\GROUPOIDS$ to $\SPECTRA$.
This implies for a fixed family $\calf$ of subgroups of $G$ and $n \in \IZ$ 
that whenever two of the three assembly maps 
\begin{eqnarray*}
A_{\calf} \colon H_n^G(\EGF{G}{\calf};\bfK_{R[t]}) & \to & K_n(R[t][G]),  \\
A_{\calf} \colon H_n^G(\EGF{G}{\calf};\bfK_{R}) & \to &  K_n(R[G]), \\
A_{\calf} \colon H_n^G ( \EGF{G}{\calf} ; \bfN \! \bfK_R ) & \to  & N\! K_n (RG)
\end{eqnarray*}
are bijective, then so is the third (compare \cite[Section 7]{Bartels-Farrell-Jones-Reich(2004a)}). 
Similarly one can define a functor 
$\bfE_R$ from the category
$\GROUPOIDS$ to $\SPECTRA$ and weak equivalences
$$\bfK_{R[t,t^{-1}]} \to \bfE_{R} \xleftarrow{} \bfK_R \vee \Sigma \bfK_R \vee  \bfN \! \bfK_R  \vee  \bfN \! \bfK_R,$$
which on homotopy groups corresponds to the Bass-Heller-Swan decomposition
(see Remark \ref{rem: Bass-Heller-Swan decomposition}). One obtains a two-out-of-three statement as above with the 
$\bfK_{R[t]}$-assembly map replaced by the $\bfK_{R[t,t^{-1}]}$-assembly map.

\subsubsection{Algebraic $K$-Theory of the Hecke Algebra}
\label{Algebraic $K$-theory of the Hecke Algebra}

In Subsection~\ref{subsec: The Baum-Connes Conjecture for Non-Discrete Groups} we mentioned the classifying space
$\EGF{G}{\calko}$ for the family of compact-open subgroups and the Baum-Connes Conjecture for a 
totally disconnected group $T$. There is an analogous conjecture dealing with the 
algebraic $K$-theory of the Hecke algebra.

Let $\calh(T)$%
\indexnotation{calh(Gamma)}
 denote the \emph{Hecke algebra}%
\index{Hecke algebra} of $T$ which consists of
locally constant functions $G \to \IC$ with compact support and
inherits its multiplicative structure from the convolution product. 
The Hecke algebra $\calh(T)$ plays the same role for $T$ as the complex group ring $\IC G$ for a discrete
group $G$ and reduces to this notion if $T$ happens to be discrete. There is a $T$-homology theory
$\calh^{T}_*$ with the property that for any open and closed subgroup $H \subseteq T$ and all $n \in \IZ$
we have $\calh^{T}_n(T/H) =  K_n(\calh(H))$, where $K_n(\calh(H))$ is the
algebraic $K$-group of the Hecke algebra $\calh(H)$. 
\begin{conjecture}[Isomorphism Conjecture for the Hecke-Algebra]
\label{Farrell-Jones Conjecture for totally disconnected groups}
\index{Conjecture!Isomorphism Conjecture for the Hecke-Algebra}
For a totally disconnected group $T$ the assembly map 
\begin{eqnarray}
A_{\calko} \colon \calh^{T}_n(\EGF{T}{\calko}) & \to &
\calh^{T}(\pt) = K_n(\calh(T))
\end{eqnarray}
induced by the projection
$\pr\colon \EGF{T}{\calko} \to \pt$ is an isomorphism for all $n \in \IZ$.
\end{conjecture}
In the case $n = 0$ this reduces to the statement that 
\begin{eqnarray}
\colim_{T/H \in \OrGF{T}{\calko}} K_0(\calh(H)) & \to & K_0(\calh(T)).
\label{Farrell-Jones Conjecture for totally disconnected groups for K_0}
\end{eqnarray}
is an isomorphism.
For $n \le -1$ one obtains the statement that $K_n(\calh(G)) = 0$. 
The group $K_0(\calh(T))$ has an interpretation in terms of the smooth
representations of $T$. The $G$-homology theory can be constructed using 
an appropriate functor $\bfE\colon \OrGF{T}{\calko} \to \SPECTRA$ 
and the recipe explained in Section~\ref{sec: Spectra over the Orbit Category}.
The desired functor $\bfE$ is given in \cite{SauerJ(2002)}.


\typeout{----------------------------  Status of the Conjectures  ----------------------------}

\section{Status of the Conjectures}
\label{chap: Status of the Conjectures}

In this section, we give the status, at the time of writing, of some of the conjectures mentioned earlier.
We begin with the Baum-Connes Conjecture.


\subsection{Status of the Baum-Connes-Conjecture}
\label{sec: Status of the Baum-Connes-Conjecture}

\subsubsection{The Baum-Connes Conjecture with Coefficients}
\label{subsec: Status of the Baum-Connes Conjecture with Coefficients}

We begin with the Baum-Connes Conjecture with Coefficients~\ref{con: Baum-Connes Conjecture with coefficients}.
It has better inheritance properties than the Baum-Connes Conjecture~\ref{con: Baum-Connes Conjecture}
itself and contains it as a special case.

\begin{theorem}{\bf (Baum-Connes Conjecture with Coefficients and a-T-menable Groups).}
 \label{the: BCC with coefficients for a-Tmenable groups}
\indextheorem{Baum-Connes Conjecture with Coefficients and a-T-menable Groups}
The discrete group $G$ satisfies the Baum-Connes Conjecture with Coefficients 
\ref{con: Baum-Connes Conjecture with coefficients} and is $K$-amenable provided that
$G$ is  a-T-menable.
\end{theorem}

This theorem is
proved in Higson-Kasparov \cite[Theorem~1.1]{Higson-Kasparov(2001)}, where more generally second countable locally
compact topological groups are treated (see also \cite{Julg(1998)}). 

A  group $G$ is \emph{a-T-menable},%
\index{group!a-T-menable}
or, equivalently, has the \emph{Haagerup property}%
\index{group!having the Haagerup property}
if $G$  admits a metrically proper isometric action on some affine Hilbert space. 
\emph{Metrically proper}%
\index{metrically proper} means that for any bounded subset $B$ the set
$\{g \in G \mid gB \cap B \not= \emptyset\}$ is finite. An extensive  treatment of
such groups is presented in 
\cite{Cherix-Cowling-Jolissaint-Julg-Valette(2001)}.  Any a-T-menable
group is countable. The class of
a-T-menable groups is closed under taking  subgroups, under extensions with finite quotients and 
under finite products. It is not closed under semi-direct products. Examples of a-T-menable
groups are countable amenable groups, countable free groups, 
discrete subgroups of $SO(n,1)$ and $SU(n,1)$, Coxeter
groups, countable groups acting properly on trees, products of trees, or simply connected CAT(0)
cubical complexes. A group $G$ has Kazhdan's \emph{property (T)}%
\index{group!having property (T)} if, whenever it acts isometrically on some affine Hilbert
space, it has a fixed point. An infinite a-T-menable group does not have property (T).
Since $SL(n,\IZ)$ for $n \ge 3$ has property (T), it cannot be a-T-menable.

Using the Higson-Kasparov result Theorem~\ref{the: BCC with coefficients for a-Tmenable groups}
and known inheritance properties of the Baum-Connes Conjecture with Coefficients (compare 
Section~\ref{sec: Inheritance Properties} and \cite{Oyono-Oyono(2001)},\cite{Oyono-Oyono(2001b)})
Mislin describes an even  larger  class of groups for which the 
conjecture is known \cite[Theorem~5.23]{Mislin-Valette(2003)}.
\begin{theorem}[The Baum-Connes Conjecture with Coefficients and the Class of Groups 
${\mathbf L}{\mathbf H}{\mathcal E}{\mathcal T}{\mathcal H}$]
 \label{the: BCC with coefficients for the class LHETH}
\indextheorem{The Baum-Connes Conjecture with Coefficients and the Class of Groups 
${\mathbf L}{\mathbf H}{\mathcal E}{\mathcal T}{\mathcal H}$}
The discrete group $G$ satisfies the Baum-Connes Conjecture with Coefficients 
\ref{con: Baum-Connes Conjecture with coefficients} provided that 
$G$  belongs to the class ${\mathbf L}{\mathbf H}{\mathcal E}{\mathcal T}{\mathcal H}$.
\end{theorem}

The class 
${\mathbf L}{\mathbf H}{\mathcal E}{\mathcal T}{\mathcal H}$ is defined as follows.
Let ${\mathbf H}{\mathcal T}{\mathcal H}$ be the smallest class of groups which contains
all a-T-menable groups and contains a group $G$ if
there is a $1$-dimensional contractible $G$-$CW$-complex whose
stabilizers belong already to ${\mathbf H}{\mathcal T}{\mathcal H}$. 
Let ${\mathbf H}{\mathcal E}{\mathcal T}{\mathcal H}$ be the smallest class of groups
containing ${\mathbf H}{\mathcal T}{\mathcal H}$ and containing a group $G$ if either $G$
is countable and admits a surjective map $p\colon G \to Q$ with $Q$ and $p^{-1}(F)$ in
${\mathbf H}{\mathcal E}{\mathcal T}{\mathcal H}$ for every finite subgroup $F \subseteq
Q$ or if $G$ admits a  $1$-dimensional contractible $G$-$CW$-complex whose
stabilizers belong already to ${\mathbf H}{\mathcal E}{\mathcal T}{\mathcal H}$.
Let ${\mathbf L}{\mathbf H}{\mathcal E}{\mathcal T}{\mathcal H}$ be the
class of groups $G$ whose finitely generated subgroups belong to 
${\mathbf H}{\mathcal E}{\mathcal T}{\mathcal H}$.

The class ${\mathbf L}{\mathbf H}{\mathcal E}{\mathcal T}{\mathcal H}$
is closed under passing to
subgroups, under extensions with torsion free quotients and under finite products.
It contains in particular one-relator groups and Haken $3$-manifold
groups (and hence all knot groups). All these facts of the class
${\mathbf L}{\mathbf H}{\mathcal E}{\mathcal T}{\mathcal H}$ and more
information can be found in \cite{Mislin-Valette(2003)}.

Vincent Lafforgue has an unpublished proof of the 
Baum-Connes Conjecture with Coefficients~\ref{con: Baum-Connes Conjecture with coefficients}
for word-hyperbolic groups.


\begin{remarknew} \label{rem: status of BC with coefficients} 
There are counterexamples to the Baum-Connes Conjecture
with (commutative) Coefficients~\ref{con: Baum-Connes Conjecture with coefficients} as soon
as the existence of finitely generated groups containing arbitrary
large expanders in their Cayley graph is shown
\cite[Section 7]{Higson-Lafforgue-Skandalis(2002)}.
The existence of such groups has been claimed by Gromov
\cite{Gromov(2000)}, \cite{Gromov(2003)}. Details of the construction
are  described  by Ghys in \cite{Ghys(2003)}.
At the time of writing no counterexample to the Baum-Connes Conjecture \ref{con: Baum-Connes Conjecture} (without coefficients)
is known to the authors.
\end{remarknew}

\subsubsection{The Baum-Connes Conjecture}
\label{subsec: Status of the Baum-Connes Conjecture}

Next we deal with the Baum-Connes Conjecture \ref{con: Baum-Connes Conjecture} itself.
Recall that all groups which satisfy the 
Baum-Connes Conjecture
with Coefficients~\ref{con: Baum-Connes Conjecture with coefficients} do  in particular
satisfy the Baum-Connes Conjecture \ref{con: Baum-Connes Conjecture}.

\begin{theorem}[Status of the Baum-Connes Conjecture]
 \label{the: BCC} 
\indextheorem{Status of the Baum-Connes Conjecture}
A group $G$ satisfies the 
Baum-Connes Conjecture \ref{con: Baum-Connes Conjecture} if it satisfies one of the 
following conditions.

\begin{enumerate}
\item \label{the: BCC: Levi-Macev}
It is a discrete subgroup of a connected Lie groups $L$, whose Levi-Malcev decomposition $L=RS$
into the radical $R$ and semisimple part $S$ is such that $S$ is locally of the form
$$S = K \times SO(n_1,1) \times \ldots \times SO(n_k,1) \times
SU(m_1,1) \times \ldots \times SU(m_l,1)$$
for a compact group $K$.

\item \label{the: BCC: RD and bolic}
The group $G$ has property (RD) and admits a proper isometric action on 
a strongly bolic weakly geodesic uniformly locally finite metric
space.

\item \label{the: BCC: subgroups of word hyperbolic groups}
$G$ is a subgroup of a word hyperbolic group.

\item \label{the. BCC: subgroups of Sp(n,1)} $G$ is a discrete subgroup of
  $Sp(n,1)$.

\end{enumerate}
\end{theorem}
\begin{proof}
The proof under condition \ref{the: BCC: Levi-Macev} is due to Julg-Kasparov
\cite{Julg-Kasparov(1995)}.
The proof under condition \ref{the: BCC: RD and bolic} is due to Lafforgue
\cite{Lafforgue(1998)} (see also \cite{Skandalis(1999)}).  
Word hyperbolic groups have property (RD)
\cite{delaHarpe(1988)}.
Any subgroup of a word hyperbolic
group satisfies the conditions appearing in the result of
Lafforgue and hence satisfies the Baum-Connes Conjecture
\ref{con: Baum-Connes Conjecture} \cite[Theorem~20]{Mineyev-Yu(2002)}.
The proof under condition \ref{the. BCC: subgroups of Sp(n,1)} is due
to Julg \cite{Julg(2002)}. 
\end{proof}

Lafforgue's result about groups satisfying condition \ref{the: BCC: RD and bolic}
yielded the first examples of infinite groups which have 
Kazhdan's property (T) and satisfy the
Baum-Connes Conjecture \ref{con: Baum-Connes Conjecture}. Here are some explanations about
condition \ref{the: BCC: RD and bolic}.

A \emph{length function}%
\index{length function}
on $G$ is a function $L\colon G \to \IR_{\ge 0}$ such that $L(1) = 0$,
$L(g) = L(g^{-1})$  for $g \in G$ and $L(g_1g_2) \le L(g_1)+ L(g_2)$ for 
$g_1,g_2 \in G$ holds. The word length metric $L_S$ associated to a finite set $S$ of
generators is an example. 
A length function $L$ on $G$ has \emph{property (RD)}%
\indexnotation{length function! having property (RD)}
(``rapid decay'')
if there exist $C,s > 0$ such that for any 
$u = \sum_{g \in G} \lambda_g \cdot g  \in \IC G$ we have
$$||\rho_G(u)||_{\infty} ~ \le ~ 
C \cdot \left(\sum_{g \in G} |\lambda_g|^2 \cdot (1 + L(g))^{2s}\right)^{1/2},$$
where $||\rho_G(u)||_{\infty}$ is the operator norm of the bounded $G$-equivariant
operator $l^2(G) \to l^2(G)$ coming from right multiplication with $u$.
A group $G$ has \emph{property (RD)}%
\index{group!having property (RD)}
if there is a length function which has property (RD).
More information about property (RD) can be found for instance in 
\cite{Chatterji-Ruane(2003)},
\cite{Lafforgue(2000a)} and \cite[Chapter 8]{Valette(2002)}. Bolicity generalizes Gromov's notion of hyperbolicity for metric spaces.
We refer to \cite{Kasparov-Skandalis(1994)} for a precise definition.


\begin{remarknew} We do not know whether all groups appearing in Theorem~\ref{the: BCC} 
satisfy also the Baum-Connes Conjecture 
with Coefficients~\ref{con: Baum-Connes Conjecture with coefficients}. 
\end{remarknew}


\begin{remarknew}
\label{Baum-Connes for SL(n,Z)} 
It is not known at the time of writing 
whether the Baum-Connes Conjecture is true for $SL(n,\IZ)$ for $n \ge 3$.
\end{remarknew}


\begin{remarknew}[The Status for Topological Groups]
We only dealt with the Baum-Connes Conjecture for discrete
groups. We already mentioned that Higson-Kasparov 
\cite{Higson-Kasparov(2001)} treat  second countable locally
compact topological groups. The Baum-Connes Conjecture for second countable almost
connected groups $G$ has been proven by Chabert-Echterhoff-Nest 
\cite{Chabert-Echterhoff-Nest(2003)} based on the work of
Higson-Kasparov \cite{Higson-Kasparov(2001)} and Lafforgue
\cite{Lafforgue(2002)}. The Baum-Connes Conjecture with 
Coefficients~\ref{con: Baum-Connes Conjecture with coefficients} has been
proved for the connected Lie groups $L$ appearing in 
Theorem~\ref{the: BCC}~\ref{the: BCC: Levi-Macev} by \cite{Julg-Kasparov(1995)}
and for $Sp(n,1)$ by Julg
\cite{Julg(2002)}. 
\end{remarknew}

\subsubsection{The Injectivity Part of the Baum-Connes Conjecture}
\label{subsec: Status of the Injectivity Part of the Baum-Connes Conjecture}

In this subsection we deal with injectivity results about the assembly map appearing in the
Baum-Connes Conjecture \ref{con: Baum-Connes Conjecture}. 
Recall that rational injectivity 
already implies the 
Novikov Conjecture~\ref{con: Novikov Conjecture} (see
Proposition~\ref{BCC implies Novikov})  and the Homological 
Stable Gromov-Lawson-Rosenberg Conjecture~\ref{con: Homological Gromov-Lawson-Rosenberg Conjecture}
(see Proposition~\ref{pro: SNC for K_*(C^*_r(G)) implies HGLR} and \ref{pro: Rational fcyc to fin}).

\begin{theorem}[Rational Injectivity of the Baum-Connes Assembly Map]
\label{the: BCC rational injectivity}
\indextheorem{Rational Injectivity of the Baum-Connes Assembly Map}
The assembly map appearing in the Baum-Connes Conjecture 
\ref{con: Baum-Connes Conjecture}  is rationally injective if $G$ belongs to one of the classes of groups below.

\begin{enumerate}

\item \label{the: BCC injectivity: action on non-positively curved manifold} 
Groups acting properly isometrically on complete manifolds
with non-positive sectional curvature.

\item \label{the: BCC injectivity: subgroups of Lie groups}
Discrete subgroups of  Lie groups with finitely many path components.

\item \label{the: BCC injectivity: p-adic groups}
Discrete subgroups of $p$-adic groups.

\end{enumerate} 
\end{theorem}
\begin{proof}
The proof of assertions 
\ref{the: BCC injectivity: action on non-positively curved manifold}  and 
\ref{the: BCC injectivity: subgroups of Lie groups}  is due to Kasparov
\cite{Kasparov(1995)}, the one of assertion \ref{the: BCC injectivity: p-adic groups} to 
Kasparov-Skandalis \cite{Kasparov-Skandalis(1991)}. 
\end{proof}

A metric space $(X,d)$ admits a \emph{uniform embedding into Hilbert space}%
\index{uniform embedding into Hilbert space}
if there exist a separable Hilbert space $H$, a map $f \colon X \to H$ and
non-decreasing functions $\rho_1$ and $\rho_2$ from
$[0,\infty) \to \IR$ such that 
$\rho_1(d(x,y)) \le ||f(x) - f(y)|| \le \rho_2(d(x,y))$ for $x,y
\in X$ and $\lim_{r \to \infty} \rho_i(r) = \infty$ for $i = 1,2$. 
A metric is \emph{proper}%
\index{proper!metric space}
if for each $r > 0$ and $x \in X$ the closed ball of radius $r$ centered at $x$ is
compact. The question whether a discrete group $G$ equipped with a proper left $G$-invariant metric $d$
admits a uniform embedding into Hilbert space is independent of the choice of $d$,
since the induced coarse structure does not depend 
on $d$ \cite[page 808]{Skandalis-Tu-Yu(2002)}. For more information
about groups admitting a uniform embedding into Hilbert space we refer to
\cite{Dranishnikov-Gong-Lafforgue-Yu(2002)},
\cite{Guentner-Higson-Weinberger(2003a)}.

The class of finitely generated groups, which embed uniformly into Hilbert space,
contains a subclass $A$, which contains all word hyperbolic groups,
finitely generated discrete subgroups of connected Lie groups and
finitely generated amenable groups and is closed under semi-direct
products  \cite[Definition~2.1, Theorem~2.2 and Proposition~2.6]{Yu(2000)}. 
Gromov \cite{Gromov(2000)},  \cite{Gromov(2003)}
has announced examples of finitely generated groups which do
not admit a uniform embedding into Hilbert space. 
Details of the construction are  described  in Ghys \cite{Ghys(2003)}.

The next theorem is proved by Skandalis-Tu-Yu~\cite[Theorem~6.1]{Skandalis-Tu-Yu(2002)} using 
ideas of Higson~\cite{Higson(2000)}.

\begin{theorem}[Injectivity of the Baum-Connes Assembly Map]
\label{the: BCC injectivity}
\indextheorem{Injectivity of the Baum-Connes Assembly Map}
Let $G$ be a countable group. Suppose that $G$ admits a $G$-invariant metric 
for which $G$ admits a uniform embedding into Hilbert space. Then
the assembly map appearing in the Baum-Connes Conjecture with 
Coefficients~\ref{con: Baum-Connes Conjecture with coefficients}  is injective.
\end{theorem}

We now discuss conditions which can be used to verify the assumption in Theorem~\ref{the: BCC injectivity}.


\begin{remarknew}[Linear Groups] \label{rem: Linear groups} 
A group $G$ is called
\emph{linear}%
\index{group!linear}
if it is a subgroup of $GL_n(F)$ for some $n$ and some field $F$.
Guentner-Higson-Weinberger \cite{Guentner-Higson-Weinberger(2003a)}  show
that every countable linear group admits a uniform embedding into Hilbert space
and hence Theorem~\ref{the: BCC injectivity} applies.
\end{remarknew}


\begin{remarknew}[Groups Acting Amenably on a Compact Space]
\label{Groups acting amenable on a compact space}
A continuous action of a discrete group $G$ on a compact space $X$ is called
\emph{amenable}
\index{amenable group action}
if there exists a sequence 
$$p_n \colon X \to M^1(G) = \{f \colon G \to [0,1] \mid \sum_{g \in G} f(g) = 1\}$$
of weak-$\ast$-continuous maps such that for each $g \in G$ one has
$$\lim_{n \to \infty} \sup_{x \in X} ||g \ast (p_n(x) - p_n(g \cdot x))||_1 ~ = ~ 0.$$
Note that a group $G$ is amenable if and only if its action on the
one-point-space is amenable.
More information about this notion can be found for instance in
\cite{Anantharaman-Renault(2001)},
\cite{Anantharaman-Delaroche-Renault(2000)}. 

Higson-Roe \cite[Theorem~1.1 and Proposition~2.3]{Higson-Roe(2000a)}
show that a finitely generated group equipped with its
word length metric admits an amenable action on a compact metric space, if and only if
it belongs to the class $A$  defined in  \cite[Definition 2.1]{Yu(2000)}, and hence admits a uniform embedding into Hilbert space.
Hence  Theorem~\ref{the: BCC injectivity} implies the result of Higson \cite[Theorem~1.1]{Higson(2000)}
that the assembly map appearing in the Baum-Connes Conjecture with 
Coefficients~\ref{con: Baum-Connes Conjecture with coefficients} is injective
if $G$ admits an amenable action on some compact space.

Word hyperbolic groups and the class of groups mentioned in Theorem~
\ref{the: BCC rational injectivity} \ref{the: BCC injectivity: subgroups of Lie groups} 
fall under the class of groups admitting an amenable action on some compact
space \cite[Section 4]{Higson-Roe(2000a)}. 
\end{remarknew}



\begin{remarknew} \label{rem: Higson's version of Carlsson-Pedersen}
Higson \cite[Theorem~5.2]{Higson(2000)}  shows that the 
assembly map appearing in the Baum-Connes Conjecture with 
Coefficients~\ref{con: Baum-Connes Conjecture with coefficients} is injective if 
$\underline{E}G$ admits an appropriate compactification. This is
a $C^*$-version of the result for $K$-and $L$-theory due to Carlsson-Pedersen
\cite{Carlsson-Pedersen(1995a)}, compare Theorem~\ref{the: Carlsson-Pedersen Rosenthal}.
\end{remarknew}


\begin{remarknew} We do not know whether the groups appearing in Theorem~
\ref{the: BCC rational injectivity} and \ref{the: BCC injectivity} satisfy
the Baum-Connes Conjecture \ref{con: Baum-Connes Conjecture}.
\end{remarknew}


Next we discuss injectivity results about the classical assembly map for topological $K$-theory.

The \emph{asymptotic dimension}%
\index{asymptotic dimension}
of a proper metric space $X$ is the infimum over all integers $n$ such that
for any $R > 0$ there exists a cover $\calu$ of $X$ with the property
that the diameter of the members of $\calu$ is uniformly bounded and
every ball of radius $R$ intersects at most $(n+1)$ elements of
$\calu$ (see \cite[page 28]{Gromov(1993)}). 

The next result is due to Yu \cite{Yu(1998a)}.

\begin{theorem} [The $C^*$-Theoretic Novikov Conjecture and Groups of Finite
  Asymptotic Dimension]
\indextheorem{The $C^*$-Theoretic Novikov Conjecture and Groups of Finite
  Asymptotic Dimension}
\label{the: Yu injectivity and asymtotic dimension}
Let $G$ be a group which possesses a finite model for $BG$ and has
finite asymptotic dimension. Then the assembly map in the Baum-Connes
Conjecture~\ref{con: BC torsion free}
\[
K_n(BG) \to K_n(C^*_r(G))
\]
is injective for all $n \in \IZ$.
\end{theorem}


\subsubsection{The Coarse Baum-Connes Conjecture}
\label{subsec: Status of the Coarse Baum-Conjecture}

The coarse Baum-Connes Conjecture was explained in Section~\ref{subsec: The Coarse Baum Connes Conjecture}. 
Recall the descent principle (Proposition~\ref{pro:descentCBC}): 
if a countable group can be equipped with a $G$-invariant metric such that the resulting metric space 
satisfies the Coarse Baum-Connes Conjecture, then the classical assembly map for topological $K$-theory is injective.

Recall that a discrete metric space has bounded geometry
if for each  $r > 0$ there exists a
$N(r)$ such that for all $x$ the ball of radius $N(r)$ centered at  $x \in X$ 
contains at most $N(r)$ elements. 

The next result is due to Yu \cite[Theorem~2.2 and Proposition~2.6]{Yu(2000)}.
\begin{theorem}[Status of the Coarse Baum-Connes Conjecture]
 \label{the: Status of the Coarse Baum-Connes Conjecture}
\indextheorem{Status of the Coarse Baum-Connes Conjecture}
The Coarse Baum-Connes Conjecture~\ref{con:CBC}
is true for a discrete metric space $X$ of bounded geometry if $X$ 
admits a uniform embedding into Hilbert space. In particular a countable group
$G$ satisfies the Coarse Baum-Connes Conjecture~\ref{con:CBC}
if $G$ equipped with a proper left
$G$-invariant metric admits a uniform embedding into Hilbert space.
\end{theorem}

Also Yu's  Theorem~\ref{the: Yu injectivity and asymtotic dimension}
is proven via a corresponding result about the Coarse Baum-Connes Conjecture.

\subsection{Status of the Farrell-Jones Conjecture}
\label{sec: Status of the Farrell-Jones Conjecture}

Next we deal with the Farrell-Jones Conjecture.

\subsubsection{The Fibered Farrell-Jones Conjecture}
\label{subsec: Status of the Fibered Farrell-Jones Conjecture}

The Fibered Farrell-Jones Conjecture~\ref{con: Fibered Farrell-Jones Conjecture} was discussed in
Subsection~\ref{subsec: Fibered Version of the Farrell-Jones Conjecture}. Recall that it has better inheritance 
properties (compare Section~\ref{sec: Inheritance Properties}) 
and contains the ordinary Farrell-Jones Conjecture~\ref{con: Farrell-Jones Conjecture} as a special case.

\begin{theorem}[Status of the Fibered Farrell-Jones Conjecture]
 \label{the: fibered Farrell-Jones and subgroups of Lie groups}
\indextheorem{Status of the Fibered Farrell-Jones Conjecture}\
 
\begin{enumerate}
\item
Let $G$ be a discrete group which satisfies one of the following conditions.
\begin{enumerate}
\item[(a)]
There is a Lie group $L$ with finitely many path components and $G$
is a cocompact discrete subgroup of $L$.
\item[(b)]
The group $G$ is virtually torsionfree and  acts properly discontinuously, cocompactly and via isometries
on a simply connected complete nonpositively curved Riemannian manifold.
\end{enumerate}
Then 
\begin{enumerate} 
\item[(1)] 
The version of the Fibered Farrell-Jones Conjecture~\ref{con: Fibered Farrell-Jones Conjecture}
for the topological and the differentiable pseudoisotopy functor is true for $G$.
\item[(2)]  
The version of the Fibered Farrell-Jones 
Conjecture~\ref{con: Fibered Farrell-Jones Conjecture} for $K$-theory and $R = \IZ$ is true for $G$
in the range $n \le 1$, i.e.\ the assembly map is bijective for $n \le 1$. 
\end{enumerate}
\end{enumerate}
Moreover we have the following statements.
\begin{enumerate}
\item[(ii)] 
The version of the Fibered Farrell-Jones 
Conjecture~\ref{con: Fibered Farrell-Jones Conjecture} for $K$-theory and $R = \IZ$ is true 
in the range $n \le 1$ for braid groups.
\item[(iii)] 
The $L$-theoretic version of the Fibered Farrell-Jones Conjecture~\ref{con: Fibered Farrell-Jones Conjecture} with $R=\IZ$ 
holds after inverting $2$ for elementary amenable groups.
\end{enumerate}
\end{theorem}
\begin{proof}
(i) For assertion (1) see~\cite[Theorem~2.1 on page 263]{Farrell-Jones(1993a)}, 
\cite[Proposition~2.3]{Farrell-Jones(1993a)} and \cite[Theorem~A]{Farrell-Roushon(2000)}.
Assertion~(2) follows from (1) by 
Remark~\ref{rem: relating K and P fibered}. \\[1mm]
(ii)
See \cite{Farrell-Roushon(2000)}. \\[1mm]
(iii) is proven in
\cite[Theorem~5.2]{Farrell-Linnell(2003a)}. For crystallographic groups see also \cite{Yamasaki(1987)}.
\end{proof}

A surjectivity result about the Fibered Farrell-Jones Conjecture for 
Pseudoisotopies appears as the last statement in Theorem~\ref{the: Vanishing of Wh(G times IZ^k) and FJC for L-th.}.

The rational comparison result between the $K$-theory and the pseudoisotopy version 
(see Proposition~\ref{pro: relating K and pseudoisotopies})
does not work in the fibered case, compare Remark~\ref{rem: relating K and P fibered}.
However, in order to exploit the good inheritance properties one can
first use the pseudoisotopy functor in the fibered set-up, then specialize to the unfibered situation
and finally do the rational comparison to $K$-theory.

\begin{remarknew} \label{rem: fibered Farrell-Jones and subgroups of Lie groups for L-theory}
The version of the Fibered Farrell-Jones Conjecture \ref{con: Fibered Farrell-Jones
Conjecture} for $L$-theory and $R = \IZ$  seems to be true if $G$ satisfies the condition (a) appearing in 
Theorem~\ref{the: fibered Farrell-Jones and subgroups of Lie groups}.
Farrell and Jones \cite[Remark 2.1.3 on page 263]{Farrell-Jones(1993a)} say that
they can also prove this  version without giving the details. 
\end{remarknew}

\begin{remarknew}
\label{rem: Wilking}
Let $G$ be a virtually poly-cyclic group. Then it contains a maximal normal finite subgroup $N$ such that
the quotient $G/N$ is a discrete cocompact subgroup of a Lie group with finitely many path components 
\cite[Theorem~3, Remark~4 on page~200]{Wilking(2000a)}.
Hence by Subsection~\ref{subsec: Extensions of Groups} and Theorem~\ref{the: fibered Farrell-Jones and subgroups of Lie groups}
the version of the 
Fibered Farrell-Jones Conjecture~\ref{con: Fibered Farrell-Jones Conjecture}
for the topological and the differentiable pseudoisotopy functor, and 
for $K$-theory and $R = \IZ$ 
in the range $n \le 1$, is true for $G$. Earlier results of this type were treated for example in 
\cite{Farrell-Hsiang(1981b)}, \cite{Farrell-Jones(1988b)}.
\end{remarknew}

\subsubsection{The Farrell-Jones Conjecture}
\label{subsec: Status of the Farrell-Jones Conjecture}

Here is a sample of some results one can deduce from Theorem~\ref{the: fibered Farrell-Jones and subgroups of Lie groups}.

\begin{theorem}[The Farrell-Jones Conjecture and Subgroups of Lie groups]
\label{the: Farrell-Jones and subgroups of Lie groups}
\indextheorem{The Farrell-Jones Conjecture and subgroups of Lie groups}
Suppose $H$ is a subgroup of $G$, where $G$ is a discrete cocompact  
subgroup of a Lie group $L$ with finitely many path components.
Then
\begin{enumerate}
\item  \label{the: Farrell-Jones and subgroups of Lie groups: K rationally}
The version of the Farrell-Jones Conjecture 
for $K$-theory and $R = \IZ$ is true for $H$ rationally, 
i.e.\ the assembly map appearing in Conjecture~\ref{con: Farrell-Jones Conjecture} 
is an isomorphism after applying $- \otimes_{\IZ} \IQ$.
\item  
\label{the: Farrell-Jones and subgroups of Lie groups: K integrally for n le 1}
The version of the Farrell-Jones Conjecture 
for $K$-theory and $R = \IZ$ is true for $H$
in the range $n \le 1$, i.e. the assembly map appearing in 
Conjecture~\ref{con: Farrell-Jones Conjecture} 
is an isomorphism for $n \le 1$. 
\end{enumerate}
\end{theorem}
\begin{proof}
The results follow from
Theorem~\ref{the: fibered Farrell-Jones and subgroups of Lie groups}, since the 
Fibered Farrell-Jones Conjecture~\ref{con: Fibered Farrell-Jones  Conjecture} 
passes to subgroups \cite[Theorem~A.8 on page 289]{Farrell-Jones(1993a)} (compare Section~\ref{subsec: Passing to Subgroups})
and implies the 
Farrell-Jones Conjecture~\ref{con: Farrell-Jones Conjecture}.
\end{proof}

We now discuss results for torsion free groups. Recall that for $R=\IZ$ 
the $K$-theoretic Farrell-Jones Conjecture in dimensions $\leq 1$ together with the 
$L$-theoretic version implies already the Borel Conjecture~\ref{con: Borel Conjecture}
in dimension $\ge 5$ (see Theorem~\ref{the: The Farrell-Jones Conjecture Implies the Borel Conjecture}).

A complete Riemannian manifold $M$ is called \emph{A-regular}%
\index{A-regular}
if there exists a sequence of positive real numbers $A_0$, $A_1$, $A_2$, $\ldots$ such
that $||\nabla^n K|| \le A_n$, where $||\nabla^nK||$ is the supremum-norm of the $n$-th covariant
derivative of the curvature tensor $K$. Every locally symmetric space is 
A-regular since $\nabla K$ is identically zero. Obviously every closed Riemannian manifold is $A$-regular.

\begin{theorem}[Status of the Farrell-Jones Conjecture for Torsionfree Groups]
 \label{the: Vanishing of Wh(G times IZ^k) and FJC for L-th.}
\indextheorem{Status of the Farrell-Jones Conjecture for Torsionfree Groups}
Consider the following conditions for the group $G$.
\begin{enumerate} 
\item \label{the: Vanishing of Wh(G times IZ^k) and FJC for L-th.: non-positively curved A-regular} 
$G = \pi_1(M)$ for a 
complete Riemannian manifold $M$ with non-positive sectional curvature which is A-regular.

\item \label{the: Vanishing of Wh(G times IZ^k) and FJC for L-th.: non-positively curved closed}
$G = \pi_1(M)$ for a 
closed Riemannian manifold $M$ with non-positive  sectional curvature.

\item \label{the: Vanishing of Wh(G times IZ^k) and FJC for L-th.: pinched negatively curved}
$G = \pi_1(M)$ for a complete Riemannian manifold with negatively pinched 
sectional curvature.

\item \label{the: Vanishing of Wh(G times IZ^k) and FJC for L-th.: subgroup of GL_n(bbR)} 
$G$ is a torsion free discrete subgroup of $GL(n,\IR)$.

\item \label{the: Vanishing of Wh(G times IZ^k) and FJC for L-th.: solvable subgroup of GL_n(bbC)} 
$G$ is a torsion free solvable discrete subgroup of $GL(n,\IC)$.

\item \label{the: Vanishing of Wh(G times IZ^k) and FJC for L-th.: negatively curved complex}
$G = \pi_1(X)$ for a non-positively curved finite simplicial complex $X$.

\item \label{the: Vanishing of Wh(G times IZ^k) and FJC for L-th.: strongly poly-free group}
$G$ is a strongly poly-free group in the sense of Aravinda-Farrell-Roushon 
\cite[Definition 1.1]{Aravinda-Farrell-Roushon(2000)}.
The pure braid group satisfies this hypothesis.
\end{enumerate}

Then 

\begin{enumerate} \label{results about FJ torsionfree under conditions}

\item[(1)] 
Suppose that $G$ satisfies one of the conditions 
\ref{the: Vanishing of Wh(G times IZ^k) and FJC for L-th.: non-positively curved A-regular}  to
\ref{the: Vanishing of Wh(G times IZ^k) and FJC for L-th.: strongly poly-free group}. Then the $K$-theoretic Farrell-Jones 
Conjecture is true for $R = \IZ$ in dimensions $n \leq 1$. In particular Conjecture~\ref{con: vanishing of lower K} holds for $G$.

\item[(2)]
Suppose that $G$ satisfies one of the conditions 
\ref{the: Vanishing of Wh(G times IZ^k) and FJC for L-th.: non-positively curved A-regular},
\ref{the: Vanishing of Wh(G times IZ^k) and FJC for L-th.: non-positively curved closed},  
\ref{the: Vanishing of Wh(G times IZ^k) and FJC for L-th.: pinched negatively curved} or
\ref{the: Vanishing of Wh(G times IZ^k) and FJC for L-th.: subgroup of GL_n(bbR)}. Then
$G$ satisfies the Farrell-Jones Conjecture for Torsion Free Groups and
L-Theory \ref{con: FJL torsion free} for $R = \IZ$.

\item[(3)] 
Suppose that $G$ satisfies 
\ref{the: Vanishing of Wh(G times IZ^k) and FJC for L-th.: non-positively curved closed}. Then
the Farrell-Jones Conjecture for Pseudoisotopies of Aspherical Spaces~\ref{con: FJ for pseudoisotopies unfibered}
holds for $G$.

\item[(4)] 
\label{jones result}
Suppose that $G$ satisfies one of the conditions
\ref{the: Vanishing of Wh(G times IZ^k) and FJC for L-th.: non-positively curved A-regular},
\ref{the: Vanishing of Wh(G times IZ^k) and FJC for L-th.: pinched negatively curved} or
\ref{the: Vanishing of Wh(G times IZ^k) and FJC for L-th.: subgroup of GL_n(bbR)}. Then
the assembly map appearing in the version of the Fibered Farrell-Jones Conjecture for 
Pseudoisotopies~\ref{con: Fibered Farrell-Jones Conjecture} is surjective, provided
that the $G$-space $Z$ appearing in Conjecture~\ref{con: Fibered Farrell-Jones Conjecture} is 
connected.

\end{enumerate}
 
\end{theorem}
\begin{proof}
Note that \ref{the: Vanishing of Wh(G times IZ^k) and FJC for L-th.: non-positively curved closed} 
is a special case of \ref{the: Vanishing of Wh(G times IZ^k) and FJC for L-th.: non-positively curved A-regular}
because every closed Riemannian manifold is A-regular. 
If $M$ is a pinched negatively
curved complete Riemannian manifold, then there is another Riemannian metric
for which $M$ is negatively curved complete and A-regular.
This fact is mentioned in \cite[page 216]{Farrell-Jones(1998)} and attributed there to
Abresch \cite{Abresch(1988)} and Shi \cite{Shi(1989)}. 
Hence also \ref{the: Vanishing of Wh(G times IZ^k) and FJC for L-th.: pinched negatively curved} can be considered as a special case of 
\ref{the: Vanishing of Wh(G times IZ^k) and FJC for L-th.: non-positively curved A-regular}.
The manifold $M = G\backslash GL(n , \IR)/O(n)$ is a non-positively curved complete locally symmetric space
and hence in particular A-regular.  So
\ref{the: Vanishing of Wh(G times IZ^k) and FJC for L-th.: subgroup of GL_n(bbR)} is a special case of
\ref{the: Vanishing of Wh(G times IZ^k) and FJC for L-th.: non-positively curved A-regular}.

Assertion~(1) under the assumption 
\ref{the: Vanishing of Wh(G times IZ^k) and FJC for L-th.: non-positively curved A-regular} is proved by Farrell-Jones  in
\cite[Proposition~0.10 and Lemma~0.12]{Farrell-Jones(1998)}. The earlier work
\cite{Farrell-Jones(1991)} treated the case \ref{the: Vanishing of Wh(G times IZ^k) and FJC for L-th.: non-positively curved closed}.
Under assumption~\ref{the: Vanishing of Wh(G times IZ^k) and FJC for L-th.: solvable subgroup of GL_n(bbC)} 
assertion~(1) is proven by Farrell-Linnell \cite[Theorem~1.1]{Farrell-Linnell(2003a)}.
The result under assumption~\ref{the: Vanishing of Wh(G times IZ^k) and FJC for L-th.: negatively curved complex}
is proved by Hu \cite{Hu(1993)},
under 
assumption~\ref{the: Vanishing of Wh(G times IZ^k) and FJC for L-th.: strongly poly-free group} it is proved 
by Aravinda-Farrell-Roushon~\cite[Theorem~1.3]{Aravinda-Farrell-Roushon(2000)}.

Assertion~(2) under assumption~\ref{the: Vanishing of Wh(G times IZ^k) and FJC for L-th.: non-positively curved A-regular}
is proven by Farrell-Jones in \cite{Farrell-Jones(1998)}. The case 
\ref{the: Vanishing of Wh(G times IZ^k) and FJC for L-th.: non-positively curved closed}
was treated earlier in \cite{Farrell-Jones(1993c)}.

Assertion~(3) is proven by Farrell-Jones in \cite{Farrell-Jones(1991)} and 
assertion~(4) by Jones in \cite{Jones(2003b)}.
\end{proof}

\begin{remarknew}
\label{rem: P and SLnIZ}
As soon as certain collapsing results (compare \cite{Farrell-Jones(1998f)},
\cite {Farrell-Jones(2003)})
are extended to orbifolds, the results
under (4) above would also apply to groups with torsion and in particular to $SL_n ( \IZ )$ for arbitrary $n$.
\end{remarknew}


\subsubsection{The Farrell-Jones Conjecture for Arbitrary Coefficients}
\label{subsec: Status of the Farrell-Jones Conjecture for $K$-theory and arbitrary coefficients}

The following result due to Bartels-Reich~\cite{Bartels-Reich(2003)} deals with algebraic $K$-theory 
for arbitrary coefficient rings $R$. It extends 
Bartels-Farrell-Jones-Reich \cite{Bartels-Farrell-Jones-Reich(2004a)}.

\begin{theorem}
\label{the: Bartels-Reich}
Suppose that $G$ is the fundamental group of a closed Riemannian
 manifold with negative sectional curvature. Then the $K$-theoretic part of the 
Farrell-Jones Conjecture \ref{con: Farrell-Jones Conjecture} is true
for any ring $R$, i.e.\ the assembly map
\[
A_{\calvcyc} \colon H_n^G( \EGF{G}{\calvcyc} ; \bfK_R ) ~ \to ~ K_n( RG ) 
\]
is an isomorphism for all $n \in \bfZ$.
\end{theorem}

Note that the assumption implies that $G$ is torsion free and hence the family $\calvcyc$
reduces to the family $\calcyc$ of cyclic subgroups.
Recall that for a regular ring $R$ the theorem above implies that the classical
assembly 
\[
A \colon H_n(BG;\bfK(R)) \to K_n(RG)
\]
is an isomorphism, compare Proposition~\ref{pro: relative assembly torsion free}~\ref{rel ass tfree i}.

\subsubsection{Injectivity Part of the Farrell-Jones Conjecture}
\label{subsec: Status of the Injectivity Part of the Farrell-Jones Conjecture}

The next result about the classical $K$-theoretic assembly map 
is due to B\"okstedt-Hsiang-Madsen \cite{Boekstedt-Hsiang-Madsen(1993)}.

\begin{theorem}[Rational Injectivity of the Classical $K$-Theoretic Assembly Map]
\label{the: Boekstedt-Hsiang-Madsen}
\indextheorem{Rational Injectivity of the Classical $K$-Theoretic Assembly Map}
Let $G$ be a group such that the integral homology $H_j(BG; \IZ )$ is finitely generated for
each $j \in \IZ$. Then the rationalized assembly map 
\[
A \colon H_n(BG;\bfK (\IZ)) \otimes_{\IZ} \IQ \cong  H_n^G(\EGF{G}{\caltr};\bfK_{\IZ} ) \otimes_{\IZ} \IQ
\to K_n(\IZ G) \otimes_{\IZ} \IQ
\]
is injective for all $n \in \bfZ$.
\end{theorem}
Because of the homological Chern character (see Remark~\ref{rem: rational computation K})
we obtain for the groups treated in Theorem~\ref{the: Boekstedt-Hsiang-Madsen} 
an injection
\begin{eqnarray}
\bigoplus_{s+t = n} H_s(BG;\IQ) \otimes_{\IQ} (K_t(\IZ) \otimes_{\IZ} \IQ )
& \to  &
K_n(\IZ G) \otimes_{\IZ} \IQ.
\label{Homological version of BHM}
\end{eqnarray}

Next we describe a generalization of Theorem~\ref{the: Boekstedt-Hsiang-Madsen} above from the trivial family
$\caltr$ to the family $\calfin$ of finite subgroups due to L\"uck-Reich-Rognes-Varisco
\cite{Lueck-Reich-Rognes-Varisco(2003)}. 
Let $\bfK_{\IZ}^{\conn} \colon \GROUPOIDS \to \SPECTRA$ be the connective version of 
the functor $\bfK_{\IZ}$ of \eqref{K^{alg}GROUPOIDS}. In particular 
$H_n( G/H ; \bfK^{\conn}_{\IZ} )$ is isomorphic to $K_n(\IZ H)$ for $n \ge 0$ and vanishes 
in negative dimensions.
For a prime $p$ we denote by $\IZ_p$%
\indexnotation{bbZ_p}
the $p$-adic integers.
Let $K_n(R;\IZ_p)$ denote the homotopy groups $\pi_n(\bfK^{\conn} (R) \widehat{_p})$ 
of the $p$-completion of the connective $K$-theory spectrum of the ring $R$.

\begin{theorem}[Rational Injectivity of the Farrell-Jones Assembly Map for Connective $K$-Theory]
 \label{the: LRRM}
\indextheorem{Rational Injectivity of the Farrell-Jones Assembly Map for Connective $K$-Theory}
Suppose that the group $G$ satisfies the following two conditions:

\begin{enumerate}
\item[(H)] \label{the: LRRM: condition (H)}
For each finite cyclic subgroup $C \subseteq G$ and all $j \geq 0$ the integral
homology group $H_j(BZ_GC; \IZ)$ of the centralizer $Z_GC$ of $C$ in $G$ is finitely generated.
\item[(K)] \label{the: LRRM: condition (K)}
There exists a prime $p$ such that for each finite cyclic subgroup $C \subseteq G$ 
and each $j \ge 1$ the map induced by the change of coefficients homomorphism
$$K_j(\IZ C;\IZ_p) \otimes_{\IZ} \IQ \to K_j(\IZ_pC;\IZ_p) \otimes_{\IZ} \IQ$$
is injective. 
\end{enumerate}
Then the rationalized assembly map
\[
A_{\calvcyc} \colon H_n^G(\EGF{G}{\calvcyc};\bfK_{\IZ}^{\conn}) \otimes_{\IZ} \IQ  \to 
K_n(\IZ G) \otimes_{\IZ} \IQ
\]
is an injection for all $n \in \IZ$.
\end{theorem}

\begin{remarknew} \label{rem: LRRJ and Chern characters}
The methods of Chapter
\ref{chap: Computations}
apply also to $\bfK_{\IZ}^{\conn}$ and yield under assumption (H) and (K) an injection 
\begin{multline*}
\bigoplus_{s+t = n, ~ t \ge 0} ~ \bigoplus_{(C) \in (\calfcyc)}
H_s(BZ_GC;\IQ) \otimes_{\IQ [W_GC]} 
\theta_C \cdot K_t(\IZ C) \otimes_{\IZ} \IQ
\\
\xrightarrow{} ~
K_n(\IZ  G) \otimes_{\IZ} \IQ.
\end{multline*}
Notice that in the index set for the direct sum appearing in the source we require
$t \ge 0$. This reflects the fact that the result deals only with 
the connective $K$-theory spectrum. 
If one drops the restriction $t \ge 0$ the 
Farrell-Jones Conjecture~\ref{con: Farrell-Jones Conjecture} predicts that the
map is an isomorphism, compare Subsection~\ref{subsec: split off nil} and
Theorem~\ref{the: Rational Computation of Algebraic K-Theory for Infinite Groups}.
If we restrict the injection to the direct sum
given by $C = 1$, we rediscover the map \eqref{Homological version of BHM} 
whose injectivity follows already from Theorem~\ref{the: Boekstedt-Hsiang-Madsen}.
\end{remarknew}

The condition (K) appearing in Theorem~\ref{the: LRRM} is conjectured to be true for all primes $p$
(compare \cite{Schneider(1979)}, \cite{Soule(1981)} and \cite{Soule(1987b)})
but no proof is known. The weaker version of condition (K), where $C$ is the trivial group
is also needed in Theorem~\ref{the: Boekstedt-Hsiang-Madsen}. But that case is known to be true and hence does not
appear in its formulation. The special case of
condition (K), where $j=1$ is implied by the Leopoldt Conjecture for abelian fields 
(compare \cite[IX, \S~3]{Neukirch-Schmidt-Wingberg(2000)}),
which is known to be true \cite[Theorem~10.3.16]{Neukirch-Schmidt-Wingberg(2000)}.
This leads to the following result.

\begin{theorem}[Rational Contribution of Finite Subgroups to $\Wh(G)$]
\label{the: Rational Contribution of Finite Subgroups to Wh(G)}
\indextheorem{Rational Contribution of Finite Subgroups to $\Wh(G)$}
 \label{the: LRRM for middle and lower K-theory}
Let $G$ be a group. Suppose that for each finite cyclic subgroup 
$C \subseteq G$ and each $j \leq 4$ the integral
homology group $H_j(BZ_GC)$ of the centralizer $Z_GC$ of $C$ in $G$ is finitely generated.
Then the map 
\[
\colim_{H \in \SubGF{G}{\calfin}} \Wh(H)  \otimes_{\IZ} \IQ
~ \to ~ \Wh(G) \otimes_{\IZ} \IQ.
\]
is injective, compare Conjecture~\ref{con: colim over rationalized Whitehead}.
\end{theorem}

The result above should be compared to the result which is  proven using
Fuglede-Kadison determinants in \cite[Section 5]{Lueck-Roerdam(1993)},
\cite[Theorem~9.38 on page 354]{Lueck(2002)}:  for \emph{every} (discrete) group $G$ and 
every finite normal subgroup $H \subseteq G$ the map
$\Wh(H) \otimes_{\IZ G} \IZ \to \Wh(G)$ induced by the inclusion $H
\to G$ is rationally injective. 

The next result is taken from Rosenthal~\cite{Rosenthal(2002)}, where the techniques and results
of Carlsson-Pedersen~\cite{Carlsson-Pedersen(1995a)} are extended from the 
trivial family $\caltr$ to the family of finite subgroups $\calfin$.
\begin{theorem} \label{the: Carlsson-Pedersen Rosenthal}
Suppose there exists a model $E$ for the classifying space $\EGF{G}{\calfin}$ which
admits a metrizable compactification $\overline{E}$ to which the group
action extends. Suppose $\overline{E}^H$ is contractible and $E^H$ is dense in $\overline{E}^H$
for every finite subgroup $H \subset G$. Suppose compact subsets of $E$ become small near $\overline{E} - E$.
Then for every ring $R$ the assembly map
\[
A_{\calfin} \colon H_n^G( \EGF{G}{\calfin} ; \bfK_R ) ~ \to ~ K_n( RG )
\]
is split injective. 
\end{theorem}
A compact subset $K \subset E$ is said to become \emph{small near} $\overline{E} - E$ if for every neighbourhood 
$U \subset \overline{E}$ of a point $x \in \overline{E} - E$
there exists a neighbourhood $V \subset \overline{E}$ such that $g \in G$ and $gK \cap V \neq \emptyset$ implies
$gK \subset U$. 
Presumably there is an analogous result for $L^{\langle - \infty \rangle}$-theory under the assumption that 
$K_{-n} ( RH)$ vanishes for finite subgroups $H$ of $G$ and $n$ large enough. 
This would extend the corresponding result for the family 
$\calf= \trivial$ which appears in Carlsson-Pedersen~\cite{Carlsson-Pedersen(1995a)}.

We finally discuss injectivity results about assembly maps for the trivial family.
The following result is due to Ferry-Weinberger \cite[Corollary 2.3]{Ferry-Weinberger(1991)} 
extending earlier work of Farrell-Hsiang \cite{Farrell-Hsiang(1981)}.
\begin{theorem}
\label{the: Ferry-Weinberger L injective}
Suppose $G= \pi_1( M)$ for a complete Riemannian manifold of non-positive sectional curvature.
Then the $L$-theory assembly map
\[
A \colon H_n ( BG ; \bfL^{\epsilon}_{\IZ} ) \to L_n^{\epsilon} ( \IZ G )
\]
is injective for $\epsilon = h,s$.
\end{theorem}
In fact Ferry-Weinberger also prove a corresponding splitting result for the classical $A$-theory assembly map. 
In \cite{Hu(1995)} Hu shows that a finite complex of non-positive curvature
is a retract of a non-positively curved $PL$-manifold and concludes split injectivity of the classical $L$-theoretic assembly
map for $R= \IZ$.

The next result due to Bartels~\cite{Bartels(2003a)} is the algebraic $K$- and $L$-theory analogue of 
Theorem~\ref{the: Yu injectivity and asymtotic dimension}.

\begin{theorem}[The $K$-and $L$-Theoretic Novikov Conjecture and Groups of Finite
  Asymptotic Dimension]
\indextheorem{The $K$-and $L$-Theoretic Novikov Conjecture and Groups of Finite
  Asymptotic Dimension}
\label{the: Bartels injectivity and asymtotic dimension}
 Let $G$ be a group which admits a
finite model for $BG$. Suppose that $G$ has finite
asymptotic dimension. Then
\begin{enumerate}
\item The assembly maps appearing in the
  Farrell-Jones Conjecture \ref{con: FJK torsion free all} 
  $$A\colon H_n(BG;\bfK (R)) \to K_n(RG)$$
  is injective for all $n \in \IZ$. 
\item   If furthermore $R$ carries an involution and $K_{-j}(R)$ vanishes for sufficiently
  large $j$, then the assembly maps appearing in the
  Farrell-Jones Conjecture \ref{con: FJL torsion free} 
  $$
  A \colon H_n(BG;\bfL^{\langle - \infty \rangle}(R)) \to L_n^{\langle - \infty \rangle}(RG)
  $$
  is injective for all $n \in \IZ$.
\end{enumerate}
\end{theorem}


Further results related to the Farrell-Jones Conjecture 
\ref{con:  Farrell-Jones Conjecture}  
can be found for instance in
\cite{Aravinda-Farrell-Roushon(1997)},
\cite{Berkhove-Farrell-Pineda-Pearson(2000)}.


\subsection{List of Groups Satisfying the Conjecture}
\label{subsec: List of Groups Satisfying the Conjecture}

In the following table we list prominent classes of groups and state whether they are known to satisfy
the Baum-Connes Conjecture \ref{con: Baum-Connes Conjecture}
(with coefficients \ref{con: Baum-Connes Conjecture with coefficients})
or the Farrell-Jones Conjecture  \ref{con: Farrell-Jones
 Conjecture} (fibered  \ref{con: Fibered Farrell-Jones Conjecture}).
 Some of the classes are redundant. A question mark means
that the authors do not know about a corresponding result.
The reader should keep in mind that there may exist results of which the authors are not aware.
\\[2em]
\begin{center}
\begin{tabular}{|| p{30mm}| p{24mm}| p{24mm}| p{24mm}||}
\hline \hline
type of group & 
Baum-Connes Conjecture~\ref{con: Baum-Connes Conjecture} 
(with coefficients \ref{con: Baum-Connes Conjecture with coefficients})
&
Farrell-Jones Conjecture \ref{con: Farrell-Jones Conjecture}
for $K$-theory for $R = \IZ$ (fibered  \ref{con: Fibered Farrell-Jones Conjecture})  
&
Farrell-Jones Conjecture \ref{con: Farrell-Jones Conjecture}
for $L$-theory for $R = \IZ$ (fibered \ref{con: Fibered Farrell-Jones Conjecture})  
\\
\hline\hline
a-T-menable groups  
& true with coefficients (see Theorem~\ref{the: BCC with coefficients for a-Tmenable groups})
& ? 
& injectivity is true after inverting $2$ (see Propositions \ref{pro: L invert 2}
and  \ref{pro: relating L and K}) 
\\
\hline
amenable groups 
& true with coefficients (see Theorem~\ref{the: BCC with coefficients for a-Tmenable groups})
& ? 
& injectivity is true after inverting $2$ (see Propositions \ref{pro: L invert 2}
and  \ref{pro: relating L and K}) 
\\
\hline
elementary amenable groups 
& true with coefficients (see Theorem~\ref{the: BCC with coefficients for a-Tmenable groups})
& ? 
& true fibered after inverting $2$ (see Theorem~\ref{the: fibered Farrell-Jones and subgroups of Lie groups})
\\
\hline
virtually poly-cyclic
& true with coefficients (see Theorem~\ref{the: BCC with coefficients for a-Tmenable groups})
& true rationally, true fibered in the range $n \le 1$ (compare Remark~\ref{rem: Wilking})
& true fibered after inverting $2$ (see Theorem~\ref{the: fibered Farrell-Jones and subgroups of Lie groups})
\\
\hline
torsion free virtually solvable subgroups of $GL(n,\IC)$ 
& true with coefficients (see Theorem~\ref{the: BCC with coefficients for a-Tmenable groups})
& true in the range $\le 1$ \cite[Theorem~1.1]{Farrell-Linnell(2003a)}
& true fibered after inverting $2$ \cite[Corollary 5.3]{Farrell-Linnell(2003a)}
\\
\hline\hline
\end{tabular}

\begin{tabular}{|| p{30mm}| p{24mm}| p{24mm}| p{24mm}||}
\hline \hline
type of group & 
Baum-Connes Conjecture~\ref{con: Baum-Connes Conjecture} 
(with coefficients \ref{con: Baum-Connes Conjecture with coefficients})
&
Farrell-Jones Conjecture \ref{con: Farrell-Jones Conjecture}
for $K$-theory for $R = \IZ$ (fibered  \ref{con: Fibered Farrell-Jones Conjecture})  
&
Farrell-Jones Conjecture \ref{con: Farrell-Jones Conjecture}
for $L$-theory for $R = \IZ$ (fibered \ref{con: Fibered Farrell-Jones Conjecture})  
\\
\hline\hline
discrete subgroups of Lie groups with finitely many path components 
& injectivity true (see Theorem~\ref{the: BCC injectivity} and 
Remark \ref{Groups acting amenable on a compact space})
& ?
& injectivity is true after inverting $2$ (see Propositions \ref{pro: L invert 2}
and  \ref{pro: relating L and K}) 
\\
\hline
subgroups of groups which are discrete cocompact subgroups of 
Lie groups with finitely many path components 
& injectivity true (see Theorem~\ref{the: BCC injectivity} and 
Remark \ref{Groups acting amenable on a compact space})
& true rationally,
  true fibered in the range $n \le 1$ 
(see Theorem~\ref{the: fibered Farrell-Jones and subgroups of Lie
  groups}) 
& probably true fibered  (see Remark 
\ref{rem: fibered Farrell-Jones and subgroups of Lie groups for L-theory}).
Injectivity is true after inverting $2$ (see Propositions \ref{pro: L invert 2}
and  \ref{pro: relating L and K}) 
\\
\hline
linear groups & injectivity is true (see Theorem~\ref{the: BCC injectivity} and 
Remark \ref{rem: Linear groups})& ? & 
injectivity is true after inverting $2$ (see Propositions \ref{pro: L invert 2}
and  \ref{pro: relating L and K})  
\\
\hline
arithmetic groups & injectivity is true (see Theorem~\ref{the: BCC injectivity} and 
Remark \ref{rem: Linear groups})& ? & 
injectivity is true after inverting $2$ (see Propositions \ref{pro: L invert 2}
and  \ref{pro: relating L and K}) 
\\

\hline
torsion free discrete subgroups of $GL(n,\IR)$ 
& injectivity is true (see Theorem~\ref{the: BCC injectivity} and 
Remark \ref{Groups acting amenable on a compact space})
& true in the range $n \le 1$ (see 
Theorem~\ref{the: Vanishing of Wh(G times IZ^k) and FJC for L-th.})
& true (see Theorem~\ref{the: Vanishing of Wh(G times IZ^k) and FJC for L-th.})
\\
\hline\hline
\end{tabular}


\begin{tabular}{|| p{30mm}| p{24mm}| p{24mm}| p{24mm}||}
\hline \hline
type of group & 
Baum-Connes Conjecture \ref{con: Baum-Connes Conjecture} 
(with coefficients \ref{con: Baum-Connes Conjecture with coefficients})
&
Farrell-Jones Conjecture \ref{con: Farrell-Jones Conjecture}
for $K$-theory for $R = \IZ$ (fibered  \ref{con: Fibered Farrell-Jones Conjecture})  
&
Farrell-Jones Conjecture \ref{con: Farrell-Jones Conjecture}
for $L$-theory for $R = \IZ$ (fibered \ref{con: Fibered Farrell-Jones Conjecture})  
\\
\hline\hline
Groups with finite BG and
finite asymptotic dimension
&
injectivity is true (see Theorem~\ref{the: Yu injectivity and asymtotic dimension})
&
injectivity is true for arbitrary coefficients $R$
(see Theorem~\ref{the: Bartels injectivity and asymtotic dimension})
&
injectivity is true for regular $R$ as coefficients
(see Theorem~\ref{the: Bartels injectivity and asymtotic dimension})
\\
\hline\hline
$G$ acts properly and isometrically on a complete Riemannian manifold $M$ with 
non-positive sectional curvature
& rational injectivity is true (see Theorem~\ref{the: BCC rational injectivity}) 
& ?
& rational injectivity is true (see Propositions \ref{pro: L invert 2}
and  \ref{pro: relating L and K}) 
\\
\hline
$\pi_1(M)$ for a complete Riemannian manifold $M$ with 
non-positive sectional curvature
& rationally injective (see Theorem~\ref{the: BCC rational injectivity}) 
& ?
& injectivity true (see Theorem~\ref{the: Ferry-Weinberger L injective})
\\
\hline
$\pi_1(M)$ for a complete Riemannian manifold $M$ with 
non-positive sectional curvature which is A-regular
& rationally injective (see Theorem~\ref{the: BCC rational injectivity}) 
& true in the range $n \le 1$, rationally surjective 
(see Theorem~\ref{the: Vanishing of Wh(G times IZ^k) and FJC for
  L-th.})
& true (see Theorem~\ref{the: Vanishing of Wh(G times IZ^k) and FJC for L-th.})
\\
\hline
$\pi_1(M)$ for a complete Riemannian manifold $M$ with 
pinched negative sectional curvature
& rational injectivity is true (see Theorem
\ref{the: BCC injectivity})
& true in the range $n \le 1$, rationally surjective 
(see Theorem~\ref{the: Vanishing of Wh(G times IZ^k) and FJC for L-th.})
& true (see Theorem~\ref{the: Vanishing of Wh(G times IZ^k) and FJC for L-th.})
\\
\hline
$\pi_1(M)$ for a closed Riemannian manifold $M$ with 
non-positive sectional curvature
& rationally injective (see Theorem~\ref{the: BCC rational injectivity}) 
& true fibered in the range $n \le 1$, true rationally
(see Theorem~\ref{the: Vanishing of Wh(G times IZ^k) and FJC for L-th.})
& true (see Theorem~\ref{the: Vanishing of Wh(G times IZ^k) and FJC for L-th.})
\\
\hline
$\pi_1(M)$ for a closed Riemannian manifold $M$ with 
 negative sectional curvature
& true for all subgroups (see Theorem~\ref{the: BCC})
& true for all coefficients $R$ (see Theorem~\ref{the: Bartels-Reich})
& true  (see Theorem~\ref{the: Vanishing of Wh(G times IZ^k) and FJC for L-th.})
\\
\hline\hline
\end{tabular}


\begin{tabular}{|| p{30mm}| p{24mm}| p{24mm}| p{24mm}||}
\hline \hline
type of group & 
Baum-Connes Conjecture \ref{con: Baum-Connes Conjecture} 
(with coefficients \ref{con: Baum-Connes Conjecture with coefficients})
&
Farrell-Jones Conjecture \ref{con: Farrell-Jones Conjecture}
for $K$-theory and $R = \IZ$ (fibered  \ref{con: Fibered Farrell-Jones Conjecture})  
&
Farrell-Jones Conjecture \ref{con: Farrell-Jones Conjecture}
for $L$-theory for $R = \IZ$ (fibered \ref{con: Fibered Farrell-Jones Conjecture}) 
\\
\hline\hline
word hyperbolic groups 
& true for all subgroups (see Theorem~\ref{the: BCC}). Unpublished proof with coefficients by V.~Lafforgue
& ? 
& injectivity is true after inverting $2$ (see Propositions \ref{pro: L invert 2}
and  \ref{pro: relating L and K}) 
\\
\hline
one-relator groups 
& true with coefficients (see Theorem~\ref{the: BCC with coefficients for the class LHETH})
& rational injectivity is true for the fibered version (see \cite{Bartels-Lueck(2004)})
& injectivity is true after inverting $2$ (see Propositions \ref{pro: L invert 2}
and  \ref{pro: relating L and K}) 
\\
\hline
torsion free one-relator groups 
& true with coefficients (see Theorem~\ref{the: BCC with coefficients for the class LHETH})
& true for $R$ regular 
\cite[Theorem~19.4 on page 249 and Theorem~19.5  on page 250]{Waldhausen(1978a)}
& true after inverting $2$ \cite[Corollary 8]{Cappell(1973)}
\\
\hline
Haken $3$-manifold groups (in particular knot groups) 
& true with coefficients (see Theorem~\ref{the: BCC with coefficients for the class LHETH})
& true in the range $n \le 1$ for $R$ regular 
\cite[Theorem~19.4 on page 249 and Theorem~19.5  on page 250]{Waldhausen(1978a)}
& true after inverting $2$ \cite[Corollary 8]{Cappell(1973)}
\\
\hline
$SL(n,\IZ), n \geq 3$ 
& injectivity is true
& compare Remark~\ref{rem: P and SLnIZ}
& injectivity is true after inverting $2$ (see Propositions \ref{pro: L invert 2}
and  \ref{pro: relating L and K}) 
\\
\hline 
Artin's braid group $B_n$
& true with coefficients \cite[Theorem~5.25]{Mislin-Valette(2003)}, \cite{Schick(2000a)} 
& true fibered in the range $n \leq 1$, true rationally \cite{Farrell-Roushon(2000)}
& injectivity is true after inverting $2$ (see Propositions \ref{pro: L invert 2}
and  \ref{pro: relating L and K}) 
\\
\hline 
pure braid group $C_n$
& true with coefficients 
& true in the range $n \le 1$ (see 
Theorem~\ref{the: Vanishing of Wh(G times IZ^k) and FJC for L-th.})
& injectivity is true after inverting $2$ (see Propositions \ref{pro: L invert 2}
and  \ref{pro: relating L and K}) 
\\
\hline
Thompson's group $F$ 
& true with coefficients \cite{Farley(2001)} 
& ? 
& injectivity is true after inverting $2$ (see Propositions \ref{pro: L invert 2}
and  \ref{pro: relating L and K}) 
\\
\hline\hline
\end{tabular}
\end{center}


\begin{remarknew} The authors have no information about the status of these conjectures for
mapping class groups of higher genus or the group of outer automorphisms of  free groups.
Since all of these spaces have finite models for $\EGF{G}{\calfin}$
Theorem~\ref{the: LRRM} applies in these cases. 
\end{remarknew}


\subsection{Inheritance Properties}
\label{sec: Inheritance Properties}

In this Subsection we list some inheritance properties of the various
conjectures.


\subsubsection{Directed Colimits}
\label{subsec: Directed Colimits}

Let $\{G_i \mid i \in I\}$ be a directed system of groups. Let $G =
\colim_{i \in I} G_i$ be the colimit. We do not require that the structure maps are
injective. If the Fibered Farrell-Jones
Conjecture \ref{con: Fibered Farrell-Jones Conjecture} is true for each $G_i$, then it is
true for $G$ \cite[Theorem~6.1]{Farrell-Linnell(2003a)}. 

Suppose that $\{G_i \mid i \in I\}$ is a system of subgroups of $G$ directed by inclusion
such that $G = \colim_{i \in I} G_i$. If each $G_i$ satisfies the 
Farrell-Jones Conjecture~\ref{con: Farrell-Jones Conjecture}, 
the Baum-Connes Conjecture \ref{con: Baum-Connes
  Conjecture} or the  Baum-Connes Conjecture with 
Coefficients~\ref{con: Baum-Connes Conjecture with coefficients}, then the same is true for $G$
\cite[Theorem~1.1]{Baum-Millington-Plymen(2002)},
\cite[Lemma 5.3]{Mislin-Valette(2003)}. We do not know a reference in Farrell-Jones case. An argument in that case 
uses Lemma~\ref{lem: G-homology theory and colimit}, the fact that $K_n (RG) = \colim_{i \in I} K_n (RG_i)$ and that
for suitable models we have $\EGF{G}{\calf} = \bigcup_{i \in I} G \times_{G_i} \EGF{G_i}{\calf \cap G_i}$.


\subsubsection{Passing to Subgroups}
\label{subsec: Passing to Subgroups}

The Baum-Connes Conjecture with Coefficients
\ref{con: Baum-Connes Conjecture with coefficients} and the
Fibered Farrell-Jones Conjecture~\ref{con: Fibered Farrell-Jones Conjecture}  
pass to subgroups, i.e. if they hold for $G$, then also for any subgroup
$H \subseteq G$. This claim for the Baum-Connes Conjecture with Coefficients
\ref{con: Baum-Connes Conjecture with coefficients}
has been stated in \cite{Baum-Connes-Higson(1994)}, a proof can be found for instance
in \cite[Theorem~2.5]{Chabert-Echterhoff(2001b)}.
For the Fibered Farrell-Jones Conjecture this is proved in 
\cite[Theorem~A.8 on page 289]{Farrell-Jones(1993a)}  for the special case $R = \IZ$, but
the proof also works for arbitrary rings $R$. 

It is not known whether the Baum-Connes Conjecture \ref{con: Baum-Connes Conjecture} 
or the Farrell-Jones Conjecture \ref{con: Farrell-Jones Conjecture}  itself
passes to subgroups. 


\subsubsection{Extensions of Groups}
\label{subsec: Extensions of Groups}

Let $p \colon G \to K$ be a surjective group homomorphism.
Suppose that the Baum-Connes Conjecture with Coefficients
\ref{con: Baum-Connes Conjecture with coefficients} or the
Fibered Farrell-Jones Conjecture \ref{con: Fibered Farrell-Jones Conjecture} respectively  holds 
for $K$ and for $p^{-1}(H)$ for any subgroup $H \subset K$ which is finite or virtually
cyclic respectively. Then the Baum-Connes Conjecture with 
Coefficients~\ref{con: Baum-Connes Conjecture with coefficients} or the
Fibered Farrell-Jones Conjecture \ref{con: Fibered Farrell-Jones Conjecture} respectively  holds 
for $G$. This is proved in \cite[Theorem~3.1]{Oyono-Oyono(2001)}
for the Baum-Connes Conjecture with
Coefficients~\ref{con: Baum-Connes Conjecture with coefficients}, 
and in \cite[Proposition~2.2 on page 263]{Farrell-Jones(1993a)} for the
Fibered Farrell-Jones Conjecture \ref{con: Fibered Farrell-Jones Conjecture}  
in the case $R = \IZ$. The same proof works for arbitrary coefficient rings.

It is not known whether the corresponding statement holds for the 
Baum-Connes Conjecture \ref{con: Baum-Connes Conjecture} 
or the Farrell-Jones Conjecture \ref{con: Farrell-Jones Conjecture} itself.

Let $H \subseteq G$ be a normal subgroup of $G$. Suppose that $H$ is a-T-menable.
Then $G$ satisfies the  
Baum-Connes Conjecture with Coefficients
\ref{con: Baum-Connes Conjecture with coefficients}
if and only if $G/H$ does
\cite[Corollary 3.14]{Chabert-Echterhoff(2001b)}. 
The corresponding statement is not known for the 
Baum-Connes Conjecture \ref{con: Baum-Connes Conjecture}.


\subsubsection{Products of Groups}
\label{subsec: Product of Groups}

The group $G_1 \times G_2$ satisfies the Baum-Connes Conjecture with Coefficients
\ref{con: Baum-Connes Conjecture with coefficients} if and only if both $G_1$ and
$G_2$ do \cite[Theorem~3.17]{Chabert-Echterhoff(2001b)},
\cite[Corollary 7.12]{Oyono-Oyono(2001)}.
The corresponding statement is not known for the 
Baum-Connes Conjecture \ref{con: Baum-Connes Conjecture}. 

Let $D_{\infty}= \IZ/ 2 \ast \IZ / 2$ denote the infinite dihedral group.
Whenever a version of the Fibered Farrell-Jones Conjecture~\ref{con: Fibered Farrell-Jones Conjecture}  
is known for 
$G= \IZ \times \IZ$, $G=\IZ \times D_{\infty}$ and $D_{\infty} \times D_{\infty}$,
then that version of the Fibered Farrell-Jones Conjecture is true for $G_1 \times G_2$
if and only if it is true for $G_1$ and $G_2$. 


\subsubsection{Subgroups of Finite Index}
\label{subsec: Subgroups of Finite Index}

It is not known whether the Baum-Connes Conjecture \ref{con: Baum-Connes Conjecture},
the Baum-Connes Conjecture with Coefficients
\ref{con: Baum-Connes Conjecture with coefficients},
the Farrell-Jones Conjecture \ref{con: Farrell-Jones Conjecture}
or the Fibered Farrell-Jones Conjecture \ref{con: Fibered Farrell-Jones Conjecture} 
is true for a group $G$ if it is true for a subgroup $H\subseteq G$ of
finite index.


\subsubsection{Groups Acting on Trees}
\label{subsec: Groups Acting on Trees}

Let $G$ be a countable discrete group acting without inversion on a tree $T$. 
Then the Baum-Connes Conjecture with Coefficients
\ref{con: Baum-Connes Conjecture with coefficients} is true for $G$ if
and only if it holds for all stabilizers of the vertices of $T$. 
This is proved by Oyono-Oyono \cite[Theorem~1.1]{Oyono-Oyono(2001b)}.
This implies that Baum-Connes Conjecture with Coefficients
\ref{con: Baum-Connes Conjecture with coefficients} is stable under
amalgamated products and HNN-extensions. Actions on trees in the context
the Farrell-Jones Conjecture \ref{con: Farrell-Jones Conjecture}
will be treated in \cite{Bartels-Lueck(2004)}.

\typeout{----------------------------- Equivariant Homology Theories --------------------}

\section{Equivariant Homology Theories}
\label{chap: Equivariant Homology Theories}

We already defined the notion of a $G$-homology theory in Subsection~\ref{subsec: G-Homology Theories}.
If $G$-homology theories for different $G$ are linked via a so called induction structure one 
obtains the notion of an equivariant homology theory. In this section we give a precise definition and we explain
how a functor from the orbit category $\Or (G)$ to the category of spectra leads to a $G$-homology theory
(see Proposition~\ref{pro: Or(G)-spectra yield a G-homology theory}) 
and how more generally a functor from the category of groupoids leads to an equivariant homology theory
(see Proposition~\ref{pro: GROUPOID-spectra and equivariant homology theories}). 
We then describe the main examples of such spectra valued functors which were already used in order to formulate the
Farrell-Jones and the Baum-Connes Conjectures in Chapter~\ref{chap: general formulation}.

\subsection{The Definition of an Equivariant Homology Theory}
\label{sec: The Definition of an Equivariant Homology Theory}

The notion of a $G$-homology theory 
$\calh_*^G$ with values in
$\Lambda$-modules  for a commutative ring $\Lambda$ was defined in 
Subsection~\ref{subsec: G-Homology Theories}.
We now recall the axioms of an equivariant homology theory from \cite[Section 1]{Lueck(2002b)}.
We will see in Section~\ref{sec: K and L-Theory Spectra over Groupoids} that the $G$-homology theories we used in
the formulation of the Baum-Connes and the Farrell-Jones Conjectures in Chapter~\ref{chap: general formulation}
are in fact the values at $G$ of suitable equivariant homology theories.

Let $\alpha\colon  H \to G$ be a group homomorphism.
Given a $H$-space $X$, define the \emph{induction of $X$ with $\alpha$}%
\index{induction!of equivariant spaces}
to be the $G$-space $\ind_{\alpha} X$ which  is the quotient of
$G \times X$ by the right $H$-action
$(g,x) \cdot h := (g\alpha(h),h^{-1} x)$
for $h \in H$ and $(g,x) \in G \times X$.
If $\alpha\colon  H \to G$ is an inclusion, we also write $\ind_H^G$ instead of
$\ind_{\alpha}$.

An \emph{equivariant homology theory $\calh_*^?$ with values in $\Lambda$-modules}%
\index{homology theory!equivariant G-homology theory@equivariant $G$-homology theory}
consists of a
$G$-homology theory $\calh^G_*$ with values in $\Lambda$-modules for each group $G$
together with the following so called \emph{induction structure}:%
\index{induction structure}
given a group homomorphism $\alpha\colon  H \to G$ and  a $H$-$CW$-pair
$(X,A)$ such that $\ker(\alpha)$ acts freely on $X$,
there are for each $n \in \IZ$
natural isomorphisms
\begin{eqnarray*}
\ind_{\alpha}\colon  \calh_n^H(X,A)
&\xrightarrow{\cong} &
\calh_n^G(\ind_{\alpha}(X,A)) 
\end{eqnarray*}
satisfying the following conditions.

\begin{enumerate}

\item Compatibility with the boundary homomorphisms\\[1mm]
$\partial_n^G \circ \ind_{\alpha} = \ind_{\alpha} \circ \partial_n^H$.

\item 
\label{functoriality}
Functoriality\\[1mm]
Let $\beta\colon  G \to K$ be another group
homomorphism such that $\ker(\beta \circ \alpha)$ acts freely on $X$.
Then we have for $n \in \IZ$
$$\ind_{\beta \circ \alpha} ~ = ~
\calh^K_n(f_1)\circ\ind_{\beta} \circ \ind_{\alpha}\colon 
\calh^H_n(X,A) \to \calh_n^K(\ind_{\beta\circ\alpha}(X,A)),$$
where $f_1\colon  \ind_{\beta}\ind_{\alpha}(X,A)
\xrightarrow{\cong} \ind_{\beta\circ \alpha}(X,A),
\hspace{3mm} (k,g,x) \mapsto (k\beta(g),x)$
is the natural $K$-homeo\-mor\-phism.

\item Compatibility with conjugation\\[1mm]
For $n \in \IZ$, $g \in G$ and a $G$-$CW$-pair $(X,A)$
the homomorphism 
\[
\ind_{c(g)\colon  G \to G}\colon  \calh^G_n(X,A)\to
\calh^G_n(\ind_{c(g)\colon  G \to G}(X,A))
\]
agrees with
$\calh_n^G(f_2)$, where the $G$-homeomorphism
\[
f_2\colon  (X,A) \to \ind_{c(g)\colon  G \to G} (X,A)
\]
sends
$x$ to $(1,g^{-1}x)$ and $c(g)\colon G \to G$ sends $g'$ to $gg'g^{-1}$.

\end{enumerate}

This induction structure links the
various homology theories for different groups $G$.

If the $G$-homology theory $\calh_*^G$%
\indexnotation{calh_*^?}
is defined or
considered only for proper $G$-$CW$-pairs $(X,A)$, we call it a
\emph{proper $G$-homology theory $\calh_*^G$ with values
in $\Lambda$-modules}.%
\index{homology theory!proper G-homology theory@proper $G$-homology theory}

\begin{examplenew}
\label{exa: homology of the quotient or Borel construction}
Let $\calk_*$ be a homology theory for (non-equivariant)
$CW$-pairs with values in $\Lambda$-modules.
Examples are singular homology, oriented bordism theory
or topological $K$-homology.
Then we obtain two equivariant homology
theories with values in $\Lambda$-modules, whose underlying $G$-homology theories for a group $G$
are given by the following  constructions
\begin{eqnarray*}
\calh^G_n(X,A) & = & \calk_n(G\backslash X,G\backslash A);
\\
\calh^G_n(X,A) & = & \calk_n(EG \times_G (X,A)).
\end{eqnarray*}
\end{examplenew}

\begin{examplenew} \label{exa: equivariant bordism}
Given a proper $G$-$CW$-pair $(X,A)$,
one can define the $G$-bordism group $\Omega_n^G(X,A)$
as the abelian group of $G$-bordism
classes of maps $f\colon  (M,\partial M) \to (X,A)$
whose sources are oriented smooth manifolds
with cocompact orientation preserving  proper smooth
$G$-actions. The definition is analogous to the one in the non-equivariant case.
This is also true for the proof
that this defines a proper $G$-homology theory.
There is an obvious induction structure coming from induction of 
equivariant spaces. Thus we obtain an equivariant proper homology theory
$\Omega_*^?$.
\end{examplenew}

\begin{examplenew} \label{Bredon homology as equivariant homology theory} 
Let $\Lambda$ be a commutative ring and let
$$M \colon \GROUPOIDS \to \Lambda\text{-}\MODULES$$
be a contravariant functor. For a group $G$ we obtain a covariant functor 
$$M^G \colon \Or(G) \to \Lambda\text{-}\MODULES$$
by its composition with the transport groupoid functor $\calg^G$ defined in \eqref{functor calg^G}.
Let $H^G_*(-;M)$ be the $G$-homology theory given by the Bredon
homology with coefficients in $M^G$ as defined in Example~\ref{Bredon homology as G-homology theory}.
There is an induction structure such that the collection of the $H^G(-;M)$ defines an
equivariant homology theory $H^?_*(-;M)$. This can be interpreted as the special case of
Proposition~\ref{pro: GROUPOID-spectra and equivariant homology theories}, where 
the covariant functor
$\GROUPOIDS \to \Omega\text{-}\SPECTRA$ is the composition of $M$ with the functor
sending a $\Lambda$-module to the associated Eilenberg-MacLane spectrum. But there is also a purely algebraic
construction.
\end{examplenew}

The next lemma was used in the proof of the Transitivity Principle~\ref{the: transitivity}.

\begin{lemma} \label{lem: EGF{G}{calf} times Z to Z}
Let $\calh^?_*$ be an equivariant homology theory with values in $\Lambda $-modules.
Let $G$ be a group and let $\calf$ a family of subgroups of $G$. 
Let $Z$ be a $G$-$CW$-complex. Consider $N \in \IZ \cup \{\infty\}$. 
For $H \subseteq G$ let $\calf \cap H$ be the family of subgroups of $H$ given by
 $\{ K \cap H \; | \; K \in \calf \}$.
Suppose for each $H \subset G$, which occurs as isotropy group in $Z$,
that the map induced by the projection
$\pr\colon \EGF{H}{\calf \cap H} \to \pt$ 
$$\calh^H_n(\pr) \colon \calh^H_n(\EGF{H}{\calf \cap H}) \to \calh^H_n(\pt)$$
is bijective for all $n \in \IZ, n \le N$. 

Then the map induced by the projection $\pr_2 \colon \EGF{G}{\calf} \times Z \to Z$
\begin{eqnarray*} \calh^G_n(\pr_2) \colon \calh^G_n(\EGF{G}{\calf} \times Z) 
& \to & \calh^G_n(Z)
\end{eqnarray*}
is bijective for $n \in \IZ, n \le N$.
\end{lemma}
\begin{proof}
We first prove the claim for 
finite-dimensional $G$-$CW$-complexes by induction over $d = \dim(Z)$.
The induction beginning $\dim(Z) = -1$, i.e. $Z = \emptyset$, is trivial.
In the induction step from $(d-1)$ to $d$ we choose a $G$-pushout
$$\comsquare{\coprod_{i \in I_d} G/H_i \times S^{d-1}}{}{Z_{d-1}}
{}{}
{\coprod_{i \in I_d} G/H_i \times D^d}{}{Z_d}
$$
If we cross it with $\EGF{G}{\calf}$, we obtain another $G$-pushout of $G$-$CW$-complexes.
The various projections induce a map
from the Mayer-Vietoris sequence of the latter $G$-pushout to the Mayer-Vietoris sequence of the
first $G$-pushout. By the Five-Lemma it suffices to prove
that the following maps
\begin{eqnarray*}
\calh_n^G(\pr_2) \colon \calh^G_n\left(\EGF{G}{\calf} \times \coprod_{i \in I_d}
G/H_i \times S^{d-1}\right) 
& \to &  
\calh^G_n\left(\coprod_{i \in I_d} G/H_i \times S^{d-1}\right);
\\
\calh_n^G(\pr_2) \colon \calh^G_n(\EGF{G}{\calf} \times Z_{d-1}) & \to &  \calh^G_n(Z_{d-1});
\\
\calh_n^G(\pr_2) \colon \calh^G_n\left(\EGF{G}{\calf} \times \coprod_{i \in I_d} G/H_i \times D^d\right) & \to &  
\calh^G_n\left(\coprod_{i \in I_d}G/H_i \times D^d\right)
\end{eqnarray*}
are bijective for $n \in \IZ, n \le N$. This follows from 
the induction hypothesis for the first two maps.
Because of the disjoint union axiom  and $G$-homotopy invariance 
of $\calh^?_*$ the claim follows for the third map
if we can show for any $H\subseteq G$ which occurs as isotropy group in $Z$ that the map
\begin{eqnarray}
\calh_n^G(\pr_2) \colon \calh^G_n(\EGF{G}{\calf} \times G/H ) & \to &  \calh^G(G/H)
\label{calh^G_n(pr_3) for G/H}
\end{eqnarray}
is bijective for $n \in \IZ, n \le N$. The $G$-map 
$$G \times_H \res_G^H \EGF{G}{\calf} \to G/H \times \EGF{G}{\calf} \hspace{5mm} (g,x) ~ \mapsto ~ (gH,gx)$$
is a $G$-homeomorphism where $\res_G^H$ denotes the restriction of the $G$-action to an $H$-action.
Obviously $\res_G^H \EGF{G}{\calf}$ is a model for $\EGF{H}{\calf \cap H}$. 
We conclude from the induction structure
that the map \eqref{calh^G_n(pr_3) for G/H}  can be identified with the map
\begin{eqnarray*}
\calh_n^G(\pr) \colon \calh^H_n(\EGF{H}{\calf \cap H}) & \to &  \calh^H(\pt)
\end{eqnarray*}
which is bijective for all $n \in \IZ, n \le N$ by assumption. This finishes the proof 
in the case  that $Z$ is finite-dimensional.
The general case follows by a colimit argument using Lemma~\ref{lem: G-homology theory and colimit}.
\end{proof}


\subsection{Constructing Equivariant Homology Theories}
\label{sec: Spectra over the Orbit Category}

Recall that a (non-equivariant) spectrum yields an associated
(non-equivariant) homology theory. In this section we explain how a spectrum
over the orbit category of a group $G$ defines a $G$-homology theory.
We would like to stress that our approach using spectra over the orbit category
should be distinguished from approaches to equivariant homology theories
using spectra with $G$-action or the more complicated notion of $G$-equivariant
spectra in the sense of \cite{Lewis-May-Steinberger(1986)}, see for example \cite{Carlsson(1992a)} for a survey. 
The latter approach leads to a much richer structure
but only works for compact Lie groups.

We briefly fix some conventions concerning spectra.
We always work in the very convenient
\emph{category} 
$\SPACES$
\emph{ of  compactly generated spaces}%
\index{category!of compactly generated spaces}
\indexnotation{SPACES} 
(see \cite{Steenrod(1967)}, \cite[I.4]{Whitehead(1978)}). 
In that category the 
adjunction homeomorphism $\map(X \times Y,Z) \xrightarrow{\cong} \map(X,\map(Y,Z))$ holds without
any further assumptions such as local compactness 
and the product of two $CW$-complexes is again a $CW$-complex.
Let $\SPACES^+$%
\indexnotation{SPACES^+}
be the category of pointed compactly generated spaces. Here the 
objects are (compactly generated) spaces $X$ with base points for
which the inclusion of the base point is a cofibration. Morphisms are
pointed maps.  If $X$ is a space, denote by $X_+$%
\indexnotation{X_+}
the pointed space obtained from $X$
by adding a disjoint base point.
For two pointed spaces $X= (X,x)$ and $Y= (Y,y)$ define their
\emph{smash product}%
\index{smash product} as the pointed space
\indexnotation{X wedge Y}
\[
X \wedge Y = X \times Y/(\{x\} \times Y \cup X \times \{y\}),
\] 
and the  \emph{reduced cone}%
\index{cone!reduced cone}
\indexnotation{cone(X)}
as
\[
\cone(X) = X \times [0,1]/ (X \times \{1\} \cup \{x\} \times [0,1]).
\]

A {\em spectrum}%
\index{spectrum}
$\mathbf{E} = \{(E(n),\sigma(n)) \mid n \in
\IZ\}$ is a sequence of pointed spaces
$\{E(n) \mid n \in \IZ\}$ together with pointed maps
called {\em structure maps}%
\index{spectrum!structure maps of a spectrum}
$\sigma(n) \colon  E(n) \wedge S^1 \longrightarrow E(n+1)$.
A {\em map}%
\index{spectrum!map of spectra}
of spectra 
$\bff \colon  \bfE \to \bfE^{\prime}$ is a sequence of maps
$f(n) \colon  E(n) \to E^{\prime}(n)$
which are compatible with the structure maps $\sigma(n)$, i.e.\
we have
$f(n+1) \circ \sigma(n) ~ = ~
\sigma^{\prime}(n) \circ \left(f(n) \wedge \id_{S^1}\right)$
for all $n \in \IZ$.
Maps of spectra are sometimes called functions in the literature, they should not be confused with the notion of a map of spectra in
the stable category (see \cite[III.2.]{Adams(1974)}). 
The category of spectra and  maps will be denoted $\SPECTRA$%
\indexnotation{SPECTRA}.
Recall
that the homotopy groups of a spectrum%
\index{spectrum!homotopy groups of a spectrum}
are defined by 
$$\pi_i(\mathbf{E})%
\indexnotation{pi_i(bfE)}
 ~ = ~
\colim_{k \to \infty} \pi_{i+k}(E(k)),$$
where the system $\pi_{i+k}(E(k))$ is given by the composition
$$\pi_{i+k}(E(k)) \xrightarrow{S} \pi_{i+k+1}(E(k)\wedge S^1)
\xrightarrow{\sigma (k)_*} \pi_{i+k+1}(E(k +1))$$
of the suspension homomorphism $S$ and the
homomorphism induced by the structure map.
A {\em weak equivalence}%
\index{weak equivalence!of spectra} 
of spectra is a map $\bff\colon \bfE \to \bfF$ 
of spectra inducing an isomorphism on all homotopy
groups.  

Given a spectrum $\bfE$ and a pointed space $X$, we can define their
smash product%
\index{smash product!of a space with a spectrum}
$X \wedge \bfE$
\indexnotation{X wedge bfE}
by $(X \wedge \bfE)(n) := X \wedge E(n)$ with the obvious structure
maps. 
It is a classical result that a spectrum $\bfE$ defines a homology
theory by setting
\[
H_n(X,A;\bfE) = \pi_n\left((X_+ \cup_{A_+} \cone(A_+)) \wedge \bfE \right).
\]
We want to extend this to $G$-homology theories. This requires the
consideration of spaces and spectra over the orbit category. Our presentation
follows \cite{Davis-Lueck(1998)}, where more details can be found.

In the sequel $\calc$ is a small category. Our main example is the orbit category
$\Or(G)$, whose objects are homogeneous $G$-spaces $G/H$
and whose morphisms are $G$-maps.

\begin{definition} \label{def:  spaces over a category}
A \emph{covariant (contravariant) $\calc$-space $X$}%
\index{space!C-space@$\calc$-space}
is a covariant (contravariant) functor
\[
X\colon  \calc ~ \to ~ \SPACES.
\]
A map between $\calc$-spaces is a natural transformation of such functors. 
Analogously a \emph{pointed $\calc$-space} is a functor from
$\calc$ to $\SPACES^+$ and a $\calc$-spectrum a functor to $\SPECTRA$.
\end{definition}

\begin{examplenew}  \label{exa: G-spaces yields Or(G)-space}
\em  Let $Y$ be a left $G$-space. Define the associated
\emph{contravariant $\Or(G)$-space}
$\map_G(-,Y)$ by
$$\map_G(-,Y)%
\indexnotation{map_G(-,Y)}
\colon ~  \Or(G) \to \SPACES, \hspace{10mm}
G/H ~ \mapsto \map_G(G/H,Y) = Y^H.$$
If $Y$ is pointed then $\map_G( - , Y)$ takes values in pointed spaces.
\end{examplenew}

Let $X$ be a contravariant and $Y$ be a covariant $\calc$-space. Define
their \emph{balanced product}%
\index{balanced product!of $\calc$-spaces}
 to be the space
$$X \times_{\calc} Y%
\indexnotation{X times_calc Y}
~ := ~ \coprod_{c \in \Ob(\calc)} X(c) \times Y(c)/\sim$$
where $\sim$ is the equivalence relation generated by
$(x\phi,y) \sim (x,\phi y)$
for all morphisms  $\phi\colon c \to d$ in $\calc$ and points
$x \in X(d)$ and $y \in Y(c)$. Here $x\phi$ stands for $X(\phi)(x)$ and $\phi y$
for $Y(\phi )(y)$. If $X$ and $Y$ are pointed, then one defines analogously
their \emph{balanced smash product}%
\index{balanced smash product!of pointed $\calc$-spaces}
to be the pointed space
\indexnotation{X wedge_calc Y}
\[
X \wedge_{\calc} Y ~ = ~
\bigvee_{c \in \Ob(\calc)} X(c) \wedge Y(c)/\sim.
\]
In \cite{Davis-Lueck(1998)} the notation $X \otimes_{\calc} Y$ was used for this space.
Doing the same construction level-wise one defines the \emph{balanced smash product}%
\index{tensor product!of pointed $\calc$-space with a $\calc$-spectrum}
$X \wedge_{\calc} \bfE$%
\indexnotation{X wedge_calc bfE}
of a contravariant pointed $\calc$-space and a covariant $\calc$-spectrum $\bfE$.

The proof of the next result is analogous to the non-equivariant case.
Details can be found in \cite[Lemma~4.4]{Davis-Lueck(1998)}, where also
cohomology theories are treated.

\begin{propositionnew}[Constructing $G$-Homology Theories]
\label{pro: Or(G)-spectra yield a G-homology theory}
Let $\bfE$ be a covariant $\Or(G)$-spectrum. It defines a $G$-homology theory
$H^G_*(-;\bfE)$%
\indexnotation{H^G_*(-;bfE)}
 by
\[
H^G_n(X,A;\bfE) ~= ~ 
\pi_n\left(\map_G\left(-,(X_+ \cup_{A_+} \cone(A_+))\right)
  \wedge_{\Or(G)} \bfE \right).
\]
In particular we have
\[
H^G_n(G/H;\bfE) ~ = ~ \pi_n(\bfE(G/H)).
\]
\end{propositionnew}

Recall that we seek an equivariant homology theory and not only a $G$-homology theory.
If the $\Or(G)$-spectrum in Proposition~\ref{pro: Or(G)-spectra yield a G-homology theory}
is obtained from a $\GROUPOIDS$-spectrum in a way we will now describe, then automatically 
we obtain the desired induction structure.

Let $\GROUPOIDS$
\indexnotation{GROUPOIDS} be the category of small groupoids with covariant functors as
morphisms. Recall that a groupoid is a category in which all morphisms are isomorphisms.
A covariant functor $f\colon \calg_0 \to \calg_1$ of groupoids
is called \emph{injective}, if for any two objects $x,y$ in $\calg_0$ the induced map
$\mor_{\calg_0}(x,y) \to \mor_{\calg_1}(f(x),f(y))$ is injective. Let
$\GROUPOIDS^{\inj}$%
\indexnotation{GROUPOIDS^inj}
 be the subcategory of $\GROUPOIDS$ with the same objects and injective
functors as morphisms. For a $G$-set $S$ we denote by $\calg^G(S)$ its associated \emph{transport groupoid}.%
\indexnotation{calg^G(S)} \index{transport groupoid}
Its objects are the elements of $S$. The set of morphisms from $s_0$ to $s_1$ consists
of those elements $g \in G$ which satisfy $g s_0 = s_1$. Composition in $\calg^G(S)$
comes from the multiplication in $G$. Thus we obtain for a group $G$ a covariant functor
\begin{eqnarray}
\calg^G\colon \Or(G) \to \GROUPOIDS^{\inj}, \hspace{7mm} G/H \mapsto \calg^G(G/H).
\label{functor calg^G}
\end{eqnarray}

A functor of small categories $F\colon \calc \to \cald$ is called an \emph{equivalence}%
\index{equivalence of categories} if there exists a
functor $G \colon \cald \to \calc$ such that both $F \circ G$ and
$G \circ F$ are naturally equivalent to the identity
functor. This is equivalent to the condition that $F$ induces a bijection on the set of isomorphisms
classes of objects and for any  objects $x,y\in \calc$ the map 
$\mor_{\calc}(x,y) \to \mor_{\cald}(F(x),F(y))$ induced by $F$ is bijective.

\begin{propositionnew}[Constructing Equivariant Homology Theories]
\label{pro: GROUPOID-spectra and equivariant homology theories}
Consider a covariant $\GROUPOIDS^{\inj}$-spectrum
$$\bfE\colon \GROUPOIDS^{\inj} \to \SPECTRA.$$
Suppose that $\bfE$ respects equivalences, i.e. it sends an
equivalence of groupoids to a weak equivalence of spectra.
Then $\bfE$ defines an equivariant homology theory
$H_*^?(-;\bfE)$,%
\indexnotation{H_*^?(-;bfE)}
whose underlying $G$-homology theory for a group $G$ is the $G$-homology theory associated to
the covariant $\Or(G)$-spectrum $\bfE \circ \calg^G \colon \Or(G) \to \SPECTRA$ 
in the previous Proposition~\ref{pro: Or(G)-spectra yield a G-homology theory}, i.e.\
\[
H^G_*(X,A;\bfE) ~ = ~ H^G_*(X,A;\bfE \circ \calg^G).
\]
In particular we have
\[
H^G_n(G/H; \bfE) ~ \cong ~ H^H_n(\pt; \bfE) ~ \cong ~ \pi_n(\bfE(I(H))),
\]
where $I(H)$ denotes $H$ considered as a groupoid with one object.
The whole construction is natural in $\bfE$.
\end{propositionnew}
\begin{proof} We have to specify the induction structure for a homomorphism
$\alpha\colon H \to G$.  We only sketch the construction in the special case
where $\alpha$ is injective and $A = \emptyset$. 
The details of the full proof can be found in \cite[Theorem~2.10 on page 21]{SauerJ(2002)}.

The functor induced by $\alpha$ on the orbit categories is denoted in the same way
$$\alpha\colon \Or(H) \to \Or(G), \hspace{7mm} H/L \mapsto \ind_{\alpha}(H/L) = G/\alpha(L).$$
There is an obvious natural equivalence of functors $\Or(H) \to \GROUPOIDS^{\inj}$
\[
T\colon \calg^H \to \calg^G \circ \alpha.
\]
Its evaluation at $H/L$ is the equivalence of groupoids
$\calg^H(H/L) \to \calg^G(G/\alpha(L))$ which sends an object $hL$ to the object
$\alpha(h)\alpha(L)$ and a morphism given by $h \in H$ to the morphism $\alpha(h) \in G$.
The desired isomorphism 
\[
\ind_{\alpha} \colon H_n^H(X;\bfE \circ \calg^H) \to H_n^G(\ind_{\alpha} X;\bfE
\circ \calg^G)
\]
is induced by the following map of spectra
\begin{multline*}
\map_H(-,X_+) \wedge_{\Or(H)} \bfE \circ \calg^H \xrightarrow{\id \wedge \bfE(T)}
\map_H(-,X_+) \wedge_{\Or(H)} \bfE \circ \calg^G \circ \alpha
\\
\xleftarrow{\simeq} 
(\alpha_*\map_H(-,X_+)) \wedge_{\Or(G)} \bfE \circ \calg^G
\xleftarrow{\simeq} \map_G(-,\ind_{\alpha}X_+) \wedge_{\Or(G)} \bfE \circ \calg^G.
\end{multline*}
Here $\alpha_*\map_H(-,X_+)$ is the pointed $\Or(G)$-space which is 
obtained from the pointed $\Or(H)$-space $\map_H(-,X_+)$ by induction, i.e.\
by taking the balanced product over $\Or(H)$ with the 
$\Or(H)$-$\Or(G)$ bimodule $\mor_{\Or(G)}(??,\alpha(?))$
\cite[Definition 1.8]{Davis-Lueck(1998)}.
Notice that $\bfE \circ \calg^G \circ \alpha$ is the same as the restriction of the
$\Or(G)$-spectrum $\bfE \circ \calg^G$ along $\alpha$ which is often denoted by
$\alpha^*(\bfE \circ \calg^G)$ in the literature
\cite[Definition 1.8]{Davis-Lueck(1998)}.
The second map is given by the adjunction homeomorphism
of induction $\alpha_*$ and restriction $\alpha^*$ (see \cite[Lemma~1.9]{Davis-Lueck(1998)}). 
The third map is the homeomorphism of $\Or(G)$-spaces
which is the adjoint of the obvious map of $\Or(H)$-spaces
$\map_H(-,X_+) \to \alpha^*\map_G(-,\ind_{\alpha}X_+)$ whose evaluation at
$H/L$ is given by $\ind_{\alpha}$. 
\end{proof}


\subsection{$K$- and $L$-Theory Spectra over Groupoids}
\label{sec: K and L-Theory Spectra over Groupoids}

Let  $\RINGS$%
\indexnotation{RINGS} be the category of associative rings with unit.
An \emph{involution}%
\index{ring with involution} on a $R$ is a map 
$R \to R, \hspace{1mm} r \mapsto \overline{r}$
satisfying $\overline{1} = 1$, $\overline{x + y} = \overline{x} + \overline{y}$
and $\overline{x \cdot y} = \overline{y} \cdot \overline{x}$ for all $x$, $y\in R$.
Let $\RINGS^{\inv}$
\indexnotation{RINGS^inv}
 be the category of rings with involution.
Let $C^{\ast}\text{-}\ALGEBRAS$%
\indexnotation{C^*-ALGEBRAS}
 be the category of $C^*$-algebras.
There are classical functors for 
$j \in {-\infty} \amalg \{j \in \IZ \mid j \le 2\}$
\begin{eqnarray}
\bfK \colon \RINGS & \to & \SPECTRA;
\label{K^{alg}RINGS}
\\
\bfL^{\langle j \rangle}\colon \RINGS^{\inv}& \to & \SPECTRA;
\label{L^{langle j rangle}RINGS^inv}
\\
\bfK^{\topo}\colon C^{\ast}\text{-}\ALGEBRAS  & \to & \SPECTRA.
\label{K^{alg}C^*ALGEBRAS}
\end{eqnarray}
The construction of such a non-connective algebraic $K$-theory functor goes back to 
Gersten \cite{Gersten(1972)} and Wagoner
\cite{Wagoner(1972)}. The spectrum for quadratic algebraic
$L$-theory is constructed by
Ranicki in \cite{Ranicki(1992)}. In a more geometric formulation it goes back to Quinn \cite{Quinn(1970)}. 
In the topological $K$-theory case a 
construction using Bott periodicity for $C^{\ast}$-algebras 
can easily be derived from the Kuiper-Mingo Theorem (see \cite[Section 2.2] {Schroeder(1993)}).
The homotopy groups of these spectra give the algebraic $K$-groups of Quillen (in high dimensions) 
and of Bass (in negative dimensions), the decorated quadratic $L$-theory groups, and the topological
$K$-groups of $C^{\ast}$-algebras.  

We emphasize again that in all three cases we need the non-connective versions of the spectra, i.e.\ the homotopy groups
in negative dimensions are non-trivial in general.
For example the version of the Farrell-Jones Conjecture where one uses connective $K$-theory 
spectra is definitely false in general, compare Remark~\ref{rem: Bass-Heller-Swan decomposition}.

Now let us fix a coefficient ring $R$ (with involution). Then sending a group $G$ to the group ring $RG$ yields functors
$R( - ) \colon \GROUPS \to \RINGS$, respectively $R( - ) \colon \GROUPS \to \RINGS^{\inv}$, where $\GROUPS$%
\indexnotation{GROUPS} denotes the category of groups.
Let $\GROUPS^{\inj}$%
\indexnotation{GROUPS^inj}
be the category of groups with injective group homomorphisms as morphisms.
Taking the reduced group $C^*$-algebra defines a functor
$C_r^{\ast} \colon \GROUPS^{\inj} \to C^{\ast}\text{-}\ALGEBRAS$.
The composition of these functors
with the functors \eqref{K^{alg}RINGS}, \eqref{L^{langle j rangle}RINGS^inv} and
\eqref{K^{alg}C^*ALGEBRAS} above yields functors
\begin{eqnarray}
\bfK R ( - )  \colon \GROUPS & \to & \SPECTRA;
\label{K^{alg}GROUPS}
\\
\bfL^{\langle j \rangle} R ( - ) \colon \GROUPS& \to & \SPECTRA;
\label{LGROUPS}
\\
\bfK^{\topo} C_r^{\ast} ( - )  \colon \GROUPS^{\inj}  & \to & \SPECTRA.
\label{K^{topo}C^*GROUPS^inj}
\end{eqnarray}
They satisfy 
\begin{eqnarray*}
\pi_n(\bfK R(G)) & = & K_n(RG);
\\
\pi_n(\bfL^{\langle j \rangle} R (G)) & = & L_n^{\langle j \rangle}(RG);
\\
\pi_n(\bfK^{\topo} C_r^{\ast} (G)) & = & K_n(C^*_r(G)),
\end{eqnarray*}
for all groups $G$ and $n \in \IZ$.
The next result essentially says that these functors can be extended to groupoids.

\begin{theorem}[$K$- and $L$-Theory Spectra over Groupoids]
\label{the: K- and L-Theory Spectra over Groupoids}
\indextheorem{K- and L-Theory Spectra over Groupoids@$K$- and $L$-Theory Spectra over Groupoids}
Let $R$ be a ring (with involution). There exist covariant functors
\index{K-theory spectrum@$K$-theory spectrum!over groupoids}
\index{L-theory spectrum@$L$-theory spectrum!over groupoids}
\begin{eqnarray}
\bfK_R%
\colon \GROUPOIDS & \to & \SPECTRA;
\label{K^{alg}GROUPOIDS}
\\
\bfL^{\langle j \rangle}_R \colon \GROUPOIDS& \to & \SPECTRA;
\label{LGROUPOIDS}
\\
\bfK^{\topo}\colon \GROUPOIDS^{\inj}  & \to & \SPECTRA
\label{K^{topo}C^*GROUPOIDS^inj}
\end{eqnarray}
with the following properties:

\begin{enumerate}

\item \label{the: K- and L-Theory Spectra over Groupoids: equivalences}
If $F\colon \calg_0 \to \calg_1$ is an equivalence of (small) groupoids, then the induced
maps $\bfK_R (F)$, $\bfL^{\langle j \rangle}_R (F)$ and
$\bfK^{\topo}(F)$ are weak equivalences of spectra.

\item \label{the: K- and L-Theory Spectra over Groupoids: groups and groupoids}
Let $I\colon \GROUPS \to \GROUPOIDS$ be the functor sending $G$ to $G$ considered as a groupoid, i.e.\ to $\calg^G(G/G)$. 
This functor restricts to a functor $\GROUPS^{\inj} \to \GROUPOIDS^{\inj}$.

There are natural transformations from $\bfK R( - )$ to $\bfK_R \circ I$, from 
$\bfL^{\langle j \rangle} R( - )$ to $\bfL^{\langle j \rangle}_R  \circ I$ and from 
$\bfK C_r^{\ast} ( - )$ to $\bfK^{\topo} \circ I$ such that the evaluation of each of these natural transformations
at a given group is an equivalence of spectra.

\item \label{the: K- and L-Theory Spectra over Groupoids: values at groups}
For every group $G$ and all $n \in \IZ$ we have
\begin{eqnarray*}
\pi_n(\bfK_R \circ I(G)) & \cong & K_n(RG);
\\
\pi_n(\bfL^{\langle j \rangle}_R \circ I^{\inv}(G)) & \cong & L_n^{\langle j \rangle}(RG);
\\
\pi_n(\bfK^{\topo} \circ I(G)) & \cong & K_n(C^*_r(G)).
\\
\end{eqnarray*}

\end{enumerate}
\end{theorem}
\begin{proof} We only sketch the strategy of the proof. More details can be
found in \cite[Section 2]{Davis-Lueck(1998)}. 

Let $\calg$ be a groupoid.
Similar to the group ring $RG$ one can define an $R$-linear category $R \calg$ by taking 
the free $R$-modules over the morphism sets of $\calg$. Composition of morphisms is extended $R$-linearly.
By formally adding finite direct sums one obtains an additive category $R \calg_{\oplus}$.
Pedersen-Weibel \cite{Pedersen-Weibel(1985)} (compare also \cite{Cardenas-Pedersen(1997)})
define a non-connective algebraic $K$-theory functor which digests additive categories and can hence be applied to
$R\calg_{\oplus}$. For the comparison result one uses that for every ring $R$ (in particular for $RG$) 
the Pedersen-Weibel functor applied 
to $R_{\oplus}$ (a small model for the category of finitely generated free $R$-modules) 
yields the non-connective $K$-theory of the ring $R$ and that it sends equivalences of 
additive categories to equivalences of spectra. 
In the $L$-theory case
$R\calg_{\oplus}$ inherits an involution
and one applies the construction of 
\cite[Example 13.6 on page 139]{Ranicki(1992)} to obtain the $L^{\langle 1 \rangle}=L^h$-version.
The versions for $j \leq 1$ can be obtained by a construction which is analogous to the Pedersen-Weibel 
construction for $K$-theory, compare \cite[Section~4]{Carlsson-Pedersen(1995a)}.
In the $C^*$-case one obtains from $\calg$ a $C^*$-category $C^*_r(\calg)$ and assigns
to it its topological $K$-theory spectrum.
There is a  construction of the topological $K$-theory spectrum
of a $C^*$-category in \cite[Section 2]{Davis-Lueck(1998)}.
However, the construction given there depends on two statements, which appeared in 
\cite[Proposition 1 and Proposition 3]{Fiedorowicz(1978)},
and those statements are  incorrect, as already pointed out by Thomason in \cite{Thomason(1980extra)}.
The construction in  \cite[Section 2]{Davis-Lueck(1998)} can easily be fixed but instead 
we recommend the reader to look at the more recent construction of Joachim \cite{Joachim(2003a)}.
\end{proof}


\subsection{Assembly Maps in Terms of Homotopy Colimits}
\label{sec: Homotopy-Theoretic Versions of the Conjectures}

In this section we describe a homotopy-theoretic formulation of the Baum-Connes and Farrell-Jones Conjectures.
For the classical assembly maps which in our set-up correspond to the trivial family such formulations
were described in \cite{Weiss-Williams(1995a)}.

For a group $G$ and a family $\calf$ of
subgroups we denote by $\OrGF{G}{\calf}$%
\indexnotation{Or(G;calf)}
the \emph{restricted orbit category}.
\index{orbit category!restricted}
Its objects are homogeneous spaces $G/H$ with $H \in \calf$. Morphisms are $G$-maps. 
If $\calf = \calall$ we get back the (full) orbit category, i.e.\  $\Or(G) =
\OrGF{G}{\calall}$. 

\begin{metaconjecture}[Homotopy-Theoretic Isomorphism Conjecture] 
\label{con: Homotopy-theoretic Isomorphism Conjecture}
\index{Conjecture!Homotopy-theoretic Isomorphism Conjecture for $(G,\bfE,\calf)$}
Let $G$ be a group and $\calf$ a family of subgroups.
Let $\bfE \colon \Or(G) \to \SPECTRA$ be a covariant functor.
Then
\[
A_{\calf} \colon \hocolim_{\OrGF{G}{\calf}} \bfE|_{\OrGF{G}{\calf}} ~ \to ~ 
\hocolim_{\Or(G)} \bfE \simeq  \bfE(G/G)
\]
is a weak equivalence of spectra.
\end{metaconjecture}

Here $\hocolim$ is the homotopy colimit of a covariant functor to spectra, which is itself
a spectrum. The map $A_{\calf}$ is induced by the obvious functor $\OrGF{G}{\calf} \to \Or(G)$.
The equivalence $\hocolim_{\Or(G)} \bfE \simeq \bfE(G/G)$ comes from the fact
that $G/G$ is a final object in $\Or(G)$. For information about homotopy-colimits we refer to 
\cite{Bousfield-Kan(1972)}, \cite[Section 3]{Davis-Lueck(1998)} and \cite{Dror-Farjoun(1987a)}.

\begin{remarknew}
If we consider the map on homotopy groups that is
induced by the map $A_{\calf}$ which appears in the 
Homotopy-Theoretic Isomorphism Conjecture above, then we obtain precisely the map with the same name
in Meta-Conjecture~\ref{metaconjecture} 
for the homology theory $H_*^G( - ; \bfE)$ associated with $\bfE$ in Proposition 
\ref{pro: Or(G)-spectra yield a G-homology theory}, compare \cite[Section 5]{Davis-Lueck(1998)}. 
In particular the Baum-Connes Conjecture~\ref{con: Baum-Connes Conjecture}  
and the Farrell-Jones Conjecture~\ref{con: Farrell-Jones Conjecture} can be seen as special cases of 
Meta-Conjecture~\ref{con: Homotopy-theoretic Isomorphism Conjecture}.
\end{remarknew}

\begin{remarknew}[Universal Property of the Homotopy-Theoretic Assembly Map]
\label{rem: Universal Property of the Homotopy-Theoretic Assembly Map}
The Homotopy-Theoretic Isomorphism Conjecture~\ref{con: Homotopy-theoretic Isomorphism Conjecture} 
is in some sense the most conceptual formulation of an Isomorphism Conjecture because
it has a universal property as the universal approximation from
the left by a (weakly) excisive $\calf$-homotopy invariant functor.
This is explained in detail in \cite[Section 6]{Davis-Lueck(1998)}. 
This universal property is important if one wants to identify different models for the assembly map, compare e.g.\
\cite[Section~6]{Bartels-Farrell-Jones-Reich(2004a)} and \cite{Hambleton-Pedersen(2004)}.

\end{remarknew}


\subsection{Naturality under Induction}
\label{sec: Naturality under Induction with Group Homomorphisms}

Consider a covariant functor $\bfE \colon \GROUPOIDS \to
\SPECTRA$ which respects equivalences.  Let $H^?_*(-;\bfE)$ be the associated equivariant
homology theory (see Proposition~\ref{pro: GROUPOID-spectra and equivariant homology theories}).
Then for a group homomorphism $\alpha \colon H \to G$ and
$H$-$CW$-pair $(X,A)$ we obtain a homomorphism 
\begin{eqnarray*}
\ind_{\alpha}\colon  H_n^H(X,A;\bfE)
&\to &
H_n^G(\ind_{\alpha}(X,A);\bfE) 
\end{eqnarray*}
which is natural in $(X,A)$. Note that we did not assume that $\ker ( \alpha)$ acts freely on $X$.
In fact the construction sketched in the proof of  
Proposition~\ref{pro: GROUPOID-spectra and equivariant homology theories} still works even though
$\ind_{\alpha}$ may not be an isomorphism as it is the case if $\ker (\alpha)$ acts freely. 
We still have functoriality as described in \ref{functoriality} towards 
the beginning of Section~\ref{sec: The Definition of an Equivariant Homology Theory}.

Now suppose that $\calh$ and $\calg$ are families of subgroups for
$H$ and $G$ such that $\alpha(K) \in \calg$ holds for all $K \in
\calh$. Then we obtain a $G$-map $f \colon \ind_{\alpha} \EGF{H}{\calh}
\to \EGF{G}{\calg}$ from the universal property of $\EGF{G}{\calg}$. 
Let $p \colon \ind_{\alpha} \pt = G/\alpha(H) \to
\pt$ be the projection. 
Let $I\colon \GROUPS \to \GROUPOIDS$ be the functor sending $G$ to $\calg^G(G/G)$. 
Then the following diagram, where the horizontal arrows are induced from the projections to the one point space, 
commutes for all $n \in \IZ$.
\[
\comsquare{H^H_n(\EGF{H}{\calh};\bfE)}
{A_{\calh}}{H^H_n(\pt;\bfE) = \pi_n(\bfE(I(H)))} 
{H^G_n(f) \circ \ind_{\alpha}}{H^G_n(p) \circ \ind_{\alpha}  = \pi_n(\bfE(I(\alpha)))}
{H^G_n(\EGF{G}{\calg};\bfE)}
{A_{\calg}}{H^G_n(\pt;\bfE) = \pi_n(\bfE(I(G))).} 
\]
If we take the special case $\bfE = \bfK_R$ and $\calh=\calg=\calvcyc$, 
we get the following commutative
diagram, where the horizontal maps are the assembly maps appearing in the
Farrell-Jones Conjecture~\ref{con: Farrell-Jones Conjecture} and $\alpha_*$ 
is the change of rings homomorphism (induction) associated to $\alpha$.
\[
\comsquare{H^H_n(\EGF{H}{\calvcyc};\bfK_R)}
{A_{\calvcyc}}{K_n(RH)}
{H^G_n(f) \circ \ind_{\alpha}}{\alpha_*}
{H^G_n(\EGF{G}{\calvcyc};\bfK_R)}
{A_{\calvcyc}}{K_n(RG).}
\]
We see that we can define a kind of induction homomorphism on the source of the assembly
maps which is compatible with the induction structure given on their
target. We get analogous diagrams for the $L$-theoretic version of the 
Farrell-Jones-Isomorphism Conjecture~\ref{con: Farrell-Jones Conjecture},
for the Bost Conjecture~\ref{con: Bost Conjecture}
and for the Baum-Connes Conjecture for maximal
group $C^*$-algebras (see \eqref{assembly map for maximal C^*-algebra} in 
Subsection~\ref{subsec: The Baum-Connes Conjecture for Maximal Group $C^*$-Algebras}).

\begin{remarknew}
\label{rem: induction for K_*(C^*_r(-)) and BCC}
The situation for the Baum-Connes Conjecture~\ref{con: Baum-Connes Conjecture} 
itself, where one has to work with reduced $C^*$-algebras, is more complicated.
Recall that not every group homomorphism $\alpha\colon H \to G$ induces
a homomorphisms of $C^*$-algebras $C^*_r(H) \to C^*_r(G)$. 
(It does  if $\ker(\alpha)$ is finite.) But it turns out that
the source $H_n^H(\EGF{H}{\calfin}; \bfK^{\topo})$ 
always admits such a homomorphism.
The point is that the isotropy groups of $\EGF{H}{\calfin}$ are all finite
and the spectra-valued functor $\bfK^{\topo}$ extends from $\GROUPOIDS^{\inj}$ to 
the category $\GROUPOIDS^{\finker}$,%
\indexnotation{GROUPOIDS^{finker}}
which has small groupoids as objects but as morphisms
only those functors $f \colon \calg_0 \to \calg_1$ with finite kernels (in the sense that
for each object $x \in \calg_0$ the group homomorphism $\aut_{\calg_0}(x) \to \aut_{\calg_1}(f(x))$ has
finite kernel).
This is enough to get for any group homomorphism $\alpha\colon H \to G$
an induced map  $\ind_{\alpha}\colon H_n^H(X,A;\bfK^{\topo}) \to H_n^G(\ind_{\alpha}(X,A);\bfK^{\topo})$
provided that $X$ is proper. Hence one can define an induction homomorphism for the source of the assembly 
map as above.

In particular the Baum-Connes Conjecture~\ref{con: Baum-Connes Conjecture}
predicts that for any group homomorphism $\alpha \colon H \to G$ there
is an induced induction homomorphism $\alpha_* \colon K_n(C^*_r(H)) \to
K_n(C^*_r(G))$ on the targets of the assembly maps although there is no induced homomorphism of
$C^*$-algebras $C^*_r(H) \to C^*_r(G)$ in general.
\end{remarknew}


\typeout{----------------------------  Summary of the Applications  ----------------------------}

\section{Methods of Proof}
\label{chap: Methods of Proof}

In Chapter~\ref{chap: general formulation}, we formulated the 
Baum-Connes Conjecture~\ref{con: Farrell-Jones Conjecture}
and the Farrell-Jones Conjecture~\ref{con: Baum-Connes Conjecture} in abstract homological terms.
We have seen that this formulation was very useful in order to understand formal properties
of assembly maps. But in order to actually prove cases of the conjectures one needs 
to interpret the assembly maps in a way that is more directly related to geometry or analysis.
In this chapter  we wish to explain such 
approaches to the assembly maps. We briefly survey some of the 
methods of proof that are used to attack the Baum-Connes Conjecture~\ref{con: Baum-Connes Conjecture}
and the Farrell-Jones Conjecture~\ref{con: Farrell-Jones Conjecture}.


\subsection{Analytic Equivariant $K$-Homology}
\label{sec: Analytic Equivariant $K$-Homology}

Recall that the covariant functor 
$\bfK^{\topo}\colon \GROUPOIDS^{\inj}  \to \SPECTRA$
introduced in \eqref{K^{topo}C^*GROUPOIDS^inj} defines 
an equivariant homology theory $H_*^?(-;\bfK^{\topo})$ in the sense of
Section \ref{sec: The Definition of an Equivariant Homology Theory} such that
$$H_n^G(G/H;\bfK^{\topo}) = H_n^H(\pt;\bfK^{\topo}) = 
\left\{\begin{array}{lll} 
R(H) & \text{ for  even  } $n$;\\
  0  & \text{ for  odd } $n$,
\end{array} 
\right.$$
holds for all groups $G$ and subgroups $H\subseteq G$ 
(see Proposition \ref{pro: GROUPOID-spectra and equivariant homology theories}).
Next we want to give for a
proper $G$-$CW$-complex $X$ an analytic definition of $H_n^G(X;\bfK^{\topo})$.

Consider a locally compact proper $G$-space $X$. Recall that a $G$-space $X$
is called \emph{proper}%
\index{proper!$G$-space}
if for each pair of points $x$ and $y$ in $X$ there are open neighborhoods
$V_x$ of $x$ and $W_y$ of $y$ in $X$ such that the subset
$\{g \in G \mid gV_x \cap W_y \not= \emptyset\}$ of $G$ is finite.
A $G$-$CW$-complex $X$ is proper if and only if all its isotropy groups are finite
\cite[Theorem 1.23]{Lueck(1989)}. 
Let $C_0(X)$%
\indexnotation{C_0(X)}
be the $C^*$-algebra of continuous functions $f \colon X \to \IC$ which vanish at
infinity. The $C^*$-norm is the supremum norm. A \emph{generalized elliptic $G$-operator}%
\index{elliptic operator!generalized elliptic $G$-operator}
is a triple $(U,\rho,F)$,%
\indexnotation{(U,rho,F)}
which consists of a unitary representation 
$U \colon G \to \calb(H)$ of $G$ on a
Hilbert space $H$, a $\ast$-representation $\rho \colon C_0(X) \to \calb(H)$ such that
$\rho(f \circ l_{g^{-1}}) = U(g) \circ \rho(f) \circ U(g)^{-1}$ holds for $g \in G$,
and a bounded selfadjoint $G$-operator $F \colon H
\to H$ such that the operators $\rho(f)(F^2-1)$ and $[\rho(f),F]$ are compact for all 
$f \in C_0(X)$. Here $\calb(H)$ is the $C^*$-algebra of bounded operators $H \to H$,
$ l_{g^{-1}}\colon H \to H$ is given by multiplication with $g^{-1}$,
and $[\rho(f),F] = \rho(f) \circ F - F \circ \rho(f)$.
We also call such a triple $(U,\rho,F)$ an \emph{odd cycle}.%
\index{cycle!odd}
If we additionally assume that $H$ comes with a $\IZ/2$-grading such that 
$\rho$ preserves the grading if we equip $C_0(X)$ with the trivial grading, and $F$
reverses it, then we call $(U,\rho,F)$ an \emph{even cycle}.%
\index{cycle!even} 
This means that we have an orthogonal decomposition $H = H_0 \oplus H_1$
such that $U$, $\rho$ and $F$ look like
\begin{eqnarray} \label{eqn: cycle as matrix}
U = \left( \begin{array}{cc} U_0 & 0 \\ 0 & U_1\end{array} \right)  \hspace{2em}
\rho = \left( \begin{array}{cc} \rho_0 & 0 \\ 0 & \rho_1\end{array} \right)  \hspace{2em}
F = \left( \begin{array}{cc} 0 & P^* \\ P & 0\end{array} \right).
\end{eqnarray}
An important example of an even cocycle is described in Section~\ref{sec: An Example of a Dirac Element}.
A cycle $(U,\rho,f)$ is called \emph{degenerate},%
\index{cycle!degenerate}
if for each $f \in C_0(X)$ we have $[\rho(f),F] = \rho(f)(F^2-1) = 0$. 
Two cycles $(U_0,\rho_0,F_0)$ and $(U_1,\rho_1,F_1)$ of the same parity are called
\emph{homotopic},%
\index{cycle!homotopic}
if $U_0 = U_1$, $\rho_0 = \rho_1$ and there exists a norm continuous path
$F_t, t \in [0,1]$ in $\calb(H)$ such that for each $t \in [0,1]$
the triple $(U_0,\rho_0,F_t)$ is again a cycle of the same parity. Two cycles
$(U_0,\rho_0,F_0)$ and $(U_1,\rho_1,F_1)$ are called \emph{equivalent},%
\index{cycle!equivalent}
if they become homotopic after taking the direct sum with degenerate cycles of the same parity.
Let $K^G_n(C_0(X))$%
\indexnotation{K^G_n(C_0(X))}
for even $n$ be the set of equivalence classes of even cycles and
$K^G_n(C_0(X))$ for odd n  be the set of equivalence classes of odd cycles. These become abelian groups
by the direct sum. The neutral element is represented by any degenerate cycle.
The inverse of an even cycle is represented by the cycle obtained by reversing the
grading of $H$. The inverse of an odd cycle $(U,\rho,F)$ is represented by $(U,\rho,-F)$.

A proper $G$-map $f \colon X \to Y$ induces a map of 
$C^*$-algebras $C_0(f) \colon C_0(Y) \to C_0(X)$ by composition and thus in the obvious
way a homomorphism of abelian groups 
$K_0^G(f) \colon K_0^G(C_0(X)) \to K_0^G(C_0(Y))$. It depends only on the proper
$G$-homotopy class of $f$. One can show that this construction defines 
a $G$-homology theory on the category of finite proper $G$-$CW$-complexes. 
It extends to a $G$-homology theory $K_*^G$ for all proper $G$-$CW$-complexes by 
\begin{eqnarray}
K_n^G(X)%
\indexnotation{K_n^G(X)}
&  =  & \colim_{Y \in I(X)} K_n^G(C_0(Y))
\label{K^G_n(X)}
\end{eqnarray}
where $I(X)$ is the set of proper finite $G$-$CW$-subcomplexes $Y \subseteq X$ directed by
inclusion. This definition is forced upon us by 
Lemma \ref{lem: G-homology theory and colimit}. The groups 
$K_n^G(X)$ and $K_n^G(C_0(X))$ agree for finite proper $G$-$CW$-complexes, in general 
they are different. 

The cycles were introduced by Atiyah \cite{Atiyah(1969a)}.
The equivalence relation, the group structure and the homological properties 
of $K_n^G(X)$ were established by Kasparov \cite{Kasparov(1993)}.  More information about
analytic $K$-homology can be found in Higson-Roe \cite{Higson-Roe(2000)}.


\subsection{The Analytic Assembly Map}
\label{sec: The Analytic Assembly Map}

For for every $G$-CW-complex $X$ the projection $\pr \colon X \to \pt$ induces
a map 
\begin{eqnarray}
H_n^G(X;\bfK^{\topo}) \to H_n^G(\pt;\bfK^{\topo}) = K_n(C^*_r(G)).
\label{homotopic assembly map for topological K-theory}
\end{eqnarray}
In the case where $X$ is the proper $G$-space $\EGF{G}{\calfin}$ we obtain the assembly map appearing in the Baum-Connes 
Conjecture~\ref{con: Baum-Connes Conjecture}. We explain its analytic analogue
\index{assembly map!analytic}
\begin{eqnarray} 
& \ind_G \colon K_n^G(X) \to K_n(C^*_r(G)).& \label{analytic assembly map ind}
\end{eqnarray}
Note that we need to assume that $X$ is a proper $G$-space since $K^G_n(X)$ was only defined for such spaces.
It suffices to define the map for a finite proper $G$-$CW$-complex $X$.
In this case it assigns to the class in $K^G_n(X) = K_n^G(C_0(X))$
represented by a cycle $(U,\rho,F)$ its $G$-index
in $K_n(C^*_r(G))$  in the sense of Mishencko-Fomenko \cite{Mishchenko-Fomenko(1979)}.
At least in the simple case, where $G$ is finite, we can give its
precise definition. The odd $K$-groups vanish in this case and $K_0(C^*_r(G))$ reduces to the complex
representation ring $R(G)$. If we write $F$ in matrix form as in \eqref{eqn: cycle as matrix} then
$P \colon H \to H$ is a $G$-equivariant Fredholm operator.
Hence its kernel and cokernel are $G$-representations and the $G$-index
of $F$ is defined as $[\ker(P)] - [\cok(P)] \in R(G)$. In the case of an infinite group
the kernel and cokernel are a priori not finitely generated projective modules over
$C^*_r(G)$, but they are after a certain pertubation. Moreover the choice of the pertubation
does not affect $[\ker(P)] - [\cok(P)] \in K_0(C^*_r(G))$.

The identification of the two assembly maps 
\eqref{homotopic assembly map for topological K-theory} and 
\eqref{analytic assembly map ind} has been carried out in
Hambleton-Pedersen \cite{Hambleton-Pedersen(2004)} using the universal
characterization of the assembly map explained in \cite[Section 6]{Davis-Lueck(1998)}. 
In particular for a proper $G$-$CW$-complex $X$ we have an identification
$H_n^G(X;\bfK^{\topo}) \cong K_n^G(X)$. Notice that
$H_n^G(X;\bfK^{\topo})$ is defined for all $G$-$CW$-complexes,
whereas $K_n^G(X)$ has only been introduced for proper $G$-$CW$-complexes.

Thus the Baum-Connes Conjecture \ref{con: Baum-Connes Conjecture}
gives an index-theoretic interpretations of elements in $K_0(C^*_r(G))$ as 
generalized elliptic operators or cycles $(U,\rho,F)$. We have explained already
in Subsection~\ref{subsec: The Trace Conjecture in the Torsion Free Case} an application
of this interpretation to the Trace Conjecture for Torsionfree Groups~\ref{con: Trace Conjecture for Torsion Free Groups}
and in Subsection~\ref{subsec: The Stable Gromov-Lawson-Rosenberg Conjecture} to  the
Stable Gromov-Lawson-Rosenberg Conjecture~\ref{con: Stable Gromov-Lawson-Rosenberg Conjecture}.


\subsection{Equivariant $KK$-theory}
\label{sec: Equivariant KK-theory}

Kasparov \cite{Kasparov(1988)} developed \emph{equivariant $KK$-theory},%
\index{KK-theory@$KK$-theory!equivariant} which we will briefly
explain next. It is one of the basic tools in the proofs of theorems about the
Baum-Connes Conjecture \ref{con: Baum-Connes Conjecture}.

A \emph{$G$-$C^*$-algebra}%
\index{C^*-algebra@$C^*$-algebra!G-C^*-algebra@$G$-$C^*$-algebra}
$A$ is a $C^*$-algebra with a $G$-action by $\ast$-automorphisms. 
To any pair of separable $G$-$C^*$-algebras $(A,B)$ Kasparov assigns abelian groups
$KK_n^G(A,B)$.%
\indexnotation{KK_n^G(A,B)}
If $G$ is trivial, we write briefly $KK_n(A,B)$.%
\indexnotation{KK_n(A,B)}
We do not give the rather complicated definition but state the main
properties. 

If we equip $\IC$ with the trivial $G$-action, then $KK^G_n(C_0(X),\IC)$ reduces
to the abelian group $K_n^G(C_0(X))$ introduced in Section
\ref{sec: Analytic Equivariant $K$-Homology}. The topological $K$-theory
$K_n(A)$ of a $C^*$-algebra coincides with $KK_n(\IC,A)$.
The equivariant $KK$-groups are covariant in the second and contravariant in the first variable under
homomorphism of $C^*$-algebras. One of the main features is the bilinear Kasparov product%
\index{Kasparov product}
\begin{eqnarray}
KK^G_i(A,B) \times KK^G_j(B,C) \to  KK_{i+j}(A,C), & 
(\alpha,\beta) \mapsto \alpha \otimes_B \beta.&
\label{Kasparov product}
\end{eqnarray}
It is associative and natural. A homomorphism $\alpha \colon A \to B$ defines an element
in $KK_0(A,B)$.  There are natural \emph{descent homomorphisms}%
\index{descent homomorphism}
\begin{eqnarray}
j_G\colon KK_n^G(A,B) \to KK_n(A\rtimes_r G,B\rtimes_r G),
\label{descent homomorphism}
\end{eqnarray}
where $A \rtimes_r G$ and $B \rtimes_r G$ denote the reduced crossed product $C^{\ast}$-algebras.


\subsection{The Dirac-Dual Dirac Method}
\label{sec: The Dirac-Dual Dirac Method}

A $G$-$C^*$-algebra $A$ is called \emph{proper}
\index{proper!G-C-algebra@$G$-$C^*$-algebra}
if there exists a locally compact proper $G$-space
$X$ and a $G$-homomorphism 
$\sigma \colon C_0(X) \to \calb(A), \hspace{1mm} f \mapsto \sigma_f$ 
satisfying $\sigma_f(ab) = a\sigma_f(b) = \sigma_f(a)b$ for $f \in C_0(X)$, $a,b \in A$ and
for every net $\{f_i \mid i \in I\}$, which converges to $1$ uniformly on compact subsets of
$X$, we have $\lim_{i \in I} \parallel \sigma_{f_i}(a) -  a \parallel ~ = ~ 0$ for all
$a \in A$. 
A locally compact $G$-space $X$ is proper  if
and only if $C_0(X)$ is proper as a $G$-$C^*$-algebra.

Given a proper $G$-$CW$-complex $X$ and a $G$-$C^*$-algebra $A$, we put
\begin{eqnarray}
KK^G_n(X;A) & = & \colim_{Y \in I(X)} KK^G_n(C_0(Y),A),
\label{KK^G_n(X;A)}
\end{eqnarray}
where $I(Y)$ is the set of proper finite $G$-$CW$-subcomplexes $Y \subseteq X$ directed by
inclusion. We have $KK^G_n(X;\IC) = K^G_n(X)$.
There is an analytic index map
\begin{eqnarray} 
& \ind_G^A \colon KK_n^G(X;A) \to K_n(A \rtimes_r G),& 
\label{analytic assembly map ind with coefficients}
\end{eqnarray}
which can be identified with the assembly map appearing in the 
Baum-Connes Conjecture with Coefficients~\ref{con: Baum-Connes Conjecture with coefficients}.
The following result is proved in Tu \cite{Tu(1999a)} 
extending results of
Kasparov-Skandalis \cite{Kasparov-Skandalis(1994)}, \cite {Kasparov-Skandalis(1991)}.

\begin{theorem} \label{the: BCC with coeff for proper A}
The Baum-Connes Conjecture with coefficients
\ref{con: Baum-Connes Conjecture with coefficients} holds for a proper $G$-$C^*$-algebra
$A$, i.e.\ $\ind_G^A \colon KK_n^G(\EGF{G}{\calfin};A) \to K_n(A \rtimes G)$
is bijective.
\end{theorem}

Now we are ready to state the \emph{Dirac-dual Dirac method}
\index{Dirac-dual Dirac method}
which is the key
strategy in many of the proofs of the Baum-Connes Conjecture 
\ref{con: Baum-Connes Conjecture} or the
 Baum-Connes Conjecture with coefficients
\ref{con: Baum-Connes Conjecture with coefficients}.

\begin{theorem}[Dirac-Dual Dirac Method] \label{the: Dirac-Dual Dirac Method}
\indextheorem{Dirac-Dual Dirac Method}
Let $G$ be a countable (discrete) group. Suppose that there exist a proper $G$-$C^*$-algebra
$A$, elements $\alpha \in KK_i^G(A,\IC)$, called the \emph{Dirac element},%
\index{Dirac element} 
and $\beta\in KK_i^G(\IC,A)$, called the \emph{dual Dirac element},%
\index{Dirac element!dual}
satisfying
$$\beta \otimes_A \alpha ~ = ~ 1 \hspace{2mm} \in KK_0^G(\IC,\IC).$$
Then the Baum-Connes Conjecture \ref{con: Baum-Connes Conjecture} is true, or, equivalently,
the analytic index map 
$$\ind_G \colon K_n^G(X) \to K_n(C^*_r(G))$$
of \ref{analytic assembly map ind} is bijective.
\end{theorem}
\begin{proof} The index map $\ind_G$ is a retract of the bijective index map $\ind_G^A$ from Theorem~\ref{the: BCC with coeff for proper A}.
This follows from the following commutative diagram 
\[
\xymatrix{
K_n^G(\EGF{G}{\calfin}) \ar[d]^-{\ind_G} \ar[r]^-{- \otimes_{\IC} \beta} & KK_n^G(\EGF{G}{\calfin};A)
\ar[d]^-{\ind_G^A} \ar[r]^-{- \otimes_{A} \alpha} & K_n^G(\EGF{G}{\calfin}) \ar[d]^-{\ind_G} \\
K_n(C_r^*(G)) \ar[r]^-{-\otimes_{C^*_r(G)} j_G(\beta)} & K_n(A \rtimes_r G) 
\ar[r]^-{- \otimes_{A \rtimes_r} j_G(\alpha)} & K_n(C^*_r(G))
}
\]
and the fact that the composition of both the top upper horizontal arrows
and lower upper horizontal arrows are bijective.
\end{proof}


\subsection{An Example of a Dirac Element}
\label{sec: An Example of a Dirac Element}

In order to give a glimpse of the basic ideas from operator theory we 
briefly describe how to define the Dirac element
$\alpha$ in the case where $G$ acts by isometries on a complete
Riemannian manifold $M$.
Let $T_{\IC}M$ be the complexified tangent bundle and let
$\Cliff(T_{\IC}M)$ be the associated Clifford bundle.
Let $A$ be the proper $G$-$C^*$-algebra given by the sections of
$\Cliff(T_{\IC}M)$ which vanish at infinity. Let $H$ be the Hilbert space $L^2(\wedge T^*_{\IC}M)$
of $L^2$-integrable differential forms on $T_{\IC}M$ with the obvious $\IZ/2$-grading
coming from even and odd forms. Let $U$ be the obvious $G$-representation on $H$
coming from the $G$-action on $M$. For a $1$-form $\omega$ on $M$ and $u \in H$ define
a $\ast$-homomorphism $\rho \colon A \to \calb(H)$ by
$$\rho_{\omega}(u) ~ := ~  \omega \wedge u + i_{\omega}(u).$$
Now $D = (d +d ^*)$ is a symmetric densely  defined operator $H \to H$ and defines
a bounded selfadjoint operator $F \colon H \to H$ by putting $F = \frac{D}{\sqrt{1 + D^2}}$.
Then $(U,\rho,F)$ is an even cocycle and defines an element 
$\alpha \in K_0^G(M) = KK_0^G(C_0(M),\IC)$.
More details of this construction and the construction of the dual Dirac element $\beta$
under the assumption that $M$ has non-positive curvature and is simply connected,
can be found for instance
in \cite[Chapter 9]{Valette(2002)}.


\subsection{Banach KK-Theory}
\label{sec: Banach KK-Theory}

Skandalis showed that the Dirac-dual Dirac method cannot work for all groups
\cite{Skandalis(1988)} as long as one works with $KK$-theory in the unitary setting. The problem is that
for a group with property (T) the trivial and the regular unitary representation cannot be
connected by a continuous path in the space of unitary representations, compare also the discussion in \cite{Julg(1997)}.
This problem can be circumvented if one drops the condition unitary and works with a variant of
$KK$-theory for Banach algebras as worked out by Lafforgue \cite{Lafforgue(1998)},
\cite{Lafforgue(2001)}, \cite{Lafforgue(2002)}.


\subsection{Controlled Topology and Algebra}
\label{sec: Controlled Topology and Algebra}

To a topological problem one can often 
associate a notion of ``size''. 
We describe a prototypical example.
Let $M$ be a Riemannian manifold. 
Recall that an $h$-cobordism $W$ over $M= \partial^{-} W$ admits
retractions $r^{\pm} \colon W \times I \to W$, $(x,t) \mapsto r_t^{\pm} ( x, t )$
which retract $W$ to $\partial^{\pm} W$, i.e.\ which satisfy $r_0^{\pm} = \id_W$ 
and $r_1^{\pm} (W) \subset \partial^{\pm} W$. Given $\epsilon >0$
we say that  $W$ is $\epsilon$-controlled if the retractions can be chosen in such a way that for every $x \in W$ the paths
(called tracks of the $h$-cobordism) $p_x^{\pm} \colon I \to M$, $t \mapsto r_1^{-} \circ r_t^{\pm} (x)$ both lie 
within an $\epsilon$-neighbourhood of their starting point.
The usefulness of this concept is illustrated by the following theorem \cite{Ferry(1977)}. 
\begin{theorem} \label{small-h-cobordism-trivial}
Let $M$ be a compact Riemannian manifold of dimension $\ge 5$. 
Then there exists an $\epsilon=\epsilon_M>0$, such that every $\epsilon$-controlled
$h$-cobordism over $M$ is trivial.
\end{theorem}
If one studies the $s$-Cobordism Theorem~\ref{the: s-cobordism theorem} and its proof one is naturally lead to algebraic
analogues of the notions above. A (geometric) $R$-module over the space $X$ is by definition a family $M=(M_x)_{x \in X}$
of free $R$-modules indexed by points of $X$ with the property that for every compact subset $K \subset X$ the module
$\oplus_{x \in K} M_x$ is a finitely generated $R$-module. A morphism $\phi$ from $M$ to $N$ is an $R$-linear map
$\phi = ( \phi_{y,x} ) \colon \oplus_{x \in X} M_x \to \oplus_{y \in X} N_y$. 
Instead of specifying fundamental group data by paths (analogues of the tracks of the $h$-cobordism) 
one can work with modules and morphisms over 
the universal covering $\widetilde{X}$, which are invariant under the operation of the fundamental group
$G = \pi_1(X)$ via deck transformations, i.e.\ we require that $M_{gx} = M_x$ and $\phi_{gy,gx} = \phi_{y,x}$. 
Such modules and morphisms form an additive category which we denote by $\calc^{G} ( \widetilde{X} ; R )$. 
In particular one can apply to it the non-connective $K$-theory functor $\bfK$ (compare \cite{Pedersen-Weibel(1985)}).
In the case where $X$ is compact the category is equivalent to the category of finitely generated free $RG$-modules and hence
$\pi_{\ast} \bfK \calc^G ( \widetilde{X} ; R ) \cong K_{\ast} ( RG)$. Now suppose $\widetilde{X}$ is equipped with
a $G$-invariant metric, then we will say that a morphism $\phi = ( \phi_{y,x} )$ 
is $\epsilon$-controlled if $\phi_{y,x}=0$, whenever $x$ and $y$ are further than 
$\epsilon$ apart. 
(Note that $\epsilon$-controlled morphisms do not form a category because the composition of two such morphisms will in general
be $2 \epsilon$-controlled.)

Theorem~\ref{small-h-cobordism-trivial} has the following algebraic 
analogue \cite{Quinn(1979a)} (see also Section~4 in \cite{Pedersen(2002)}).
\begin{theorem}
Let $M$ be a compact Riemannian manifold with fundamental group $G$. There exists an 
$\epsilon= \epsilon_M >0$ with the following property.
The $K_1$-class of every $G$-invariant automorphism of modules over $\widetilde{M}$
which together with its inverse is $\epsilon$-controlled
lies in the image of the classical assembly map
\[
H_1 ( BG ; \bfK R )
\to K_1( RG ) \cong  K_1 ( \calc^G ( \widetilde{M} ; R ) ).
\]
\end{theorem}
To understand the relation to Theorem~\ref{small-h-cobordism-trivial} note that for $R= \IZ$ such an $\epsilon$-controlled
automorphism represents the trivial element in the Whitehead group which is in bijection with the $h$-cobordisms over $M$, compare
Theorem~\ref{the: s-cobordism theorem}.

There are many variants to the  simple concept of ``metric $\epsilon$-control'' we used above. 
In particular it is very useful to not measure size directly in $M$ but instead use
a map $p \colon M \to X$ to measure size in some auxiliary space $X$. (For example we have seen in 
Subsection~\ref{subsec: bounded h-cobordisms} and \ref{subsec: negative K pseudo} that ``bounded'' 
control over $\IR^k$ may be used in order to define or describe negative $K$-groups.)

Before we proceed we would like to mention that there are 
analogous control-notions for pseudoisotopies and homotopy equivalences.
The tracks of a pseudoisotopy $f \colon M \times I \to M \times I$
are defined as the paths in $M$ which are given by the composition
\[
\xymatrix{
p_x \colon I = \{ x \} \times I \subset M \times I \ar[r]^-f & M \times I \ar[r]^-p & M
         }
\]
for each $x \in M$, where the last map is the projection onto the $M$-factor. 
Suppose $f \colon N \to M$ is a homotopy equivalence, $g \colon M \to N$ its inverse and $h_t$ and $h^{\prime}_t$ are homotopies from
$f \circ g$ to $\id_M$ respectively from $g \circ f$ to $\id_N$ then the tracks are defined to be the paths in $M$ that are given
by $t \mapsto h_t ( x )$ for $x \in M$ and $t \mapsto f \circ h^{\prime}_t(y)$ for $y \in N$. 
In both cases, for pseudoisotopies and for homotopy equivalences, the tracks can be used to define 
$\epsilon$-control.

\subsection{Assembly as Forget Control}
\label{sec:Assembly as forget control}

If instead of a single problem over $M$ one defines a family of problems over $M \times [ 1 , \infty)$ and requires
the control to tend to zero for $t \to \infty$ in a suitable sense, then one obtains something which is a homology theory in $M$.
Relaxing the control to just bounded families yields the classical assembly map.
This idea appears in \cite{Quinn(1982a)} in the context of pseudoisotopies 
and in a more categorical fashion suitable for higher algebraic $K$-theory in 
\cite{Carlsson-Pedersen(1995a)} and 
\cite{Pedersen-Weibel(1989)}.
We spell out some details in the case of algebraic $K$-theory, i.e.\ for geometric modules. 

Let $M$ be a Riemannian manifold with fundamental group $G$ and let $\cals(1/t)$ be the space of 
all functions $[ 1 , \infty) \to [0 , \infty)$, $t \mapsto \delta_t$ such that $t \mapsto t \cdot \delta_t$ is bounded. 
Similarly let $\cals( 1 )$ be the space of all functions $t \mapsto \delta_t$ which are bounded.
Note that $\cals ( 1/t ) \subset \cals ( 1 )$.
A $G$-invariant morphism $\phi$ over $\widetilde{M} \times [1 , \infty)$ is $\cals$-controlled for 
$\cals= \cals(1)$ or $\cals( 1/t)$
if there exists an $\alpha>0$ and a 
$\delta_t \in \cals$ (both depending on the morphism) such that $\phi_{(x,t),(x',t')} \neq 0$ implies that $|t-t'| \leq \alpha$ and 
$d_{\widetilde{M}} ( x , x') \leq \delta_{\min \{ t,t' \} }$. We denote by 
$\calc^G ( \widetilde{M} \times [ 1 , \infty ) , \cals ; R )$ the category of all $\cals$-controlled morphisms. 
Furthermore $\calc^G ( \widetilde{M} \times [ 1 , \infty ) , \cals ; R )^{\infty}$ denotes the quotient category which has the 
same objects, but where two morphisms are identified, if their difference factorizes 
over an object which lives over $\widetilde{M} \times [1, N]$
for some large but finite number $N$. This passage to the quotient category is called ``taking germs at infinity''. 
It is a special case of a Karoubi quotient, compare \cite{Cardenas-Pedersen(1997)}.

\begin{theorem}[Classical Assembly as Forget Control]
\label{the: assembly as forget control for K}
Suppose $M$ is aspherical, i.e.\ $M$ is a model for $BG$, then for all $n \in \IZ$ the map
\[
\pi_{n} ( \bfK  \calc^{G} ( \widetilde{M} \times [ 1 , \infty ) , \cals( 1/t ) ; R )^{\infty} ) \to
\pi_{n} ( \bfK  \calc^{G} ( \widetilde{M} \times [ 1 , \infty ) , \cals(1) ; R )^{\infty} )
\]
can be identified up to an index shift with the classical assembly map that appears in 
Conjecture~\ref{con: FJK torsion free all}, i.e.\ with
\[
H_{n-1} ( BG ; \bfK (R) ) \to K_{n-1} ( RG ).
\]
\end{theorem}
Note that the only difference between the left and the right hand side is that on the left morphism are required to become smaller
in a $1/t$-fashion, whereas on the right hand side they are only required to stay bounded in the $[1, \infty)$-direction.

Using so called equivariant continuous control 
(see \cite{Anderson-Connolly-Ferry-Pedersen(1994)} and \cite[Section~2]{Bartels-Farrell-Jones-Reich(2004a)} for the equivariant version)
one can define an equivariant homology theory which applies to 
arbitrary $G$-CW-complexes. This leads to a ``forget-control description'' for the generalized assembly maps that appear in the 
Farrell-Jones Conjecture~\ref{con: Farrell-Jones Conjecture}. 
Alternatively one can use stratified spaces and stratified Riemannian manifolds in order to describe generalized
assembly maps in terms of metric control. Compare \cite[3.6 on p.270]{Farrell-Jones(1993a)}
and \cite[Appendix]{Quinn(1982a)}.

\subsection{Methods to Improve Control}
\label{sec: Methods to improve control}

From the above description of assembly maps we learn that the problem of proving surjectivity results translates into the problem
of improving control.
A combination of many different techniques is used in order to 
achieve such control-improvements. We discuss some prototypical arguments
which go back to \cite{Farrell-Hsiang(1978b)} and  \cite{Farrell-Jones(1986a)}
and again restrict attention to $K$-theory.
Of course this can only give a very incomplete impression of the whole
 program which is mainly due to Farrell-Hsiang and Farrell-Jones.
The reader should consult \cite{Farrell(2002)} and \cite{Jones(2002)} for a more detailed survey.

We restrict to the basic case, where $M$ is a compact Riemannian manifold with negative sectional curvature. 
In order to explain a contracting property of the 
geodesic flow $\Phi \colon \IR \times S \widetilde{M} \to S \widetilde{M}$ on the 
unit sphere bundle $S \widetilde{M}$,
we introduce the notion of foliated control. We think of $S \widetilde{M}$ as 
a manifold equipped with the one-dimensional foliation by the 
flow lines of $\Phi$ and equip it with its natural Riemannian 
metric. Two vectors $v$ and $w$ in $S \widetilde{M}$ 
are called foliated $(\alpha , \delta)$-controlled if there 
exists a path of length $\alpha$ inside one flow line such that 
$v$ lies within distance $\delta/2$ of the starting point of 
that path and $w$ lies within distance $\delta/2$ of its endpoint.

Two vectors $v$ and $w \in S \widetilde{M}$ are called asymptotic if the distance between their associated 
geodesic rays is bounded. These rays will then determine the same point on the sphere at infinity which can be 
introduced to compactify $\widetilde{M}$ to a disk.
Recall that the universal covering  of a negatively curved manifold is diffeomorphic to $\IR^n$.
Suppose $v$ and $w$ are $\alpha$-controlled asymptotic vectors, i.e.\ their distance is smaller than
$\alpha > 0$. As a consequence of negative sectional curvature the vectors $\Phi_t (v)$ and $\Phi_t(w)$ are 
foliated $(C \alpha , \delta_t )$-controlled, where $C>1$ is a 
constant and $\delta_t>0$ tends to zero when $t$ tends to $\infty$.
So roughly speaking the flow contracts the directions transverse to the flow lines and leaves the 
flow direction as it is, at least if we only apply it to asymptotic vectors.

This property can be used in order to find foliated $(\alpha , \delta)$-controlled 
representatives of $K$-theory classes 
with arbitrary small $\delta$
if one is able to define a suitable 
transfer from $M$ to $S \widetilde{M}$, which yields representatives 
whose support is in an asymptotic starting position for the flow. 
Here one needs to take care of the additional problem that in general 
such a transfer may not induce an isomorphism in $K$-theory.

Finally one is left with the problem of improving foliated control 
to ordinary control. Corresponding statements  are called
``Foliated Control Theorems''. Compare 
\cite{Bartels-Farrell-Jones-Reich(2003b)},
\cite{Farrell-Jones(1986)}, 
\cite{Farrell-Jones(1987c)}, 
\cite{Farrell-Jones(1988a)} and
\cite{Farrell-Jones(1991b)}.

If such an improvement were possible without further hypothesis,
we could prove that the classical assembly map, i.e.\ the assembly map  
with respect to the trivial family is surjective.
We know however that this is not true in general. 
It fails for example in the case of topological pseudoisotopies or for algebraic $K$-theory with 
arbitrary coefficients. In fact the geometric arguments  
that are involved in a ``Foliated Control Theorem'' need
to exclude flow lines in $S \widetilde{M}$ which correspond to 
``short'' closed geodesic loops in $S M$. 
But the techniques mentioned above can be used in order to achieve $\epsilon$-control 
for arbitrary small $\epsilon>0$ outside of a suitably chosen neighbourhood of ``short'' closed geodesics. 
This is the right kind of control for the source of the 
assembly map which involves the family of cyclic subgroups.
(Note that a closed a loop in $M$ determines the conjugacy class of a maximal infinite cyclic subgroup  
inside $G= \pi_1(M)$.) We see that  even in the torsionfree case the class 
of cyclic subgroups of $G$ naturally appears during the proof of 
a surjectivity result.

Another source for processes which improve control are 
expanding self-maps.
Think for example of an $n$-torus $\IR^n / \IZ^n$ and the 
self-map $f_s$ which is induced by $m_s \colon \IR^n \to \IR^n$, $x \to sx$
for a large positive integer $s$. If one pulls an automorphism back along such a map one can improve control, but unfortunately
the new automorphism describes a different $K$-theory class. 
Additional algebraic arguments nevertheless make this technique very 
successful.
Compare for example \cite{Farrell-Hsiang(1978b)}.
Sometimes a clever mixture between flows and expanding self-maps is needed in order to achieve the goal, 
compare \cite{Farrell-Jones(1988b)}. Recent work of Farrell-Jones (see 
\cite{Farrell-Jones(1998f)}, 
\cite{Farrell-Jones(1998)},
\cite{Farrell-Jones(2003)} and  
\cite{Jones(2003)}) makes use of a variant of the Cheeger-Fukaya-Gromov collapsing theory.

\begin{remarknew}[Algebraicizing the Farrell-Jones Approach]
\label{rem:Algebraicizing the Farrell-Jones Approach}
In this Subsection we sketched some of the geometric ideas which are used in order to 
obtain control over an $h$-cobordism, a pseudisotopy or an automorphism of a geometric module 
representing a single class in $K_1$. In Subsection \ref{sec:Assembly as forget control} we
used families over the cone $M \times [1, \infty)$ in order to 
described the whole algebraic K-theory assembly map at once in categorical terms without ever referring to 
a single K-theory element. 
The recent work \cite {Bartels-Reich(2003)} shows that the geometric 
ideas can be adapted to this more categorical set-up, at least in the case where the group is the 
fundamental group of a 
Riemannian manifold with strictly negative curvature. However serious difficulties
 had to be overcome in order to achieve this. 
One needs to formulate and prove a Foliated Control Theorem in this context and also
construct a transfer map to the sphere bundle for higher K-theory which is in a 
suitable sense compatible
 with the control structures.
\em
\end{remarknew}

\subsection{The Descent Principle}
\label{sec: The descent principle}

In Theorem~\ref{the: assembly as forget control for K} we described the classical assembly map as a forget control map
using $G$-invariant geometric modules over $\widetilde{M} \times [1 , \infty)$. 
If in that context  one does not require the modules and morphisms
to be invariant under the $G$-action one nevertheless obtains a forget control functor between additive categories for which we 
introduce the notation 
\[
\cald (1/t) = \calc ( \widetilde{M} \times [ 1 , \infty ) , \cals( 1/t ) ; R )^{\infty}  \to
\cald ( 1 ) = \calc ( \widetilde{M} \times [ 1 , \infty ) , \cals(1) ; R )^{\infty} .
\]
Applying $K$-theory yields a version of a ``coarse'' assembly map which is the algebraic $K$-theory analogue of the map
described in Section~\ref{subsec: The Coarse Baum Connes Conjecture}. 
A crucial feature of such a construction is that the 
left hand side can be interpreted as a locally finite homology theory 
evaluated on $\widetilde{M}$. It is hence an invariant of the proper homotopy 
type of $\widetilde{M}$. Compare 
\cite{Anderson-Connolly-Ferry-Pedersen(1994)} and \cite{Weiss(2002)}.
It is usually 
a lot easier to prove that this coarse assembly map is an equivalence.
Suppose for example that $M$ has non-positive curvature, choose a point $x_0 \in M$ (this will destroy the $G$-invariance)
and with increasing $t \in [1, \infty)$ move the modules along geodesics towards $x_0$. In this way one can show that the
coarse assembly map is an isomorphism. Such coarse assembly maps exist also in the context of algebraic $L$-theory and topological
$K$-theory, compare \cite{Higson-Pedersen-Roe(1997)}, \cite{Roe(1996)}.

Results about these maps 
(compare e.g.\ \cite{Bartels(2003a)}, \cite{Carlsson-Pedersen(1995a)}, \cite{Yu(1998)}, \cite{Yu(2000)})
lead to injectivity results for the classical assembly map by the ``descent principle'' 
(compare \cite{Carlsson(2004)}, \cite{Carlsson-Pedersen(1995a)}, \cite{Roe(1996)}) 
which we will now briefly describe in the context of algebraic $K$-theory.
(We already stated an analytic version in Section~\ref{subsec: The Coarse Baum Connes Conjecture}.)
For a spectrum $\bfE$ with $G$-action we denote by $\bfE^{hG}$ the homotopy fixed points.
Since there is a natural map from fixed points to homotopy fixed points we obtain a commutative diagram
\[
\xymatrix{
\bfK (\cald(1/t))^G  \ar[r] \ar[d] &
\bfK (\cald(1))^G \ar[d] . \\
\bfK (\cald(1/t))^{hG}  \ar[r] & 
\bfK (\cald(1))^{hG}.
         }
\]
If one uses a suitable model $K$-theory commutes with taking fixed points and hence the upper horizontal map
can be identified with the classical assembly map by Theorem~\ref{the: assembly as forget control for K}.
Using that $K$-theory commutes with infinite products \cite{Carlsson(1995)}, one can show by an 
induction over equivariant cells, that the vertical map on the left is an equivalence.
Since we assume that the map $\bfK (\cald(1/t)) \to \bfK (\cald(1))$ is an equivalence, a standard property of 
the homotopy fixed point construction implies that the lower horizontal map is an equivalence. It follows that
the upper horizontal map and hence the classical assembly map is split injective.
A version of this argument which involves the assembly map for the family of finite subgroups can be found in
\cite{Rosenthal(2002)}.

\subsection{Comparing to Other Theories}
\label{sec: Comparing to Other Theories}

Every natural transformation of $G$-homology theories  leads to a comparison between the associated 
assembly maps. For example one can compare topological $K$-theory
to periodic cyclic homology \cite{Connes-Moscovici(1990)}, i.e.\ for every Banach algebra 
completion $\cala(G)$ of $\IC G$ inside $C_r^{\ast} (G)$ there exists a commutative diagram
\[
\xymatrix{
K_{\ast}( BG ) \ar[d] \ar[r] & K_{\ast} ( \cala (G)) \ar[d] \\
H_{\ast} ( BG ; HP_{\ast}( \IC )) \ar[r] & HP_{\ast} ( \cala (G) ) .
}
\]
This is used in \cite{Connes-Moscovici(1990)} to prove injectivity results for word hyperbolic groups. 
Similar diagrams exist for other cyclic theories (compare for example \cite{Puschnigg(2002)}).

A suitable model for the cyclotomic trace $trc \colon K_n ( RG ) \to TC_n (RG )$ from (connective)
algebraic $K$-theory to topological cyclic homology \cite{Boekstedt-Hsiang-Madsen(1993)} leads for 
every family $\calf$ to a commutative diagram
\[
\xymatrix{
H_n ( \EGF{G}{\calf} ; \bfK^{con}_{\IZ} ) \ar[d] \ar[r] & K^{con}_n ( \IZ G ) \ar[d] \\
H_n ( \EGF{G}{\calf} ; \bfT \bfC_{\IZ} ) \ar[r] & TC_n ( \IZ G ).
         }
\]
Injectivity results about the left hand and the lower horizontal 
map lead to injectivity results about the upper horizontal map. This is the principle behind
Theorem~\ref{the: Boekstedt-Hsiang-Madsen} and \ref{the: LRRM}.


\typeout{----------------------------  Computations ----------------------------}
\section{Computations}
\label{chap: Computations}

Our ultimate goal is to compute $K$- and $L$-groups such as $K_n(RG)$,
$L_n^{\langle -\infty \rangle}(RG)$ and $K_n(C^*_r(G))$. 
Assuming that the Baum-Connes Conjecture~\ref{con: Baum-Connes Conjecture} 
or the Farrell-Jones Conjecture~\ref{con: Farrell-Jones Conjecture}
is true, this reduces to  the computation of the left hand side of the corresponding assembly map, i.e.\ to 
$H_n^G(\EGF {G}{\calfin};\bfK^{\topo})$, 
$H_n^G(\EGF{G}{\calvcyc};\bfK_R)$
and
$H_n^G (\EGF{G}{\calvcyc};\bfL_R^{\langle - \infty \rangle})$.
This is much easier since here we can use standard methods from algebraic topology such as
spectral sequences, Mayer-Vietoris sequences and Chern characters.
Nevertheless such computations can be pretty hard. Roughly speaking, one can obtain a
general reasonable answer after rationalization, but integral computations
have only been done case by case and no general pattern is known.


\subsection{$K$- and $L$- Groups for Finite Groups}
\label{sec: K- and L- Groups for Finite Groups}

In all these computations the answer is given in terms of the
values of $K_n(RG)$,
$L_n^{\langle -\infty \rangle}(RG)$ and $K_n(C^*_r(G))$ for finite groups $G$. Therefore
we briefly recall some of the results known for finite groups focusing on the case
$R = \IZ$

\subsubsection{Topological $K$-Theory for Finite Groups}
\label{subsec: Topological K-Theory for Finite Groups}

Let $G$ be a finite group. By 
$r_F (G)$,%
\indexnotation{r_F(G)}
we denote the number of isomorphism classes of irreducible 
representations of $G$ over the field $F$. By
\indexnotation{r_{IR}(G;IR)}
$r_{\IR}(G;\IR)$,
\indexnotation{r_{IR}(G;IC)}
$r_{\IR}(G;\IC)$,
respectively
\indexnotation{r_{IR}(G;IH)}
$r_{\IR}(G;\IH)$
we denote the number of isomorphism classes of irreducible real $G$-representations $V$,
which are of real, complex respectively of quaternionic type,
i.e.\ $\aut_{\IR G}(V)$ is isomorphic to the field of real numbers $\IR$,
complex numbers $\IC$ or quaternions $\IH$.
Let $RO(G)$%
\indexnotation{RO(G)}
respectively $R(G)$%
\indexnotation{R(G)}
be the \emph{real} respectively the \emph{complex representation ring}.%
\index{representation ring} 


Notice that 
$\IC G = l^1(G) = C^*_r(G) = C^*_{\max}(G)$ holds for a finite group, and analogous
for the real versions. 

\begin{propositionnew}
\label{pro: Topological K-theory for finite groups}

Let $G$ be a finite group.

\begin{enumerate}
\item \label{pro: Topological K-theory for finite groups: complex}
We have
$$K_n(C^*_r(G)) ~ \cong ~ 
\left\{\begin{array}{lllll} R(G) & \cong & \IZ^{r_{\IC}(G)} & & \text{ for } n \text{ even};
\\
0 & & & & \text{ for } n \text{ odd}.
\end{array}
\right.
$$ 

\item \label{pro: Topological K-theory for finite groups: real}
There is an isomorphism of topological $K$-groups
$$K_n(C_r^*(G;\IR)) ~ \cong ~ 
K_n(\IR)^{r_{\IR}(G; \IR)} \times K_n(\IC)^{r_{\IR}(G;\IC)} \times  K_n(\IH)^{r_{\IR}(G ; \IH)}.$$
Moreover $K_n(\IC)$ is $2$-periodic with values $\IZ$, $0$ for $n = 0,1$,
$K_n(\IR)$ is $8$-periodic with values $\IZ$, $\IZ/2$, $\IZ/2$, $0$, $\IZ$, $0$, $0$, $0$
for $n = 0,1, \ldots , 7$ and $K_n(\IH) = K_{n + 4}(\IR)$ for $n \in
\IZ$. 

\end{enumerate}
\end{propositionnew}
\begin{proof} One gets isomorphisms of $C^*$-algebras
$$C^*_r(G) ~ \cong ~ \prod_{j=1}^{r_{\IC}(G)} M_{n_i}(\IC)$$ 
and
$$C^*_r(G;\IR)~ \cong ~ 
\prod_{i=1}^{r_{\IR}(G; \IR)} M_{m_i}(\IR) \times \prod_{i=1}^{r_{\IR}(G; \IC)} M_{n_i}(\IC) \times 
\prod_{i=1}^{r_{\IR}(G;\IH)} M_{p_i}(\IH)$$ 
from \cite[Theorem~7 on page 19, Corollary 2 on page 96, page 102, page 106]{Serre(1977)}.
Now the claim follows from Morita invariance and the well-known
values for $K_n(\IR)$, $K_n(\IC)$ and $K_n(\IH)$
(see for instance \cite[page 216]{Switzer(1975)}).
\end{proof}

To summarize, the values of $K_n(C^*_r(G))$ and $K_n(C^*_r(G;\IR))$ are explicitly known for
finite groups $G$ and are in the complex case in contrast to the real case always torsion
free.


\subsubsection{Algebraic $K$-Theory for Finite Groups}
\label{subsec: Algebraic K-Theory for Finite Groups}

Here are some facts about the algebraic $K$-theory of integral group rings of finite groups. 

\begin{propositionnew} 
\label{pro: Algebraic K-theory for finite groups}
Let $G$ be a finite group.
\begin{enumerate}

\item \label{pro: Algebraic K-theory for finite groups: K_q(ZG):q le -2}
$K_n(\IZ G) = 0$ for $n \le -2$.

\item \label{pro: Algebraic K-theory for finite groups: K_q(ZG):q = -1}
We have
\[
K_{-1}(\IZ G) \cong \IZ^r  \oplus (\IZ/2)^s,
\]
where 
\[
r = 1  - r_{\IQ}(G) + \sum_{p \; \mid \; |G|}  r_{\IQ_p}(G)  - r_{\IF_p}(G)
\]
and the sum runs over all primes dividing the order of $G$.
(Recall that $r_F (G)$ denotes the number of isomorphism classes of irreducible representations of $G$ over the field $F$.)
There is an explicit description of the integer $s$ in terms of global and local Schur indices \cite{Carter(1980)}.
If $G$ contains a normal abelian subgroup of odd index, then $s = 0$.

\item \label{pro: Algebraic K-theory for finite groups: K_0(ZG)}
The group $\widetilde{K}_0(\IZ G)$ is finite. 

\item \label{pro: Algebraic K-theory for finite groups: K_q(ZG): q=1}
The group $\Wh(G)$ is a finitely generated abelian group and its rank is $r_{\IR}(G) - r_{\IQ}(G)$.

\item \label{pro: Algebraic K-theory for finite groups: K_q(ZG): all q}
The groups $K_n(\IZ G)$ are finitely generated for all $n \in \IZ$.

\item \label{pro: Algebraic K-theory for finite groups: K_q(ZG): special G}
We have $K_{-1}(\IZ G )=0$, $\widetilde{K}_0( \IZ G)=0$ and $\Wh(G) = 0$ for the following finite groups
$G = $
$\{1\}$, $\IZ/2$, $\IZ/3$, $\IZ/4$, $\IZ/2 \times \IZ/2$, $D_6$, $D_8$,
where $D_m$ is the dihedral group of order $m$.

If $p$ is a prime, then
$K_{-1}(\IZ[\IZ/p]) = K_{-1}(\IZ[\IZ/p \times \IZ/p]) = 0$.

We have
\begin{eqnarray*}
K_{-1}( \IZ[\IZ /6] ) \cong \IZ, \quad \widetilde{K}_0( \IZ [\IZ/6] ) = 0, \quad \Wh(\IZ/6) = 0 \\
K_{-1}( \IZ[D_{12}] ) \cong \IZ, \quad \widetilde{K}_0( \IZ [D_{12}] ) = 0, \quad \Wh(D_{12}) = 0.
\end{eqnarray*}

\item \label{pro: Algebraic K-theory for finite groups: K_q(ZG): Wh_2}
Let $\Wh_2(G)$ denote the cokernel of the assembly map 
\[
H_2( BG ; \bfK(\IZ)) \to K_2( \IZ G).
\]
We have $\Wh_2(G)=0$ for $G = \{1\}$, $\IZ/2$, $\IZ/3$ and $\IZ/4$. Moreover
$|\Wh_2(\IZ/6)| \le  2$,
$|\Wh_2(\IZ/2 \times \IZ/2)| \ge 2$ and 
$\Wh_2(D_6) = \IZ/2$.
\end{enumerate}

\end{propositionnew}
\begin{proof}
\ref{pro: Algebraic K-theory for finite groups: K_q(ZG):q le -2} and 
\ref{pro: Algebraic K-theory for finite groups: K_q(ZG):q = -1} are proved 
in \cite{Carter(1980)}. 
\\[1mm]
\ref{pro: Algebraic K-theory for finite groups: K_0(ZG)} is proved in
\cite[Proposition~9.1 on page 573]{Swan(1960a)}.
\\[1mm]
\ref{pro: Algebraic K-theory for finite groups: K_q(ZG): q=1}
This is proved for instance in \cite{Oliver(1989)}.
\\[1mm]
\ref{pro: Algebraic K-theory for finite groups: K_q(ZG): all q}
See \cite{Kuku(1986)}, \cite{Quillen(1973a)}.
\\[1mm]
\ref{pro: Algebraic K-theory for finite groups: K_q(ZG): special G} and
\ref{pro: Algebraic K-theory for finite groups: K_q(ZG): Wh_2}
The computation  $K_{-1}(\IZ G)=0$ for $G =\IZ/p$ or $\IZ/p
\times\IZ/p$ can be found in \cite[Theorem~10.6, p. 695]{Bass(1968)}
and is a special case of \cite{Carter(1980)}.

The vanishing of $\widetilde{K}_0(\IZ G)$ is proven for $G = D_6$ in
\cite[Theorem~8.2]{Reiner(1976)} 
and for $G =D_8$ in
\cite[Theorem~6.4]{Reiner(1976)}. The cases $G=\IZ/2, \IZ/3, \IZ/4,
\IZ/6$, and $(\IZ/2)^2$ are treated in \cite[Corollary 5.17]{Curtis-Reiner(1987)}. 
Finally, $\widetilde{K}_0(\IZ D_{12})=0$
follows from \cite[Theorem~50.29 on page 266]{Curtis-Reiner(1987)} and
the fact that $\IQ D_{12}$ as a $\IQ$-algebra splits into copies of
$\IQ$ and matrix algebras over $\IQ$, so its maximal order has
vanishing class group by Morita equivalence.

The claims about $\Wh_2(\IZ/n)$  for $n= 2,3,4,6$ and for
$\Wh_2((\IZ/2)^2)$ are taken from
\cite[Proposition~5]{Dennis-Keating-Stein(1976)},
\cite[p.482]{Dunwoody(1975)} and
\cite[p. 218 and 221]{Stein(1980)}.

We get $K_2(\IZ D_6) \cong (\IZ/2)^3$ from
\cite[Theorem~3.1]{Stein(1980)}.
The assembly map $H_2(B\IZ/2;\bfK(\IZ)) \to K_2(\IZ[\IZ/2])$
is an isomorphism by
\cite[Theorem~on p. 482]{Dunwoody(1975)}.
Now construct a commutative diagram
$$\comsquare
{H_2(B\IZ/2; \bfK(\IZ))}{\cong }{H_2(BD_6; \bfK(\IZ))}
{\cong}{}
{K_2(\IZ[\IZ/2])}{}{K_2(\IZ D_6)}
$$
whose lower horizontal arrow is split injective
and whose upper horizontal arrow is an isomorphism
by the Atiyah-Hirzebruch spectral sequence.
Hence the right vertical arrow is split injective and $\Wh_2(D_6) = \IZ/2$.
\end{proof}

Let us summarize. We already mentioned that a complete computation of
$K_n(\IZ)$ is not known. Also a complete computation of
$\widetilde{K}_0(\IZ[\IZ/p])$ for arbitrary primes $p$ is out of reach 
(see \cite[page 29,30]{Milnor(1971)}). 
There is a complete formula for $K_{-1}(\IZ G)$
and $K_n(\IZ G) = 0$ for $n \le -2$ and 
one has a good understanding of $\Wh(G)$
(see \cite{Oliver(1989)}).  We have already mentioned Borel's formula for
$K_n(\IZ) \otimes_{\IZ} \IQ$ for all $n \in \IZ$
(see Remark~\ref{rem: coefficients K}). For more rational information see also \ref{rem: already interesting for finite}.


\subsubsection{Algebraic $L$-Theory for Finite Groups}
\label{subsec: Algebraic L-Theory for Finite Groups}

Here are some facts about $L$-groups of finite groups.
\begin{propositionnew} Let $G$ be a finite group. Then
\label{pro: Algebraic L-theory for finite groups}
\begin{enumerate}
\item \label{pro: Algebraic L-theory for finite groups: fin. gen 2-tors}
For each $j \le 1$ the groups $L_n^{\langle j \rangle }(\IZ G)$ are finitely generated as
abelian groups and contain no $p$-torsion for odd primes $p$. Moreover, they are finite for odd $n$.
\item \label{pro: Algebraic L-theory for finite groups: after inverting 2}
For every decoration $\langle j \rangle$ we have
\begin{eqnarray*}
L_n^{\langle j \rangle}(\IZ G)[1/2] ~ \cong ~ L_n^{\langle j \rangle} (\IR G)[1/2] & \cong &
\left\{
\begin{array}{lcl}
\IZ[1/2]^{r_{\IR}(G)} & & n \equiv 0 \hspace{2mm} (4);
\\
\IZ[1/2]^{r_{\IC}(G)} & & n \equiv 2 \hspace{2mm} (4);
\\
0 & & n \equiv 1,3 \hspace{2mm} (4).
\end{array}
\right.
\end{eqnarray*}
\item \label{pro: Algebraic L-theory for finite groups: odd order}
If $G$ has odd order and $n$ is odd, then $L_n^{\epsilon}(\IZ G) = 0$
for $\epsilon = p,h,s$.
\end{enumerate}
\end{propositionnew}
\begin{proof}
\ref{pro: Algebraic L-theory for finite groups: fin. gen 2-tors} 
See for instance \cite{Hambleton-Taylor(2000)}.
\\[1mm]
\ref{pro: Algebraic L-theory for finite groups: after inverting 2}
See \cite[Proposition~22.34 on page 253]{Ranicki(1992)}.
\\[1mm]
\ref{pro: Algebraic L-theory for finite groups: odd order}
See \cite{Bak(1975)} or \cite{Hambleton-Taylor(2000)}.
\end{proof}

The $L$-groups of $\IZ G$ are pretty well understood for finite groups $G$.
More information about them can be found in \cite{Hambleton-Taylor(2000)}.


\subsection{Rational Computations for Infinite Groups}
\label{sec: Rational Computations for Infinite Groups}

Next we state what is known rationally about the $K$- and $L$-groups
of an infinite (discrete) group, provided the Baum-Connes Conjecture 
\ref{con: Baum-Connes Conjecture} or the relevant version of the
Farrell-Jones Conjecture~\ref{con: Farrell-Jones Conjecture} is known.

In the sequel let $(\calfcyc)$%
\indexnotation{(calfcyc)}
be the set of conjugacy classes $(C)$ for finite cyclic subgroups $C \subseteq G$.
For $H \subseteq G$ let $N_GH = \{g \in G \mid gHg^{-1} = H\}$%
\indexnotation{N_GH}
be its \emph{normalizer},%
\index{normalizer}
let $Z_GH = \{g \in G \mid ghg^{-1} = h \text{ for } h \in H\}$%
\indexnotation{Z_GH}
be its \emph{centralizer},%
\index{centralizer}
and put 
$$W_GH%
\indexnotation{W_GH}
:= N_GH/(H \cdot Z_GH),$$
where $H \cdot Z_GH$ is the normal subgroup of $N_GH$ consisting of elements of
the form $hu$ for $h \in H$ and $u \in Z_GH$.
Notice that $W_GH$ is finite if $H$ is finite.

Recall that the \emph{Burnside ring}%
\index{Burnside ring}
$A(G)$%
\indexnotation{A(G)}
of a finite group is the Grothendieck group associated to the abelian monoid of isomorphism
classes of finite $G$-sets with respect to the disjoint union.  The
ring multiplication comes from the cartesian product. The zero element
is represented by the empty set, the unit is represented by
$G/G = \pt$. For a finite group $G$ the abelian groups
$K_q(C^*_r(G))$, $K_q(RG)$ 
and $L^{\langle - \infty \rangle}(RG)$ become modules
over $A(G)$ because these
functors come with a Mackey structure and $[G/H]$ acts by $\ind_H^G
\circ \res_G^H$.   

We obtain a ring homomorphism
$$\chi^G \colon  A(G) \to \prod_{(H) \in \calfin} \IZ, \hspace{5mm}
[S] \mapsto (|S^H|)_{(H) \in \calfin}$$
which sends the class of a finite $G$-set $S$ to the element given by
the cardinalities of the $H$-fixed point sets. This is an injection with
finite cokernel. This leads to an isomorphism of $\IQ$-algebras
\begin{eqnarray}
\chi_{\IQ}^G :=  \chi^G \otimes_{\IZ} \id_{\IQ} \colon
A(G) \otimes_{\IZ} \IQ & \xrightarrow{\cong} & \prod_{(H) \in (\calfin)} \IQ.
\label{chi_{bbQ}^G}
\end{eqnarray}
For a finite cyclic group $C$ let 
\begin{eqnarray}
& \theta_C ~ \in ~ A(C) \otimes_{\IZ} \IZ[1/|C|] &
\label{theta_C}
\end{eqnarray}
be the element which is sent under
the isomorphism $\chi^C_{\IQ}\colon A(C) \otimes_{\IZ} \IQ \xrightarrow{\cong}
\prod_{(H)\in \calfin} \IQ$ of \eqref{chi_{bbQ}^G}
to the element, whose entry  is one if $(H) = (C)$ and is zero if $(H) \not= (C)$.
Notice that $\theta_C$ is an idempotent. In particular 
we get a direct summand
$\theta_C \cdot K_q(C^*_r(C))\otimes_{\IZ} \IQ$ in
$K_q(C^*_r(C))\otimes_{\IZ} \IQ$ and analogously for  $K_q(RC) \otimes_{\IZ} \IQ$ 
and $L^{\langle - \infty \rangle}(RG)\otimes_{\IZ} \IQ$.


\subsubsection{Rationalized Topological $K$-Theory for Infinite Groups}
\label{subsec: Rationalized Topological K-Theory for Ininite Groups}

The next result is taken from \cite[Theorem 0.4 and page
127]{Lueck(2002d)}. Recall that $\Lambda^G$ is the ring $\IZ \subseteq \Lambda^G \subseteq \IQ$
which is obtained from $\IZ$ by inverting the orders of the finite
subgroups of $G$.

\begin{theorem}[Rational Computation of Topological $K$-Theory for Infinite Groups]
\label{the: Rational Computation of Topological K-Theory for Infinite Groups}
\indextheorem{Rational Computation of Topological $K$-Theory for Infinite Groups}
Suppose that the group $G$ satisfies the 
Baum-Connes Conjecture ~\ref{con: Baum-Connes Conjecture}. 
Then there is an isomorphism 
\begin{multline*}
\bigoplus_{p+q = n} ~ \bigoplus_{(C) \in (\calfcyc)}
K_p(BZ_GC) \otimes_{\IZ[W_G C]} \theta_C \cdot
K_q(C^*_r(C))\otimes_{\IZ} \Lambda^G 
\\
\xrightarrow{\cong} ~
K_n(C^*_r(G)) \otimes_{\IZ} \Lambda^G.
\end{multline*}

If we tensor with $\IQ$, we get an isomorphism
\begin{multline*}
\bigoplus_{p+q = n} ~ \bigoplus_{(C) \in (\calfcyc)}
H_p(BZ_GC;\IQ) \otimes_{\IQ [W_GC]} \theta_C \cdot K_q(C^*_r(C)) \otimes_{\IZ} \IQ.
\\
\xrightarrow{\cong} ~
K_n(C^*_r(G)) \otimes_{\IZ} \IQ.
\end{multline*}
\end{theorem}

\subsubsection{Rationalized Algebraic $K$-Theory for Infinite Groups}
\label{subsec: Rationalized Algebraic K-Theory for Ininite Groups}

Recall that for algebraic $K$-theory of the integral group ring we know 
because of Proposition~\ref{pro: rational vanishing of lower nil} that in the Farrell-Jones Conjecture
we can reduce to the family of finite subgroups. A reduction to the family of finite subgroups also works
if the coefficient ring is a regular $\IQ$-algebra, compare \ref{pro: K for special coefficients}.
The next result is a variation of \cite[Theorem 0.4]{Lueck(2002b)}.

\begin{theorem}[Rational Computation of Algebraic $K$-Theory] 
\label{the: Rational Computation of Algebraic K-Theory for Infinite Groups}
\indextheorem{Rational Computation of Algebraic $K$-Theory for Infinite Groups}
Suppose 
that the group $G$ satisfies the 
Farrell-Jones Conjecture~\ref{con: Farrell-Jones Conjecture} 
for the algebraic $K$-theory of $RG$, where
either $R=\IZ$ or $R$ is a regular ring with $\IQ \subset R$. Then we get an isomorphism
\begin{multline*}
\bigoplus_{p+q = n} ~ \bigoplus_{(C) \in (\calfcyc)}
H_p(BZ_GC;\IQ) \otimes_{\IQ [W_GC]} 
\theta_C \cdot K_q(RC) \otimes_{\IZ} \IQ
\\
\xrightarrow{\cong} ~
K_n(RG) \otimes_{\IZ} \IQ.
\end{multline*}
\end{theorem}

\begin{remarknew} \label{rem: Theorem refRational Computation of Algebraic $K$-Theory ... for arbitrary R}
If in Theorem~\ref{the: Rational Computation of Algebraic K-Theory for Infinite Groups}
we assume the Farrell-Jones Conjecture for the algebraic $K$-theory of $RG$ but make no 
assumption on the coefficient ring $R$, then we still obtain 
that the map appearing there is split injective. 
\end{remarknew}

\begin{examplenew}[The Comparison Map from Algebraic to Topological $K$-theory]
\label{exa: K^{alg}(CG) and K(C^*_r(G))} 
If we consider $R = \IC$ as coefficient ring and apply $- \otimes_{\IZ} \IC$ 
instead of $- \otimes_{\IZ} \IQ$ , the formulas
simplify. Suppose that $G$ satisfies the Baum-Connes Conjecture 
\ref{con: Baum-Connes Conjecture} and the Farrell-Jones Conjecture 
\ref{con: Farrell-Jones Conjecture} for algebraic $K$-theory with $\IC$ as coefficient ring.
Recall that $\con(G)_f$ is the set of conjugacy
classes $(g)$ of elements $g \in G$ of finite order. We denote for $g \in G$ by
$\langle g \rangle$%
\indexnotation{langle g rangle}
the cyclic subgroup generated by $g$.

Then we get the following commutative square, whose horizontal maps are isomorphisms and
whose vertical maps are induced by the obvious change of theory homorphisms
(see \cite[Theorem~0.5]{Lueck(2002b)}) 
$$\comsquare{\bigoplus_{p+q=n} \bigoplus_{(g) \in \con(G)_f}
H_p(Z_G\langle g \rangle;\IC) \otimes_{\IZ}  K_q (\IC)}{\cong}
{K_n(\IC G) \otimes_{\IZ} \IC}{}{}
{\bigoplus_{p+q= n} \bigoplus_{(g) \in \con(G)_f}
H_p(Z_G\langle g \rangle;\IC) \otimes_{\IZ}  K_q^{\topo}(\IC)}{\cong}
{K_n(C_r^*(G)) \otimes_{\IZ} \IC}
$$
The Chern character appearing in the lower row of the commutative square above has already
been constructed by different methods in \cite{Baum-Connes(1988a)}. The construction in \cite{Lueck(2002b)}
works also for $\IQ$ (and even smaller rings) and other theories like algebraic $K$- and
$L$-theory. This is important for the proof of 
Theorem~\ref{the: Baum Connes Conjecture implies modified Trace Conjecture} and to get the
commutative square above.
\end{examplenew}

\begin{examplenew}[A Formula for $K_0(\IZ G) \otimes_{\IZ} \IQ$]
\label{exa: Formula for K_0(IZ G)}
Suppose that the Farrell-Jones Conjecture is true rationally for $K_0( \IZ G )$, i.e.\ the assembly map 
\[
A_{\calvcyc} \colon H_0^G(\EGF{G}{\calvcyc};\bfK_{\IZ}) \otimes_{\IZ} \IQ  \to 
K_0(\IZ G)\otimes_{\IZ} \IQ
\]
is an isomorphism. Then we obtain 
\begin{eqnarray*}
& K_0(\IZ G) \otimes_{\IZ} \IQ  \cong  &\\
& K_0(\IZ) \otimes_{\IZ} \IQ  \oplus \bigoplus_{(C) \in (\calfcyc)}
H_1(BZ_GC;\IQ) \otimes_{\IQ [W_GC]} 
\theta_C \cdot K_{-1}(RC) \otimes_{\IZ} \IQ & .
\end{eqnarray*}

Notice that $\widetilde{K}_0(\IZ G) \otimes_{\IZ} \IQ$ contains only contributions
from $K_{-1}(\IZ C) \otimes_{\IZ} \IQ$ for finite cyclic subgroups $C \subseteq G$.
\end{examplenew}

\begin{remarknew}
\label{rem: already interesting for finite}
Note that these statements are interesting already for finite groups.
For instance Theorem~\ref{the: Rational Computation of Topological K-Theory for Infinite Groups}
yields for a finite group $G$ and $R = \IC$ an isomorphism
$$\bigoplus_{(C) \in (\calfcyc)}
\Lambda_G \otimes_{\Lambda_G [W_GC]} 
\theta_C \cdot R(C) \otimes_{\IZ} \Lambda_G 
~ \cong ~ R(G) \otimes_{\IZ} \Lambda_G.$$
which in turn implies Artin's Theorem discussed in
Remark \ref{rem: isomorphisms conjectures as induction theorem}.
\end{remarknew}


\subsubsection{Rationalized Algebraic $L$-Theory for Infinite Groups}
\label{subsec: Rationalized Algebraic L-Theory for Infinite Groups}

Here is the $L$-theory analogue of the results above. Compare \cite[Theorem 0.4]{Lueck(2002b)}.

\begin{theorem}[Rational Computation of Algebraic $L$-Theory for Infinite Groups] 
\label{the: Rational Computation of Algebraic L-Theory for Infinite Groups}
\indextheorem{Rational Computation of Algebraic $K$-Theory for Infinite Groups}
Suppose that the group $G$ satisfies the Farrell-Jones Conjecture  
\ref{con: Farrell-Jones Conjecture} for $L$-theory. Then we get for all
$j \in \IZ, j \le 1$ an isomorphism
\begin{multline*}
\bigoplus_{p+q = n} ~ \bigoplus_{(C) \in (\calfcyc)}
H_p(BZ_GC;\IQ) \otimes_{\IQ [W_GC]} 
\theta_C \cdot L_q^{\langle j \rangle}(RC) \otimes_{\IZ} \IQ
\\
\xrightarrow{\cong} ~
L^{\langle j \rangle}_n(RG) \otimes_{\IZ} \IQ.
\end{multline*}
\end{theorem}

\begin{remarknew}[Separation of Variables] 
\label{rem: Separation of Variables for calfin}
Notice that 
in Theorem~\ref{the: Rational Computation of Topological K-Theory for Infinite Groups}, 
\ref{the: Rational Computation of Algebraic K-Theory for Infinite Groups} and 
\ref{the: Rational Computation of Algebraic L-Theory for Infinite Groups} we see again 
the principle we called
\emph{separation of variables}%
\index{principle!separation of variables}
in Remark~\ref{rem: separation of variables for K-theory in the torsion free case}.
There is a group homology part which is independent of the coefficient ring 
$R$ and the $K$- or $L$-theory under consideration and a part
depending only on the values of the theory under consideration on $RC$ or $C^*_r(C)$
for all finite cyclic subgroups $C \subseteq G$.
\end{remarknew}


\subsection{Integral Computations for Infinite Groups}
\label{sec: Integral Computations for Infinite Groups}

As mentioned above, no general pattern for integral calculations is known or expected.
We mention at least one situation where a certain class of groups can be treated
simultaneously.
Let $\calmfin$%
\indexnotation{calmfin}
be the subset of $\calfin$ consisting
of elements in $\calfin$ which are maximal in $\calfin$. 
Consider the following assertions on the group $G$.\\[1mm]
\begin{description}

\item[(M)] $M_1, M_1 \in \calmfin, M_1 \cap M_2 \not= 1 ~ \Rightarrow ~ M_1 = M_2$;

\item[(NM)] $M \in \calmfin ~ \Rightarrow ~ N_GM = M$;

\item[(VCL$_I$)] If $V$ is an infinite virtually cyclic subgroup of $G$, then
$V$ is of type  I (see Lemma~\ref{lem: virtually cyclic});

\item[(FJK$_N$)]
The Isomorphism Conjecture of Farrell-Jones for algebraic $K$-theory
\ref{con: Farrell-Jones Conjecture}
is true for $\IZ G$ in the range $n \le N$ for a 
fixed element $N \in \IZ \amalg \{\infty\}$,
i.e.\ the assembly map
$
A \colon
\calh^G_n(\EGF{G}{\calvcyc};\bfK_R) ~ \xrightarrow{\cong} ~ K_n(RG)
$
is bijective for $n \in \IZ$ with $n \le N$.

\end{description}

Let $\widetilde{K}_n(C^*_r(H))$%
\indexnotation{widetilde{K}_n(C^*_r(H))}
be the cokernel of the map 
$K_n(C^*_r(\{1\})) \to K_n(C^*_r(H))$ and 
$\overline{L}^{\langle j \rangle}_n(RG)$%
\indexnotation{overline{L}^{langle j rangle}_n(RG)}
be the cokernel of the map 
$L_n^{\langle j \rangle}(R) \to L_n^{\langle j \rangle}(RG)$. 
This coincides with 
$\widetilde{L}^{\langle j \rangle}_n(R)$, 
\indexnotation{widetilde{L}^{langle j rangle}_n(R)}
which is the cokernel of the map 
$L_n^{\langle j \rangle}(\IZ) \to L_n^{\langle j \rangle}(R)$ if $R=\IZ$ but not in general.
Denote by $\Wh_n^R(G)$%
\indexnotation{Wh_n^R(G)}
the $n$-th Whitehead group%
\index{Whitehead group!higher}
of $RG$ which is the $(n-1)$-th homotopy group
of the homotopy fiber of the assembly map 
$BG_+ \wedge \bfK (R) \to \bfK (RG)$. It agrees with the previous defined
notions if $R=\IZ$.
The next result is taken from \cite[Theorem 4.1]{Davis-Lueck(2002a)}.

\begin{theorem} \label{the: computations based on IC}
Let $\IZ \subseteq \Lambda \subseteq \IQ$ be a ring  such that the order of any 
finite subgroup of $G$ is invertible in $\Lambda $.  Let $(\calmfin)$ be the set
of conjugacy classes (H) of subgroups of $G$ such that $H$ belongs to $\calmfin$. Then: 
\begin{enumerate} 

\item \label{the: computations based on IC:  (M), (NM) and (BC)}
If $G$ satisfies (M), (NM) and the Baum-Connes Conjecture 
\ref{con: Baum-Connes Conjecture}, then for $n \in \IZ$ there is an exact sequence of topological $K$-groups
$$
0 \to \bigoplus_{(H) \in (\calmfin)} \widetilde{K}_n(C^*_r(H))
 \to K_n(C_r^*(G)) \to K_n(G\backslash \EGF{G}{\calfin}) \to 0,
$$
which splits after applying $- \otimes_{\IZ} \Lambda$.

\item \label{the: computations based on IC:  (M), (NM), (VCL_I) and (FJL)}
If $G$ satisfies (M), (NM), (VCL$_I$) and the $L$-theory part of the
Farrell-Jones Conjecture~\ref{con: Farrell-Jones Conjecture}, then for all $n \in \IZ$ there is an exact sequence
\begin{multline*}
\ldots \to H_{n+1}(G\backslash \EGF{G}{\calfin};\bfL^{\langle -\infty\rangle}(R)) 
\to \bigoplus_{(H) \in (\calmfin)} \overline{L}_n^{\langle -\infty\rangle}(RH)
\\
\to L_n^{\langle -\infty\rangle}(RG) \to 
H_n(G\backslash \EGF{G}{\calfin};\bfL^{\langle -\infty\rangle}(R)) \to \ldots
\end{multline*}
It splits after applying $- \otimes_{\IZ} \Lambda$, more precisely
$$L_n^{\langle -\infty\rangle}(RG) \otimes_{\IZ} \Lambda 
\to H_n(G\backslash \EGF{G}{\calfin};\bfL^{\langle -\infty\rangle}(R)) \otimes_{\IZ} \Lambda$$
is a split-surjective map of $\Lambda$-modules.

\item \label{the: computations based on IC:  (M), (NM) and [(FJL[1/2])} If $G$ satisfies (M), (NM) and the Farrell-Jones Conjecture 
\ref{con: Farrell-Jones Conjecture} for $L_n(RG){[1/2]}$, then the conclusion of
assertion~\ref{the: computations based on IC:  (M), (NM), (VCL_I) and (FJL)} still holds
if we invert $2$ everywhere. Moreover, in the case $R = \IZ$ 
the sequence reduces to a short exact sequence
\begin{multline*}
0 \to \bigoplus_{(H) \in (\calmfin)} \widetilde{L}^{\langle j \rangle}_n(\IZ H)[\frac{1}{2}]
 \to L_n^{\langle j \rangle}(\IZ G)[\frac{1}{2}]
\\
\to H_n(G\backslash \EGF{G}{\calfin};\bfL(\IZ)[\frac{1}{2}] \to 0,
\end{multline*}
which splits after applying 
$- \otimes_{\IZ[\frac{1}{2}]}\Lambda[\frac{1}{2}]$.

\item \label{the: computations based on IC:  (M), (NM), and (FJK_N)}
If $G$ satisfies (M), (NM),  and (FJK$_N$), then there is for
$n \in \IZ, n \le N$ an isomorphism 
$$H_n( \EGF{G}{\calvcyc} , \EGF{G}{\calfin} ; \bfK_R ) \oplus 
\bigoplus_{(H) \in (\calmfin)} \Wh_n^R(H) \xrightarrow{\cong} \Wh_n^R(G),$$
where $\Wh_n^R(H) \to \Wh_n^R(G)$ is induced by the inclusion $H \to G$.
\end{enumerate}

\end{theorem}

\begin{remarknew}[Role of $G\backslash\EGF{G}{\calfin}$] 
\label{rem: Role of E_FIN G}
Theorem~\ref{the: computations based on IC} illustrates that for 
such computations a good understanding of the geometry of 
the orbit space $G\backslash \EGF{G}{\calfin}$ is
necessary.
\end{remarknew}

\begin{remarknew} \label{rem: examples for groups with (MFIN),..}
In \cite{Davis-Lueck(2002a)} it is explained that the following classes of groups
do satisfy the assumption appearing in Theorem~\ref{the: computations based on IC} 
and what the conclusions are in the case $R = \IZ$.  Some of these cases have been treated
earlier in \cite{Berkhove-Juan-Pineda-Pearson(2001)}, 
\cite{Lueck-Stamm(2000)}.

\begin{itemize}
\item Extensions $1 \to \IZ^n \to G \to F \to 1$ for finite $F$ such that the conjugation
  action of $F$ on $\IZ^n$ is free outside $0 \in \IZ^n$;
\item Fuchsian groups $F$;
\item One-relator groups $G$.
\end{itemize}

Theorem~\ref{the: computations based on IC} is generalized in
\cite{Lueck(2003a)}  
in order to treat for instance the semi-direct product of the discrete 
three-dimensional Heisenberg group by $\IZ/4$. For this group
$G\backslash \EGF{G}{\calfin}$ is $S^3$.

A calculation for $2$-dimensional crystallographic groups and more general
cocompact NEC-groups is presented in \cite{Lueck-Stamm(2000)}
(see also \cite{Pearson(1998)}). For these groups the orbit spaces
$G\backslash \EGF{G}{\calfin}$ are compact surfaces possibly with boundary.
\end{remarknew}

\begin{examplenew} Let $F$ be a cocompact Fuchsian group with presentation
\begin{multline*}
F= \langle a_1,b_1,\ldots,a_g,b_g,c_1,\ldots,c_t \mid
\\
c_1^{\gamma_1}=\ldots =c_t^{\gamma_t}=c_1^{-1}\cdots c_t^{-1}[a_1,b_1]
\cdots [a_g,b_g] =1 \rangle
\end{multline*}
for integers $g,t \ge 0$ and $\gamma_i > 1$. Then $G\backslash\EGF{G}{\calfin}$
is a closed orientable surface of genus $g$. The following 
is a consequence of Theorem~\ref{the: computations based on IC}
(see \cite{Lueck-Stamm(2000)} for more details).

\begin{itemize}

\item 
There are isomorphisms
$$K_n(C^*_r(F)) \cong \left\{
\begin{array}{lcl}
\left(2 + \sum_{i=1}^t (\gamma_i - 1)\right) \cdot \IZ
& & n = 0;
\\
(2g) \cdot \IZ
& & n = 1.
\end{array}\right.
$$

\item 
The inclusions of the maximal subgroups $\IZ/\gamma_i = \langle c_i \rangle$
induce an isomorphism
$$\bigoplus_{i=1}^t \Wh_n(\IZ/\gamma_i) \xrightarrow{\cong} \Wh_n(F)$$
for $n \le 1$. 

\item 
There are isomorphisms
$$L_n(\IZ F)[1/2] \cong \left\{
\begin{array}{lcl}
\left(1 + \sum_{i=1}^t \left[\frac{\gamma_i}{2}\right]\right) \cdot
\IZ[1/2]
& & n \equiv 0 \hspace{2mm} (4);
\\
(2g) \cdot \IZ[1/2]
& & n \equiv 1 \hspace{2mm} (4);
\\
\left(1 + \sum_{i=1}^t \left[\frac{\gamma_i - 1}{2}\right]\right) \cdot
\IZ[1/2]
& & n \equiv 2 \hspace{2mm} (4);
\\
0
& & n \equiv 3 \hspace{2mm} (4),
\end{array}\right.
$$
where $[r]$ for $r \in \IR$ denotes the largest integer less
than or equal to $r$.

From now on suppose that each $\gamma_i$ is odd.
Then the number $m$ above is odd and we get for 
for $\epsilon =p$ and $s$
$$L_n^{\epsilon}(\IZ F) \cong \left\{
\begin{array}{lcl}
\IZ/2 \bigoplus \left(1 + \sum_{i=1}^t \frac{\gamma_i-1}{2}\right) \cdot \IZ
& & n \equiv 0 \hspace{2mm} (4);
\\
(2g) \cdot \IZ
& & n \equiv 1 \hspace{2mm} (4);
\\
\IZ/2 \bigoplus  \left(1 + \sum_{i=1}^t \frac{\gamma_i - 1}{2}\right) \cdot \IZ
& & q \equiv 2 \hspace{2mm} (4);
\\
(2g) \cdot \IZ/2
& & n \equiv 3 \hspace{2mm} (4).
\end{array}\right.
$$
For $\epsilon =h$ we do not know an explicit formula.
The problem is that no general formula is known for the $2$-torsion
contained in $\widetilde{L}^h_{2q}(\IZ[\IZ/m])$, for  $m$ odd,
since it is given by the term
$\widehat{H}^2(\IZ/2; \widetilde{K}_0(\IZ[\IZ/m]))$,
see \cite[Theorem~2]{Bak(1976)}.

\end{itemize}
\end{examplenew}

Information about the left hand side of the Farrell-Jones assembly map for algebraic $K$-theory
in the case where $G$ is $SL_3( \IZ)$ can be found in \cite{Upadhyay(1996)}.


\subsection{Techniques for Computations}
\label{sec: Techniques for Computations}

We briefly outline some methods that are fundamental for computations and for the proofs
of some of the theorems above. 


\subsubsection{Equivariant Atiyah-Hirzebruch Spectral Sequence}
\label{subsec: Equivariant Atiyah-Hirzebruch Spectral Sequence}

Let $\calh^G_*$ be a $G$-homology theory with values in
$\Lambda$-modules.  Then there are two
spectral sequences which can be used to compute it. The first
one is the rather obvious equivariant version of the
\emph{Atiyah-Hirzebuch spectral sequence.}%
\index{spectral sequence!equivariant Atiyah-Hirzebruch spectral sequence}
It converges to $\calh_n^G(X)$ and its $E^2$-term is given in terms of
Bredon homology
$$E^2_{p,q} ~ = ~ H^G_p(X;\calh^G_q(G/H))$$
of $X$ with respect to the coefficient system, which is given by the covariant functor
$\Or(G) \to \Lambda\text{-}\MODULES, \; G/H \mapsto \calh^G_q(G/H)$.
More details can be found for instance in \cite[Theorem~4.7]{Davis-Lueck(1998)}.


\subsubsection{$p$-Chain Spectral Sequence}
\label{subsec: p-Chain Spectral Sequence}

There is another spectral sequence, the \emph{$p$-chain spectral sequence}%
\index{spectral sequence!p-chain spectral sequence@$p$-chain spectral sequence}
\cite{Davis-Lueck(2002a)}. Consider a covariant functor 
$\bfE \colon \Or(G) \to \SPECTRA$. It defines a $G$-homology theory $\calh^G_*(-;\bfE)$
(see Proposition~\ref{pro: Or(G)-spectra yield a G-homology theory}).
The $p$-chain spectral sequence converges to $\calh^G_n(X)$ 
but has a different setup and in particular a different $E^2$-term than the equivariant
Atiyah-Hirzebruch spectral sequence.  We describe the $E^1$-term for simplicity
only for a proper $G$-$CW$-complex.

A \emph{$p$-chain}%
\index{p-chain@$p$-chain}
is a sequence of conjugacy classes of finite subgroups
$$(H_0) < \ldots  <  (H_p)$$ 
where $(H_{i-1}) < (H_i)$ means that
$H_{i-1}$ is subconjugate, but not
conjugate to $(H_i)$. Notice for the sequel
that the group of automorphism of $G/H$ in $\Or(G)$ is isomorphic to
$N\!H/H$. To such a $p$-chain there is associated the
$N\!H_p/H_p$-$N\!H_0/H_0$-set
\begin{multline*}
S((H_0) < \ldots < (H_p)) ~ = ~ 
\map(G/H_{p-1},G/H_p)^G \times_{N\!H_{p-1}/H_{p-1}}
\\ 
\ldots \times_{N\!H_1/H_1} \map(G/H_0,G/H_1)^G. 
\end{multline*}
The $E^1$-term $E^1_{p,q}$ of the $p$-chain spectral sequence is
$$
\bigoplus_{(H_0) < \ldots <  (H_p)}~
\pi_q\left(\left(X^{H_p} \times_{N\!H_p/H_p} S((H_0) <\ldots <  (H_p))\right)_+ \wedge_{N\!H_0/H_0}
  \bfE(G/H_0)\right)
$$
where $Y_+$ means the pointed space obtained from $Y$ by adjoining an extra base point.
There are many situations where the $p$-chain spectral sequence is much more useful
than the equivariant Atiyah-Hirzebruch spectral sequence. Sometimes a combination of
both is necessary to carry through the desired calculation.


\subsubsection{Equivariant Chern Characters}
\label{subsubsec: Equivariant Chern Characters}

Equivariant Chern characters have been studied in \cite{Lueck(2002b)}
and \cite{Lueck(2002d)} and allow to compute equivariant homology
theories for proper $G$-$CW$-complexes. 
The existence of the equivariant Chern character  says that under
certain conditions
the Atiyah-Hirzebruch spectral
sequence collapses and, indeed,  the source of the equivariant Chern character is canonically
isomorphic to $\bigoplus_{p+q} E^2_{p,q}$, where $E^2_{p,q}$ is the $E^2$-term 
of the equivariant Atiyah-Hirzebruch spectral sequence.

The results of Section \ref{sec: Rational Computations for Infinite Groups}
are essentially proved by applying the equivariant Chern character to
the source of the assembly map for the family of finite subgroups.

\typeout{-------------------- References -------------------------------}

\addcontentsline{toc}{section}{References}
\bibliographystyle{abbrv}
\bibliography{dbdef,dbpub,dbpre,dbextrabcsfinal}

\typeout{-------------------- Notation ---------------------------------}
\twocolumn
\section*{Notation}
\addcontentsline{toc}{section}{Notation}


\noindent
\entry{$A(G)$}{A(G)}\\
\entry{$\bfA (X)$}{bfA(X)}\\
\entry{$A_{\calf}$}{A_calf}\\
\entry{$\widehat{A}_x(M,u)$}{widehat{A}_x(M)}\\
\entry{$BG$}{BG}\\
\entry{$\class_0(G)$}{class_0(G)}\\
\entry{$\class_0(G)_f$}{class_0(G)_f}\\
\entry{$\con(G)$}{con(G)}\\
\entry{$\con(G)_f$}{con(G)_f}\\
\entry{$\cone(X)$}{cone(X)}\\
\entry{$C^{\ast}\text{-}\ALGEBRAS$}{C^*-ALGEBRAS}\\
\entry{$C^*_{\max}(G)$}{C^*_{max}(G)}\\
\entry{$C^*_r(G)$}{C^*_r(G)}\\
\entry{$C^*_r(G;\IR)$}{C^*_r(G;IR)}\\
\entry{$C_0(X)$}{C_0(X)}\\
\entry{$EG$}{EG}\\
\entry{$\EGF{G}{\calf}$}{E(G;calf)}\\
\entry{$\underline{E}G$}{underline{E}G}\\
\entry{$G_n(R)$}{G_n(R)}\\
\entry{$(g)$}{(g)}\\
\entry{$\langle g \rangle$}{langle g rangle}\\
\entry{$\GROUPS$}{GROUPS}\\
\entry{$\GROUPS^{\inj}$}{GROUPS^inj}\\
\entry{$\GROUPOIDS$}{GROUPOIDS}\\
\entry{$\GROUPOIDS^{\finker}$}{GROUPOIDS^{finker}}\\
\entry{$\GROUPOIDS^{\inj}$ }{GROUPOIDS^inj}\\
\entry{$H^G_*(-;\bfE)$}{H^G_*(-;bfE)}\\
\entry{$H_*^?(-;\bfE)$}{H_*^?(-;bfE)}\\
\entry{$H_*^G(X;M)$}{H_*^G(X;M)}\\
\entry{$\HS_{\IC G}$}{HS_IC G}\\
\entry{$\inn(K)$}{inn(K)}\\
\entry{$K_n(A)$}{K_n(A)}\\
\entry{$KK_n(A,B)$}{KK_n(A,B)}\\
\entry{$KK_n^G(A,B)$}{KK_n^G(A,B)}\\
\entry{$\widetilde{K}_n(C^*_r(H))$}{widetilde{K}_n(C^*_r(H))}\\
\entry{$K_n(\caln \! il (R))$}{K_n(calnil(R))}\\
\entry{$K_n(R)$}{K_n(R)}\\
\entry{$\widetilde{K}_n (R)$}{widetildeK_n(R)}\\
\entry{$K^G_n(C_0(X))$}{K^G_n(C_0(X))}\\
\entry{$K_n^G(X)$}{K_n^G(X)}\\
\entry{$K_*(Y)$}{K_*(Y)}\\
\entry{$KO_n(C^*_r(G;\IR))$}{KO_n(C^*_r(G;R)}\\
\entry{$KO_*(Y)$}{KO_n(Y)}\\
\entry{$L_n^{\langle - \infty \rangle}(R)$}{L_n^j(R)}\\
\entry{$\widetilde{L}^{\langle j \rangle}_n(R)$}{widetilde{L}^{langle j rangle}_n(R)}\\
\entry{$\overline{L}^{\langle j \rangle}_n(RG)$}{overline{L}^{langle j rangle}_n(RG)}\\
\entry{$\bfL^{\langle j  \rangle}(R)$}{bfL^j(R)}\\
\entry{$\map_G(-,Y)$}{map_G(-,Y)}\\
\entry{$N_GH$}{N_GH}\\
\entry{$N\!K_n(R)$}{NK_n(R)}\\
\entry{$\Nil_{n-1} (R)$}{Nil_n(R)}\\
\entry{$\Or(G)$}{Or(G)}\\
\entry{$\OrGF{G}{\calf}$}{Or(G;calf)}\\
\entry{$\widetilde{o}(X)$}{widetilde o(X)}\\
\entry{$P ( M )$}{P(M)}\\
\entry{$\bfP (M)$}{bfP(M)}\\
\entry{$P_b( M ; \IR^k )$}{P_b(M;R^k)}\\
\entry{$\pt$}{pt}\\
\entry{$\RINGS$}{RINGS}\\
\entry{$\RINGS^{\inv}$}{RINGS^inv}\\
\entry{$r_F(G)$}{r_F(G)} \\
\entry{$r_{\IR}(G;\IR)$}{r_{IR}(G;IR)}\\
\entry{$r_{\IR}(G;\IC)$}{r_{IR}(G;IC)}\\
\entry{$r_{\IR}(G;\IH)$}{r_{IR}(G;IH)}\\
\entry{$R(G)$}{R(G)}\\
\entry{$RO(G)$}{RO(G)}\\
\entry{$\sign_x(M)$}{sign_x(M)}\\
\entry{$\SPACES$}{SPACES}\\
\entry{$\SPACES^+$}{SPACES^+}\\
\entry{$\SPECTRA$}{SPECTRA}\\
\entry{$\SubGF{G}{\calf}$}{sub_FG}\\
\entry{$\tr_{C^*_r(G)}$}{tr_{C^*_r(H)}}\\
\entry{$\tr^u_{\IC G}$}{tr^u_IC G}\\
\entry{$(U,\rho,F)$}{(U,rho,F)}\\
\entry{$W_GH$}{W_GH}\\
\entry{$\Wh(G)$}{Wh(G)}\\
\entry{$\Wh_0^R(G)$}{Wh_0^R(G)}\\
\entry{$\Wh_1^R(G)$}{Wh_1^R(G)}\\
\entry{$\Wh_n^R(G)$}{Wh_n^R(G)}\\
\entry{$\Wh^{\Diff}( X )$}{Wh^Diff(X)}\\
\entry{$\bfWh^{PL}( X )$}{bfWh^PL(X)}\\
\entry{$X_+$}{X_+}\\
\entry{$X \wedge \bfE$}{X wedge bfE}\\
\entry{$X \wedge Y$}{X wedge Y}\\
\entry{$X \times_{\calc} Y$}{X times_calc Y}\\
\entry{$X \wedge_{\calc} Y$}{X wedge_calc Y}\\
\entry{$X \wedge_{\calc} \bfE$}{X wedge_calc bfE}\\
\entry{$Z_GH$}{Z_GH}\\
\entry{$\IZ_p$}{bbZ_p}\\
\entry{$\eta(M)$}{eta(M)}\\
\entry{$\eta^{(2)}(\widetilde{M})$}{eta^(2)(widetilde M)}\\
\entry{$\Lambda^G$}{Lambda^G}\\
\entry{$\pi_i(\mathbf{E})$}{pi_i(bfE)}\\
\entry{$\rho^{(2)}(M)$}{rho^(2)(M)}\\
\entry{$\tau^{(2)}(M) \in \IR$}{tau^(2)}\\
\entry{$\widehat{\cala}(M)$}{widehat cala(M)}\\
\entry{$\calall$}{familycalall}\\
\entry{$\calcyc$}{calcyc}\\
\entry{$\calfcyc$}{familycalfcyc}\\
\entry{$(\calfcyc)$}{(calfcyc)}\\
\entry{$\calfin$}{familycalfin}\\
\entry{$\calg^G(S)$}{calg^G(S)}\\
\entry{$\calh(G)$}{calh(Gamma)}\\
\entry{$\calh_*^G$}{calh_*^G}\\
\entry{$\calh_*^?$}{calh_*^?}\\
\entry{$\call(M)$}{call(M)}\\
\entry{$\calmfin$}{calmfin}\\
\entry{$\calp ( M )$}{calp(M)}\\
\entry{$\calp^{\Diff} ( M )$}{calp^Diff(M)}\\
\entry{$\calvcyc$}{familycalvcyc}\\
\entry{$\calvcyc_{I}$}{calvcyc_I}\\
\entry{$\caltr$}{trivial family}

\onecolumn
 
\typeout{-------------------- Index ---------------------------------}

\flushbottom
\addcontentsline{toc}{section}{Index}
\printindex                                  

\end{document}